\documentclass[11pt]{article}

\usepackage{amsmath,amssymb,amsthm,mathtools}
\usepackage{dsfont}
\usepackage[margin=1.25in]{geometry}
\usepackage{graphicx}
\usepackage{microtype}
\usepackage{longtable}
\usepackage{array}
\usepackage{tikz}
\usetikzlibrary{arrows.meta,matrix,fit,positioning,backgrounds}
\usepackage{tikz-cd}
\usepackage{hyperref}

\DeclareMathOperator{\tr}{tr}

\DeclareMathOperator{\Stab}{Stab}
\DeclareMathOperator{\stab}{\mathfrak{stab}}

\DeclareMathOperator{\End}{End}
\DeclareMathOperator{\vol}{vol}

\DeclareMathOperator{\pr}{pr}
\DeclareMathOperator{\Id}{id}

\DeclareMathOperator{\diag}{diag}

\DeclareMathOperator{\Span}{span}

\DeclareMathOperator{\Sym}{Sym}

\DeclareMathOperator{\Lie}{Lie}
\renewcommand{\Im}{\mathrm{Im}}

\newcommand{\transp}{\top}
\newcommand{\ii}{\mathrm{i}}

\newcommand{\tadj}{{T_W}}

\newcommand{\R}{\mathbb{R}}
\newcommand{\C}{\mathbb{C}}

\renewcommand{\O}{\mathbb{O}}

\newcommand{\VC}{V_{\mathbb C}}
\newcommand{\WC}{W_{\mathbb C}}
\newcommand{\WdualC}{W_{\mathbb C}^{*}}
\newcommand{\gtwo}{\mathfrak g_2}
\newcommand{\esix}{\mathfrak e_6}

\newcommand{\gtwoC}{\mathfrak g_2^{\mathbb C}}
\newcommand{\liegl}{\mathfrak{gl}}

\newcommand{\lieso}{\mathfrak{so}}

\newtheorem{theorem}{Theorem}[section]
\newtheorem{proposition}[theorem]{Proposition}
\newtheorem{lemma}[theorem]{Lemma}
\newtheorem{corollary}[theorem]{Corollary}

\theoremstyle{definition}
\newtheorem{definition}[theorem]{Definition}
\newtheorem{example}[theorem]{Example}

\theoremstyle{remark}
\newtheorem{remark}[theorem]{Remark}

\numberwithin{equation}{section}

\title{Exceptional cubics on 7-manifolds}
\author{Sergey Grigorian\\University of Texas Rio Grande Valley}
\date{}

\begin{document}
\maketitle

\begin{abstract}
A $G_2$-structure $\varphi$ on a $7$-manifold $\mathcal M$ determines a two-dimensional space of invariant totally symmetric cubic forms on each fiber of the rank-$27$ bundle $S^2_0T^*\mathcal M$ of traceless symmetric $2$-tensors. This space is spanned by the trace cubic $C_0$ and a cubic $C_\varphi$ that, up to a constant factor, is the pointwise integrand of the Yukawa coupling on the $G_2$ moduli space. We determine the exact infinitesimal stabilizers over $\mathbb C$ of all nonzero cubics in this family. For $C_0+\lambda C_\varphi$, the stabilizer has dimension $78$ precisely at the two non-real values $\lambda_\pm=(1\pm\ii\sqrt7)/8$. Its dimension is $21$ for $C_0$ and $14$ in all other cases. At $\lambda_\pm$, the stabilizer is simple of type $E_6$, and the corresponding cubic is linearly equivalent to the complex Albert algebra cubic norm.

The two cubics at $\lambda_\pm$, which we call critical, define complex-conjugate $E_6$-reductions of the complexified bundle. At any element with nonzero critical norm, the corresponding cubic determines a complex Albert algebra with basepoint stabilizer $F_4$; a line subbundle of the tangent bundle selects a smooth bundle of these algebras with tensorial operations and a Euclidean real form. The induced Levi-Civita connection preserves either critical cubic exactly when the $G_2$-structure is torsion-free.
Its orthogonal projections onto the two exceptional reductions
coincide with the connection induced by the minimal
$G_2$-connection.

\medskip
\noindent\textit{Keywords:} $G_2$-structures; invariant cubic forms; Yukawa coupling; exceptional Lie algebras; $E_6$; $F_4$; Albert algebra; exceptional Jordan algebras; cubic norms; reductions of the structure group; torsion-free $G_2$-structures.

\smallskip
\noindent\textit{Mathematics Subject Classification (2020):} Primary 53C10; Secondary 17B25, 17C40, 53C05, 53C29.
\end{abstract}
\tableofcontents

\section{Introduction}
\label{sec:introduction}

On a compact torsion-free $G_2$-manifold, Joyce showed that the moduli space of torsion-free $G_2$-structures is locally modeled on harmonic $3$-forms \cite{Joycebook}. These forms decompose into irreducible $G_2$-types. Through the identification with traceless symmetric tensors, the harmonic $\mathbf{27}$-sector corresponds to traceless metric deformations, as described by Grigorian and Yau \cite{GrigorianYau1} and Karigiannis and Leung \cite{karigiannis-2007a}. Lee and Leung introduced a natural symmetric cubic on the moduli space, called the Yukawa coupling \cite{Lee:2002fa}. On harmonic $\mathbf{27}$-directions, it is obtained by integrating a pointwise cubic in the corresponding traceless symmetric tensors \cite{GrigorianYau1,karigiannis-2007a}.

In their 2002 preprint, Lee and Leung also proposed an exceptional algebraic interpretation of the Yukawa coupling. In \cite[\S1.2.1, Remark]{Lee:2002fa}, they suggested that the cubic structure on the $27$-dimensional $G_2$-module $\Lambda^3_{27}$ should arise from the exceptional Jordan algebra, also known as the Albert algebra, with an underlying symmetry group of type $E_6$. They further proposed that when the holonomy reduces from $G_2$ to $SU(3)$, as on $\mathcal M=\mathcal Y\times S^1$ for a Calabi--Yau threefold $\mathcal Y$, the distinguished direction associated with the Kähler form gives a splitting $\mathbf{27}=\mathbf1\oplus\mathbf{26}$ and a corresponding reduction from $E_6$ to $F_4$. In his thesis, the author also revisited this proposal and asked whether the natural
product on $\Lambda^3_{27}$ could be modified to satisfy the Jordan identity
\cite[\S4.4.1, pp.~129--130]{GrigorianPhD}. The Albert algebra and the Jordan identity are recalled below. In Section~\ref{subsec:invariant-cubic-plane} we describe the two $G_2$-equivariant commutative products on \(\mathbf{27}\).

Almost 25 years after their proposal, we resolve this question at the level of pointwise cubics, although the answer differs from the original proposal in one key respect. The pointwise Yukawa cubic is not itself an Albert norm (Proposition~\ref{prop:yukawa-not-albert}). The exceptional norms arise instead from precisely two non-real linear combinations of the Yukawa and trace cubics. We determine which members of the natural two-dimensional space of $G_2$-invariant cubics possess exceptional symmetry and describe the associated Albert algebras and bundle reductions. These cubics are defined for every $G_2$-structure, without compactness or a torsion hypothesis. Our main result is stated in Theorem~\ref{thm:intro-exceptional-pencil}, while Theorem~\ref{thm:intro-albert-bundles} describes the Albert algebra bundles associated with a choice of smooth line subbundle \(\ell\subset T\mathcal M\). Theorem~\ref{thm:intro-compatible-connections} describes the implications for connections on the induced bundles.

Let $\mathcal M$ be an oriented smooth $7$-manifold with a $G_2$-structure, not necessarily torsion-free, given by
a positive $3$-form $\varphi$. The form
induces a metric $g$ and the rank-$27$ bundle
\begin{equation*}
  \mathcal W_\varphi=S^2_0T^*\mathcal M\cong\Lambda^3_{27}T^*\mathcal M.
\end{equation*}
At $x\in\mathcal M$, write $V=T_x\mathcal M$ and
$W=S^2_0V^*\cong\Lambda^3_{27}V^*$. We work in the symmetric-tensor model and write
$\mathcal W_{\varphi,\mathbb C}=\mathcal W_\varphi\otimes_{\mathbb R}\mathbb C$
for its complexification. Symmetric tensors describe pointwise metric
variations and the corresponding metric-changing variations of $\varphi$.
We explain this deformation interpretation in
Section~\ref{subsec:deformation-module}.

Throughout this paper we use the Einstein summation convention: an index repeated once upstairs
and once downstairs is summed from $1$ to $7$, over all ordered values,
unless another range is specified.
Indices are raised and lowered with $g$. Explicit summation signs are
retained for restricted sums and matrix-entry formulas.
For a traceless symmetric tensor $h$, with indices raised by $g$, we define the
two cubic polynomials
\begin{equation}
\label{eq:intro-cubic-generators}
\begin{aligned}
  C_0(h)&=\tr(h^3),& \text{and}\qquad
  C_\varphi(h)&=\varphi_{abc}\varphi_{def}h^{ad}h^{be}h^{cf}.
\end{aligned}
\end{equation}
Complexification extends all tensors complex-bilinearly, without
complex conjugation. We also use the standard facts on
$G_2$-structures from Bryant \cite{bryant-2003}; the conventions and
the $G_2$-contraction identities are collected in
Section~\ref{sec:deformation-representation}.

The two cubics span the $G_2$-invariant cubic space on $W$, over both
$\mathbb R$ and $\mathbb C$. This dimension count is standard and follows
from the tensor-product decomposition in
\cite[Table~A.121]{feger-2020-lieart}.
Both generators already occur together in
\cite[(4.30b)]{GrigorianG2Review} and
\cite[(5.45)]{GrigorianYau1}, and have natural geometric
interpretations. Up to normalization constants, $C_0$ and
$C_\varphi$ arise as third affine variations of the pointwise
volume on the spaces of metrics and positive $3$-forms,
respectively (Proposition~\ref{prop:affine-volume-cubics}). The cubic $C_\varphi$ also describes the
variation of the metric pairing of fixed $3$-form directions
(Proposition~\ref{prop:form-metric-response}).
The same cubic plane occurs in the expansion of the $G_2$ Hitchin
functional \cite[\S2.2]{deBoer:2007zu}. Via the metric pairing,
it also describes the two-dimensional space of $G_2$-equivariant
commutative products on $W$
\cite[Proposition~3.8]{desmet-2025-nonassociative}.
On the moduli space of torsion-free $G_2$-structures,
the integrated cubic $C_\varphi$ gives, up to normalization,
the Yukawa coupling. Its relationship with the global
volume functional was established by Karigiannis and Leung
for irreducible $G_2$-manifolds \cite{karigiannis-2007a}.
More recently, Langlais showed that, on harmonic
$\mathbf{27}$-directions and under the hypothesis $b^1(\mathcal M)=0$,
the Yukawa coupling agrees, up to a constant, with the third
derivative of the logarithm of volume in period coordinates
\cite[Proposition~5.5]{langlais-2025-g2},
\cite[\S3.2, Proposition~3.4]{langlais-2025-geometry}. We revisit the relationship between these cubics and the
volume functional in Section~\ref{subsec:volume-cubics}.

We define the homogeneous pencil as
\begin{equation*}
  N_{a,b}=aC_0+bC_\varphi,\qquad
  \text{where }[a:b]\in\mathbb P^1(\mathbb C)
  \quad\text{and}\quad \lambda=b/a\quad\text{if }a\ne0.
\end{equation*}
Fixing $a=1$, we also define the affine pencil
\begin{equation*}
      N_\lambda=C_0+\lambda C_\varphi \qquad \text{for} \qquad \lambda \in \C.
\end{equation*}
We show that precisely two members of this pencil have
exceptional $E_6$-symmetry and are linearly equivalent to the
complex Albert norm, and we determine the infinitesimal
stabilizers of all remaining members.
Here the exact stabilizer preserves the cubic itself, rather than only
its one-dimensional span, and $\stab(N)$ denotes its Lie algebra
(Definition~\ref{def:exact-stabilizer}). We write $E_6(\mathbb C)$ for
the simply connected complex group of type $E_6$. Its two mutually dual
minuscule representations of dimension $27$ are irreducible, and the weights of
each form a single Weyl-group orbit; these are the representations appearing below.
Set $N_\pm=N_{\lambda_\pm}$, where
$\lambda_\pm=(1\pm\ii\sqrt7)/8$.

\begingroup
\renewcommand{\thetheorem}{A}
\begin{theorem}[The exceptional pencil]
\label{thm:intro-exceptional-pencil}
Let $(\mathcal M,\varphi)$ be a smooth $G_2$-structure on a
$7$-manifold, with no torsion assumption, and take the cubics
$C_0$ and $C_\varphi$ in this order, with the normalization of
\eqref{eq:intro-cubic-generators}. Then the following hold:
\begin{enumerate}
  \renewcommand{\labelenumi}{(\arabic{enumi})}
  \renewcommand{\theenumi}{\arabic{enumi}}

  \item\label{item:intro-pencil-dimensions}
  On each complexified fiber $W_{\mathbb C}$, the exact
  infinitesimal linear stabilizers of
  $N_{a,b}=aC_0+bC_\varphi$ are
  \[
    \stab(N_{a,b})\cong
    \begin{cases}
      \mathfrak{so}(7,\mathbb C),
        & [a:b]=[1:0]
          \quad(\dim_{\mathbb C}=21),\\
      \mathfrak e_6,
        & a^2-2ab+8b^2=0
          \quad(\dim_{\mathbb C}=78),\\
      \mathfrak g_2^{\mathbb C},
        & \text{otherwise}
          \quad(\dim_{\mathbb C}=14).
    \end{cases}
  \]
  The two critical affine parameters are
  \[
    \lambda_\pm=\frac{1\pm\ii\sqrt7}{8}, \qquad \text{and we write} \qquad
    \qquad N_\pm=C_0+\lambda_\pm C_\varphi.
  \]
  The critical Lie algebras are simple of type $E_6$, and each acts
  on $W_{\mathbb C}$ through one of the two mutually dual $27$-dimensional
  minuscule representations.

  \item\label{item:intro-pencil-albert}
  For either sign, the full exact linear stabilizer is
  \[
    G_\pm:=\Stab_{GL(W_{\mathbb C})}(N_\pm)
    \cong E_6(\mathbb C),
  \]
  the simply connected group of type $E_6$ in its faithful
  minuscule representation. Moreover, $N_\pm$ is linearly
  equivalent over $\mathbb C$ to the cubic norm of the
  complex Albert algebra.

  \item\label{item:intro-pencil-reductions}
  The critical cubics $N_\pm$ are global smooth sections of
$S^3\mathcal W^*_{\varphi,\mathbb C}$. Each defines a canonical
reduction of the complex frame bundle of $\mathcal W_{\varphi,\mathbb C}$ to
  $G_\pm\cong E_6(\mathbb C)$, and complex conjugation
  exchanges the two reductions.
\end{enumerate}
\end{theorem}
\endgroup

Theorem~\ref{thm:intro-exceptional-pencil} follows from the
critical-parameter characterization in Proposition~\ref{prop:double-adjoint},
the stabilizer computation in Proposition~\ref{prop:block-kernels}
and Theorem~\ref{thm:stabilizer-trichotomy}, and
the identification of the critical stabilizers with $E_6$
in Theorems~\ref{thm:critical-type} and~\ref{thm:critical-group}. Corollary~\ref{cor:stabilizer-types}
collects the isomorphism types of all infinitesimal stabilizers. The Albert norm equivalence is shown in Theorem~\ref{thm:albert-norm},
and the global reductions are constructed in
Theorem~\ref{thm:exceptional-reductions}.

Theorem~\ref{thm:intro-exceptional-pencil} has the following bundle consequences, established
in Sections~\ref{subsec:exceptional-reductions} and~\ref{subsec:connections-pointings}:
\begin{itemize}
  \item The complexified bundle \(S^3\mathcal W^*_{\varphi,\C}\) admits two complex-conjugate trivial complex line subbundles generated by \(N_\pm\).
  \item The infinitesimal stabilizers of \(N_\pm\) define adjoint \(\esix\)-bundles associated with the \(G_2\)-structure bundle. With respect to the action of \(G_2\), we have a decomposition \(\esix\cong\mathbf{14}\oplus\mathbf{64}\), and hence the \(\mathbf{64}\)-component provides the additional fiberwise infinitesimal symmetries of \(N_\pm\) beyond those of \(\varphi\).
  \item The full stabilizer groups of \(N_\pm\) define concrete, complex-conjugate \(E_6(\C)\)-principal bundles \(\mathcal F_\pm\) on \(\mathcal M\), which canonically extend the original \(G_2\) structure group.
\end{itemize}

The metric also determines compact forms of the exceptional groups.
Let $\langle h,k\rangle_W=\tr(hk)$. In each reduction, the frames that
are unitary for the Hermitian form $\langle h,\bar k\rangle_W$ form a
principal bundle under the compact simply connected group of type
$E_6$, and this bundle contains the adapted $G_2$ frames
(Proposition~\ref{prop:compact-form} and
Section~\ref{subsec:compact-reductions}).

The exceptional embedding of $G_2$ in $E_6$ and the corresponding
branching of the adjoint representation are classical.
Over $\mathbb C$, Dynkin identified two conjugacy classes of irreducible
$G_2$ subgroups in $E_6$ acting irreducibly on the $\mathbf{27}$,
exchanged by outer automorphisms, and for each of them,
$\mathbf{78}=\mathbf{14}\oplus\mathbf{64}$
\cite[Theorem~11.1, Tables~24--26]{dynkin-1952-semisimple}.
In positive characteristic, Testerman constructed the corresponding
embeddings over fields containing $\sqrt{-7}$, under the characteristic
restrictions of \cite[Theorems~1(a), 2(a)]{testerman-1989-construction}.
She also reports Aschbacher's result that, under suitable characteristic
assumptions, exactly two invariant cubics up to scalar are of $E_6$-type
if and only if the field contains $\sqrt{-7}$
\cite[p.~300]{testerman-1989-construction}.
In this paper, we locate these cubics explicitly in the plane generated by the geometric
cubics $C_0,C_\varphi$, determine the infinitesimal
stabilizers throughout the pencil, and express the additional
symmetries through tensor operations on traceless symmetric $2$-tensors in \(7\) dimensions. This realization then supplies the exceptional
bundle structures and compatible connections described below.

The complex Albert algebra $\mathbb A=H_3(\mathbb O_{\mathbb C})$ is the
$27$-dimensional exceptional simple Jordan algebra of $3\times3$
Hermitian matrices over the complexified octonions
\cite[Example~6.6 and Definition~6.7]{garibaldi-2023-albert}. Its product is
$x\circ y=\frac12(xy+yx)$, with $xy$ the octonionic matrix product.
It is commutative and satisfies the Jordan identity
$x\circ((x\circ x)\circ y)=(x\circ x)\circ(x\circ y)$,
although it is not associative. Its cubic norm is the Albert determinant.
The exact stabilizer of this norm is $E_6(\mathbb C)$, and the subgroup
fixing the unit is $F_4(\mathbb C)$
\cite[\S15]{garibaldi-2023-albert}.

The critical cubic determines an Albert norm, but does not
specify a Jordan product without a choice of unit.
For either sign, an element $e\in W_{\mathbb C}$ with
$N_\pm(e)\ne0$ determines a unital complex Albert algebra
$\mathfrak J_e$ through the normalized cubic
$N_\pm/N_\pm(e)$ and the basepoint $e$
(Theorem~\ref{thm:based-algebras}).
This is the classical construction of a Jordan algebra
from a cubic norm structure
\cite{springer-1962-characterization,
mccrimmon-1969-freudenthal}.
Its automorphism group is
$\Stab_{G_\pm}(e)\cong F_4(\mathbb C)$
(Proposition~\ref{prop:f4-stabilizer};
\cite[\S15]{garibaldi-2023-albert}).
Moreover, each nonzero norm level is a single
$G_\pm\cong E_6(\mathbb C)$-orbit, isomorphic to
$E_6(\mathbb C)/F_4(\mathbb C)$
(Theorem~\ref{thm:based-algebras};
\cite[Proposition~12(i)]{krutelevich-2007-jordan}).
Thus the value of $N_\pm$ distinguishes the orbits of
basepoints under the exact norm stabilizer.
Allowing a scalar adjustment, every basepoint can be
transported to a fixed model basepoint, and all the
resulting Albert algebras are isomorphic
(Corollary~\ref{cor:basepoint-transport};
\cite[Corollary~31.25]{garibaldi-2024-albert}).

A tangent line provides a natural geometric choice of
unit. Let $\ell\subset T\mathcal M$ be a smooth real
line subbundle, and let $P_\ell$ be its orthogonal
projection tensor. Define
\[
E_\ell=7P_\ell-g.
\]
This tensor depends only on the unoriented line and
satisfies
\[
N_\pm(E_\ell)=n_\pm
=42(5+2\lambda_\pm)\ne0.
\]
This follows from \eqref{eq:su3-norm-trace}.
Consequently, $\ell$ determines an
$F_4(\mathbb C)$-reduction of each canonical
$E_6(\mathbb C)$-reduction, together with a smooth
bundle of unital complex Albert algebras.
Their Jordan products are given explicitly in terms
of $g$, $\varphi$, and $P_\ell$
(Theorem~\ref{thm:albert-bundle}).
The basepoint $E_\ell$ has residual $G_2$-stabilizer
$SU(3)\rtimes\mathbb Z/2$
(Proposition~\ref{ex:structured-basepoints}), realizing
the $E_6$-to-$F_4$ reduction anticipated by
Lee and Leung \cite[\S1.2.1, Remark]{Lee:2002fa}.
The construction also gives an explicit Euclidean
real form of each based Albert algebra
(Propositions~\ref{prop:euclidean-form}
and~\ref{prop:albert-real-form}).

Moreover, the residual stabilizer of \emph{any} nonzero basepoint
in $W_{\mathbb C}$ has dimension at most $8$, and the
line basepoints attain this bound
(Proposition~\ref{prop:g2-at-basepoint}).
Since $W_{\mathbb C}$ has no nonzero $G_2$-invariant
element, the original $G_2$-structure alone does not
select a Jordan unit. On the other hand, the choice of a basepoint has no topological
obstruction: $E_6(\mathbb C)/F_4(\mathbb C)$ is $8$-connected, and
every $F_4(\mathbb C)$-reduction is homotopic to the one supplied by
a line field (Proposition~\ref{prop:no-obstruction}).

\begingroup
\renewcommand{\thetheorem}{B}
\begin{theorem}[Albert algebra bundles]
\label{thm:intro-albert-bundles}
For either critical parameter and any smooth line subbundle $\ell\subset T\mathcal M$, not necessarily orientable,
the following hold:
\begin{enumerate}
\renewcommand{\labelenumi}{(\arabic{enumi})}
\renewcommand{\theenumi}{\arabic{enumi}}
\item\label{item:intro-bundle-reduction} $E_\ell$ is a smooth global section of $\mathcal W_\varphi$,
independent of the orientation of $\ell$, and $N_\pm(E_\ell)=n_\pm$ is a
nonzero constant. Fixing $E_\ell$ reduces the $E_6(\mathbb C)$-reduction
of Theorem~\ref{thm:intro-exceptional-pencil} to a principal
$F_4(\mathbb C)$-bundle.
\item\label{item:intro-bundle-operations} The based algebras at $E_\ell$ form a smooth bundle of complex
Albert algebras on $\mathcal W_{\varphi,\mathbb C}$, with norm
$N_\pm/n_\pm$ and unit $E_\ell$. Their trace, trace form, based
adjoint and Jordan product are globally defined and are invariant under $u\mapsto-u$ for any local unit section $u$ of $\ell$.
\item\label{item:intro-bundle-real-form} There is an explicit smooth antilinear involutive
automorphism $\varsigma_{\ell,\pm}$ of the Albert algebra bundle.
Its fixed-point bundle is a bundle of Euclidean Albert algebras,
with each fiber isomorphic to $H_3(\mathbb O)$, the Euclidean Albert algebra of $3\times3$ Hermitian octonionic matrices.
\end{enumerate}
\end{theorem}
\endgroup

This theorem combines the construction~\eqref{eq:El-def} of $E_\ell$ with
Lemma~\ref{lem:line-field-basepoint}, Theorem~\ref{thm:albert-bundle}
and Proposition~\ref{prop:albert-real-form}. The reduction uses the
standard associated-bundle construction. The algebra operations and
the Euclidean real form realize classical Albert structures explicitly
in the deformation bundle, in terms of $g,\varphi$ and $P_\ell$.

The critical cubics also capture the first-order geometry
of the $G_2$-structure. Let
$T\in\Omega^1(\mathcal M,T\mathcal M)$ denote its intrinsic
torsion, defined by
\[
  D^{LC}_X\varphi=\iota_{T(X)}\psi,
\]
where $\psi=*_{g}\varphi$ and \(D^{LC}\) denotes the Levi-Civita connection.
Write $\Theta_v(w)=w\times v$, with the cross product of~\eqref{eq:cross-product}, and let
\[
  D^{\min}=D^{LC}-\tfrac13\Theta_{T(\cdot)}
\]
be the minimal $G_2$-connection on $T\mathcal M$
\cite[\S3.1]{CleytonIvanovCurv}.

We use the same notation for connections on $T\mathcal M$
and their induced connections on $\mathcal W_{\varphi,\mathbb C}$,
as well as, when compatible, the corresponding principal
connections on the $E_6(\C)$-bundles $\mathcal F_\pm$.
Write $\rho:\mathfrak{so}(V)\to\End(W)$ for the induced infinitesimal action on symmetric tensors, defined in \eqref{eq:rho-action}.

To understand connections on \(\mathcal{W}_{\varphi,\C}\) with respect to the \(E_6(\C)\)-structure, we need to understand the corresponding Lie algebra. Indeed, with $\mathfrak g_\pm=\stab(N_\pm)$, we have the following decomposition, orthogonal for the complex-bilinear pairing $(X,Y)\mapsto\tr_W(XY)$:
\[
  \End(W_{\mathbb C})
  =\mathbb C\Id_W\oplus\mathfrak g_\pm
   \oplus\mathfrak c_\pm,
  \qquad
  \mathfrak g_\pm\cong\mathfrak e_6,\quad
  \mathfrak c_\pm\cong\mathbf{650}.
\]

\begingroup
\renewcommand{\thetheorem}{C}
\begin{theorem}[Torsion reconstruction and exceptional projection]
\label{thm:intro-compatible-connections}
Let $(\mathcal M,\varphi)$ be a smooth $G_2$-structure,
with no torsion assumption. Then:

\begin{enumerate}
\item[(i)]
For either sign,
\[
  D^{LC}_XN_\pm=\lambda_\pm L_\varphi(T(X)),
\]
where
\[
  L_\varphi:T\mathcal M\longrightarrow
  S^3\mathcal W_\varphi^*
\]
is a fiberwise injective \(G_2\)-equivariant map. Thus the Levi-Civita derivative of either critical cubic uniquely determines the full intrinsic torsion of the \(G_2\)-structure, and in particular,
$D^{LC}N_\pm=0$ if and only if $\varphi$ is torsion-free.
More generally, a metric connection on $T\mathcal M$
preserves $\varphi$ if and only if its induced connection
preserves either critical cubic, or even its line.

\item[(ii)]
For either sign,
\[
  D^{LC}-D^{\min}
   =\tfrac13\rho(\Theta_{T(\cdot)})
   \in\Omega^1(\mathcal M,\mathfrak c_\pm).
\]
The induced minimal connection is the unique connection
on $\mathcal F_\pm$ satisfying this
$\mathfrak c_\pm$-valued difference condition.
Equivalently, it is the orthogonal projection of the
induced Levi-Civita connection onto $\mathcal F_\pm$,
and the unique compatible connection pointwise closest to the
induced Levi-Civita connection for the induced Hermitian norm.
Both exceptional reductions give the same projection,
which also preserves their compact reductions.
\end{enumerate}
\end{theorem}
\endgroup

The torsion reconstruction follows from
Lemma~\ref{lem:phi-cubic-variation} and
Corollary~\ref{cor:critical-connection-criteria}.
The exceptional projection follows from
Lemma~\ref{lem:exceptional-torsion-orthogonality}
and Proposition~\ref{prop:projected-levi-civita}.

The connection $D^{\min}$ is the minimal connection
of the original tangent-bundle $G_2$-structure in
the sense of Cleyton and Swann
\cite[Definition~3.1 and Proposition~3.2]
{cleyton-2004-einstein}.
Their construction applies to a Riemannian
$G$-structure on $T\mathcal M$, with
$G\subset SO(n)$, and projects the Levi-Civita
connection using
$\mathfrak{so}(n)=\mathfrak g\oplus\mathfrak g^\perp$.
The exceptional reductions considered here are
instead structures on the auxiliary bundle
$\mathcal W_{\varphi,\mathbb C}
=S^2_0T^*\mathcal M\otimes\mathbb C$.
Their orthogonal projections take place in
$\End(W_{\mathbb C})$, with complement
$\mathbb C\Id_W\oplus\mathbf{650}$.

Theorem~\ref{thm:intro-compatible-connections}
shows that these two constructions are compatible:
the connection induced by the minimal
$G_2$-connection agrees with the projection of the
induced Levi-Civita connection onto either exceptional
reduction. The underlying reason is the
representation-theoretic inclusion
\[
  \rho(\Lambda^2_7)\subset\mathfrak c_\pm
  \cong\mathbf{650}.
\]
Thus the original intrinsic torsion, which takes
values in $\Lambda^2_7$, becomes a distinguished
$\mathbf7$-valued defect inside the exceptional
$\mathbf{650}$.

This paper develops the pointwise algebra and global bundle
geometry of the critical cubic pair for arbitrary smooth
$G_2$-structures. Tangent line subbundles give explicit Albert
products and Euclidean real forms, while the connection
results relate these structures to intrinsic torsion and
parallelism. In particular, the metric selects
the same minimal $G_2$-connection on both exceptional
reductions by projection of the induced Levi-Civita
connection. Further work will investigate curvature
conditions for compatible exceptional connections and the
relation between this pointwise exceptional algebra and the
integrated Yukawa coupling on moduli spaces of torsion-free
$G_2$-structures.

\begin{figure}[tp]
\centering

\begin{tikzpicture}[
  >=Latex,
  font=\small,
  box/.style={
    draw,
    rounded corners,
    align=center,
    text width=5cm,
    inner sep=5pt
  },
  inputbox/.style={
    box,
    dashed
  },
  thmbox/.style={
    draw,
    rounded corners,
    thick,
    densely dotted,
    inner xsep=13pt,
    inner ysep=20pt
  },
  arr/.style={->,thick}
]

\node[inputbox,text width=4.2cm] (assumptions) at (-4.5,0.25) {%
\textbf{Starting assumptions}\\[2pt]
A smooth $7$-manifold $\mathcal M$
with a $G_2$-structure $\varphi$
and induced metric $g=g_\varphi$.
};

\node[box,text width=7.8cm] (cubics) at (2.2,0.25) {%
\textbf{Two natural cubics}\\[2pt]
The trace cubic $C_0$ and pointwise Yukawa cubic $C_\varphi$
span the invariant cubic plane on
$\mathcal W_\varphi=S^2_0T^*\mathcal M$.\\[3pt]
{\footnotesize\itshape
Lemma~\ref{lem:invariant-cubic-plane};
Propositions~\ref{prop:affine-volume-cubics}
and~\ref{prop:form-metric-response}}
};

\node[box] (classify) at (-3.1,-4.1) {%
\textbf{Exceptional members of the pencil}\\[2pt]
Critical parameters $\lambda_\pm$,
complete infinitesimal stabilizer classification,
and identification of the critical algebras
as $\mathfrak e_6$.\\[3pt]
{\footnotesize\itshape
Proposition~\ref{prop:double-adjoint};
Theorems~\ref{thm:stabilizer-trichotomy}
and~\ref{thm:critical-type}}
};

\node[box] (reductions) at (3.1,-4.1) {%
\textbf{Albert norms and exceptional reductions}\\[2pt]
The critical cubics are complex Albert norms
and define canonical conjugate
$E_6(\C)$-reductions
$\mathcal F_\pm$ of
$\mathcal W_{\varphi,\C}$.\\[3pt]
{\footnotesize\itshape
Theorems~\ref{thm:critical-group},
\ref{thm:albert-norm}
and~\ref{thm:exceptional-reductions}}
};

\node[inputbox] (linefield) at (-3.1,-9.9) {%
\textbf{Additional geometric choice}\\[2pt]
A tangent line field
$\ell\subset T\mathcal M$ supplies
a basepoint of nonzero critical norm:
\[
  E_\ell=7P_\ell-g.
\]
{\footnotesize\itshape
Lemma~\ref{lem:line-field-basepoint}}
};

\node[box] (albert) at (3.1,-9.9) {%
\textbf{Albert algebra bundles}\\[2pt]
Explicit Jordan products,
$F_4(\C)$-reductions,
Euclidean real forms, and
residual $G_2$-symmetry.\\[3pt]
{\footnotesize\itshape
Theorem~\ref{thm:based-algebras};
Proposition~\ref{prop:g2-at-basepoint};
Theorem~\ref{thm:albert-bundle};
Proposition~\ref{prop:albert-real-form}}
};

\node[box] (torsion) at (-3.1,-15.5) {%
\textbf{Torsion reconstruction}\\[2pt]
The Levi-Civita derivative of either
critical cubic determines the full
intrinsic torsion:
\[
  D^{LC}N_\pm=\lambda_\pm L_\varphi(T).
\]
{\footnotesize\itshape
Lemma~\ref{lem:phi-cubic-variation};
Corollary~\ref{cor:critical-connection-criteria}}
};

\node[box] (connection) at (3.1,-15.5) {%
\textbf{Exceptional connection projection}\\[2pt]
Projection of the induced Levi-Civita
connection onto either exceptional
reduction gives the induced minimal
$G_2$-connection $D^{\min}$,
with $\mathbf{650}$-valued defect.\\[3pt]
{\footnotesize\itshape
Lemma~\ref{lem:exceptional-torsion-orthogonality};
Proposition~\ref{prop:projected-levi-civita}}
};

\begin{scope}[on background layer]

\node[
  thmbox,
  fit=(classify)(reductions)
] (thmA) {};

\node[
  thmbox,
  fit=(linefield)(albert)
] (thmB) {};

\node[
  thmbox,
  fit=(torsion)(connection)
] (thmC) {};

\end{scope}

\node[
  anchor=north west,
  font=\bfseries,
  inner sep=0pt
] at ([xshift=11pt,yshift=-4pt]thmA.north west)
  {Theorem~\ref{thm:intro-exceptional-pencil}};

\node[
  anchor=north west,
  font=\bfseries,
  inner sep=0pt
] at ([xshift=11pt,yshift=-4pt]thmB.north west)
  {Theorem~\ref{thm:intro-albert-bundles}};

\node[
  anchor=north west,
  font=\bfseries,
  inner sep=0pt
] at ([xshift=11pt,yshift=-4pt]thmC.north west)
  {Theorem~\ref{thm:intro-compatible-connections}};

\coordinate (leftRailTop) at (-7.2,-0.9);
\coordinate (leftRailMid) at (-7.2,-15.5);

\draw[arr] (assumptions.east) -- (cubics.west);

\draw[arr]
  (assumptions.south) -- (leftRailTop)
  -- (leftRailMid)
  -- (torsion.west);

\draw[arr]
  (cubics.south)
  to[out=-120,in=80]
  (classify.north);

\draw[arr]
  (classify.east) -- (reductions.west);

\draw[arr]
  (reductions.south) -- (albert.north);

\draw[arr]
  (linefield.east) -- (albert.west);

\draw[arr]
  (torsion.east) -- (connection.west);

\draw[arr]
  (classify.west) -- ++(-1.15,0)
  |- (torsion.west);

\draw[arr]
  (reductions.east) -- ++(1.15,0)
  |- (connection.east);

\end{tikzpicture}

\caption{%
Overall construction of the paper.
}
\label{fig:overall-construction}
\end{figure}
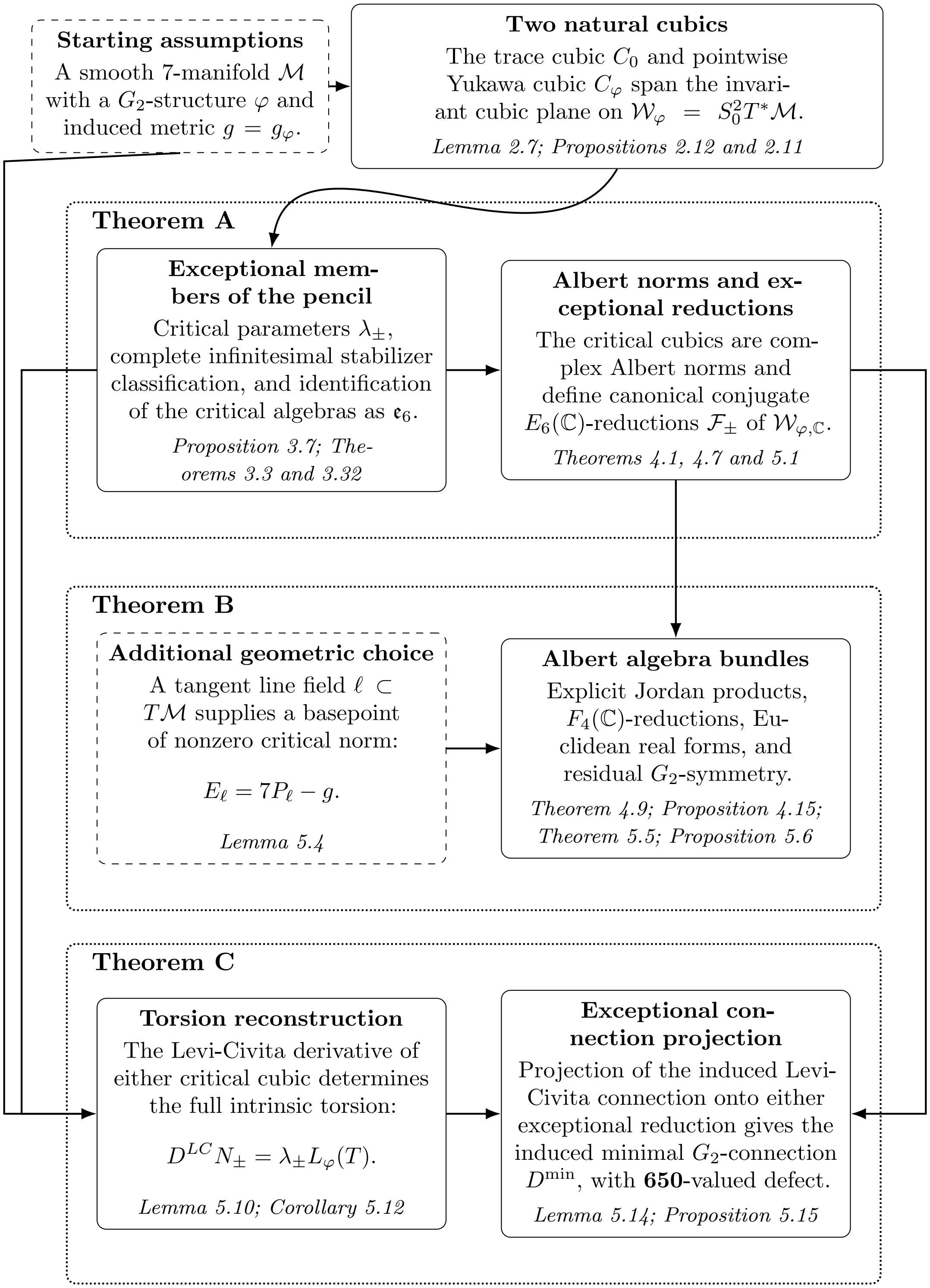

Figure \ref{fig:overall-construction} summarizes the overall construction in the paper.

We begin the algebraic argument with contraction
identities for the two cubics $C_0$ and $C_\varphi$. From these,
we derive two universal quartic identities and use them to
characterize the critical cubics through a double-adjoint identity
(Proposition~\ref{prop:double-adjoint}).
For \(\lambda\in\C\), define the quadratic adjoint maps $F_\lambda:W_{\mathbb C}\to W_{\mathbb C}$ by
\begin{equation}
\label{eq:intro-quadratic-adjoint}
\langle F_\lambda(h),k\rangle_W=N_\lambda(h,h,k),
\end{equation}
for \(h,k\in W_\C\), where $\langle h,k\rangle_W=\tr(hk)$ is the $W$-pairing, extended complex-bilinearly to \(W_\C\). We use the same symbol $N_\lambda$ for its symmetric trilinear polarization,
normalized by $N_\lambda(h,h,h)=N_\lambda(h)$ as in
\eqref{eq:cubic-polarization}.

Using quartic contraction identities for \(C_0\) and \(C_\varphi\) (Lemma \ref{lem:universal-quartic-identities}), we reduce the double-adjoint condition
$F_\mu(F_\lambda(h))=\kappa N_\lambda(h)h$, for every $h$ and some $\kappa\in\C$, to
$\lambda+\mu=\frac14$ and $\lambda\mu=\frac18$.
Thus $\lambda,\mu$ must be the two roots of
$8\lambda^2-2\lambda+1$, which explains both the critical quadratic
and the appearance of $\sqrt{-7}$
(Proposition~\ref{prop:double-adjoint} and~\eqref{eq:composition-formula}).
Writing $F_\pm=F_{\lambda_\pm}$, the identity is
\begin{equation}
\label{eq:intro-double-adjoint}
F_\mp(F_\pm(h))=\frac23 N_\pm(h)h.
\end{equation}
A key feature is that the critical parameters are identified before the stabilizers of \(N_\pm\) are shown to be \(E_6\). The double-adjoint identity is analogous to the classical identity
$\operatorname{adj}(\operatorname{adj} A)=\det(A)A$
for $3\times3$ matrices. Moreover, \eqref{eq:intro-quadratic-adjoint} gives $F_\lambda=\frac13\nabla N_\lambda$ with respect to
the \(W\)-pairing, in analogy with the classical adjugate identity
$\operatorname{adj}(A)=\nabla\det(A)$ for symmetric matrices. Our quadratic
adjoints pair the conjugate critical cubics $N_+$ and $N_-$.

The double-adjoint identity also has an \emph{a posteriori}
explanation through the representation theory of $E_6(\C)$.
This group has two inequivalent, mutually dual minuscule
representations, $\mathbf{27}$ and
$\overline{\mathbf{27}}=\mathbf{27}^*$
\cite[\S3, Table~1]{podder-2023-k-ring}.
Their tensor products decompose as
\[
\begin{aligned}
S^2\mathbf{27}
  &=\overline{\mathbf{27}}\oplus\mathbf{351}',\\
\mathbf{27}\otimes\overline{\mathbf{27}}
  &=\mathbf1\oplus\mathbf{78}\oplus\mathbf{650},\\
S^3\mathbf{27}
  &=\mathbf1\oplus\mathbf{650}\oplus\mathbf{3003}.
\end{aligned}
\]
See \cite{feger-2020-lieart} for the representation data.
Neither minuscule representation admits an invariant bilinear
form, but the two carry a natural invariant pairing.
The unique trivial summand in the symmetric cube corresponds
to the Albert norm
(Proposition~\ref{prop:unique-cubic};
\cite[12.4(ii)]{springer-1973-jordan}),
while the occurrence of $\overline{\mathbf{27}}$ in
$S^2\mathbf{27}$ gives, up to scale, a unique equivariant
quadratic adjoint into the dual representation, and
similarly in the reverse direction.

The composition of these adjoints is a quartic equivariant
self-map of $\mathbf{27}$. It must be proportional to
$h\mapsto N(h)h$: at a point of nonzero norm, the stabilizer
$F_4(\C)$ fixes precisely the line through that point, and
$E_6(\C)$ acts transitively on each nonzero norm level
(Proposition~\ref{prop:f4-stabilizer} and
Theorem~\ref{thm:based-algebras}).
In our $G_2$ realization, the metric-induced pairing identifies the
two dual modules as vector spaces, but not as $E_6$-modules.
The corresponding dual actions are realized by the two
critical stabilizers, with quadratic adjoints $F_+$ and
$F_-$ (Remark~\ref{rem:koecher-tkk}).
Thus representation theory explains the form of
\eqref{eq:intro-double-adjoint} once the exceptional
structures have been identified. By contrast,
Proposition~\ref{prop:double-adjoint} determines the
critical parameters and the coefficient $\tfrac23$
directly from $G_2$ tensor identities, without assuming
exceptional symmetry.

The double-adjoint identity \eqref{eq:intro-double-adjoint} also allows us to construct
infinitesimal symmetries through the Koecher operators, which are
linear endomorphisms associated with the Tits--Kantor--Koecher
construction (Remark~\ref{rem:koecher-tkk}). Explicitly, with
$F_\pm(k,h)$ denoting the halved bilinear polarization of $F_\pm$, they are
\begin{equation*}
 S_\pm(k,l)h=\langle h,l\rangle_Wk
 -2F_\mp(F_\pm(k,h),l)+\tfrac13\langle k,l\rangle_Wh.
\end{equation*}
By the double-adjoint identity, these operators are infinitesimal
symmetries (Corollary~\ref{cor:koecher}).
Separately, Proposition~\ref{prop:five-blocks}
excludes five $G_2$-isotypic components of $\End\WC$ from the
stabilizer of every nonzero member of the pencil. Namely, it excludes a \(\mathbf{27}\oplus\mathbf{27}\)-block, two blocks corresponding to the two inequivalent \(77\)-dimensional representations of \(G_2\), as well as the higher-dimensional \(\mathbf{182}\) and \(\mathbf{189}\) blocks.
The next key step is to evaluate the span of the critical Koecher operators $S_\pm(k,l)$ at a test tensor
$E_u=7u^\flat\otimes u^\flat-g$ for a unit vector $u$.
This evaluation maps onto $E_u^\perp$, of dimension $26$.
After the five exclusions, the available irreducible types other
than the $64$-dimensional module $V(1,1)$ have total dimension
$1+7+14=22$. Hence a copy of $V(1,1)$ must occur.
Together with $\mathfrak g_2^{\mathbb C}$, it gives the lower
bound $78$ for the dimension of the critical stabilizer
(Proposition~\ref{prop:koecher-jump}).

Using this bound and the stabilizers of the two generators, we then identify
the $\mathbf{64}$-parts of the two critical stabilizers as explicit graphs of $G_2$-isomorphisms
between the two copies of $V(1,1)$ in the $\mathbf{64}\oplus\mathbf{64}$ component
of $\End\WC$ (Lemma~\ref{lem:64-block}).
The count over the remaining isotypic components gives the exact stabilizer dimensions
throughout the pencil (Proposition~\ref{prop:block-kernels}).
Finally, we prove that the critical stabilizer Lie algebras are
simple and identify their type as $E_6$ using the classification
of simple Lie algebras and their irreducible representations
(Theorem~\ref{thm:critical-type}).
By simplicity of the Lie algebra, the Koecher operators in fact span the whole critical
stabilizer (Corollary~\ref{cor:koecher-span}), giving an explicit
realization of $\mathfrak e_6$ on $W_{\mathbb C}$ from the quadratic
maps $F_\pm$. The double-adjoint characterization and the
$78$-dimensional lower bound precede, and do not use, the
classification or the Albert identification.

Theorem~\ref{thm:intro-exceptional-pencil} is proved in
Sections~\ref{sec:deformation-representation}--\ref{sec:reductions},
and Theorems~\ref{thm:intro-albert-bundles}
and~\ref{thm:intro-compatible-connections} in
Section~\ref{sec:reductions}.
In Section~\ref{sec:deformation-representation}, we introduce the
deformation module, its invariant cubic plane and the description of
the two generators as variations of the volume.
In Section~\ref{sec:exceptional-pencil}, we locate the critical pair
and determine every stabilizer through the double-adjoint and Koecher
argument just described.
In Section~\ref{sec:e6-albert-norm}, we determine the full stabilizer groups, identify the critical cubics
with the Albert norm, discuss basepoints and their residual $G_2$
symmetry, compare the result with the Yukawa cubic, and construct
the compact and Euclidean real forms.
In Section~\ref{sec:reductions}, we return to the manifold: we
construct the conjugate exceptional reductions and their compact forms,
prove Theorem~\ref{thm:intro-albert-bundles} and the homotopy
uniqueness of the $F_4(\mathbb C)$-reductions, and conclude with compatible connections
and parallel pointings, that is, parallel sections of nonzero critical norm.

In Appendix~\ref{app:contraction-proofs}, we derive the one- and
two-index contraction identities for \(C_0\) and \(C_\varphi\) (c-contraction identities) in detail, compute the contractions used in the stabilizer exclusion, and prove the corresponding quartic identities. Appendix \ref{subsec:notation} gives the full notation inventory for this paper.

\subsection*{Acknowledgments}

Codex (most recently, with GPT-6.1 Sol and GPT 6 Astra models) and Claude Code (most recently, with Opus 5.5 and Fable 5.1/5.2 models) were used to organize and verify the research
artifacts behind this paper and to assist with typesetting. ChatGPT was
used for assistance with ideation. All results have been verified by the
author, who takes full responsibility for the content.
\section{The traceless symmetric tensors and their natural cubics}
\label{sec:deformation-representation}

The critical parameters of the pencil in Section~\ref{sec:exceptional-pencil},
and the fields they generate, depend on the specified generators
$C_0,C_\varphi$ and on the ground field. In this section we therefore fix
explicitly the scale of the positive form, of the insertion map from
symmetric tensors to $3$-forms, and of the two cubics.

\subsection{The standard form and its \texorpdfstring{$G_2$}{G2}-contraction identities}
\label{subsec:standard-form-contractions}

Let $V=\mathbb R^7$, with ordered basis $(e_1,\ldots,e_7)$ and dual
coframe $(e^1,\ldots,e^7)$. The Euclidean metric $g$, the orientation
form $\vol$ and the standard positive $3$-form $\varphi$ are
\begin{equation}
\begin{gathered}
  g=\sum_{i=1}^7e^i\otimes e^i,\qquad\text{and}\qquad \vol=e^{1234567},\\
  \varphi=e^{123}+e^{145}+e^{167}+e^{246}-e^{257}-e^{347}-e^{356},
\end{gathered}
  \label{eq:standard-positive-form}
\end{equation}
where $e^{i_1\cdots i_k}=e^{i_1}\wedge\cdots\wedge e^{i_k}$.

The Hodge star is determined by
\begin{equation}
  \alpha\wedge\ast\alpha'=\langle\alpha,\alpha'\rangle\vol,
  \label{eq:hodge-convention}
\end{equation}
where the increasing wedge monomials are orthonormal. In these conventions,
\begin{equation}
\begin{aligned}
  \psi=\ast\varphi
    &=e^{4567}+e^{2367}+e^{2345}+e^{1357}
       -e^{1346}-e^{1256}-e^{1247},\\
  |\varphi|^2&=|\psi|^2=7,\qquad\text{and}\qquad \ast\psi=\varphi.
\end{aligned}
\label{eq:standard-dual-form}
\end{equation}
The signs agree with Bryant's conventions
\cite[(2.1)--(2.8)]{bryant-2003}.

We use the index convention of Section~\ref{sec:introduction}; in this
adapted frame, $g_{ij}=\delta_{ij}$ and $g^{ij}=\delta^{ij}$. The alternating
components are $\varphi_{ijk}=\varphi(e_i,e_j,e_k)$ and
$\psi_{ijkl}=\psi(e_i,e_j,e_k,e_l)$, so that
\begin{equation*}
  \varphi=\frac1{3!}\varphi_{ijk}e^{ijk},\qquad\text{and}\qquad
  \psi=\frac1{4!}\psi_{ijkl}e^{ijkl},\qquad\text{with }\varphi_{123}=1.
\end{equation*}
The factorials appear in the expansion of each form, not in its
components.

The stabilizer
\begin{equation}
  G_2=\{O\in GL(V):O^*\varphi=\varphi\}
  \label{eq:model-g2}
\end{equation}
preserves $g$ and $\vol$ \cite[(2.2)]{bryant-2003}. All the algebraic
objects extend to $V_{\mathbb C}=V\otimes_{\mathbb R}\mathbb C$ by complex
multilinearity. In particular, no complex conjugation enters any bilinear
or cubic formula, including the complex critical cubics of
Section~\ref{sec:exceptional-pencil}. These conventions are in
characteristic zero, and the comparisons with positive characteristic
cited in the text hold under their own stated hypotheses.

For $u,v\in V$, the cross product $u\times v\in V$ is defined by
\begin{equation}
  \iota_v\iota_u\varphi=(u\times v)^\flat,\qquad\text{equivalently}\qquad
  (u\times v)_i=\varphi_{iab}u^av^b,
  \label{eq:cross-product}
\end{equation}
so that $g(u\times v,w)=\varphi(u,v,w)$ for $w\in V$.

The following contraction identities of $\varphi$ and $\psi$ are used
throughout Sections~\ref{sec:deformation-representation}
and~\ref{sec:exceptional-pencil} and
Appendix~\ref{app:contraction-proofs}.

\begin{lemma}
\label{lem:standard-contractions}
For $u,v\in V$,
\begin{equation}
  (\iota_u\varphi)\wedge(\iota_v\varphi)\wedge\varphi
     =6g(u,v)\vol.
  \label{eq:metric-volume-contraction}
\end{equation}
The component contractions are
\begin{equation}
  \varphi_{iab}\varphi_{j}{}^{ab}=6\delta_{ij},\qquad\text{and}\qquad
  \varphi_{ija}\varphi_{kl}{}^{a}
    =\delta_{ik}\delta_{jl}-\delta_{il}\delta_{jk}+\psi_{ijkl}.
  \label{eq:phi-phi-contractions}
\end{equation}
The mixed contraction of $\varphi$ with $\psi$ is
\begin{equation}
  \varphi_{iab}\psi_{jk}{}^{ab}=4\varphi_{ijk}.
  \label{eq:mixed-contraction}
\end{equation}
\end{lemma}
\begin{proof}
The identities \eqref{eq:metric-volume-contraction} and
\eqref{eq:phi-phi-contractions} follow by a componentwise check from
\eqref{eq:standard-positive-form} and \eqref{eq:standard-dual-form}.
Substituting the second identity of \eqref{eq:phi-phi-contractions} into
the mixed contraction and then applying the first gives
\begin{equation*}
  \varphi_{iab}\psi_{jk}{}^{ab}
     =6\varphi_{jki}-2\varphi_{ijk}=4\varphi_{ijk}.
\end{equation*}
\end{proof}
In Bryant's notation, after relabeling, the identities
\eqref{eq:phi-phi-contractions} and
\eqref{eq:mixed-contraction} are given in
\cite[(2.6)--(2.8)]{bryant-2003}.
For exterior-algebra identities, see also
\cite[Lemma~2.2.3]{karigiannis-2005-57}.
Component contractions are collected in
\cite[Lemmas~A.12--A.14]{karigiannis-2007},
where the orientation convention gives the opposite sign for the
$\psi$ terms.

\subsection{Symmetric tensors and pointwise deformations}
\label{subsec:deformation-module}

A positive $3$-form $\varphi$ on a $7$-manifold $\mathcal M$ is pointwise
linearly equivalent to the standard model \eqref{eq:standard-positive-form}
and therefore determines a
$G_2$-structure, a metric $g=g_\varphi$, and an orientation
\cite[Section~3.1.1]{bryant-2003}. On both $V$ and $\mathcal M$, we write
$\psi=*\varphi$.

Infinitesimal deformations of $\varphi$ split into the irreducible
summands $\Lambda^3_1$, $\Lambda^3_7$, and $\Lambda^3_{27}$.
The first gives conformal metric variations, the second leaves
the metric unchanged to first order, and the third gives traceless
metric variations. Our cubics are sections of
$S^3\mathcal W_\varphi^*$, where
\[
  \mathcal W_\varphi
  =S^2_0T^*\mathcal M\cong\Lambda^3_{27}T^*\mathcal M.
\]

We first specify the deformation bundle and its frames, before describing
its fibers as pointwise metric and form deformations.

\begin{definition}
  \label{def:deformation-bundle}
  Let $(\mathcal M,\varphi)$ be a $7$-manifold with a $G_2$-structure, with no assumption on its
  torsion, and let $W=S^2_0(\mathbb R^7)^*$ be the model fiber.
  \begin{enumerate}
    \item[(i)] $\mathrm{Fr}_\varphi\to \mathcal M$ is the principal $G_2$-bundle of $\varphi$-adapted
      coframes, which identify $\varphi$ pointwise with the standard model~\eqref{eq:standard-positive-form}, and
      $\mathrm{Fr}_{\varphi,\mathbb C}=\mathrm{Fr}_\varphi\times_{G_2}G_2^{\mathbb C}$ is its
      complexification, with $G_2$ acting on $G_2^{\mathbb C}$ by left
      multiplication. Here $G_2^{\mathbb C}$ is the complex Lie group of type $G_2$,
      the complexification of $G_2$ (recalled in Section~\ref{sec:e6-albert-norm}).
    \item[(ii)] The \emph{deformation bundle} of $\varphi$ is
      \begin{equation*}
        \mathcal W_\varphi=S^2_0T^*\mathcal M=\mathrm{Fr}_\varphi\times_\rho W,
        \qquad\text{where}\qquad
        \rho(O)h=(O^{-1})^*h\quad(O\in G_2),
      \end{equation*}
      and its complexification is
      $\mathcal W_{\varphi,\mathbb C}=\mathcal W_\varphi\otimes_{\mathbb R}\mathbb C
      =\mathrm{Fr}_\varphi\times_\rho W_{\mathbb C}$, a complex vector bundle of rank $27$.
    \item[(iii)] For a coframe $p\in\mathrm{Fr}_\varphi$ over $x\in \mathcal M$, the induced frame of the
      fiber is
      \begin{equation*}
        f_p\colon W_{\mathbb C}\to\mathcal W_{\varphi,\mathbb C,x},
        \qquad\text{with}\qquad
        f_p(h)=[p,h].
      \end{equation*}
  \end{enumerate}
  Fiberwise, the insertion map $\mathcal I_\varphi$ defined below identifies $\mathcal W_\varphi$ with
  $\Lambda^3_{27}T^*\mathcal M$.
\end{definition}

We use the tensor model $S^2_0T^*\mathcal M$ throughout; the form model
$\Lambda^3_{27}T^*\mathcal M$ enters through $\mathcal I_\varphi$. Note that $\mathcal{W}_{\varphi}$ actually only depends on the metric \(g\), since trace removal in \(S^2T^*\mathcal M\) requires the metric. Hence,  $\mathcal{W}_{\varphi}$ can also be regarded as an associated bundle to the orthonormal frame bundle. In our case, since we assume a \(G_2\)-structure, which induces the metric \(g=g_\varphi\), we regard  $\mathcal{W}_{\varphi}$ as a bundle associated to \(\mathrm{Fr}_\varphi\).

The cubic $C_0$ depends only on the metric, while $C_\varphi$ also depends
on $\varphi$, as in~\eqref{eq:natural-cubics} below.

We now fix $x\in \mathcal M$ and identify $V=T_x\mathcal M$ with the model $\mathbb R^7$
by an adapted coframe. Throughout the rest of this section and
Sections~\ref{sec:exceptional-pencil}--\ref{sec:e6-albert-norm} we work
in this fiber, with $W=\mathcal W_{\varphi,x}\cong S^2_0V^*$.

Elements of $S^2V^*$ are covariant symmetric $2$-tensors, and we multiply
covariant $2$-tensors by contracting through the metric:
\begin{equation}
(hk)_{ab}=h_{ac}g^{cd}k_{db}.
\label{eq:contracted-product}
\end{equation}
This product is associative and has $g$ as its identity, but $hk$ is in
general neither equal to $kh$ nor symmetric. We write $h^n$ for the
$n$-fold contracted product, with $h^0=g$, and
$\tr(h)=g^{ab}h_{ab}$, so that $\tr(h^n)=g^{ab}(h^n)_{ab}$ and
$\tr(hk)=\tr(kh)$. The transpose of a covariant $2$-tensor is
$(h^{\transp})_{ab}=h_{ba}$. The metric identifies $h$ with the
endomorphism
\begin{equation}
g^{-1}h\in\End(V),\qquad (g^{-1}h)^a{}_b=g^{ac}h_{cb},
\label{eq:metric-identification}
\end{equation}
which is self-adjoint with respect to $g$ exactly when $h$ is symmetric.
This identification carries contracted products to compositions,
$g^{-1}(hk)=(g^{-1}h)(g^{-1}k)$, preserves traces, and sends $g$ to the
identity endomorphism $I$ of $V$. We use it when compositions,
commutators or eigenspaces are needed, without a separate symbol for
$g^{-1}h$, and write $\Id_W$ for the identity endomorphism of $W$.
In a $g$-orthonormal frame, $h$ and $g^{-1}h$ have the same matrix, $g$ is
represented by the identity matrix, the contracted product by the matrix
product and $h^{\transp}$ by the matrix transpose; coordinate calculations
below use this matrix notation.

\begin{definition}
\label{def:deformation-module}
The space of traceless covariant symmetric $2$-tensors and its pairing
induced by $g$ are
\begin{equation}
W=S^2_0V^*=\{h\in S^2V^*:\tr(h)=0\},
\qquad
\langle h,k\rangle_W=g^{ac}g^{bd}h_{ab}k_{cd}=\tr(hk).
\label{eq:deformation-module}
\end{equation}
We call $\langle\cdot,\cdot\rangle_W$ the \emph{$W$-pairing}, and use the
same notation for its complex-bilinear extension to $W_{\mathbb C}$; no
complex conjugation is involved.
\end{definition}

We use the following tensor conventions throughout the paper.
For $h\in S^2 V^*$, define its traceless part $h_0$ by
\begin{equation}
\label{eq:traceless}
h_0=h-\frac17 \tr(h) g.
\end{equation}
Let $\Sigma\subseteq V$ be a subspace with orthogonal projection
$\pr_\Sigma\in\End(V)$. For $u\in V$, or $u\in V_{\mathbb C}$ with $\flat$
extended complex-linearly, define
\begin{equation}
\label{eq:projection-tensors}
P_u=u^\flat\otimes u^\flat,
\qquad
P_\Sigma(v,w)=g(\pr_\Sigma v,w),
\end{equation}
so that $g^{-1}P_\Sigma=\pr_\Sigma$ and $P_{\Sigma^\perp}=g-P_\Sigma$.
For a unit vector $u$, $g^{-1}P_u=\pr_{\mathbb Ru}$, so that $P_u$ depends
only on the line $\ell=\mathbb Ru$ and equals $P_\ell$; we write $P_\ell$
when a line is specified and $P_u$ when only a vector is. For a general $u$,
$P_u$ is only a rank-one tensor, and for a complex null vector
$g^{-1}P_u$ is nilpotent. The calligraphic letter $\mathcal P$ with two
arguments is reserved for the bilinear operation $\mathcal P(h,k)$ defined
below.
The traceless parts of these tensors, rescaled by $7$, serve as test
tensors and basepoints throughout the paper. For $v\in V$ or
$v\in V_{\mathbb C}$, and for a subspace $\Sigma\subseteq V$, define
\begin{equation}
\label{eq:test-tensors}
E_v=7(P_v)_0=7P_v-g(v,v)\,g,
\qquad
E_\Sigma=7(P_\Sigma)_0=7P_\Sigma-(\dim\Sigma)\,g.
\end{equation}
For a unit vector $u$ this is $E_u=7P_u-g$. It depends only on the line
$\ell=\mathbb Ru$, and we write $E_\ell=7P_\ell-g$ when the line is
specified. For a $3$-plane $\Sigma$,
$E_\Sigma=7P_\Sigma-3g=4P_\Sigma-3P_{\Sigma^\perp}$. For the basis vectors
$e_i$ of \eqref{eq:standard-positive-form} we abbreviate
$P_i=P_{e_i}$ and $E_i=E_{e_i}$.

For $h\in S^2 V^*$, define the insertion map by
\begin{equation}
\begin{aligned}
\mathcal I_\varphi(h)
 &=\left.\frac{d}{dt}\right|_{t=0}(I+t\,h)^*\varphi,\\
(\mathcal I_\varphi(h))_{ijk}
 &=h_i{}^a\varphi_{ajk}
   +h_j{}^a\varphi_{iak}
   +h_k{}^a\varphi_{ija}.
\end{aligned}
\label{eq:model-map}
\end{equation}
We write $\chi_h=\mathcal I_\varphi(h)$. For $h=g$, \eqref{eq:model-map}
gives $\mathcal I_\varphi(g)=3\varphi$. Note that in \cite{GrigorianG2Review,GrigorianYau1}, the averaged insertion $i_{\varphi}(h)=\mathcal I_\varphi(h)/3$ is used. In \cite{dwivedi-2025-flows}, the notation $h\diamond \varphi$ is used, and is extended to not necessarily symmetric $h$.
The standard $G_2$-decomposition of $3$-forms
\cite{bryant-2003}, \cite[(2.19)--(2.21)]{karigiannis-2005-57} is the
orthogonal sum
\begin{equation}
\Lambda^3V^*
 =\mathbb R\varphi
  \oplus\{\iota_v\psi:v\in V\}
  \oplus\Lambda^3_{27}V^*,
\qquad\text{where}\qquad
\Lambda^3_{27}V^*
 =\{\chi:\chi\wedge\varphi=0,\ \chi\wedge\psi=0\}.
\label{eq:three-form-splitting}
\end{equation}
We write $\pr_1$, $\pr_7$ and $\pr_{27}$ for the projections onto the three
summands of \eqref{eq:three-form-splitting}, in this order; for another
positive form, they denote the projections onto the summands of its
splitting.
The following contractions, with $\chi_h=\mathcal I_\varphi(h)$ as in
\eqref{eq:model-map}, fix the normalization of $\mathcal I_\varphi$.

\begin{lemma}
\label{lem:model-map-identities}
For $h,k\in S^2V^*$, in the inner product of forms
\eqref{eq:hodge-convention},
\begin{subequations}
\begin{align}
  (\chi_h)_{iab}\varphi_{j}{}^{ab}
     &=4h_{ij}+2(\tr h)\delta_{ij},
       \label{eq:model-map-contraction}\\
  \langle\mathcal I_\varphi(h),\mathcal I_\varphi(k)\rangle
     &=2\tr(hk)+(\tr h)(\tr k).
       \label{eq:model-map-inner-product}
\end{align}
\end{subequations}
In particular, for $h\in W$,
\begin{equation}
|\chi_h|^2=2\tr(h^2),
\qquad\text{and}\qquad
h_{ij}=\frac14(\chi_h)_{iab}\varphi_j{}^{ab}.
\label{eq:insertion-normalization}
\end{equation}
\end{lemma}
\begin{proof}
For \eqref{eq:model-map-contraction}, expand the three terms of
$\chi_h$ and apply \eqref{eq:phi-phi-contractions}. The first term
contributes $6h_{ij}$; each of the other two contributes
$(\tr h)\delta_{ij}-h_{ij}$. The contraction with $\psi$ vanishes
because $h$ is symmetric. For the inner product, we use the component
pairing with its factor $1/6$. By alternation, the three terms of
$\mathcal I_\varphi(k)$ contribute equally. Thus we obtain one half of
\eqref{eq:model-map-contraction} paired with $k$, which proves
\eqref{eq:model-map-inner-product}. The identities
\eqref{eq:insertion-normalization} on the traceless subspace follow.
\end{proof}

For $h\in W$, the trace of \eqref{eq:model-map-contraction} gives
$\langle\chi_h,\varphi\rangle=\tfrac16(\chi_h)_{iab}\varphi^{iab}
=3\tr h=0$. The three summands of $\chi_h$ in \eqref{eq:model-map} have
the same contraction with $\psi_{mijk}$ after cyclic relabeling, so
\eqref{eq:mixed-contraction} gives
\begin{equation*}
  (\chi_h)_{ijk}\psi_{m}{}^{ijk}
    =3h_{ia}\varphi^{a}{}_{jk}\psi_{m}{}^{ijk}
    =12h_{ia}\varphi^{a}{}_{m}{}^{i}=0,
\end{equation*}
by symmetry of $h$. Thus $\chi_h$ is orthogonal to $\iota_{e_m}\psi$ for
every $m$. By the orthogonal splitting \eqref{eq:three-form-splitting},
$\chi_h\in\Lambda^3_{27}V^*$. Directly from the explicit forms,
$\varphi,\iota_{e_1}\psi,\ldots,\iota_{e_7}\psi$ are mutually orthogonal
with squared norms $7,4,\ldots,4$; in particular,
$|\iota_v\psi|^2=4|v|^2$, and $\Lambda^3_{27}V^*$ has dimension
$35-8=27$. The identity \eqref{eq:insertion-normalization} gives
injectivity, and $W$ also has dimension $27$, so
$\mathcal I_\varphi\colon W\to\Lambda^3_{27}V^*$ is an isomorphism.
We use the pullback derivative itself, rather than the isometric
rescaling $\mathcal I_\varphi/\sqrt2$.

The interpretation of $W$ in terms of deformations follows from the same
pullback curve.
For $h\in S^2V^*$, the induced metric is
$(I+t\,h)^*g$, whose derivative at $t=0$ is $2h$.
Thus $h\in W$ describes a pointwise traceless metric variation, and
$\chi_h$ is the corresponding metric-changing variation of $\varphi$.
We call $W$ the deformation representation (or module) in this sense,
and $\mathcal W_\varphi$ its associated deformation bundle. Their elements
are pointwise directions, and the terminology involves no differential
equation or harmonic condition.
For finite deformation formulas for the induced metric and torsion, see
\cite[\S4; Lemma~5.1, equation~(5.31)]{GrigorianG2Torsion1}.

For a Lie-algebra element $M\in\mathfrak{gl}(V)$, define $\theta_M$, the
infinitesimal pullback of $\varphi$ by $M$, by
\begin{subequations}
  \begin{align}
\theta_M(v_1,v_2,v_3)
 &=\varphi(Mv_1,v_2,v_3)+\varphi(v_1,Mv_2,v_3)
   +\varphi(v_1,v_2,Mv_3), \text{ and} \label{eq:insert-M}\\
\mathfrak g_2=\operatorname{Lie}(G_2)
 &=\{M\in\mathfrak{gl}(V):\theta_M=0\} \subseteq\mathfrak{so}(7). \label{eq:annihilator-map}
\end{align}
\end{subequations}

The orthogonal group acts on covariant $2$-tensors by the natural left
action, which preserves $W$:
\begin{equation}
\rho(O)h=(O^{-1})^*h,
\qquad\text{so that}\qquad
g^{-1}(\rho(O)h)=O(g^{-1}h)O^{-1}\quad(O\in O(7)).
\label{eq:rho-group}
\end{equation}
Its derivative is
\begin{equation}
\rho:\mathfrak{so}(7)\longrightarrow\End(W),
\qquad\text{given by}\qquad
\rho(M)h=g\bigl([M,g^{-1}h]\,\cdot,\cdot\bigr),
\label{eq:rho-action}
\end{equation}
that is, $(\rho(M)h)(u,v)=-h(Mu,v)-h(u,Mv)$, the commutator being taken
in $\End(V)$. In a $g$-orthonormal frame $\rho(M)h=[M,h]$ as matrices,
and we use this form in matrix calculations.
Each $\rho(M)$ is skew-adjoint for the $W$-pairing ($W$-skew-adjoint in the terminology introduced after~\eqref{eq:traceTw} below), since cyclicity of trace gives
\begin{equation*}
\langle\rho(M)h,k\rangle_W+\langle h,\rho(M)k\rangle_W
 =\tr\bigl([M,g^{-1}h]\,g^{-1}k+g^{-1}h\,[M,g^{-1}k]\bigr)=0.
\end{equation*}
This representation is injective, and its restriction to $\mathfrak g_2$ is also an injective representation on $W$.
We use the standard identification of $\mathfrak g_2$ as the compact
simple Lie algebra of type $G_2$, of dimension $14$, and write
$\mathfrak g_2^{\mathbb C}$ for its complexification.

Using $g$ we identify $\mathfrak{so}(7)$ with $\Lambda^2V^*$. The standard
splitting of $2$-forms is
\begin{equation}
\begin{gathered}
\Lambda^2V^*=\Lambda^2_7\oplus\Lambda^2_{14},\qquad
\Lambda^2_7=\{\iota_v\varphi:v\in V\}
 =\{\alpha:\psi_{abcd}\alpha^{cd}=4\alpha_{ab}\},\\
\Lambda^2_{14}=\{\alpha:\varphi_{abc}\alpha^{bc}=0\}
 =\{\alpha:\psi_{abcd}\alpha^{cd}=-2\alpha_{ab}\}=\mathfrak g_2,
\end{gathered}
\label{eq:two-form-splitting}
\end{equation}
where $(\iota_v\varphi)_{bc}=v^a\varphi_{abc}$ \cite{bryant-2003}. The
identification of $\Lambda^2_{14}$ with the Lie algebra of $G_2$ is given
in \cite[(3.33), (3.37)]{GrigorianYau1}, and the two eigenvalues of
$\alpha\mapsto(\psi_{abcd}\alpha^{cd})_{ab}$ follow from the
$G_2$-contraction identities \eqref{eq:phi-phi-contractions} and
\eqref{eq:mixed-contraction}. In particular, the contraction
$\alpha\mapsto(\varphi_{abc}\alpha^{bc})_a$ maps $\iota_v\varphi$ to $6v$,
so $\dim\Lambda^2_7=7$ and $\dim\mathfrak g_2=21-7=14$.

Write $V(j_1,j_2)$ for the irreducible
$\mathfrak g_2^{\mathbb C}$-module of highest weight
$j_1\omega_1+j_2\omega_2$, where the first simple root is short.
The fundamental weights $\omega_1,\omega_2$ give the defining
seven-dimensional and adjoint representations, respectively, so
$V(0,1)=\mathfrak g_2^{\mathbb C}$ and $W_{\mathbb C}=V(2,0)$.
We also use Bryant's identification $V(j,0)=S^j_0V_{\mathbb C}$,
the traceless symmetric $j$-tensors \cite{bryant-2003}.
The relevant decompositions are recorded in the tensor-product tables
\cite[Table~A.121]{feger-2020-lieart},
\cite[Table~1500]{yamatsu-2020-finite}.
A boldface number $\mathbf n$ denotes an irreducible representation of
dimension $n$ of the group under consideration.

\begin{proposition}
\label{prop:deformation-squares}
The symmetric and alternating squares of $W_{\mathbb C}$ are
\begin{subequations}
\begin{align}
S^2W_{\mathbb C}
 &=V(0,0)\oplus2V(2,0)\oplus V(1,1)\oplus V(0,2)\oplus V(4,0),
 \label{eq:symmetric-square}\\
\Lambda^2W_{\mathbb C}
 &=V(1,0)\oplus V(0,1)\oplus V(1,1)\oplus V(3,0)\oplus V(2,1).
 \label{eq:alternating-square}
\end{align}
\end{subequations}
Here $V(1,1)$ has dimension $64$, while the two modules of dimension
$77$ are $V(0,2)$ in the symmetric square and $V(3,0)$ in the
alternating square. In particular, $V(1,0)$ and $V(0,1)$ occur only in
the alternating square, each with multiplicity one. Moreover,
\begin{equation}
V(1,0)\otimes V(0,1)=V(1,1)\oplus V(2,0)\oplus V(1,0).
\label{eq:seven-adjoint-product}
\end{equation}
\end{proposition}

The \emph{$W$-adjoint} of $X\in\End(W_{\mathbb C})$ is its adjoint
$X^{\tadj}$ with respect to the $W$-pairing:
\begin{equation}
\label{eq:traceTw}
\langle Xh,k \rangle_W=\langle h,X^{\tadj}k \rangle_W
\qquad\text{for all }h,k\in W_{\mathbb C}.
\end{equation}
We call $X$ \emph{$W$-self-adjoint} if $X^{\tadj}=X$ and
\emph{$W$-skew-adjoint} if $X^{\tadj}=-X$. In a $W$-orthonormal basis,
$X\mapsto X^{\tadj}$ is matrix transposition, with no complex conjugation.
The $W$-pairing on $W_{\C}$ identifies an endomorphism $X\in\End(W_{\mathbb C})$ with the bilinear
form $(h,k)\mapsto\langle Xh,k\rangle_W$.
Taking its $W$-self-adjoint and $W$-skew-adjoint parts gives the equivariant
identification
\begin{equation}
\End(W_{\mathbb C})\cong S^2W_{\mathbb C}\oplus\Lambda^2W_{\mathbb C}.
\label{eq:endomorphism-inventory}
\end{equation}
The two copies of \(V(1,1)\) in \(\operatorname{End}(W_{\mathbb C})\) are
distinguished by their $W$-adjoint parity: the copy in
\(S^2W_{\mathbb C}\) consists of $W$-self-adjoint operators, while the copy in
\(\Lambda^2W_{\mathbb C}\) consists of $W$-skew-adjoint operators.

Table~\ref{tab:endomorphism-inventory} records the irreducible modules
and their multiplicities in the $W$-self-adjoint and $W$-skew-adjoint
components.
\begin{table}[htbp]
\centering
\begin{tabular}{lrrr}
Type & Dimension & $S^2W$ & $\Lambda^2W$\\
\hline
$V(0,0)$ & $1$ & $1$ & $0$\\
$V(1,0)$ & $7$ & $0$ & $1$\\
$V(0,1)$ & $14$ & $0$ & $1$\\
$V(2,0)$ & $27$ & $2$ & $0$\\
$V(1,1)$ & $64$ & $1$ & $1$\\
$V(0,2)$ & $77$ & $1$ & $0$\\
$V(3,0)$ & $77$ & $0$ & $1$\\
$V(4,0)$ & $182$ & $1$ & $0$\\
$V(2,1)$ & $189$ & $0$ & $1$
\end{tabular}
\caption{Irreducible $G_2^{\mathbb C}$-modules in $\End(W_{\mathbb C})$.
The $S^2W$ and $\Lambda^2W$ columns give their multiplicities in the
$W$-self-adjoint and $W$-skew-adjoint components, respectively, after complexification.}
\label{tab:endomorphism-inventory}
\end{table}

\subsection{The invariant cubic plane}
\label{subsec:invariant-cubic-plane}

We show that the trace cubic and the $\varphi$-cubic span the invariant
cubics on $W$, and introduce a product in terms of which we write the
$\varphi$-cubic.

\begin{definition}
\label{def:natural-cubics}
The trace cubic and the $\varphi$-cubic on $W$ are, respectively,
\begin{equation}
\begin{aligned}
C_0(h_1,h_2,h_3)
 &=\tr(h_1h_2h_3),\\
C_\varphi(h_1,h_2,h_3)
 &=\varphi_{a_1a_2a_3}\varphi_{b_1b_2b_3}
   h_1^{a_1b_1}h_2^{a_2b_2}h_3^{a_3b_3}.
\end{aligned}
\label{eq:natural-cubics}
\end{equation}
We identify each symmetric trilinear form with its associated
homogeneous cubic polynomial, writing
$C_0(h)=C_0(h,h,h)=\tr(h^3)$ and
$C_\varphi(h)=C_\varphi(h,h,h)$.
\end{definition}

The contraction in $C_\varphi$ is over all ordered index values, with no
prefactor. Both cubics are symmetric and $G_2$-invariant. By cyclicity
of the trace, $C_0$ takes the same value on three of the six orderings of
its arguments, and since the arguments are symmetric, transposition shows
that it takes this value on the other three as well.
Permuting the slots of $C_\varphi$ applies the same sign twice
to the alternating tensors. Trace invariance and direct contraction give
invariance under $\rho(O)$ \eqref{eq:rho-group} for $O\in G_2$.
The trace cubic involves only $g$ and is invariant under the action of $O(7)$:
$\tr((\rho(O)h)^3)=\tr(h^3)$.

We identify a homogeneous cubic $N(h)=N(h,h,h)$ on $W$ or $W_{\mathbb C}$
with its symmetric trilinear polarization by
\begin{equation}
\begin{aligned}
  N(h_1,h_2,h_3)=\tfrac16\bigl(&N(h_1+h_2+h_3)-N(h_1+h_2)-N(h_1+h_3)\\
                        &-N(h_2+h_3)+N(h_1)+N(h_2)+N(h_3)\bigr),
\end{aligned}
\label{eq:cubic-polarization}
\end{equation}
so $D^3N(0)[h_1,h_2,h_3]=6N(h_1,h_2,h_3)$.

Using the representation counts given above, we show in
Lemma \ref{lem:invariant-cubic-plane} below that these two cubics exhaust
the invariant cubic forms.

\begin{lemma}
\label{lem:invariant-cubic-plane}
The invariant cubic spaces are
\begin{equation}
(S^3W^*)^{G_2}=\mathbb RC_0\oplus\mathbb RC_\varphi,
\qquad\text{and}\qquad
(S^3W_{\mathbb C}^*)^{G_2^{\mathbb C}}
 =\mathbb CC_0\oplus\mathbb CC_\varphi.
\label{eq:invariant-cubic-plane}
\end{equation}
\end{lemma}

\begin{proof}
We identify $W_{\mathbb C}$ with its dual via the
complex-bilinear $W$-pairing.
As shown in Table \ref{tab:endomorphism-inventory} above, $W_{\mathbb C}\cong V(2,0)$ occurs with multiplicity two in $S^2W_{\mathbb C}$ and does not occur in $\Lambda^2W_{\mathbb C}$. Consequently, we can read off the dimensions of $G_2^{\C}$-invariant subspaces:
\[
\dim (W_{\mathbb C}^{*\otimes3})^{G_2^{\mathbb C}}=2,
\qquad
(W_{\mathbb C}^*\otimes
\Lambda^2W_{\mathbb C}^*)^{G_2^{\mathbb C}}=0.
\]
The second equality implies that every invariant
trilinear form is symmetric in any pair of arguments.
Hence
\[
\dim (S^3W_{\mathbb C}^*)^{G_2^{\mathbb C}}=2.
\]

To see that $C_0$ and $C_\varphi$ are independent,
recall that $C_0$ is $SO(7)$-invariant. For any orthonormal triple $u,v,w$,
applying the $G_2$-contraction identities \eqref{eq:phi-phi-contractions}
to the tensors \eqref{eq:test-tensors} gives
\[
C_0(E_u,E_v,E_w)=14,
\qquad
C_\varphi(E_u,E_v,E_w)
=343\,\varphi(u,v,w)^2-63.
\]
Since $\varphi(u,v,w)=g(u\times v,w)$, with $u\times v$ the cross product
\eqref{eq:cross-product}, the second
expression varies with the orthonormal triple,
whereas the first does not. Thus $C_\varphi$
is not $SO(7)$-invariant and is not proportional
to $C_0$.

The two cubics therefore span the complex invariant
space. Since both are real, the corresponding real
invariant space is their real span.
\end{proof}

With the generators and normalization of this section, we write the cubic
pencil as
\begin{equation}
\{N_{a,b}=aC_0+bC_\varphi:[a:b]\in\mathbb P^1(\mathbb C)\},
\qquad\text{with}\qquad
N_\lambda=C_0+\lambda C_\varphi,\quad\text{where }\lambda=b/a.
\label{eq:invariant-pencil}
\end{equation}
In Section~\ref{sec:exceptional-pencil}, we determine which members of the pencil have infinitesimal stabilizers larger than $\rho(\mathfrak g_2^{\mathbb C})$.

\begin{remark}[Change of generators]
\label{rem:pencil-coordinates}
An invertible change of generators of the cubic plane induces a
projective linear change of coordinates on $\mathbb P^1$.
Thus the members with larger stabilizer, classified in
Theorem~\ref{thm:stabilizer-trichotomy}, and their stabilizer dimensions
are independent of the generators, whereas their parameter values
depend on the chosen normalization. Real changes preserve real
lines and conjugate pairs.
\end{remark}

Under the standard identification $V\cong\Im\O$, Desmet's model
$S^2_0\Im\O$ identifies with $W$. The cubic plane corresponds to his
two-dimensional space of $G_2$-equivariant commutative products
\cite[Proposition~3.8]{desmet-2025-nonassociative}; see also
\cite[Section~9]{chayet-2021-class}.
We return to the product interpretation in
Section~\ref{sec:exceptional-pencil}. Interestingly, in positive characteristic, Craven gives a two-dimensional space
of invariant symmetric trilinear forms for $G_2(11)$ on its
$27$-dimensional module \cite[\S5.2.4, p.~683]{craven-2023-maximal}.

The invariant $W$-pairing identifies $C_0$ and $C_\varphi$ with two
$G_2$-equivariant symmetric bilinear products $A$ and $B$.
For symmetric tensors $h,k\in S^2V^*$, define
\begin{equation}
  A(h,k)=\frac12(hk+kh)_0,
  \label{eq:trace-product}
\end{equation}
and
\begin{equation}
  \mathcal P(h,k)_{ij}
  =\varphi_{iab}\varphi_{jcd}h^{ac}k^{bd},
  \qquad
  B(h,k)=\mathcal P(h,k)
  +\frac17\bigl(\tr(hk)-\tr(h)\tr(k)\bigr)g.
  \label{eq:phi-product}
\end{equation}
Here the subscript $0$ denotes the traceless part \eqref{eq:traceless},
so $B(h,k)=\mathcal P(h,k)_0$ by \eqref{eq:p-trace-rule} below. We will use $\mathcal P$ only where the trace part is needed.
For $h,k,l\in W$, these products satisfy
\begin{equation}\label{eq:quadratic-adjoint-pairings}
\langle A(h,k),l\rangle_W=C_0(h,k,l),
  \qquad
  \langle B(h,k),l\rangle_W=C_\varphi(h,k,l).
\end{equation}
Their quadratic maps $A(h)=A(h,h)$ and $B(h)=B(h,h)$
will be used to construct the quadratic adjoints
of the cubic pencil in Section~\ref{sec:exceptional-pencil}. As needed, we may also restrict the above products to traceless symmetric tensors, in which case, we get \(A,B:S^2W\to W\).

Our product $\mathcal P$ coincides
with $(\varphi hh\varphi)$ of Grigorian--Yau
\cite[(5.39)]{GrigorianYau1} and the operation
$A\circledcirc B$ of Dwivedi--Gianniotis--Karigiannis
\cite[(2.46)]{dwivedi-2025-flows}. Its contraction identities
follow from the standard $\varphi\varphi$ identity
\cite{bryant-2003},
\cite[Prop.~2.48]{dwivedi-2025-flows}.

The product $\mathcal P$ satisfies the following contraction identities.

\begin{lemma}
\label{lem:product-contractions}
For $v\in V$ let $P_v=v^\flat\otimes v^\flat$ be the rank-one tensor of
\eqref{eq:projection-tensors}, let $v\times v'$ be the cross product
\eqref{eq:cross-product}, and let
$h,h_1,h_2,h_3$ be symmetric $2$-tensors, not necessarily traceless. The
contraction
$\mathcal P(h_1,h_2)_{ij}
=\varphi_{ia_1a_2}\varphi_{jb_1b_2}h_1^{a_1b_1}h_2^{a_2b_2}$ of
\eqref{eq:phi-product} satisfies
\begin{subequations}
\begin{align}
\mathcal P(P_v,P_{v'})&=P_{v\times v'},
\qquad \mathcal P(P_v,P_v)=0,\label{eq:p-rank-one-rule}\\
\mathcal P(g,h)&=(\tr h)g-h,
\label{eq:p-identity-rule}\\
\tr\mathcal P(h_1,h_2)
&=(\tr h_1)(\tr h_2)-\tr(h_1h_2).
\label{eq:p-trace-rule}
\end{align}
\end{subequations}
The product is symmetric in its arguments, and
$\tr(h_3\mathcal P(h_1,h_2))$ is symmetric in all
three symmetric tensors.
\end{lemma}
\begin{proof}
For rank-one arguments the two alternating factors are
$(v\times v')_i$ and $(v\times v')_j$. Permuting the arguments in the
trilinear contraction permutes both copies of $\varphi$, so their signs
cancel. Contracting \eqref{eq:phi-phi-contractions} with $h$ gives
\begin{equation*}
\mathcal P(g,h)_{ij}
 =(\delta_{bd}\delta_{ij}-\delta_{bj}\delta_{id}+\psi_{bidj})h^{bd}
 =(\tr h)\delta_{ij}-h_{ij},
\end{equation*}
since the $\psi$ term vanishes. Finally, trilinear symmetry gives
$\tr\mathcal P(h_1,h_2)
=\tr(h_1\mathcal P(h_2,g))$, proving the trace identity.
\end{proof}

Using $\mathcal P$, the cubic can then be written
\begin{equation}
C_\varphi(h_1,h_2,h_3)
 =\tr(\mathcal P(h_1,h_2)h_3)
 =\frac12h_3^{ij}(\chi_{h_1})_{iab}(\chi_{h_2})_j{}^{ab},
\label{eq:phi-cubic-contraction}
\end{equation}
where $\chi_{h_1}=\mathcal I_\varphi(h_1)$ and $\chi_{h_2}=\mathcal I_\varphi(h_2)$, as in \eqref{eq:model-map}. These expressions thus allow the cubic to be written in terms of $3$-forms.
The second equality in \eqref{eq:phi-cubic-contraction} follows from the
quadratic insertion identity of Lemma~\ref{lem:quadratic-insertion}.

\begin{lemma}
\label{lem:quadratic-insertion}
For any $h\in S^2V^*$, not necessarily traceless,
\begin{equation}
\label{eq:chi-chi-contract1}
(\chi_h)_{iab}(\chi_h)_j{}^{ab}
 =2\mathcal P(h,h)_{ij}
  +4(\tr h)h_{ij}+2\tr(h^2)g_{ij}.
\end{equation}
\end{lemma}
\begin{proof}
To prove \eqref{eq:chi-chi-contract1}, diagonalize $h$ in an orthonormal frame, with eigenvalues
$d_i$, allowing $\varphi$ to have general components in this frame.
The components are
\begin{equation*}
(\chi_h)_{iab}=(d_i+d_a+d_b)\varphi_{iab}
\qquad\text{(no sum over }i,a,b\text{)}.
\end{equation*}
We expand the contracted product by grouping the terms according to
their eigenvalue factors. The terms involving only $d_i,d_j$ contribute
$6d_i^2\delta_{ij}$, those linear in $d_a+d_b$ contribute
$4d_i(\tr h-d_i)\delta_{ij}$, and those involving $d_a^2+d_b^2$
contribute $2(\tr(h^2)-d_i^2)\delta_{ij}$ (no sum over $i$).
The remaining $2d_ad_b$ contribution is $2\mathcal P(h,h)_{ij}$.
The $d_i^2$ terms cancel.
\end{proof}

For $h,k\in W$, polarization of \eqref{eq:chi-chi-contract1} yields
\begin{equation}
\frac12\bigl((\chi_h)_{iab}(\chi_k)_j{}^{ab}
           +(\chi_k)_{iab}(\chi_h)_j{}^{ab}\bigr)
 =2\mathcal P(h,k)_{ij}+2\tr(hk)g_{ij}.
\label{eq:polarized-insertion}
\end{equation}
Contracting \eqref{eq:polarized-insertion} with a traceless symmetric
$h_3^{ij}$ gives the second equality of \eqref{eq:phi-cubic-contraction}.

Written in terms of $3$-forms, the cubic $C_\varphi$ is, up to a constant,
the pointwise integrand of the Yukawa coupling of Lee and Leung.
In our notation, their parallel $3$-form and its Hodge dual are
$\varphi$ and $\psi$. The operation in
\cite[Definition~1.1]{Lee:2002fa} sends a $3$-form $\alpha$ to the
vector-valued one-form
\begin{equation*}
  \widehat\alpha
    =g^{ab}\,{\ast}\bigl(\alpha\wedge\iota_{e_a}\psi\bigr)\otimes e_b
    =\widehat\alpha_i{}^{b}\,e^i\otimes e_b,
\end{equation*}
and the Yukawa integrand is
$\varphi(\widehat\alpha_1,\widehat\alpha_2,\widehat\alpha_3)\wedge\psi$.
We read the evaluation of $\varphi$ on vector-valued one-forms as
evaluation on the vector components and the wedge product of the one-form
components:
$\varphi(\theta_1\otimes v_1,\theta_2\otimes v_2,\theta_3\otimes v_3)
=\varphi(v_1,v_2,v_3)\,\theta_1\wedge\theta_2\wedge\theta_3$.
Since $\psi=\ast\varphi$, the Hodge convention \eqref{eq:hodge-convention}
gives $\eta\wedge\psi=\langle\eta,\varphi\rangle\vol$ for every $3$-form
$\eta$. We wedge the Leibniz rule
$\alpha\wedge\iota_v\psi=\iota_v\alpha\wedge\psi-\iota_v(\alpha\wedge\psi)$
with a one-form and apply this identity to both resulting terms; we also
apply it to $\eta=e^i\wedge e^j\wedge e^k$. In components,
this gives
\begin{equation}
\label{eq:yukawa-hat-components}
\begin{aligned}
  \widehat\alpha_i{}^{a}
    &=\tfrac12\alpha^{a}{}_{jk}\varphi_i{}^{jk}
      -\langle\alpha,\varphi\rangle\delta_i{}^{a},\\
  \varphi(\widehat\alpha_1,\widehat\alpha_2,\widehat\alpha_3)\wedge\psi
    &=\varphi_{abc}\varphi^{ijk}\,
      (\widehat\alpha_1)_i{}^{a}(\widehat\alpha_2)_j{}^{b}
      (\widehat\alpha_3)_k{}^{c}\,\vol.
\end{aligned}
\end{equation}
For $\alpha=\chi_h$ with $h\in W$, we have $\langle\chi_h,\varphi\rangle=0$
by \eqref{eq:three-form-splitting}, and the second identity in
\eqref{eq:insertion-normalization} gives
$(\widehat{\chi_h})_i{}^{a}=2h_i{}^{a}$. Therefore
\begin{equation}
\label{eq:yukawa-integrand}
\begin{aligned}
  \varphi\bigl(\widehat{\chi_{h_1}},\widehat{\chi_{h_2}},
               \widehat{\chi_{h_3}}\bigr)\wedge\psi
    &=\tfrac18\varphi_{a_1a_2a_3}\varphi^{b_1b_2b_3}\,
      (\chi_{h_1})^{a_1}{}_{pq}\varphi_{b_1}{}^{pq}\,
      (\chi_{h_2})^{a_2}{}_{rs}\varphi_{b_2}{}^{rs}\,
      (\chi_{h_3})^{a_3}{}_{mn}\varphi_{b_3}{}^{mn}\,\vol\\
    &=8\varphi_{a_1a_2a_3}\varphi_{b_1b_2b_3}
      h_1^{a_1b_1}h_2^{a_2b_2}h_3^{a_3b_3}\,\vol
     =8C_\varphi(h_1,h_2,h_3)\,\vol.
\end{aligned}
\end{equation}
The pointwise Yukawa integrand, pulled back to $W$ by $\mathcal I_\varphi$, is
therefore $8C_\varphi$ under this evaluation convention. We also call $C_\varphi$ itself the \emph{pointwise Yukawa cubic}; it is the $\varphi$-cubic of Definition~\ref{def:natural-cubics}. A different
evaluation convention changes the scalar factor but not the line
$\C C_\varphi$, the endpoint $[0:1]$ of the pencil
\eqref{eq:invariant-pencil}. We return to this comparison in
Proposition~\ref{prop:yukawa-not-albert}.

We use these product and insertion identities to relate the cubics to
volume.

\subsection{Volume and the two cubics}
\label{subsec:volume-cubics}

The two generators of the invariant cubic plane arise, up to constant
factors, as third derivatives of volume along two different affine
paths. Fix a positive $3$-form $\varphi$ on $V$, with metric
$g=g_\varphi$, put $\psi=\ast_\varphi\varphi$, and write
\begin{equation}
\nu=\vol_g=\frac17\varphi\wedge\psi.
\label{eq:form-volume-definition}
\end{equation}
The factor $\frac17$ comes from $|\varphi|_g^2=7$. The form $\nu$ is the
pointwise density in Hitchin's volume functional
\cite[\S7.3]{Hitchin:2000jd-arxiv}.
For $h\in W$, consider the paths
\begin{equation}
\begin{aligned}
q_h(t)&=g+2th,\qquad\varphi_h(t)=\varphi+t\chi_h,
\qquad \chi_h=\mathcal I_\varphi(h),\\
g_h(t)&=g_{\varphi_h(t)},
\qquad \psi_h(t)=\ast_{g_h(t)}\varphi_h(t),
\qquad \vol_{g_h(t)}=\frac17\varphi_h(t)\wedge\psi_h(t).
\end{aligned}
\label{eq:affine-volume-paths}
\end{equation}
For sufficiently small $t$, $q_h(t)$ is a metric and $\varphi_h(t)$
is positive. The paths $q_h(t)$ and $g_h(t)$ have the same initial
metric velocity $2h$, but $g_h(t)$ is generally not
affine. All insertions and index identifications in these paths are
taken at the basepoint $(\varphi,g)$; in particular, $\chi_h$ is held fixed.
We first compute the second derivatives and then the third derivatives.

We use the pointwise deformation formulas of Grigorian--Yau
\cite[Section~5]{GrigorianYau1}. Their averaged insertion map satisfies
$i_\varphi=\mathcal I_\varphi/3$, so the symmetric tensor that they
associate with our direction $\chi_h$ is $3h$. With this substitution, the second-order
terms of their (5.19), (5.35), (5.43), and (5.48) give
\begin{subequations}
\begin{align}
g_h(t)
 &=g+2th+t^2\left(2h^2-\frac12\mathcal P(h,h)
                   -\frac16\tr(h^2)g\right)+\mathcal O(t^3),
 \label{eq:affine-form-metric-expansion}\\
\psi_h(t)
 &=\psi-t\ast\chi_h
   +t^2\ast\left(\frac12\mathcal I_\varphi(B(h,h))
                  -\frac1{21}\tr(h^2)\varphi\right)
   +\mathcal O(t^3).
 \label{eq:affine-form-duality-expansion}
\end{align}
\end{subequations}
Here and below, a Hodge star without a subscript is the star of the
basepoint, and $h^2$ is the contracted square
\eqref{eq:contracted-product}.
The expansions are Taylor expansions in $t$, so each second derivative is
twice the coefficient of $t^2$. Polarizing gives the following.

\begin{proposition}
\label{prop:form-metric-response}
For $h_1,h_2,h_3\in W$, with the ambient form directions fixed,
\begin{subequations}
\begin{align}
D_\varphi^2g[\chi_{h_1},\chi_{h_2}]
 &=2(h_1h_2+h_2h_1)-\mathcal P(h_1,h_2)-\frac13\tr(h_1h_2)g,
 \label{eq:form-metric-hessian}\\
\ast D_\varphi^2(\ast_\varphi\varphi)[\chi_{h_1},\chi_{h_2}]
 &=\mathcal I_\varphi(B(h_1,h_2))
   -\frac{2}{21}\tr(h_1h_2)\varphi,
 \label{eq:duality-hessian}\\
D_\varphi\langle\chi_{h_1},\chi_{h_2}\rangle_{g_\varphi}[\chi_{h_3}]
 &=-2C_\varphi(h_1,h_2,h_3).
 \label{eq:form-pairing-response}
\end{align}
\end{subequations}
\end{proposition}
\begin{proof}
The first two identities are the polarized expansions
\eqref{eq:affine-form-metric-expansion} and
\eqref{eq:affine-form-duality-expansion}. For the
last identity, $D_\varphi g[\chi_{h_3}]=2h_3$, so differentiating the three
inverse metrics in the inner product of the fixed forms gives
\begin{equation*}
D_\varphi\langle\chi_{h_1},\chi_{h_2}\rangle_{g_\varphi}[\chi_{h_3}]
 =-h_3^{ij}(\chi_{h_1})_{iab}(\chi_{h_2})_j{}^{ab}
 =-2C_\varphi(h_1,h_2,h_3),
\end{equation*}
by \eqref{eq:phi-cubic-contraction}.
\end{proof}

We now differentiate the volume itself. We regard $\vol$ both as the
function $q\mapsto\vol_q$ on metrics and as the function
$\varphi\mapsto\vol_{g_\varphi}$ on positive forms; the subscript on $D$
indicates which is differentiated, and $\nu$ always denotes the value at the
basepoint. For a metric velocity $2h$, the determinant formula gives
\begin{equation}
\log\frac{\vol_{q_h(t)}}{\nu}
 =\frac12\tr\log(I+2t\,g^{-1}h)
 =-t^2\tr(h^2)+\frac43t^3\tr(h^3)+\mathcal O(t^4).
\label{eq:metric-volume-log-expansion}
\end{equation}
Thus the first derivative vanishes, and the Hessian of the metric volume is
\begin{equation}
\frac{D_g^2\vol[2h,2k]}{\nu}=-2\tr(hk).
\label{eq:metric-volume-hessian}
\end{equation}

For positive forms, we call $\varphi\mapsto\ast_\varphi\varphi$ the duality
map. The
first-variation identities, valid at every positive form, are as follows;
for the derivative of the duality map, see Hitchin
\cite[Lemma~20]{Hitchin:2000jd-arxiv}.
\begin{subequations}
\begin{align}
D_\varphi\vol[\xi]&=\frac13\xi\wedge\psi,
 \label{eq:volume-first-variation}\\
D_\varphi(\ast_\varphi\varphi)[\xi]
 &=\ast\left(\frac43\pr_1\xi+\pr_7\xi-\pr_{27}\xi\right).
 \label{eq:duality-first-variation}
\end{align}
\end{subequations}
The projections are those of the form at which the derivative is
taken. To see the factor $1/3$, differentiate
$\vol_{g_\varphi}=\varphi\wedge(\ast_\varphi\varphi)/7$: the second identity gives
$\varphi\wedge D_\varphi(\ast_\varphi\varphi)[\xi]
=\frac43\xi\wedge\psi$. Along the affine path in
\eqref{eq:affine-volume-paths}, these identities yield
\begin{equation*}
\frac{d}{dt}\vol_{g_h(t)}=\frac13\chi_h\wedge\psi_h(t),
\qquad
\left.\frac{d^2}{dt^2}\right|_{t=0}\vol_{g_h(t)}
 =\frac13\chi_h\wedge\dot\psi_h(0)
 =-\frac13|\chi_h|^2\nu.
\end{equation*}
Since $|\chi_h|^2=2\tr(h^2)$, polarization gives
\begin{equation}
\frac{D_\varphi^2\vol[\chi_h,\chi_k]}{\nu}
 =-\frac13\langle\chi_h,\chi_k\rangle
 =-\frac23\tr(hk).
\label{eq:form-volume-hessian}
\end{equation}
In particular, the two paths with the same initial metric velocity
have different volume Hessians.

We can now take the third derivatives. For mixed derivatives the paths
are replaced by the affine families
\begin{equation*}
q(s,t,r)=g+2(sh_1+th_2+rh_3),
\qquad
\varphi(s,t,r)=\varphi+s\chi_{h_1}+t\chi_{h_2}+r\chi_{h_3}.
\end{equation*}
The direction forms are fixed in $\Lambda^3V^*$, rather than transported
into the moving $\Lambda^3_{27}$ subspaces.

\begin{proposition}
\label{prop:affine-volume-cubics}
With flat affine derivatives on the cones of metrics and positive
$3$-forms, for $h_1,h_2,h_3\in W$,
\begin{equation}
\frac{D_g^3\vol[2h_1,2h_2,2h_3]}{\nu}=8C_0(h_1,h_2,h_3),
\qquad
\frac{D_\varphi^3\vol[\chi_{h_1},\chi_{h_2},\chi_{h_3}]}{\nu}
 =\frac23C_\varphi(h_1,h_2,h_3).
\label{eq:affine-volume-cubics}
\end{equation}
\end{proposition}
\begin{proof}
For the metric path, the vanishing first derivative means that the
third derivatives of the volume ratio and its logarithm agree at zero.
The determinant expansion \eqref{eq:metric-volume-log-expansion} gives $8\tr(h^3)$, and polarization
gives the first identity.

For the form path, differentiate \eqref{eq:volume-first-variation}
twice in fixed ambient directions and use \eqref{eq:duality-hessian}:
\begin{align*}
D_\varphi^3\vol[\chi_{h_1},\chi_{h_2},\chi_{h_3}]
 &=\frac13\chi_{h_3}\wedge D_\varphi^2(\ast_\varphi\varphi)[\chi_{h_1},\chi_{h_2}]\\
 &=\frac13\left\langle\chi_{h_3},
       \mathcal I_\varphi(B(h_1,h_2))
       -\frac{2}{21}\tr(h_1h_2)\varphi\right\rangle\nu\\
 &=\frac23\tr\bigl(h_3B(h_1,h_2)\bigr)\nu
  =\frac23C_\varphi(h_1,h_2,h_3)\nu.
\end{align*}
Here $\chi_{h_3}\perp\varphi$, the insertion map multiplies the
$W$-pairing by $2$, and $h_3$ is traceless. Equivalently, the one-parameter
volume expansions are
\begin{align*}
\frac{\vol_{q_h(t)}}{\nu}
 &=1-t^2\tr(h^2)+\frac43t^3C_0(h)+\mathcal O(t^4),\\
\frac{\vol_{g_h(t)}}{\nu}
 &=1-\frac13t^2\tr(h^2)+\frac19t^3C_\varphi(h)
   +\mathcal O(t^4).
\end{align*}
Thus the second-order term of $\psi_h(t)$ in Grigorian--Yau's (5.48)
already determines the third derivative of the volume density of positive
forms.
\end{proof}

Thus $C_\varphi$ measures both the third affine derivative of the
positive-form volume and the change of the inner product on fixed
ambient $3$-forms in \eqref{eq:form-pairing-response}. For the metric
volume, the corresponding cubic is $C_0$.

An arbitrary nonlinear metric curve with initial velocity $2h$ need
not have these volume derivatives. For instance, along a curve in the fixed-volume
hypersurface of metrics the volume is constant, so all its derivatives vanish.

The second derivative identities are the polarized second-order terms
of Grigorian--Yau's expansions \cite[Section~5]{GrigorianYau1}.
The contraction formula \eqref{eq:phi-cubic-contraction} appears in
Karigiannis--Leung \cite[Lemma~3.12]{karigiannis-2007a}, and the third
derivative identity appears in integrated form in
\cite[Theorem~3.10]{karigiannis-2007a}.
A cubic $C_s$ on $W$, with $s\in\{0,\varphi\}$, is transported to
$\Lambda^3_{27}V^*$ through $\mathcal I_\varphi^{-1}$. Replacing
$\mathcal I_\varphi$ by a nonzero multiple, with the cubics on $W$ fixed,
multiplies both transported cubics by the inverse cube of the multiple.
For the isometry $\mathcal I_\varphi/\sqrt2$, this factor is $2\sqrt2$.
Since both cubics are scaled by the same factor, the pencil parameter is
unchanged.

In Section~\ref{sec:exceptional-pencil}, we use the invariant cubic
plane, its pencil, the product $B$ and the decomposition of
$\End(\WC)$ into irreducible representations to locate the parameters at
which the stabilizer jumps.
\section{The exceptional cubic pencil}
\label{sec:exceptional-pencil}

\subsection{The question and the theorem}
\label{subsec:pencil-question-theorem}

In Section~\ref{sec:deformation-representation} we identified the $G_2$-invariant
cubics on $W_{\mathbb C}$ as the pencil $N_{a,b}=aC_0+bC_\varphi$.
Every member therefore has infinitesimal symmetries
$\rho(\mathfrak g_2^{\mathbb C})$. We determine the parameters
$[a:b]\in\mathbb P^1(\mathbb C)$ for which the stabilizer is larger and
identify the additional symmetries.

We use the splitting \eqref{eq:two-form-splitting} (with
$\Lambda^2_{14}=\mathfrak g_2$), the $G_2$-contraction identities of
Lemma~\ref{lem:standard-contractions}, and the contraction rules of
Lemma~\ref{lem:product-contractions} from
Section~\ref{sec:deformation-representation}.
The contraction data are stated where they enter the proof; their detailed
derivations are collected in Appendix~\ref{app:contraction-proofs}.

We work over $\mathbb C$ unless a statement explicitly allows
$\mathbb F\in\{\mathbb R,\mathbb C\}$; then $\gtwo^{\mathbb F}$ denotes
$\gtwo$ for $\mathbb F=\mathbb R$ and $\gtwoC$ for $\mathbb F=\mathbb C$. As in
Section~\ref{sec:deformation-representation}, we identify a cubic polynomial
with its symmetric trilinear polarization, normalized by $N(h)=N(h,h,h)$
as in \eqref{eq:cubic-polarization}.
All pairings and adjoints on $W_{\mathbb C}$ use the complex-bilinear
extension of the $W$-pairing \eqref{eq:deformation-module}; in particular
$X^{\tadj}$ is the $W$-adjoint \eqref{eq:traceTw}. In the stabilizer formulas, $X$ (and in Section~\ref{subsec:explicit-64}, $L$) denotes
an endomorphism of $W_{\mathbb F}$.

\begin{definition}
\label{def:exact-stabilizer}
Let $N$ be a cubic on $W_{\mathbb F}$. Its \emph{exact infinitesimal stabilizer} is
\begin{equation*}
\stab(N)=\bigl\{X\in\End(W_{\mathbb F}):\ N(Xh,h,h)=0\ \text{for all }h\in W_{\mathbb F}\bigr\}.
\end{equation*}
Its \emph{action map} is $\Phi_N\colon\End(W_{\mathbb F})\to S^3W_{\mathbb F}^{*}$,
\begin{equation}
\Phi_N(X)(h_1,h_2,h_3)=N(Xh_1,h_2,h_3)+N(h_1,Xh_2,h_3)+N(h_1,h_2,Xh_3).
\label{eq:action-map}
\end{equation}
\end{definition}

Since $N$ is symmetric, $\Phi_N(X)$ is a symmetric trilinear form, and
$\Phi_N(X)(h,h,h)=3N(Xh,h,h)$. A symmetric trilinear form vanishes when its cubic does, so
\begin{equation*}
\stab(N)=\ker\Phi_N .
\end{equation*}
For the pencil members we write $\Phi_{a,b}=\Phi_{N_{a,b}}=a\Phi_{C_0}+b\Phi_{C_\varphi}$.
The \emph{exact} stabilizer preserves the cubic itself, whereas preserving
its one-dimensional span $\mathbb F N$ requires only
$\Phi_N(X)\in\mathbb F N$. We return to stabilizers of cubic lines in Sections~\ref{subsec:explicit-64} and~\ref{sec:e6-albert-norm}.

The $G_2$-invariance of the generators in \eqref{eq:natural-cubics}
has two immediate consequences for the action map.

\begin{lemma}
\label{lem:common-g2}
For every $[a:b]\in\mathbb P^1(\mathbb F)$, with
$\mathbb F\in\{\mathbb R,\mathbb C\}$,
\begin{equation*}
\rho(\mathfrak g_2^{\mathbb F})\subseteq\ker\Phi_{a,b}=\stab(N_{a,b}).
\end{equation*}
Moreover, $\Phi_{a,b}$ is $\mathfrak g_2^{\mathbb F}$-equivariant:
\begin{equation*}
\Phi_{a,b}([\rho(M),X])=M\cdot\Phi_{a,b}(X).
\end{equation*}
Here $M\in\mathfrak g_2^{\mathbb F}$ acts on endomorphisms by commutators
with $\rho(M)$ and on cubic tensors by the induced action.
Consequently, every stabilizer is a $\mathfrak g_2^{\mathbb F}$-submodule
of $\End(W_{\mathbb F})$.
\end{lemma}

\begin{proof}
Differentiating the invariance of $N_{a,b}$ using \eqref{eq:rho-action}
and \eqref{eq:action-map} gives the containment. By the same
invariance, the action map is equivariant, so its kernel is a submodule.
Both statements extend complex-linearly.
\end{proof}

We show in Section~\ref{sec:e6-albert-norm} that the stabilizers
of any two distinct members intersect precisely in
$\rho(\mathfrak g_2^{\mathbb F})$.

The main theorem of this section determines the dimension of the
stabilizer of every member of the pencil, including the parameters at which the stabilizer is larger
than $\rho(\mathfrak g_2^{\mathbb C})$. We use the normalization of
$C_0$ and $C_\varphi$ fixed in Definition~\ref{def:natural-cubics};
Remark~\ref{rem:pencil-coordinates} describes how the parameters
depend on this choice.

\begin{theorem}
\label{thm:stabilizer-trichotomy}
For $[a:b]\in\mathbb P^1(\mathbb C)$, let $N_{a,b}=aC_0+bC_\varphi$ on $W_{\mathbb C}$. Then
\begin{equation}
\dim_{\mathbb C}\stab(N_{a,b})=
\begin{cases}
21, & [a:b]=[1:0],\\
78, & a^2-2ab+8b^2=0,\\
14, & \text{otherwise, including }[0:1].
\end{cases}
\label{eq:stabilizer-trichotomy}
\end{equation}
The quadratic $a^2-2ab+8b^2$ has two distinct non-real zeros in $\mathbb P^1(\mathbb C)$. In the affine
coordinate $\lambda=b/a$ they are the roots of $8\lambda^2-2\lambda+1=0$,
\begin{equation}
\lambda_{\pm}=\frac{1\pm\ii\sqrt7}{8}.
\label{eq:critical-parameters}
\end{equation}
\end{theorem}

We use the critical cubics $N_\pm$ and critical lines $\mathbb CN_\pm$
introduced in Section~\ref{sec:introduction}, with the parameters
$\lambda_\pm$ of \eqref{eq:critical-parameters}.
Since $C_0$ and $C_\varphi$ are real and linearly
independent, complex conjugation exchanges $N_+$
and $N_-$ and their stabilizers. Neither critical
line contains a nonzero real cubic. Consequently,
by \eqref{eq:stabilizer-trichotomy}, the real pencil has exactly one member with
stabilizer dimension greater than $14$: the trace
cubic $C_0$, whose stabilizer has dimension $21$.

The proof of Theorem~\ref{thm:stabilizer-trichotomy} has four stages.
\begin{enumerate}
\item We develop the contraction identities of the two cubic tensors on
$W$, the c-contraction identities, including two universal quartic
identities. Their symbolic proof in
Section~\ref{subsec:quartic-identities} uses the evaluated one-index
c-contraction identity, a formal resolvent and a triple-form contraction (Lemma~\ref{lem:triple-form-reduction}).
\item We characterize the affine members $N_\lambda=C_0+\lambda C_\varphi$
whose quadratic adjoints satisfy a double-adjoint identity
(Proposition~\ref{prop:double-adjoint}). Evaluation on elementary test
tensors proves necessity, and the quartic identities prove sufficiency
and produce the parameters $\lambda_\pm$.
\item Using the cubic contractions, we exclude five isotypic components
from every kernel. We treat the double $2V(2,0)$ component with the
four-map Gram matrix $\mathcal G$, and the multiplicity-one components
with four $2\times2$ matrices $\Gamma$. Their contraction proofs are
given in Section~\ref{subsec:stabilizer-contractions}; the invertibility
argument remains in Section~\ref{subsec:stabilizer-exclusion}.
These exclusions and the Koecher operators
constructed from the double adjoints give $\dim\stab(N_\pm)\ge78$
(Proposition~\ref{prop:koecher-jump}).
\item Finally, the trace-cubic stabilizer and the Koecher lower bound
determine the two critical kernels inside the double $V(1,1)$ component.
A single explicit computation determines their graphs in the parity
model. From the resulting relative operator we obtain the exact kernel
dimension throughout the projective pencil
(Proposition~\ref{prop:block-kernels}), including the endpoint
$C_\varphi$, and rule out all other jumps.
\end{enumerate}

An isotypic component is the sum of all irreducible submodules of a fixed
type. By equivariance, components of different types do not mix, and by
Schur's lemma, maps between repeated copies reduce to small matrices.
This approach has recently been applied to the octonionic structure operator,
see \cite{grigorian-2026-structure-operator,grigorian-2026-two-step}. In
this paper, the main example is the pair of copies of the
$64$-dimensional module. We use the decompositions from Section~\ref{sec:deformation-representation},
Schur's lemma and complete reducibility throughout these module arguments.
We use the classification of complex simple Lie algebras only in
Theorem~\ref{thm:critical-type}, after the dimensions have been established,
to identify the critical stabilizer as $\mathfrak e_6$.

\subsection{Quadratic adjoints and c-contraction identities}
\label{subsec:cubic-contraction-calculus}

Recall the $G_2$-equivariant symmetric bilinear products
$A,B:S^2W\to W$ from \eqref{eq:trace-product}
and \eqref{eq:phi-product}. We write
$A(h)=A(h,h)$ and $B(h)=B(h,h)$, and for $\lambda\in\mathbb C$  set
\[
  F_\lambda=A+\lambda B.
\]
We use the same symbol $F_\lambda$ for its halved polarization:
\begin{equation}
\begin{aligned}
F_\lambda(h_1,h_2)
&=\tfrac12\bigl(F_\lambda(h_1+h_2)-F_\lambda(h_1)-F_\lambda(h_2)\bigr)\\
&=A(h_1,h_2)+\lambda B(h_1,h_2).
\end{aligned}
\label{eq:F-lambda-def}
\end{equation}
Thus $F_\lambda(h,h)=F_\lambda(h)$, and
\eqref{eq:quadratic-adjoint-pairings} gives
\begin{equation}
\langle F_\lambda(h_1,h_2),h_3\rangle_W=N_\lambda(h_1,h_2,h_3).
\label{eq:F-lambda-pairing}
\end{equation}
These identities extend complex-bilinearly to $W_\C$.

\begin{definition}
\label{def:quadratic-adjoint}
The map $F_\lambda$ is the \emph{quadratic adjoint} of $N_\lambda$
relative to the $W$-pairing. It is characterized by
\begin{equation}
N_\lambda(h,h,k)=\langle F_\lambda(h),k\rangle_W.
\label{eq:F-adjoint}
\end{equation}
In particular, $dN_\lambda|_h(k)=3\langle F_\lambda(h),k\rangle_W$,
so $F_\lambda=\tfrac13\nabla N_\lambda$ is one third of the $W$-gradient
of $N_\lambda$, where the $W$-gradient $\nabla N$ of a cubic $N$ is defined by
$\langle\nabla N(h),k\rangle_W=dN|_h(k)$.
For the critical parameters we write $F_\pm=F_{\lambda_\pm}$.
\end{definition}

The terminology comes from the adjoint, or sharp, map of a cubic norm
structure. However, in this definition we do not assume a Jordan algebra structure.

Fix a real orthonormal basis $(k_\alpha)_{\alpha=1}^{27}$ of $W$ for the
$W$-pairing of Definition~\ref{def:deformation-module}
and use it also as a basis of $W_{\mathbb C}$, with the pairing extended
complex-bilinearly. Greek indices run from $1$ to $27$ and
are raised and lowered with $\delta_{\alpha\beta}$; Latin indices still
run from $1$ to $7$ on $V$. All repeated indices are contracted with
these bilinear metrics, without complex conjugation.

For such a basis we have the \emph{projection identity}
\begin{equation}
\sum_{\alpha=1}^{27}(k_\alpha)_{ij}(k_\alpha)^{mn}
 =\frac12(\delta_i^m\delta_j^n+\delta_i^n\delta_j^m)
   -\frac17g_{ij}g^{mn}.
\label{eq:projection-identity}
\end{equation}
Indeed, adjoining $g/\sqrt7$ gives a basis of $S^2V^*$ that is orthonormal
for the $W$-pairing $\tr(h_1h_2)$, extended to $S^2V^*$ by the same formula;
subtracting the rank-one projection onto $g/\sqrt7$ from the symmetrization gives
\eqref{eq:projection-identity}. Thus \eqref{eq:projection-identity}
represents, on covariant $2$-tensors $\mathsf X$, the projection
$\mathsf X\mapsto\tfrac12(\mathsf X+\mathsf X^{\transp})
-\tfrac17(\tr\mathsf X)g$ onto $W$, and it is independent of the chosen basis.

The superscripts \(0\) and \(\varphi\) indicate that a quantity is
associated with \(C_0\) or \(C_\varphi\), respectively. For
$s,t\in\{0,\varphi\}$ define
\begin{equation}
c^s_{\alpha\beta\gamma}=C_s(k_\alpha,k_\beta,k_\gamma),
\qquad
\mathsf T^{st}_{\alpha\beta\gamma\delta}
=c^s_{\alpha\beta\sigma}(c^t)^\sigma{}_{\gamma\delta}.
\label{eq:cubic-tensor-coordinates}
\end{equation}
Each $c^s$ is fully symmetric. The pairings \eqref{eq:quadratic-adjoint-pairings} between the products \(A\) and \(B\) and the two cubics give
their coordinate description:
\begin{equation*}
\begin{aligned}
A(h_1,h_2)^\alpha&=(c^0)^\alpha{}_{\beta\gamma}h_1^\beta h_2^\gamma,
&B(h_1,h_2)^\alpha&=(c^\varphi)^\alpha{}_{\beta\gamma}h_1^\beta h_2^\gamma,\\
F_\lambda(h)^\alpha
&=(c^0+\lambda c^\varphi)^\alpha{}_{\beta\gamma}h^\beta h^\gamma.
\end{aligned}
\end{equation*}
We call the four-index tensors $\mathsf T^{st}$, which contract one
index of $c^s$ with one of $c^t$, the \emph{one-index c-contractions};
contracting two indices gives the \emph{two-index c-contractions} below.
Formulas that evaluate them are the \emph{c-contraction identities}.
The one-index c-contractions are symmetric within each pair and satisfy
\begin{equation*}
\mathsf T^{st}_{\alpha\beta\gamma\delta}
=\mathsf T^{ts}_{\gamma\delta\alpha\beta}.
\end{equation*}
Thus interchanging the two cubic types also interchanges the two pairs
of indices. The two mixed one-index c-contractions $\mathsf T^{0\varphi}$
and $\mathsf T^{\varphi0}$ need not agree when all four indices are kept
in place.

Contracting two indices produces an invariant bilinear form on \(W\), hence a scalar multiple of the \(W\)-pairing. The resulting $2\times2$ matrix $\mathsf a=(\mathsf a_{st})$ of scalars, indexed by $s,t\in\{0,\varphi\}$, computes the pairings of the multiplication operators associated with \(A\) and \(B\). Its nonsingularity will distinguish the two copies of \(W\) in \(\operatorname{End}(W)\). The two-index c-contraction identity is
\begin{equation}
c^s_{\alpha\beta\gamma}(c^t)_\delta{}^{\beta\gamma}
=\mathsf a_{st}\delta_{\alpha\delta},
\qquad
\mathsf a=\begin{pmatrix}55/28&11/14\\11/14&159/7\end{pmatrix},
\qquad \det\mathsf a=44.
\label{eq:two-index-c-contraction}
\end{equation}
We derive all three entries in Section~\ref{subsec:cubic-tensor-contractions}
from the projection identity \eqref{eq:projection-identity} and the
one-index c-contraction identities. Identity
\eqref{eq:two-index-c-contraction} is tensorial and does not depend on
the choice of orthonormal basis.

In terms of the products, the
one-index c-contractions are
\[
\mathsf T^{st}(h_1,h_2;h_3,h_4)
 =\langle F^{(s)}(h_1,h_2),F^{(t)}(h_3,h_4) \rangle_W,
 \qquad F^{(0)}=A,\quad F^{(\varphi)}=B,
\]
where $s,t\in\{0,\varphi\}$ and $h_1,h_2,h_3,h_4\in W$.
Expanding these contractions gives
\begin{equation}
\begin{aligned}
\mathsf T^{00}(h_1,h_2;h_3,h_4)
 &=\tfrac12\tr(h_1h_2h_3h_4+h_1h_2h_4h_3)
   -\tfrac17 \langle h_1,h_2 \rangle_W \langle h_3,h_4 \rangle_W,\\
\mathsf T^{0\varphi}(h_1,h_2;h_3,h_4)
 &=C_\varphi(h_3,h_4,\tfrac12(h_1h_2+h_2h_1))
   +\tfrac17 \langle h_1,h_2 \rangle_W \langle h_3,h_4 \rangle_W,\\
\mathsf T^{\varphi0}(h_1,h_2;h_3,h_4)
 &=C_\varphi(h_1,h_2,\tfrac12(h_3h_4+h_4h_3))
   +\tfrac17 \langle h_1,h_2 \rangle_W \langle h_3,h_4 \rangle_W.
\end{aligned}
\label{eq:one-index-c-contraction-mixed}
\end{equation}
Here the contraction defining $C_\varphi$ is extended to
all symmetric matrices when its third argument is not
traceless.

To evaluate the fourth contraction, define the vector obtained by
contracting the commutator $[h_1,h_2]$ with $\varphi$,
\[
\eta(h_1,h_2)_i
 =\varphi_{iab}(h_1h_2)^{ab}
 =\tfrac12\varphi_{iab}[h_1,h_2]^{ab},
\]
so that $\eta(h_2,h_1)=-\eta(h_1,h_2)$, and let $\Xi$ be the sum of the six
mixed cubic terms,
\begin{equation*}
\begin{aligned}
\Xi(h_1,h_2,h_3,h_4)={}&C_\varphi(h_1,h_2,A(h_3,h_4))
 +C_\varphi(h_1,h_3,A(h_2,h_4))
 +C_\varphi(h_1,h_4,A(h_2,h_3))\\
 &+C_\varphi(h_2,h_3,A(h_1,h_4))
 +C_\varphi(h_2,h_4,A(h_1,h_3))
 +C_\varphi(h_3,h_4,A(h_1,h_2)).
\end{aligned}
\end{equation*}
Then
\begin{equation}
\begin{aligned}
\mathsf T^{\varphi\varphi}(h_1,h_2;h_3,h_4)={}&
 \tfrac27\langle h_1,h_2 \rangle_W \langle h_3,h_4 \rangle_W
 +\tfrac{10}{7}\langle h_1,h_3 \rangle_W \langle h_2,h_4 \rangle_W
 +\tfrac{10}{7}\langle h_1,h_4 \rangle_W \langle h_2,h_3 \rangle_W\\
 &-8\tr(h_1h_3h_2h_4)+2\Xi(h_1,h_2,h_3,h_4)\\
 &-2\bigl(\eta(h_1,h_3)\cdot\eta(h_2,h_4)
          +\eta(h_1,h_4)\cdot\eta(h_2,h_3)\bigr).
\end{aligned}
\label{eq:one-index-c-contraction-phiphi}
\end{equation}
The dot denotes the complex-bilinear metric pairing
on $V_\mathbb C$. We prove these identities in
Section~\ref{subsec:cubic-tensor-contractions}
by contracting the cubic products using the
projection identity \eqref{eq:projection-identity} and evaluating the
resulting $4$-form contraction.

We obtain a scalar quartic identity by setting
$h_1=h_2=h_3=h_4=h$ in
\eqref{eq:one-index-c-contraction-phiphi}.
On the left-hand side,
\[
\mathsf T^{\varphi\varphi}(h,h;h,h)
 =\langle B(h),B(h)\rangle_W.
\]
Since $\eta(h,h)=0$ and
\[
\Xi(h,h,h,h)
 =6\langle A(h),B(h)\rangle_W,
\]
the contraction identity gives
\[
\begin{aligned}
\langle B(h),B(h)\rangle_W
 ={}&\tfrac{22}{7}\langle h,h\rangle_W^2
     -8\tr(h^4)\\
    &+12\langle A(h),B(h)\rangle_W.
\end{aligned}
\]
Using
\[
\tr(h^4)
 =\langle A(h),A(h)\rangle_W
  +\tfrac17\langle h,h\rangle_W^2,
\]
we therefore obtain
\begin{equation}
\begin{aligned}
8\langle A(h),A(h)\rangle_W
 -12\langle A(h),B(h)\rangle_W
 +\langle B(h),B(h)\rangle_W
 =2\langle h,h\rangle_W^2.
\end{aligned}
\label{eq:scalar-cubic-contraction}
\end{equation}
Equivalently, polarization gives
\begin{equation*}
\operatorname{Sym}_{\alpha\beta\gamma\delta}
 \bigl(8\mathsf T^{00}_{\alpha\beta\gamma\delta}
       -12\mathsf T^{0\varphi}_{\alpha\beta\gamma\delta}
       +\mathsf T^{\varphi\varphi}_{\alpha\beta\gamma\delta}\bigr)
 =2\operatorname{Sym}_{\alpha\beta\gamma\delta}
 \bigl(\delta_{\alpha\beta}\delta_{\gamma\delta}\bigr).
\end{equation*}
Here symmetrization is normalized to average over
all permutations of the displayed indices.

Although this scalar identity yields further relations
by polarization and differentiation, it does not by
itself establish the $W$-valued quartic identities
needed below. Their proofs also use the triple-form reduction of Lemma~\ref{lem:triple-form-reduction} in
Section~\ref{subsec:cubic-tensor-contractions} and the
formal resolvent of Section~\ref{subsec:quartic-identities}.

The double-adjoint argument will use the two general identities of
Lemma~\ref{lem:universal-quartic-identities} below. Their proof is given in
Section~\ref{subsec:quartic-identities}. Composing two quadratic adjoints produces six quartic expressions formed from \(A\) and \(B\). The following two identities relate those expressions to \(C_0(h)h\) and \(C_\varphi(h)h\).

\begin{lemma}
\label{lem:universal-quartic-identities}
For every $h\in\WC$,
\begin{gather}
  8A(A(h))-A(B(h))+2B(A(h))+2B(A(h),B(h))=\tfrac{16}{3}C_0(h)h,
  \label{eq:quartic-identity-one}\\
  16A(A(h),B(h))+2A(B(h))-8B(A(h))+B(B(h))=\tfrac{16}{3}C_\varphi(h)h.
  \label{eq:quartic-identity-two}
\end{gather}
\end{lemma}

These two identities supply the sufficiency argument in
Proposition~\ref{prop:double-adjoint}: at the critical parameters, the
composition of the quadratic adjoints is a linear combination of their
left-hand sides.

\subsection{Double-adjoint characterization of the critical cubics}
\label{subsec:critical-members-jordan}

For a cubic norm structure, the classical adjoint identity
$(x^\sharp)^\sharp=N(x)x$ says that applying the adjoint twice returns
the algebra element $x$ multiplied by its norm. Cubic norm structures and their adjoints are recalled in
Section~\ref{subsec:albert-norm}, see \eqref{eq:cubic-norm-structure}.
This identity motivates the following question.
Here we allow two different quadratic adjoints $F_\lambda$ and $F_\mu$
(Definition~\ref{def:quadratic-adjoint}) and ask when
\begin{equation}
F_\mu\bigl(F_\lambda(h)\bigr)=\kappa N_\lambda(h)h
\qquad\text{for every }h\in W_{\mathbb C}.
\label{eq:double-adjoint}
\end{equation}
The parameters $\lambda,\mu,\kappa\in\mathbb C$ are initially unrestricted.

\begin{definition}
\label{def:associative-plane}
A real $3$-plane $\Sigma\subset V$ is \emph{associative} if
$\lvert\varphi(a,b,c)\rvert=1$ for an orthonormal basis $(a,b,c)$ of
$\Sigma$. The calibration $\varphi$ then specifies its orientation.
\end{definition}

\begin{proposition}
\label{prop:double-adjoint}
Let $\lambda,\mu,\kappa\in\mathbb C$. Then, \eqref{eq:double-adjoint} holds if and only if $(\lambda,\mu,\kappa)=(\lambda_\pm,\lambda_\mp,\tfrac23)$. In particular, identically on
$W_{\mathbb C}$,
\begin{equation*}
F_-\bigl(F_+(h)\bigr)=\tfrac23N_+(h)\,h,
\qquad\text{and}\qquad
F_+\bigl(F_-(h)\bigr)=\tfrac23N_-(h)\,h .
\end{equation*}
\end{proposition}

\begin{proof}
\emph{Necessity.}
Assume \eqref{eq:double-adjoint} holds. We will use three test tensors in \eqref{eq:double-adjoint} to determine $(\lambda,\mu,\kappa)$. Recall from
\eqref{eq:projection-tensors} that $P_i=e_i^\flat\otimes e_i^\flat\in S^2V^*$
is the projection tensor of $\mathbb R e_i$, so that
$g^{-1}P_i=\pr_{\mathbb Re_i}$.  For the standard $G_2$ $3$-form, the $3$-plane
$\Sigma=\Span(e_1,e_2,e_3)$ is associative, with $e_1\times e_2=e_3$ and its cyclic
variants. Its projection tensor is $P_\Sigma=P_1+P_2+P_3$. We then use the test tensors
\eqref{eq:test-tensors} of the line $\mathbb Re_1$ and of $\Sigma$, and also define $\Delta$:
\[
E_1=7P_1-g,\qquad E_\Sigma=7P_\Sigma-3g,\qquad \Delta=P_1-P_2.
\]
Clearly, the tensors $E_1,E_\Sigma,\Delta$ are traceless. Geometrically, $E_1=7(P_1)_0$ and $E_\Sigma=7(P_\Sigma)_0$ are multiples of the traceless parts of the projection tensors of a line and of the associative $3$-plane, respectively, while $\Delta$ has rank $2$.  We will use the identities
\eqref{eq:p-rank-one-rule} and \eqref{eq:p-identity-rule}:
\[
\mathcal P(P_u,P_v)=P_{u\times v},\qquad
\mathcal P(g,h)=(\tr h)g-h.
\]
Here $P_u=u^\flat \otimes u^{\flat}$ as in \eqref{eq:projection-tensors}, and $h$ is any symmetric $2$-tensor. The preceding identities give
\[
\mathcal P(P_1,P_1)=0,\qquad
\mathcal P(P_\Sigma,P_\Sigma)=2P_\Sigma,\qquad
\mathcal P(g,g)=6g.
\]
For example, the six ordered pairs of distinct basis vectors in the
associative $3$-plane give two copies of each $P_i$ in $\mathcal P(P_\Sigma,P_\Sigma)$.
We substitute $h=E_1$ and $h=E_\Sigma$ into $A(h)=(h^2)_0$ and
$B(h)=\mathcal P(h,h)_0$, and then into $F_\lambda(h)=A(h)+\lambda B(h)$
and $N_\lambda(h)=\langle F_\lambda(h),h\rangle_W$. This gives the
following values.
\begin{equation}
\label{eq:line-plane-values}
\begin{array}{c|cccc}
h&A(h)&B(h)&F_\lambda(h)&N_\lambda(h)\\ \hline
E_1&5E_1&2E_1&(5+2\lambda)E_1&42(5+2\lambda)\\
E_\Sigma&E_\Sigma&20E_\Sigma&(1+20\lambda)E_\Sigma&84(1+20\lambda)
\end{array}
\end{equation}

We then evaluate \eqref{eq:double-adjoint} at $E_1$ and $E_\Sigma$, respectively. Since $F_\mu$ is quadratic, we obtain the following equations, where we do not yet cancel any parameter-dependent factors:
\begin{subequations}
\begin{align}
(5+2\lambda)^2(5+2\mu)&=42\kappa(5+2\lambda),
\label{eq:line-test-equation}\\
(1+20\lambda)^2(1+20\mu)&=84\kappa(1+20\lambda).
\label{eq:plane-test-equation}
\end{align}
\end{subequations}
For $\Delta=P_1-P_2$, the same identities give
\[
A(\Delta)=P_1+P_2-\tfrac27g,
\qquad B(\Delta)=-2P_3+\tfrac27g.
\]
Both are orthogonal to $\Delta$ under the $W$-pairing, so
$N_\lambda(\Delta)=0$. Thus \eqref{eq:double-adjoint} forces
$F_\mu(X)=0$, where
\[
X=F_\lambda(\Delta)=zg+R-2\lambda P_3,
\qquad R=P_1+P_2,\qquad z=\frac{2(\lambda-1)}7.
\]
Only the nonscalar terms of $X^2+\mu\mathcal P(X,X)$ matter, since
$F_\mu(X)$ is its traceless part. Orthogonality of the projection tensors,
$P_iP_j=\delta_{ij}P_i$ for the contracted product \eqref{eq:contracted-product}, gives
\[
(X^2)_0=\bigl((2z+1)R+(4\lambda^2-4\lambda z)P_3\bigr)_0.
\]
Also,
\[
\mathcal P(R,R)=2P_3,\qquad
\mathcal P(R,P_3)=R,\qquad
\mathcal P(P_3,P_3)=0.
\]
Bilinear expansion, using \eqref{eq:p-identity-rule} to discard scalar terms, now gives
\[
\mathcal P(X,X)_0
=\bigl((-2z-4\lambda)R+(2+4\lambda z)P_3\bigr)_0.
\]
Since the tensors $g,R,P_3$ are linearly independent,
$F_\mu(X)=0$ holds exactly when the coefficients of $R$ and $P_3$ vanish.
Rewriting $z$ in terms of $\lambda$, we obtain
\begin{subequations}
\begin{align}
3+4\lambda+(4-32\lambda)\mu&=0,
\label{eq:block-equation-one}\\
8\lambda+20\lambda^2+(14-8\lambda+8\lambda^2)\mu&=0.
\label{eq:block-equation-two}
\end{align}
\end{subequations}
These two equations say that the matrix
\[
\begin{pmatrix}
3+4\lambda&4-32\lambda\\
8\lambda+20\lambda^2&14-8\lambda+8\lambda^2
\end{pmatrix}
\]
annihilates the nonzero vector $(1,\mu)^{\transp}$. Its determinant must
therefore vanish, which gives the following equation
\[
(2\lambda+1)(8\lambda^2-2\lambda+1)=0.
\]
If $\lambda=-\tfrac12$, equation~\eqref{eq:block-equation-one} gives
$\mu=-\tfrac1{20}$. The line test \eqref{eq:line-test-equation} then forces
$\kappa=\tfrac7{15}$, whereas the $3$-plane test
\eqref{eq:plane-test-equation} forces $\kappa=0$, giving a contradiction.

Hence $8\lambda^2-2\lambda+1=0$. In particular,
$4-32\lambda\ne0$, and thus solving
\eqref{eq:block-equation-one} and using the quadratic relation yields
\[
\mu=\frac14-\lambda.
\]
This immediately implies that $\lambda$ and $\mu$ are the two distinct roots $\lambda_\pm$.
Finally, $5+2\lambda\ne0$, so we may now
solve \eqref{eq:line-test-equation} for $\kappa$. Using $\lambda+\mu=\frac14$ and $\lambda\mu=\frac18$ from the quadratic equation, we obtain
\[
\kappa=\frac{(5+2\lambda)(5+2\mu)}{42}
=\frac{25+10(\lambda+\mu)+4\lambda\mu}{42}=\frac23.
\]

\emph{Sufficiency.}
Conversely, let $\lambda,\mu$ be the two roots of
$8t^2-2t+1$, so $\lambda+\mu=1/4$ and $\lambda\mu=1/8$.
Bilinearity of the halved polarizations gives
\begin{equation*}
\begin{aligned}
F_\mu(F_\lambda(h))={}&A(A(h))+2\lambda A(A(h),B(h))+\lambda^2A(B(h))\\
 &+\mu B(A(h))+2\lambda\mu\,B(A(h),B(h))
   +\lambda^2\mu B(B(h)).
\end{aligned}
\end{equation*}
Use $\mu=1/4-\lambda$, $\lambda^2=\lambda/4-1/8$ and
$\lambda^2\mu=\lambda/8$ to collect this as
\begin{equation}
\begin{aligned}
8F_\mu(F_\lambda(h))={}&
 \bigl[8A(A(h))-A(B(h))+2B(A(h))
             +2B(A(h),B(h))\bigr]\\
 &+\lambda\bigl[16A(A(h),B(h))+2A(B(h))
                        -8B(A(h))+B(B(h))\bigr].
\end{aligned}
\label{eq:critical-composition-expansion}
\end{equation}
Lemma~\ref{lem:universal-quartic-identities} identifies the brackets as
$16C_0(h)h/3$ and $16C_\varphi(h)h/3$, respectively. Dividing by $8$,
we obtain $F_\mu(F_\lambda(h))=\tfrac23N_\lambda(h)h$ on all of $\WC$.
\end{proof}

\begin{lemma}
\label{lem:critical-adjoint-dominance}
The image of $F_\pm$ contains $\{u\in\WC:N_\mp(u)\ne0\}$.
In particular, $F_\pm$ is dominant, and a polynomial identity holding on
its image holds on all of $\WC$.
\end{lemma}

\begin{proof}
Proposition~\ref{prop:double-adjoint} gives
$F_\pm(F_\mp(u))=\tfrac23N_\mp(u)u$.
If $N_\mp(u)\ne0$, choose $t\in\C$ with
$t^2\tfrac23N_\mp(u)=1$. By homogeneity,
$F_\pm(tF_\mp(u))=u$.
The nonvanishing set of the nonzero cubic $N_\mp$ is Zariski-dense,
which proves the remaining assertions.
\end{proof}

The bilinear expansion of $F_\mu(F_\lambda(h))$ in the proof of
Proposition~\ref{prop:double-adjoint} gives an identity for unrestricted
$\lambda,\mu\in\C$. Put
\begin{equation*}
\begin{aligned}
\mathcal R^{(1)}_\lambda(h)&=B(A(h))+\lambda A(B(h)),\\
\mathcal R^{(2)}_\lambda(h)&=2B(A(h),B(h))-A(B(h))
                                 +\lambda B(B(h)).
\end{aligned}
\end{equation*}
Then, as a polynomial identity on $\WC$,
\begin{equation}
\begin{aligned}
F_\mu(F_\lambda(h))={}&\tfrac23N_\lambda(h)h\\
 &+\bigl(\lambda+\mu-\tfrac14\bigr)\mathcal R^{(1)}_\lambda(h)
  +\bigl(\lambda\mu-\tfrac18\bigr)\mathcal R^{(2)}_\lambda(h).
\end{aligned}
\label{eq:composition-formula}
\end{equation}
Indeed, subtracting the two remainder terms from that expansion
leaves one eighth of the left-hand side of \eqref{eq:quartic-identity-one}
plus $\lambda/8$ times that of \eqref{eq:quartic-identity-two}. For example, the coefficient of $A(B(h))$ becomes
$\lambda^2-\lambda(\lambda+\mu-\tfrac14)+(\lambda\mu-\tfrac18)
=\lambda/4-1/8$. Thus \eqref{eq:composition-formula} specializes to the identity
$F_\mu(F_\lambda(h))=\tfrac23N_\lambda(h)h$ precisely when the two remainder
coefficients vanish.

\begin{remark}
\label{rem:double-adjoint}
In Proposition~\ref{prop:double-adjoint}, the parameters $\lambda_\pm$ are determined without computing a stabilizer.
The later stabilizer calculation shows that these same parameters are
exactly the two $78$-dimensional cases in
Theorem~\ref{thm:stabilizer-trichotomy}. Double-adjoint identities are classical in the theory of cubic Jordan
norms \cite{springer-1962-characterization}. Their relation to polynomial
multiplicative Legendre transforms and semisimple rank-three Jordan
algebras is developed in \cite[Theorem~3.10]{etingof-2002-fourier},
\cite[Theorem~2 and Corollary~4]{chaput-2013-homaloidal}, and
\cite[Corollary~5.9]{pirio-2011-extremal}, under the hypotheses stated there. McCrimmon gives a coordinate proof of $x^{\sharp\sharp}=N(x)x$ in the
Hermitian-matrix model \cite[pp.~502--503]{mccrimmon-1969-freudenthal}.
In the proof of Proposition~\ref{prop:double-adjoint}, we use neither Jordan theory nor representation theory. The parameters are determined directly from the $G_2$ model.
\end{remark}

The double-adjoint identity of Proposition~\ref{prop:double-adjoint} will
produce stabilizer elements through the Koecher operators of
Section~\ref{subsec:koecher-jump}. We first exclude the five isotypic
blocks that can never contribute, so that the resulting lower bound
can be compared with the contribution of the remaining blocks.

\subsection{Stabilizer exclusion}
\label{subsec:stabilizer-exclusion}

We will show that five isotypic blocks of $\End(\WC)$ never enter the infinitesimal
stabilizer $\ker\Phi_{a,b}$ for $(a,b)\ne0$. We prove this exclusion
here because it shortens the determination of the stabilizer: after the
scalar component is excluded for a nonzero cubic, only $V(1,0)$,
$V(0,1)$ and $V(1,1)$ remain as possible components in
Table~\ref{tab:endomorphism-inventory}.
The two copies of $V(2,0)$ are handled together by the Gram matrix of
four maps. On each of the four multiplicity-one blocks, a $2\times2$
matrix of the scalars by which the compositions $\Phi_s^\top\Phi_t$ of
\eqref{eq:back-contraction} act gives the exclusion.

With the cubic coordinates already fixed, write
$\Phi_s=\Phi_{C_s}$ for $s\in\{0,\varphi\}$. The action
\eqref{eq:action-map} has components
\begin{equation}
\Phi_s(X)_{\alpha\beta\gamma}
 =X^\sigma{}_\alpha c^s_{\sigma\beta\gamma}
  +X^\sigma{}_\beta c^s_{\alpha\sigma\gamma}
  +X^\sigma{}_\gamma c^s_{\alpha\beta\sigma}.
\label{eq:action-map-ortho}
\end{equation}
Use the full ordered pairing on cubic tensors and the trace-adjoint
pairing on endomorphisms,
\begin{equation}
\label{eq:endomorphism-pairing}
\langle U,U'\rangle_3=U_{\alpha\beta\gamma}U'^{\alpha\beta\gamma},
\qquad
\langle X,Y\rangle_{\End}=\tr_W(X^{\tadj}Y).
\end{equation}
The cubic pairing has no factor $1/3!$. The algebraic transpose
$\Phi_s^\top$ is characterized by
$\langle\Phi_s^\top U,X\rangle_{\End}
=\langle U,\Phi_s(X)\rangle_3$, and therefore
\begin{equation}
(\Phi_s^\top U)^\sigma{}_\alpha
 =3(c^s)^{\sigma\beta\gamma}U_{\alpha\beta\gamma}.
\label{eq:action-map-adjoint-ortho}
\end{equation}
Define the operators of multiplication by $h$ in the two symmetric
products $A$ and $B$,
\begin{equation}
m_h=A(h,\cdot),\qquad n_h=B(h,\cdot),
\qquad m_h^{(0)}=m_h,\quad m_h^{(\varphi)}=n_h.
\label{eq:mn-operators}
\end{equation}

\begin{lemma}
\label{lem:two-v20}
The operators $m_h^{(s)}$ are $W$-self-adjoint and satisfy
\begin{equation}
\langle m_{h_1}^{(s)},m_{h_2}^{(t)}\rangle_{\End}
 =\mathsf a_{st}\langle h_1,h_2\rangle_W.
\label{eq:m-trace-pairing}
\end{equation}
The map $(h_1,h_2)\mapsto m_{h_1}+n_{h_2}$ is a $G_2$-equivariant
isomorphism from $\WC\oplus\WC$ onto the full
$2V(2,0)$-isotypic component of $\End(\WC)$.
\end{lemma}
\begin{proof}
Self-adjointness follows from the full symmetry of $c^s$, and
\eqref{eq:m-trace-pairing} is exactly the two-index c-contraction identity
\eqref{eq:two-index-c-contraction}. If $m_{h_1}+n_{h_2}=0$, pair it with $m_k$ and
$n_k$ for arbitrary $k\in\WC$. The resulting system has coefficient
matrix $\mathsf a$, whose determinant is $44$, so the nondegenerate
$W$-pairing gives $h_1=h_2=0$. The map is \(G_2\)-equivariant, and by
Table~\ref{tab:endomorphism-inventory}, $\End(\WC)$ contains exactly two
copies of $W=V(2,0)$. The image of the map, of dimension $54$, is
therefore the whole isotypic component.
\end{proof}

To apply the c-contraction identities to endomorphisms, we define
$\mathcal Z_{st}$ by the first line below; the second line then holds:
\begin{equation}
\begin{aligned}
\mathcal Z_{st}(X)^\alpha{}_\beta
 &=(c^s)^{\alpha\gamma\delta}c^t_{\beta\sigma\delta}
      X^\sigma{}_\gamma,\\
\Phi_s^\top\Phi_t(X)&=3\mathsf a_{st}X+6\mathcal Z_{st}(X).
\end{aligned}
\label{eq:back-contraction}
\end{equation}
Indeed, substitute the three terms of \eqref{eq:action-map-ortho} into
\eqref{eq:action-map-adjoint-ortho}. The first contracts to
$3\mathsf a_{st}X$ by \eqref{eq:two-index-c-contraction}; the other two coincide
by cubic symmetry and give $6\mathcal Z_{st}(X)$. Equivalently,
\begin{equation*}
\mathcal Z_{st}(X)=
 \sum_\alpha m_{k_\alpha}^{(s)}X^{\tadj}m_{k^\alpha}^{(t)}.
\end{equation*}
Applying $\Phi_0^\top$ and $\Phi_\varphi^\top$ to
$\Phi_{a,b}(X)=a\Phi_0(X)+b\Phi_\varphi(X)=0$ shows that every
$X\in\ker\Phi_{a,b}$ satisfies
\begin{subequations}
\begin{align}
a\Phi_0^\top\Phi_0(X)+b\Phi_0^\top\Phi_\varphi(X)&=0,
 \label{eq:kerPhi-eq00}\\
a\Phi_\varphi^\top\Phi_0(X)+b\Phi_\varphi^\top\Phi_\varphi(X)&=0.
 \label{eq:kerPhi-eqpp}
\end{align}
\end{subequations}

For the $2V(2,0)$ block, consider the four maps
\[
\mathcal B_{st}\colon\WC\longrightarrow S^3\WdualC,
\qquad
\mathcal B_{st}(h)=\Phi_s(m_h^{(t)}).
\]
Enumerate the four maps by
\[
\mathcal B_1=\mathcal B_{00},\qquad
\mathcal B_2=\mathcal B_{\varphi0},\qquad
\mathcal B_3=\mathcal B_{0\varphi},\qquad
\mathcal B_4=\mathcal B_{\varphi\varphi}.
\]
Let $\mathcal B_q^\top\colon S^3\WdualC\to\WC$ denote the transpose of $\mathcal B_q$ with respect to
$\langle\cdot,\cdot\rangle_3$ and the $W$-pairing.
Write $\mathcal G$ for the $4\times4$ matrix of the scalars by which the
compositions $\mathcal B_p^\top\mathcal B_q$ act.

\begin{lemma}
\label{lem:four-map-gram}
The four maps satisfy
$\mathcal B_p^\top\mathcal B_q=\mathcal G_{pq}\Id_W$, where
\begin{equation}
\mathcal G=\frac1{196}
 \begin{pmatrix}
 4455&1782&1782&-9636\\
 1782&16236&-9636&27192\\
 1782&-9636&16236&27192\\
 -9636&27192&27192&456816
 \end{pmatrix},
\qquad
\det\mathcal G=\frac{179397504}{7}\ne0.
\label{eq:four-map-gram-matrix}
\end{equation}
\end{lemma}
The contraction proof is given in Section~\ref{subsec:stabilizer-contractions},
\hyperref[proof:four-map-gram]{the proof of Lemma~\ref{lem:four-map-gram}}.

On a multiplicity-one block, write $\Gamma_{st}$ for the scalar
by which $\Phi_s^\top\Phi_t$ acts, with $s,t\in\{0,\varphi\}$, and
let $\Gamma=(\Gamma_{st})$ in the order $0,\varphi$.
The following table records all four matrices.

\begin{lemma}
\label{lem:four-block-scalars}
On the four multiplicity-one blocks, the compositions satisfy
$\Phi_s^\top\Phi_t=\Gamma_{st}\Id$, with the following scalars and
determinants:
{\renewcommand{\arraystretch}{1.35}
\begin{center}
\begin{tabular}{lcccc}
  Component & $\Gamma_{00}$
    & $\Gamma_{0\varphi}=\Gamma_{\varphi0}$
    & $\Gamma_{\varphi\varphi}$ & $\det\Gamma$\\ \hline
  $V(4,0)$ & $\tfrac{225}{28}$ & $\tfrac{45}{14}$ & $\tfrac{513}{7}$
    & $\tfrac{4050}{7}$\\
  $V(0,2)$ & $\tfrac{99}{28}$ & $\tfrac{45}{14}$ & $\tfrac{765}{7}$
    & $\tfrac{5265}{14}$\\
  $V(3,0)$ & $\tfrac{27}{4}$ & $\tfrac{27}{2}$ & $63$ & $243$\\
  $V(2,1)$ & $\tfrac{27}{4}$ & $\tfrac32$ & $87$ & $585$
\end{tabular}
\end{center}
}

\end{lemma}
The contraction proof is given in Section~\ref{subsec:stabilizer-contractions},
\hyperref[proof:four-block-scalars]{the proof of Lemma~\ref{lem:four-block-scalars}}.

\begin{proposition}
\label{prop:five-blocks}
For every $(a,b)\in\C^2\setminus\{0\}$, the restriction of the action map $\Phi_{a,b}=a\Phi_{C_0}+b\Phi_{C_{\varphi}}$
to each of the following isotypic components is injective:
\[
\begin{array}{c|c}
\text{$W$-self-adjoint part }S^2\WC&2V(2,0),\quad V(0,2),\quad V(4,0)\\
\text{$W$-skew-adjoint part }\Lambda^2\WC&V(3,0),\quad V(2,1).
\end{array}
\]
Here injectivity on $2V(2,0)$ means injectivity on the sum of
both copies. None of these five types occurs in $\ker\Phi_{a,b}$.
\end{proposition}

\begin{proof}
By Lemma~\ref{lem:two-v20}, an element of the $2V(2,0)$ block has
the unique form $X=m_{h_1}+n_{h_2}$. If $\Phi_{a,b}(X)=0$, then
\[
\mathcal B_1(ah_1)+\mathcal B_2(bh_1)
+\mathcal B_3(ah_2)+\mathcal B_4(bh_2)=0.
\]
Applying the four transposes gives
$\mathcal G(ah_1,bh_1,ah_2,bh_2)^{\transp}=0$.
By Lemma~\ref{lem:four-map-gram}, $\mathcal G$ is invertible, so
$ah_1=bh_1=ah_2=bh_2=0$. Since $(a,b)\ne0$, we have $h_1=h_2=0$.

On any of the four multiplicity-one blocks, applying $\Phi_0^\top$
and $\Phi_\varphi^\top$ to $\Phi_{a,b}(X)=0$ gives
$\Gamma(aX,bX)^{\transp}=0$ by
\eqref{eq:kerPhi-eq00} and \eqref{eq:kerPhi-eqpp}.
Every determinant in Lemma~\ref{lem:four-block-scalars} is nonzero.
Thus $aX=bX=0$, and again $X=0$.
\end{proof}

Once the scalar endomorphisms are also excluded (see the proof of
Proposition~\ref{prop:block-kernels}), only the $\mathbf7$,
$\mathbf{14}$ and $\mathbf{64}$ components can occur in the stabilizer
of a nonzero cubic.
We next construct enough infinitesimal symmetries at the critical
parameters to force the $\mathbf{64}$ component to occur.

\subsection{Koecher operators and a conceptual lower bound}
\label{subsec:koecher-jump}

Differentiating the double-adjoint identity produces endomorphisms that
preserve the cubic. We call these Koecher operators because they are
related to the Tits--Kantor--Koecher construction in the theory of Jordan
algebras \cite{demedts-2026-cubic,garibaldi-2024-albert,springer-1973-jordan}.
Lemma~\ref{lem:koecher} applies to any pair of cubics with quadratic
adjoints. We first need a precise definition.

\begin{definition}
\label{def:koecher}
Let $\mathcal X$ be a finite-dimensional vector space over a field $\mathbb F$ of
characteristic other than $2$ and $3$, with a nondegenerate symmetric bilinear
form $\langle\cdot,\cdot\rangle$. Let $N$ and $N'$ be cubic forms on $\mathcal X$
with the polarization convention $N(h,h,h)=N(h)$ and
$N'(h,h,h)=N'(h)$. Define symmetric bilinear maps
$F,F'\colon\mathcal X\times\mathcal X\to\mathcal X$ by
\[
\langle F(h_1,h_2),h_3\rangle=N(h_1,h_2,h_3),\qquad
\langle F'(h_1,h_2),h_3\rangle=N'(h_1,h_2,h_3).
\]
Write $F(h)=F(h,h)$ and $F'(h)=F'(h,h)$, and define $\stab(N)$ as in
Definition~\ref{def:exact-stabilizer}. For
$0\ne\kappa\in\mathbb F$ and $h_1,h_2\in\mathcal X$ define the \emph{Koecher operator}
$S(h_1,h_2)\in\End(\mathcal X)$ by
\begin{equation}
\label{eq:def-koecher}
  S(h_1,h_2)h=\langle h,h_2\rangle h_1-\frac{4}{3\kappa}F'\bigl(F(h_1,h),h_2\bigr)
  +\frac13\langle h_1,h_2\rangle h.
\end{equation}
\end{definition}

\begin{lemma}
\label{lem:koecher}
With notation as in Definition~\ref{def:koecher}, the following are equivalent:
\begin{enumerate}
\renewcommand{\labelenumi}{(\roman{enumi})}
\renewcommand{\theenumi}{\roman{enumi}}
  \item\label{item:koecher-composition} $F'(F(h))=\kappa N(h)h$ for all $h\in\mathcal X$;
  \item\label{item:koecher-stabilizers} $S(h_1,h_2)\in\stab(N)$ for all $h_1,h_2\in\mathcal X$.
\end{enumerate}
\end{lemma}

For the pencil, Lemma~\ref{lem:koecher} gives the following corollary immediately.

\begin{corollary}
\label{cor:koecher}
Take $\mathcal X=\WC$ with the $W$-pairing, $N=N_\lambda$ and
$N'=N_\mu$, so that $F=F_\lambda$ and $F'=F_\mu$ by \eqref{eq:F-lambda-pairing}.
For $\lambda,\mu\in\C$ and $\kappa\ne0$, the identity \eqref{eq:double-adjoint}
holds at $(\lambda,\mu,\kappa)$ if and only if every Koecher operator
$S(h_1,h_2)=S_{\lambda\mu\kappa}(h_1,h_2)$ lies in $\stab(N_\lambda)$.
\end{corollary}

We now proceed with the proof of Lemma~\ref{lem:koecher}.

\begin{proof}[Proof of Lemma~\ref{lem:koecher}]
Regarding (\ref{item:koecher-composition}) as a quartic polynomial equation in $h$, we differentiate it in the direction $h_1$. Since
$F(h+th_1)=F(h)+2tF(h,h_1)+t^2F(h_1)$, the coefficient of $t$ in $F'(F(h+th_1))$ is
$2F'(F(h),2F(h,h_1))$, while that of $N(h+th_1)(h+th_1)$ is
$N(h)h_1+3N(h,h,h_1)h=N(h)h_1+3\langle F(h),h_1\rangle h$. Hence (\ref{item:koecher-composition}) gives
\begin{equation}
  4F'\bigl(F(h),F(h,h_1)\bigr)=\kappa\bigl(N(h)h_1+3\langle F(h),h_1\rangle h\bigr).
  \label{eq:koecher-polarized}
\end{equation}
Next pair $S(h_1,h_2)h$ with $F(h)$. By the full symmetry of
$\langle F'(\cdot,\cdot),\cdot\rangle=N'$,
\[
\langle F(h),F'(F(h_1,h),h_2)\rangle
=\langle F'(F(h),F(h,h_1)),h_2\rangle.
\]
Also $N(h_1,h,h)=\langle F(h),h_1\rangle$, so
\begin{equation}
\begin{aligned}
  N\bigl(S(h_1,h_2)h,h,h\bigr)&=\langle F(h),S(h_1,h_2)h\rangle\\
  &=\langle h,h_2\rangle\langle F(h),h_1\rangle
   -\frac{4}{3\kappa}\bigl\langle F'\bigl(F(h),F(h,h_1)\bigr),h_2\bigr\rangle
   +\frac13\langle h_1,h_2\rangle N(h),
\end{aligned}
\label{eq:koecher-pairing}
\end{equation}
and (\ref{item:koecher-stabilizers}) is the vanishing of \eqref{eq:koecher-pairing} for all $h,h_1,h_2$. If (\ref{item:koecher-composition})
holds, then by \eqref{eq:koecher-polarized} the middle term equals
$-\tfrac13\bigl(N(h)\langle h_1,h_2\rangle+3\langle F(h),h_1\rangle\langle h,h_2\rangle\bigr)$,
and \eqref{eq:koecher-pairing} vanishes identically. Conversely, if
\eqref{eq:koecher-pairing} vanishes for all $h_2$, nondegeneracy gives
\eqref{eq:koecher-polarized}, and setting $h_1=h$ in it gives
$4F'(F(h))=4\kappa N(h)h$.
\end{proof}

Each $O\in G_2$ acts on $\WC$ by $\rho(O)$, preserving the
$W$-pairing and the maps $F_\lambda$
and $F_\mu$. Consequently,
\[
S(\rho(O)h_1,\rho(O)h_2)=\rho(O)S(h_1,h_2)\rho(O)^{-1},
\]
and the span of the Koecher operators is a $G_2$-submodule of $\End(\WC)$. At $\kappa=\tfrac23$ put
\begin{gather*}
  D(h_1,h_2)=S(h_1,h_2)+\tfrac23\langle h_1,h_2\rangle_W\Id_W,\\
  D(h_1,h_2)h=\langle h_1,h_2\rangle_Wh+\langle h,h_2\rangle_Wh_1-2F_\mu(F_\lambda(h_1,h),h_2).
\end{gather*}
When the double-adjoint identity holds, so that Lemma~\ref{lem:koecher}(\ref{item:koecher-stabilizers}) holds as well, equation \eqref{eq:koecher-pairing} gives
$3\langle F_\lambda(h),D(h_1,h_2)h\rangle_W=2\langle h_1,h_2\rangle_WN_\lambda(h)$, that
is, $\Phi_{N_\lambda}(D(h_1,h_2))=2\langle h_1,h_2\rangle_WN_\lambda$: each $D(h_1,h_2)$ is an
infinitesimal similarity of $N_\lambda$ with multiplier $2\langle h_1,h_2\rangle_W$.

\begin{remark}
\label{rem:koecher-literature}
It is classical that, in the Jordan algebra setting, the operators
$D(h_1,h_2)$ act as infinitesimal similarities of the norm with
multiplier $2\langle h_1,h_2\rangle_W$
\cite[\S\S4.3--4.4 and 14.22]{springer-1973-jordan},
\cite[(33a.30) and Proposition~50.5]{garibaldi-2024-albert},
and De Medts--Wiedemann also establish this property in
\cite[Lemma~3.3]{demedts-2026-cubic} (arXiv v1, February 2026).
Lemma~\ref{lem:koecher} also establishes the converse
(\ref{item:koecher-stabilizers})$\Rightarrow$(\ref{item:koecher-composition}) and requires only the one-sided composition
$F_\mu\circ F_\lambda$. We have not found these two features
explicitly stated in the literature.
\end{remark}

We now use the Koecher operators to obtain a lower bound for the
dimension of the critical stabilizers.
The argument has two steps. First, evaluation
at a suitable element $E_1$ shows that the Koecher operators span
a space of dimension at least $26$, containing symmetries outside of $\rho(\lieso(7,\C))$. We then use Proposition~\ref{prop:five-blocks}
and the $G_2$-module decomposition of $\End(\WC)$ to show that
their span must contain a copy of $V(1,1)$, giving the lower bound
$14+64=78$.

Let $E_1=7P_1-g$, as in \eqref{eq:test-tensors}, be the test
tensor in $W$ used for the line test in the proof of
Proposition~\ref{prop:double-adjoint}, so that
\[
F_\lambda(E_1)=(5+2\lambda)E_1,
\qquad
\langle E_1,E_1\rangle_W=42.
\]
For a real unit vector $u\in V$, define $\Theta_u\in\mathfrak{so}(V)$ by
$\Theta_u(v)=v\times u$, and put $\Theta_1=\Theta_{e_1}$, so that
\[
\Theta_1(v)=v\times e_1,
\qquad
(\Theta_1)_{ia}=\varphi_{ia1}.
\]
In the standard frame of \eqref{eq:standard-positive-form}, the
terms of $\varphi$ containing the index $1$ are
$e^{123}$, $e^{145}$ and $e^{167}$, giving
\[
\Theta_1=0\oplus J_0\oplus J_0\oplus J_0,
\qquad
J_0=\begin{pmatrix}0&1\\-1&0\end{pmatrix},
\]
on the blocks $(e_1)$, $(e_2,e_3)$, $(e_4,e_5)$ and $(e_6,e_7)$.
In particular, $\Theta_1e_1=0$ and
$\Theta_1^2=-\pr_{e_1^\perp}$, so its restriction to $e_1^\perp$
is a complex structure. The choice of $e_1$ reduces $G_2$ to its stabilizer
$SU(3)$. The complex structure $\Theta_1$ on
$e_1^\perp$ identifies this space with $\C^3$.
We use the corresponding standard $SU(3)$-module
decomposition
\[
\left.\WC\right|_{SU(3)}
\cong
\mathbf1\oplus
(\mathbf3\oplus\overline{\mathbf3})
\oplus\mathbf8\oplus
(\mathbf6\oplus\overline{\mathbf6}),
\]
and thus define, for a real unit vector $u\in V$, the following real
subspaces of $W$, which we call the $SU(3)$ \emph{sectors}:
\begin{gather*}
W_6(u)
=\{u^\flat\otimes v^\flat+v^\flat\otimes u^\flat:v\perp u\},\\
W_8(u)
=\{h\in W:h(u,\cdot)=0,\ \Theta_u h=h \Theta_u\},\\
W_{12}(u)
=\{h\in W:h(u,\cdot)=0,\ \Theta_u h=-h \Theta_u\},
\end{gather*}
where indices are contracted with the metric, so that \((\Theta_u h)_{ab}=(\Theta_u)_{ai}h_{jb}g^{ij}\).
The subscript $\C$ denotes complexification:
$W_r(u)_\C=W_r(u)\otimes\C$ for $r=6,8,12$.
At $u=e_1$ these give the decomposition
\[
\WC=\C E_1\oplus W_6(e_1)_\C\oplus
W_8(e_1)_\C\oplus W_{12}(e_1)_\C,
\qquad
E_1^\perp=W_6(e_1)_\C\oplus W_8(e_1)_\C
\oplus W_{12}(e_1)_\C,
\]
where $E_1^\perp$ is the $W$-orthogonal complement of $E_1$.
The subscripts $6,8,12$ indicate the dimensions of the summands. The components $W_6(e_1)_\C$ and $W_{12}(e_1)_\C$ are sums of the conjugate irreducible $SU(3)$-modules $\mathbf3\oplus\overline{\mathbf3}$ and $\mathbf6\oplus\overline{\mathbf6}$, respectively.

By Proposition~\ref{prop:double-adjoint}, the double-adjoint
identity \eqref{eq:double-adjoint} holds precisely for
\[
(\lambda,\mu,\kappa)
=(\lambda_+,\lambda_-,\tfrac23)
\quad\text{or}\quad
(\lambda_-,\lambda_+,\tfrac23).
\]
For either choice, let $S=S_{\lambda\mu\kappa}$ denote
the Koecher operators of Corollary~\ref{cor:koecher}.

By Definition~\ref{def:exact-stabilizer} and \eqref{eq:F-adjoint},
every $X\in\stab(N_\lambda)$ satisfies
\[
0=N_\lambda(Xh,h,h)=\langle F_\lambda(h),Xh\rangle_W
\qquad(h\in\WC).
\]
Taking $h=E_1$ and using $F_\lambda(E_1)=(5+2\lambda)E_1$ shows that
$XE_1\in E_1^\perp$ whenever $5+2\lambda\ne0$, in particular at either
critical parameter. We therefore have the evaluation map
\begin{equation}
\operatorname{ev}_{E_1}\colon\stab(N_\lambda)\longrightarrow E_1^\perp,
\qquad X\longmapsto XE_1.
\label{eq:stabilizer-evaluation}
\end{equation}

\begin{proposition}
\label{prop:koecher-jump}
Let $\lambda\in\{\lambda_+,\lambda_-\}$, let $\mu$ be the other critical
parameter, and let $S=S_{\lambda\mu\kappa}$ with $\kappa=\tfrac23$. For the
critical cubic $N_\lambda$, the evaluation map
$\operatorname{ev}_{E_1}$ in \eqref{eq:stabilizer-evaluation} is surjective.
Moreover, for every
$0\ne h\in W_8(e_1)_\C\oplus W_{12}(e_1)_\C$,
the operator $S(h,E_1)$ belongs to $\stab(N_\lambda)$
but not to $\rho(\lieso(7,\C))$.

The span of the Koecher operators contains a
$G_2$-submodule $\tau\cong V(1,1)$, and
\[
\rho(\gtwoC)\oplus\tau\subseteq\stab(N_\lambda),
\qquad
\dim\stab(N_\lambda)\ge14+64=78.
\]
\end{proposition}

Equality is established in
Section~\ref{subsec:nothing-else-jumps}.

\begin{proof}
By Corollary~\ref{cor:koecher}, the Koecher operators
$S(h_1,h_2)$ of each critical cubic $N_\lambda$ belong to
$\stab(N_\lambda)$. Since the evaluation map
$\operatorname{ev}_{E_1}$ takes values in $E_1^\perp$ by
\eqref{eq:stabilizer-evaluation}, it suffices to
show that the elements $S(h,E_1)E_1$, for $h\in\WC$,
span $E_1^\perp$.

We begin by computing the action of $F_\lambda(\,\cdot\,,E_1)$
on the four summands of the $SU(3)$-decomposition
of $\WC$. These are preserved by
$A(\,\cdot\,,E_1)$ and $B(\,\cdot\,,E_1)$,
which act by the scalars listed in Table~\ref{tab:koecher-sector-response}.
\begin{table}[htbp]
\centering
\caption{Eigenvalues of $A(\,\cdot\,,E_1)$, $B(\,\cdot\,,E_1)$ and
$F_\lambda(\,\cdot\,,E_1)$ on the $SU(3)$-summands of $\WC$, at
$E_1=7P_1-g$.}
\label{tab:koecher-sector-response}
\begin{tabular}{lccc}
  Sector & $A(h,E_1)$ & $B(h,E_1)$
        & $F_\lambda(h,E_1)$\\ \hline
  $\C E_1$          & $5h$          & $2h$
                 & $(5+2\lambda)h$\\
  $W_6(e_1)_\C$  & $\tfrac52h$    & $h$
                 & $(\tfrac52+\lambda)h$\\
  $W_8(e_1)_\C$  & $-h$          & $8h$
                 & $(-1+8\lambda)h$\\
  $W_{12}(e_1)_\C$ & $-h$        & $-6h$
                 & $(-1-6\lambda)h$
\end{tabular}
\end{table}
The first row follows from \eqref{eq:line-plane-values}. For
$h\in E_1^\perp$, we have $h_{11}=0$, and
\[
A(h,E_1)=\tfrac72(hP_1+P_1h)-h,
\qquad
B(h,E_1)
=7\mathcal P(h,P_1)-\mathcal P(h,g).
\]
The remaining entries follow from $\mathcal P(h,g)=-h$, by
\eqref{eq:p-identity-rule}, and from
$\mathcal P(h,P_1)=h(\Theta_1\cdot,\Theta_1\cdot)$, which in the
standard orthonormal frame is the matrix $-\Theta_1h\Theta_1$.

On each sector, we thus have $F_\lambda(h,E_1)=f(\lambda)h$, where $f(\lambda)$ is the entry of the last column for that sector. Applying $F_\mu$ gives $F_\mu(F_\lambda(h,E_1),E_1)=f(\mu)f(\lambda)h$. At a triple satisfying \eqref{eq:double-adjoint},
we have $\lambda+\mu=\tfrac14$ and
$\lambda\mu=\tfrac18$. Direct computation of the eigenvalues of $F_\mu(\,\cdot\,,E_1)\circ F_\lambda(\,\cdot\,,E_1)$ at these values therefore gives $7$ on each of
$W_6(e_1)_\C$, $W_8(e_1)_\C$ and $W_{12}(e_1)_\C$,
and $28$ on $\C E_1$. Consequently,
\[
F_\mu(F_\lambda(h,E_1),E_1)
=
\begin{cases}
7h, & h\in E_1^\perp,\\
28h, & h\in\C E_1.
\end{cases}
\]

Substituting this into the definition of the Koecher
operators \eqref{eq:def-koecher}, with $\kappa=\tfrac23$ and
$\langle E_1,E_1\rangle_W=42$, gives for $h\in E_1^\perp$
\[
S(h,E_1)E_1=42h-2(7h)=28h
\]
while
\[
S(E_1,E_1)E_1=(42-2\cdot28+14)E_1=0.
\]
By Corollary~\ref{cor:koecher}, every $S(h,E_1)$
belongs to $\stab(N_\lambda)$. Hence the evaluation map has image exactly $E_1^\perp$,
of dimension $26$.

The same computation also yields symmetries outside
$\rho(\lieso(7,\C))$. Indeed, for
$M\in\lieso(7,\C)$,
\[
\rho(M)E_1
=7\bigl((Me_1)^\flat\otimes e_1^\flat+e_1^\flat\otimes(Me_1)^\flat\bigr)
\in W_6(e_1)_\C.
\]
On the other hand, $S(h,E_1)E_1=28h$ for
$h\in W_8(e_1)_\C\oplus W_{12}(e_1)_\C$.
Thus $S(h,E_1)\notin\rho(\lieso(7,\C))$
whenever $h\ne0$.

\emph{The $G_2$-module count.}
Let $\mathcal K$ be the linear span of all Koecher
operators $S(h_1,h_2)$. The map $(h_1,h_2)\mapsto S(h_1,h_2)$ is
$G_2$-equivariant, so $\mathcal K$ is a
$G_2$-submodule of
$\stab(N_\lambda)=\ker\Phi_{1,\lambda}$.
Since evaluation at $E_1$ maps $\mathcal K$
onto $E_1^\perp$, we have
\[
\dim\mathcal K\ge26.
\]

We now use the decomposition of $\End(\WC)$.
By Proposition~\ref{prop:five-blocks}, none of the
five isotypic types
\[
V(2,0),\quad V(0,2),\quad V(3,0),\quad
V(4,0),\quad V(2,1)
\]
can occur in the stabilizer.
Of the remaining types, the three smaller ones are
$V(0,0)$, $V(1,0)$ and $V(0,1)$. Each occurs with multiplicity one in
Table~\ref{tab:endomorphism-inventory}, so their isotypic components have
total dimension
\[
1+7+14=22<26.
\]
They cannot account for all of $\mathcal K$.
The only remaining isotypic type is $V(1,1)$,
of dimension $64$. By complete reducibility,
$\mathcal K$ therefore contains an irreducible
submodule $\tau\cong V(1,1)$.

Finally, $\rho(\gtwoC)\cong V(0,1)$ is contained
in every stabilizer by Lemma~\ref{lem:common-g2}.
Since $V(0,1)$ and $V(1,1)$ are inequivalent,
$\rho(\gtwoC)\cap\tau=0$. Hence
\[
\rho(\gtwoC)\oplus\tau
\subseteq\stab(N_\lambda)
\qquad \text{and }
\dim\stab(N_\lambda)\ge14+64=78.
\]

\end{proof}

\begin{remark}
\label{rem:structure-algebra}
For the Albert norm, the Lie algebra of infinitesimal similarities
has dimension $78$ modulo scalar operators
\cite[14.24--14.25]{springer-1973-jordan},
\cite[\S1.8]{petersson-2019-survey}.
The surjectivity of evaluation at $E_1$ also has a classical counterpart
at the unit of the Albert algebra
\cite[proof of~5.5 and 14.12]{springer-1973-jordan}.

Although the dimension bound $78$ is classical for the Albert norm,
we obtain the lower bound $78$ directly from the double-adjoint
identity, without first identifying the critical cubic with the
Albert norm and without any classification result.
The $G_2$-isotypic count also identifies a
$64$-dimensional constituent of type $V(1,1)$
alongside $\rho(\gtwoC)$.
\end{remark}

\subsection{Generator stabilizers and a model of the \texorpdfstring{$\mathbf{64}$}{64}}
\label{subsec:explicit-64}

The exact count needs two further ingredients: the stabilizers of
$C_0$ and $C_\varphi$, and a concrete realization of the two copies of
$V(1,1)$ in $\End(\WC)$. Although \(W=S^2_0V^*\) consists of symmetric tensors on \(V\), a general endomorphism \(L\in\operatorname{End}(W)\) need not arise from an endomorphism of \(V\). In particular, when determining the full infinitesimal stabilizer of a cubic on \(W\), we cannot initially assume that \(L=\rho(M)\) for some \(M\in \End V\). Lemma \ref{lem:rank-one-recognition} gives a criterion for this reduction. In Propositions \ref{prop:trace-cubic-stabilizer} and \ref{prop:phi-cubic-stabilizer}, we first show that the stabilizer condition implies this criterion, and only then identify the corresponding skew-adjoint endomorphism of \(V\).

Let $\mathbb F\in\{\R,\C\}$. As in \eqref{eq:test-tensors}, for $u\in V$ let
\[
E_u=7P_u-g(u,u)\,g\in W.
\]
For $L\in\End(W_{\mathbb F})$, let
$l(u)=g^{-1}(LE_u)\in\End(V_{\mathbb F})$ be the endomorphism
corresponding to $LE_u$ under the metric identification \eqref{eq:metric-identification}.
Thus
\[
g(l(u)v,w)=(LE_u)(v,w),
\qquad v,w\in V_{\mathbb F}.
\]
For $u,v\in V$ we define the test form $\beta$ by
\[
\beta(u,v)
=\tfrac17\langle LE_u,E_v\rangle_W
=(LE_u)(v,v)
=g(l(u)v,v),
\]
where $\tfrac17\langle LE_u,E_v\rangle_W=(LE_u)(v,v)$ because $LE_u$ is symmetric and traceless.
By Lemma~\ref{lem:rank-one-recognition}, $\beta$ vanishes on orthogonal
pairs exactly for $L\in\rho(\lieso(7,\mathbb F))$.

The following lemma characterizes the endomorphisms of $W_{\mathbb F}$
that arise from the action of $\lieso(7,\mathbb F)$. Recall from
\eqref{eq:rho-action} that this Lie algebra acts on $W_{\mathbb F}$
by $\rho$ and that the operators $\rho(M)$ are $W$-skew-adjoint.
Although stated for dimension $7$, the lemma and its
proof apply to $S^2_0(\mathbb F^n)$ for every $n\ge2$.
\begin{lemma}[Rank-one recognition]
\label{lem:rank-one-recognition}
\leavevmode
\begin{enumerate}
\renewcommand{\labelenumi}{(\roman{enumi})}
\renewcommand{\theenumi}{\roman{enumi}}
\item\label{item:rank-one-recognition-span} The elements $E_u$, for $u\in V$, span $W_{\mathbb F}$.
\item\label{item:rank-one-recognition-criterion} An endomorphism $L\in\End(W_{\mathbb F})$ lies in
$\rho(\lieso(7,\mathbb F))$ if and only if
$\beta(u,v)=0$ for all $u,v\in V$ with $g(u,v)=0$.
In that case, $L=\rho(M)$ for a unique
$M\in\lieso(7,\mathbb F)$.
\end{enumerate}
The test vectors $u,v$ are real even when $\mathbb F=\C$.
\end{lemma}

\begin{proof}
For (\ref{item:rank-one-recognition-span}), the identities
\[
\begin{aligned}
E_i-E_j
&=7\bigl(e_i^\flat\otimes e_i^\flat-e_j^\flat\otimes e_j^\flat\bigr),\\
E_{e_i+e_j}-E_i-E_j
&=7\bigl(e_i^\flat\otimes e_j^\flat+e_j^\flat\otimes e_i^\flat\bigr)
\end{aligned}
\]
for $i\ne j$ give a spanning set of $W_{\mathbb F}$.

For (\ref{item:rank-one-recognition-criterion}), if $L=\rho(M)$ with
$M\in\lieso(7,\mathbb F)$, then
\[
LE_u=7\bigl(u^\flat\otimes(Mu)^\flat+(Mu)^\flat\otimes u^\flat\bigr),
\qquad
\beta(u,v)=14g(u,v)g(Mu,v),
\]
so $\beta$ vanishes on orthogonal pairs.

Conversely, suppose that $\beta(u,v)=0$ whenever $g(u,v)=0$. For fixed
real $u\ne0$, by polarization, $LE_u$ vanishes as a bilinear form on
$u^\perp$.
Since $LE_u$ is traceless, there is a unique $w(u)\in V_{\mathbb F}$,
orthogonal to $u$, such that
\[
LE_u=7\bigl(u^\flat\otimes w(u)^\flat+w(u)^\flat\otimes u^\flat\bigr).
\]
Each diagonal entry $d_i(u)=(LE_u)(e_i,e_i)$ is a homogeneous quadratic
polynomial vanishing when $u_i=0$, since then $e_i\in u^\perp$. Dividing by $u_i$ therefore gives a
linear polynomial, so there is a matrix $M$, independent of $u$, with
$d_i(u)=14u_i(Mu)_i$. Comparing with $d_i(u)=14u_iw_i(u)$ gives
$w(u)=Mu$ whenever all coordinates of $u$ are nonzero. Both sides of
\[
LE_u=7\bigl(u^\flat\otimes(Mu)^\flat+(Mu)^\flat\otimes u^\flat\bigr)
\]
are polynomial in $u$, so the identity holds for every real $u$.
Taking traces gives $14g(Mu,u)=0$ for all real $u$, hence $M$ is
$g$-skew-adjoint. Thus $LE_u=\rho(M)E_u$; since the $E_u$ span
$W_{\mathbb F}$ by (\ref{item:rank-one-recognition-span}), $L=\rho(M)$,
and the injectivity of $\rho$ gives uniqueness. The polynomial argument
applies over either coefficient field with the test vectors kept real.
\end{proof}

The stabilizer of the trace cubic is already known. Garibaldi and
Guralnick showed that, for $n=3$ and $n\ge5$, the identity component of
the stabilizer of $\tr(h^3)$ on $S^2_0(\C^n)$ is the image of
$SO(n,\C)$ \cite[proof of Theorem~8.1]{garibaldi-2015-simple}.
Their proof uses Seitz's classification. We give a direct proof of the
corresponding infinitesimal statement for $n=7$.

Both generator stabilizers can be found by the same rank-one test.
For $L$ annihilating the cubic, evaluation on $(E_u,E_v,E_v)$ first
forces $\beta(v,v)=0$. Evaluations on orthonormal pairs, supplemented
for $C_\varphi$ by their cross products, then force $\beta(u,v)=0$
whenever $u\perp v$. Lemma~\ref{lem:rank-one-recognition} gives
$L=\rho(M)$ for a skew endomorphism $M$ of $V$.
Every such $M$ preserves $C_0$; for $C_\varphi$, the evaluation formula of
Lemma~\ref{lem:rank-one-phi-recovery} will identify those that also preserve $\varphi$.
In Propositions~\ref{prop:trace-cubic-stabilizer} and~\ref{prop:phi-cubic-stabilizer} below, the line stabilizer is obtained from the kernel by
adding multiples of $\Id_W$: for any nonzero cubic $N$,
\begin{equation}
\Phi_N(\Id_W)=3N,
\qquad
\{X:\Phi_N(X)\in\mathbb FN\}
=\ker\Phi_N\oplus\mathbb F\,\Id_W.
\label{eq:generator-line-stabilizer}
\end{equation}
Indeed, $\Phi_N(X)=cN$ if and only if
$X-\tfrac c3\Id_W\in\ker\Phi_N$.

\begin{proposition}
\label{prop:trace-cubic-stabilizer}
For $\mathbb F\in\{\R,\C\}$, $\ker\Phi_{C_0}=\rho(\lieso(7,\mathbb F))$,
of dimension $21$ over $\mathbb F$, and the stabilizer
$\{X\in\End(W_{\mathbb F}):\Phi_{C_0}(X)\in\mathbb FC_0\}$ of the line of
$C_0$ has dimension $22$.
\end{proposition}

\begin{proof}
For $M\in\lieso(7,\mathbb F)$, the action is the matrix commutator
$\rho(M)h=[M,h]$ by \eqref{eq:rho-action}. Cyclicity of the trace gives
\[
\tfrac13\Phi_{C_0}(\rho(M))(h,h,h)
=\tr([M,h]h^2)=\tr(Mh^3)-\tr(hMh^2)=0.
\]
Polarization therefore gives
$\rho(\lieso(7,\mathbb F))\subseteq\ker\Phi_{C_0}$.

Conversely, let $L\in\ker\Phi_{C_0}$. Using the contracted product
\eqref{eq:contracted-product}, we have
\[
C_0(h_1,h_2,h_3)
=\bigl\langle h_1,\tfrac12(h_2h_3+h_3h_2)_0\bigr\rangle_W,
\]
and the definition of $E_u$ gives
\[
\begin{aligned}
E_v^2&=35g(v,v)P_v+g(v,v)^2g,\\
\tfrac12(E_uE_v+E_vE_u)
&=\tfrac{49}2g(u,v)\bigl(u^\flat\otimes v^\flat+v^\flat\otimes u^\flat\bigr)
-7g(v,v)P_u\\
&\quad-7g(u,u)P_v+g(u,u)g(v,v)g.
\end{aligned}
\]
Contracting with the traceless tensors $LE_u$ and $LE_v$ yields
\begin{equation}
\begin{aligned}
0&=\Phi_{C_0}(L)(E_u,E_v,E_v)\\
&=C_0(LE_u,E_v,E_v)+2C_0(LE_v,E_u,E_v)\\
&=35g(v,v)\beta(u,v)+98g(u,v)(LE_v)(u,v)\\
&\quad-14g(v,v)\beta(v,u)-14g(u,u)\beta(v,v).
\end{aligned}
\label{eq:trace-cubic-evaluation}
\end{equation}
At $u=v$, this becomes $105g(v,v)\beta(v,v)=0$.
Since $g$ is positive definite, $\beta(v,v)=0$ for every real $v$.
For orthonormal $u,v$, \eqref{eq:trace-cubic-evaluation} and the same
identity with $u$ and $v$ exchanged give
\[
5\beta(u,v)=2\beta(v,u),
\qquad
5\beta(v,u)=2\beta(u,v),
\]
so both values vanish. By homogeneity, $\beta$ vanishes on all real
orthogonal pairs. Lemma~\ref{lem:rank-one-recognition} now gives
$L\in\rho(\lieso(7,\mathbb F))$, proving the equality.
The kernel has dimension $21$ by the injectivity of $\rho$, and by
\eqref{eq:generator-line-stabilizer}, the stabilizer of the line has
dimension $22$.
\end{proof}

By \eqref{eq:two-form-splitting} and the equivariance and injectivity of $\rho$,
$\ker\Phi_{C_0}=\rho(\Lambda^2_7)\oplus\rho(\Lambda^2_{14})\cong\mathbf7\oplus\mathbf{14}$
as a $G_2$-module.

The $\varphi$-cubic has an equally explicit value on the same test tensors $E_u$.

\begin{lemma}
\label{lem:rank-one-phi-recovery}
Let $\varphi$ be the standard $G_2$ $3$-form, or more generally
any positive $3$-form on $V$ compatible with $g$, obtained by
transporting the standard form by an isometry of either orientation.
Let $C_\varphi$ be defined as in \eqref{eq:natural-cubics}, and
let $E_u=7P_u-g(u,u)\,g$ be as in \eqref{eq:test-tensors}.
\begin{enumerate}
\renewcommand{\labelenumi}{(\roman{enumi})}
\renewcommand{\theenumi}{\roman{enumi}}
\item\label{item:rank-one-phi-recovery-evaluation} For all $u,v,w\in V$,
\begin{equation}
\begin{aligned}
C_\varphi(E_u,E_v,E_w)
={}&343\varphi(u,v,w)^2\\
&+49\bigl(
g(u,u) g(v,w)^2
+g(v,v) g(u,w)^2
+g(w,w) g(u,v)^2\bigr)\\
&-63g(u,u)g(v,v)g(w,w).
\end{aligned}
\label{eq:rank-one-phi-recovery}
\end{equation}

\item\label{item:rank-one-phi-recovery-form} The cubic $C_\varphi$ determines $\varphi$ up to sign.
More precisely, if $\varphi'$ is another such form with
$C_{\varphi'}=C_\varphi$, then $\varphi'=\pm\varphi$.
If the two forms induce the same orientation, then
$\varphi'=\varphi$.

\item\label{item:rank-one-phi-recovery-stabilizer} The subgroup of $O(7)$ whose induced action on $W$
preserves $C_\varphi$ is $\{\pm I\}G_2$, where $G_2$ is the
stabilizer of $\varphi$. Since $\rho(-I)=\Id_W$, its image
in $GL(W)$ is precisely $\rho(G_2)$.
\end{enumerate}
\end{lemma}

For orthonormal $u,v,w$, \eqref{eq:rank-one-phi-recovery}
reduces to
\[
C_\varphi(E_u,E_v,E_w)
=343\varphi(u,v,w)^2-63,
\]
recovering the identity used in the proof of
Lemma~\ref{lem:invariant-cubic-plane}.

\begin{proof}
\emph{(\ref{item:rank-one-phi-recovery-evaluation}) Evaluation on test tensors.}
Both sides of \eqref{eq:rank-one-phi-recovery} are invariant
under simultaneous orthogonal transport of $\varphi$ and
$u,v,w$, so it suffices to consider the standard form. Recall,
from \eqref{eq:projection-tensors} and \eqref{eq:test-tensors},
\[
P_u=u^\flat\otimes u^\flat,
\qquad
E_u=7P_u-g(u,u)\,g.
\]
We extend $C_\varphi$ to all symmetric $2$-tensors using
the contraction formula \eqref{eq:natural-cubics}. The rank-one and metric contraction rules
\eqref{eq:p-rank-one-rule} and \eqref{eq:p-identity-rule}
give
\[
\begin{aligned}
C_\varphi(P_u,P_v,P_w)
&=\varphi(u,v,w)^2,\\
C_\varphi(g,P_v,P_w)
&=g(v,v)g(w,w)-g(v,w)^2,\\
C_\varphi(g,g,P_w)&=6g(w,w),\\
C_\varphi(g,g,g)&=42.
\end{aligned}
\]
The first identity follows from
$g(u\times v,w)=\varphi(u,v,w)$. Expanding
$C_\varphi(7P_u-g(u,u)g,
7P_v-g(v,v)g,
7P_w-g(w,w)g)$
using these identities gives
\eqref{eq:rank-one-phi-recovery}.

\emph{(\ref{item:rank-one-phi-recovery-form}) Recovery of the $3$-form.}
Suppose $C_{\varphi'}=C_\varphi$.
Since both forms are compatible with $g$,
the metric terms in \eqref{eq:rank-one-phi-recovery}
agree. Hence
\[
\varphi'(u,v,w)^2=\varphi(u,v,w)^2
\qquad(u,v,w\in V).
\]
Equivalently,
\[
[(\varphi'-\varphi)(u,v,w)]
[(\varphi'+\varphi)(u,v,w)]=0.
\]
As an identity in the polynomial ring
$\mathbb R[u,v,w]$, which is an integral domain,
this implies $\varphi'=\varphi$ or
$\varphi'=-\varphi$.

By the normalization in
Lemma~\ref{lem:standard-contractions},
$\varphi$ and $-\varphi$ induce opposite orientations.
Thus the orientation determines the sign.

\emph{(\ref{item:rank-one-phi-recovery-stabilizer}) Orthogonal stabilizer.}
For $O\in O(7)$, \eqref{eq:rho-group} gives
$\rho(O)E_u=E_{Ou}$.
By \eqref{eq:rank-one-phi-recovery},
$\rho(O)$ preserves $C_\varphi$ if and only if
\[
(O^*\varphi)(u,v,w)^2=\varphi(u,v,w)^2
\qquad(u,v,w\in V).
\]
By (\ref{item:rank-one-phi-recovery-form}), this is equivalent to
$O^*\varphi=\pm\varphi$.
The orientation normalization forces the sign
to equal $\det O$. Therefore the subgroup
in question is
\[
G_2\cup(-I)G_2=\{\pm I\}G_2.
\]
Finally, $\rho(-I)=\Id_W$, so the image of $\{\pm I\}G_2$
in $GL(W)$ is $\rho(G_2)$.
\end{proof}

Lemma \ref{lem:rank-one-phi-recovery} concerns the orthogonal group only. In the next proposition, we consider the full infinitesimal stabilizer of $C_\varphi$ in $\End(W)$.

\begin{proposition}
\label{prop:phi-cubic-stabilizer}
Let $\mathbb F\in\{\R,\C\}$. Recall that $\theta_M$ denotes the
infinitesimal pullback of $\varphi$ by $M$ from \eqref{eq:insert-M},
and that
\[
\gtwo^{\mathbb F}
=\{M\in\liegl(V_{\mathbb F}):\theta_M=0\}
\]
is the infinitesimal stabilizer of $\varphi$, extended
complex-linearly when $\mathbb F=\C$. Then
\[
\ker\Phi_{C_\varphi}=\rho(\gtwo^{\mathbb F}),
\qquad
\dim_{\mathbb F}\ker\Phi_{C_\varphi}=14.
\]
The infinitesimal stabilizer of the line
$\mathbb F C_\varphi$ is
\[
\rho(\gtwo^{\mathbb F})
\oplus\mathbb F\,\Id_W,
\]
of dimension $15$.
\end{proposition}

\begin{proof}
\emph{Reduction to skew endomorphisms of $V$.}
Let $L\in\ker\Phi_{C_\varphi}$. For $h\in W_{\mathbb F}$, the
identities \eqref{eq:p-rank-one-rule}, \eqref{eq:p-identity-rule}
and \eqref{eq:phi-cubic-contraction} give
\[
C_\varphi(h,E_u,E_v)
=49h(u\times v,u\times v)+7g(v,v)h(u,u)+7g(u,u)h(v,v),
\]
where \(h(u,u)=h_{ij}u^iu^j\). Substituting into
$0=C_\varphi(LE_u,E_v,E_v)+2C_\varphi(LE_v,E_u,E_v)$ and dividing by
$14$ gives
\begin{equation}
0=g(v,v)\bigl[\beta(u,v)+\beta(v,u)\bigr]
+7\beta(v,u\times v)+g(u,u)\beta(v,v).
\label{eq:phi-cubic-evaluation}
\end{equation}
Taking $u=v$ yields $3g(v,v)\beta(v,v)=0$, hence $\beta(v,v)=0$
for every real $v$.

For orthonormal $u,v$, put $w=u\times v$. By
\eqref{eq:phi-phi-contractions}, $(u,v,w)$ is orthonormal with
$v\times w=u$ and $w\times u=v$. Equation~\eqref{eq:phi-cubic-evaluation}
and the same equation with $u$ and $v$ exchanged give
\[
\begin{aligned}
\beta(u,v)+\beta(v,u)+7\beta(v,w)&=0,\\
\beta(u,v)+\beta(v,u)+7\beta(u,w)&=0,
\end{aligned}
\]
where $\beta(u,-w)=\beta(u,w)$. Thus $\beta(u,w)=\beta(v,w)$.
Applying the same argument cyclically, we may write
\[
\begin{aligned}
b_u&=\beta(v,u)=\beta(w,u),\\
b_v&=\beta(w,v)=\beta(u,v),\\
b_w&=\beta(u,w)=\beta(v,w).
\end{aligned}
\]
The three cyclic evaluations of \eqref{eq:phi-cubic-evaluation} become
\[
b_u+b_v+7b_w=0,\qquad
b_v+b_w+7b_u=0,\qquad
b_w+b_u+7b_v=0.
\]
Subtracting shows $b_u=b_v=b_w$, so all three vanish.
Homogeneity now gives $\beta(u,v)=0$ for all real orthogonal pairs,
and Lemma~\ref{lem:rank-one-recognition} gives
$L=\rho(M)$ for some $M\in\lieso(7,\mathbb F)$.

\emph{Identification of $\gtwo^{\mathbb F}$.}
For real skew $M$, put $O_t=e^{tM}$. Since $\rho(O_t)E_u=E_{O_tu}$,
differentiating \eqref{eq:rank-one-phi-recovery} gives
\begin{equation}
\begin{aligned}
\Phi_{C_\varphi}(\rho(M))(E_u,E_v,E_w)
&=343\left.\frac{d}{dt}\right|_{t=0}
\varphi(O_tu,O_tv,O_tw)^2\\
&=686\varphi(u,v,w)\theta_M(u,v,w).
\end{aligned}
\label{eq:skew-phi-cubic}
\end{equation}
The metric terms have zero derivative because $O_t$ is orthogonal.
By linearity in $M$, the identity extends to complex skew endomorphisms.

If $\rho(M)$ annihilates $C_\varphi$, then
$\varphi(u,v,w)\theta_M(u,v,w)=0$. Hence $\theta_M$ vanishes on the set
where $\varphi\ne0$, which is dense in $V^3$, and therefore vanishes
everywhere by continuity. Conversely, if $\theta_M=0$, then
\eqref{eq:skew-phi-cubic} vanishes on all triples of test tensors $(E_u,E_v,E_w)$;
since the $E_u$ span $W_{\mathbb F}$ by Lemma~\ref{lem:rank-one-recognition}\,(\ref{item:rank-one-recognition-span}), $\rho(M)\in\ker\Phi_{C_\varphi}$.
This proves $\ker\Phi_{C_\varphi}=\rho(\gtwo^{\mathbb F})$.
Its dimension is $14$ by \eqref{eq:two-form-splitting} and injectivity of
$\rho$; \eqref{eq:generator-line-stabilizer} gives the line stabilizer
and its dimension $15$.
\end{proof}

We next construct the two copies of $V(1,1)$ in $\End \WC$ explicitly. The same
$64$-dimensional tensor space will map to both copies, one consisting
of $W$-self-adjoint endomorphisms and the other of $W$-skew-adjoint endomorphisms. The tensor model $\mathcal Q_{64}$ as used below is standard, and can be found in
\cite[Lemma~2.1]{CleytonIvanovCurv} and \cite[(4.27)]{dwivedi-2025-flows}.

\begin{definition}
\label{def:q64}
Let $\mathcal Q_{64}=\ker\bigl(\operatorname{Alt}\colon V^*\otimes\Lambda^2_{14}\to\Lambda^3V^*\bigr)$;
equivalently, $\mathcal Q_{64}$ is the space of real tensors $Q_{abc}$ with
\begin{subequations}
\begin{align}
Q_{abc}&=-Q_{acb},\label{eq:q64-skew}\\
Q_{abc}\varphi_d{}^{bc}&=0,\label{eq:q64-phi}\\
Q_{abc}+Q_{bca}+Q_{cab}&=0,\label{eq:q64-cyclic}
\end{align}
\end{subequations}
and $\mathcal Q_{64,\C}=\mathcal Q_{64}\otimes\C$. For $Q\in\mathcal Q_{64}$
define the bilinear form $\varepsilon_Q$ on covariant $2$-tensors
$h,k\in V^*\otimes V^*$ by
\[
\varepsilon_Q(h,k)=Q_{abc}\varphi^c{}_{de}h^{ad}k^{be}.
\]
This is the form through which $Q$ acts on $W$.
Its symmetric and alternating parts are
\[
\varepsilon^{\mathrm s}_Q(h,k)
=\tfrac12\bigl(\varepsilon_Q(h,k)+\varepsilon_Q(k,h)\bigr),
\qquad
\varepsilon^{\mathrm a}_Q(h,k)
=\tfrac12\bigl(\varepsilon_Q(h,k)-\varepsilon_Q(k,h)\bigr).
\]
On $W$, let $J^{\mathrm s}_Q,J^{\mathrm a}_Q \in \End W$ be the corresponding endomorphisms under the $W$-pairing:
\begin{equation}
\label{eq:JQ-epsQ}
\langle J^{\mathrm s}_Qh,k\rangle_W=\varepsilon^{\mathrm s}_Q(h,k),
\qquad
\langle J^{\mathrm a}_Qh,k\rangle_W=\varepsilon^{\mathrm a}_Q(h,k).
\end{equation}
We call $Q\mapsto J^{\mathrm s}_Q$ and $Q\mapsto J^{\mathrm a}_Q$ the \emph{parity embeddings}; the labels $\mathrm s$ and $\mathrm a$ refer
to symmetry and skew-symmetry under the $W$-adjoint. All these definitions
extend complex-linearly in $Q$ and complex-bilinearly in $h,k$.
Equivalently, for $h\in W$, the operators act by
\begin{equation}
\begin{aligned}
(J_Q^{\mathrm s} h)_{ij}
&=\frac14\Bigl[
(Q_{aic}- Q_{iac})\varphi^c{}_{dj}
+(Q_{ajc}- Q_{jac})\varphi^c{}_{di}
\Bigr]h^{ad},\\
(J_Q^{\mathrm a} h)_{ij}
&=\frac14\Bigl[
(Q_{aic}+ Q_{iac})\varphi^c{}_{dj}
+(Q_{ajc}+ Q_{jac})\varphi^c{}_{di}
\Bigr]h^{ad}.
\end{aligned}
\label{eq:q64-parity-action}
\end{equation}
The right-hand side is symmetric and traceless in
$i,j$, as follows from the defining identities
of $\mathcal Q_{64}$.
\end{definition}

Conditions \eqref{eq:q64-skew} and
\eqref{eq:q64-phi} ensure that each slice
$Q_a=Q(e_a,\cdot,\cdot)$ belongs to $\Lambda^2_{14}$.
Since $Q$ is skew in its last two indices,
\[
\operatorname{Alt}(Q)_{abc}
=\tfrac13(Q_{abc}+Q_{bca}+Q_{cab}),
\]
so \eqref{eq:q64-cyclic} is equivalent to
$\operatorname{Alt}(Q)=0$. This establishes the equivalence
of the two descriptions of $\mathcal Q_{64}$.
Finally, the nondegenerate $W$-pairing
uniquely determines $J^{\mathrm s}_Q$ from $\varepsilon^{\mathrm s}_Q$
and $J^{\mathrm a}_Q$ from $\varepsilon^{\mathrm a}_Q$.

The contraction identities of the model and the symmetries of its forms are as follows.

\begin{lemma}
\label{lem:q64-model}
Let $Q\in\mathcal Q_{64}$.
\begin{enumerate}
\renewcommand{\labelenumi}{(\roman{enumi})}
\renewcommand{\theenumi}{\roman{enumi}}
\item\label{item:q64-model-eigenvalue} Every $\alpha\in\Lambda^2_{14}$ satisfies
$\psi_{rsbc}\alpha^{bc}=-2\alpha_{rs}$.
In particular,
\begin{equation}
\psi_{rsbc}Q_a{}^{bc}=-2Q_{ars}.
\label{eq:q64-eigenvalue}
\end{equation}

\item\label{item:q64-model-contractions} The following contractions vanish:
\begin{subequations}
\begin{align}
Q^a{}_{ab}&=0,
\label{eq:q64-trace}\\
Q_{abc}\varphi^{ab}{}_{d}&=0.
\label{eq:q64-phi-contraction}
\end{align}
\end{subequations}

\item\label{item:q64-model-traceless} For all symmetric $2$-tensors $h,k$,
\[
\varepsilon_Q(h,g)=\varepsilon_Q(g,k)=0.
\]
Thus $\varepsilon_Q$, $\varepsilon^{\mathrm s}_Q$ and $\varepsilon^{\mathrm a}_Q$
depend only on the traceless parts of their arguments.

\item\label{item:q64-model-equivariance} The maps
$J^{\mathrm s},J^{\mathrm a}:\mathcal Q_{64}\to\End(W)$ are
$G_2$-equivariant, with $O\in G_2$ acting on
$\End(W)$ by $X\mapsto\rho(O)X\rho(O)^{-1}$. Moreover, $J^{\mathrm s}_Q$ is $W$-self-adjoint
and $J^{\mathrm a}_Q$ is $W$-skew-adjoint.
\end{enumerate}
\end{lemma}

\begin{proof}
Part (\ref{item:q64-model-eigenvalue}) follows from \eqref{eq:two-form-splitting},
since each slice $Q_a$ belongs to $\Lambda^2_{14}$.

For (\ref{item:q64-model-contractions}), equation~\eqref{eq:q64-eigenvalue}
and the cyclic identity give
\[
Q^a{}_{as}
=\tfrac12\psi_{sabc}Q^{abc}
=\tfrac16\psi_{sabc}
  (Q^{abc}+Q^{bca}+Q^{cab})=0.
\]
Similarly, contracting the cyclic identity
with $\varphi^{ab}{}_{d}$ gives
\[
Q_{abc}\varphi^{ab}{}_{d}
=-Q_{bca}\varphi^{ab}{}_{d}
=-Q_{abc}\varphi^{ab}{}_{d},
\]
where the $Q_{cab}$ term vanishes by
\eqref{eq:q64-phi}, and the last equality uses
\eqref{eq:q64-skew}. This proves (\ref{item:q64-model-contractions}).

For (\ref{item:q64-model-traceless}), inserting $g$ into the second
argument of $\varepsilon_Q$ gives zero by
\eqref{eq:q64-phi}; inserting it into the
first gives zero by
\eqref{eq:q64-phi-contraction} and
\eqref{eq:q64-skew}.

Finally, $\mathcal Q_{64}$ is $G_2$-invariant
because $\Lambda^2_{14}$ and the alternation
map are equivariant. The defining contraction
of $\varepsilon_Q$ and the $W$-pairing are
also $G_2$-invariant, proving equivariance
of $J^{\mathrm s}$ and $J^{\mathrm a}$. Their adjointness properties
follow directly from the symmetry of
$\varepsilon^{\mathrm s}_Q$ and the skew-symmetry of $\varepsilon^{\mathrm a}_Q$.
\end{proof}

Differentiating the equivariance in
Lemma~\ref{lem:q64-model}\,(\ref{item:q64-model-equivariance}) gives
\[
J^{\mathrm s}_{M\cdot Q}=[\rho(M),J^{\mathrm s}_Q],
\qquad
J^{\mathrm a}_{M\cdot Q}=[\rho(M),J^{\mathrm a}_Q],
\qquad M\in\gtwo,
\]
where $M\cdot Q$ denotes the natural infinitesimal
action on $\mathcal Q_{64}$.

\begin{corollary}
\label{cor:q64-is-64}
The complexification $\mathcal Q_{64,\C}$ is
the irreducible $\gtwoC$-module $V(1,1)$.
Consequently, $\mathcal Q_{64}$ has real
dimension $64$ and is irreducible under
both $G_2$ and $\gtwo$.

Moreover, under the metric identification
$\Lambda^2_{14}\cong\gtwo$,
\[
\VC\otimes\gtwoC
\cong\mathcal Q_{64,\C}\oplus\WC\oplus\VC.
\]
In particular, $\mathcal Q_{64,\C}$ is the
unique submodule of type $V(1,1)$ in this
tensor product.
\end{corollary}

\begin{proof}
The decomposition
\eqref{eq:seven-adjoint-product} gives
\[
(V^*\otimes\Lambda^2_{14})_{\C}
\cong V(1,1)\oplus V(2,0)\oplus V(1,0).
\]
By \cite[Lemma~2.1]{CleytonIvanovCurv},
the kernel of alternation is precisely the
$V(1,1)$ summand. The other two summands
are $\WC$ and $\VC$, respectively.
The dimension and uniqueness follow, while
irreducibility over $\R$ follows by
complexification.
\end{proof}

\begin{example}
\label{ex:q64-test-tensor}
In the standard frame of \eqref{eq:standard-positive-form}, choose
\[
Q_1=-e^{14}+e^{27},\qquad Q_2=e^{17}+e^{24},\qquad
Q_3=\cdots=Q_7=0.
\]
The contractions of the two nonzero slices with $\varphi$ cancel
in the $e_5$ and $e_6$ components, respectively. The cyclic identity
holds on repeated indices by skew-symmetry, and the only triple of
distinct indices on which it is nontrivial is $\{1,2,7\}$, where
$Q_{127}+Q_{271}+Q_{712}=1-1+0=0$.
Thus $Q\in\mathcal Q_{64}$.
\end{example}

In Section~\ref{sec:deformation-representation} we identified the
$W$-self-adjoint and $W$-skew-adjoint endomorphisms of $\WC$ equivariantly with $S^2\WC$ and
$\Lambda^2\WC$; see \eqref{eq:endomorphism-inventory},
\eqref{eq:symmetric-square} and \eqref{eq:alternating-square}.
Since $\mathcal Q_{64,\C}$ is irreducible, each parity map is
injective whenever it is nonzero. A single pair of test elements shows
that both maps are nonzero.

\begin{lemma}
\label{lem:parity-embeddings-injective}
For the maps $Q\mapsto J_Q^{\mathrm s}$ and $Q\mapsto J_Q^{\mathrm a}$ defined by \eqref{eq:q64-parity-action}, we have the following.
\begin{enumerate}
\renewcommand{\labelenumi}{(\roman{enumi})}
\renewcommand{\theenumi}{\roman{enumi}}
\item\label{item:parity-embeddings-injectivity} The maps $Q\mapsto J^{\mathrm s}_Q$ and $Q\mapsto J^{\mathrm a}_Q$ are injective on $\mathcal Q_{64}$
  and on $\mathcal Q_{64,\C}$; equivalently, for every nonzero $Q$, real or complex, neither
  $\varepsilon^{\mathrm s}_Q$ nor $\varepsilon^{\mathrm a}_Q$ vanishes identically on $W$, respectively $\WC$.
\item\label{item:parity-embeddings-images} The maps $Q\mapsto J_Q^{\mathrm s}$ and $Q\mapsto J_Q^{\mathrm a}$ are $G_2$-equivariant isomorphisms of $\mathcal Q_{64,\C}$ onto the unique
  copies of the $64$-dimensional module $V(1,1)$ in $S^2\WC$ and in $\Lambda^2\WC$, respectively.
\end{enumerate}
\end{lemma}

\begin{proof}
Use the tensor $Q$ of Example~\ref{ex:q64-test-tensor} and take
\[
h=e^1\otimes e^1-e^7\otimes e^7,\qquad
k=e^1\otimes e^5+e^5\otimes e^1.
\]
Both tensors belong to $W$. The defining contraction gives
\[
\begin{aligned}
\varepsilon_Q(h,k)&=Q_{114}\varphi_{415}=1,\\
\varepsilon_Q(k,h)&=Q_{114}\varphi_{451}-Q_{172}\varphi_{257}=-2.
\end{aligned}
\]
Hence $\varepsilon_Q^{\mathrm s}(h,k)=-\tfrac12$ and
$\varepsilon_Q^{\mathrm a}(h,k)=\tfrac32$, so both parity maps are nonzero.
By Lemma~\ref{lem:q64-model}\,(\ref{item:q64-model-equivariance}),
their complexifications are equivariant. Their kernels are therefore
submodules of the irreducible $\mathcal Q_{64,\C}$, so both maps are
injective, as are their real restrictions.

The images lie in $S^2\WC$ and $\Lambda^2\WC$, respectively,
by the adjointness properties in
Lemma~\ref{lem:q64-model}\,(\ref{item:q64-model-equivariance}).
Each image is a copy of $V(1,1)$, which occurs with multiplicity one
in each target by \eqref{eq:symmetric-square} and
\eqref{eq:alternating-square}.
\end{proof}

The two generator stabilizers and the parity model of the
$\mathbf{64}$ now reduce the remaining count to the kernel restricted
to the isotypic component of $\mathbf{64}$ in $\End(\WC)$.

\subsection{The exact stabilizer count}
\label{subsec:nothing-else-jumps}

The remaining component is
\[
\mathcal H_{64}=\End(\WC)_{V(1,1)}
 =J^{\mathrm s}(\mathcal Q_{64,\C})\oplus J^{\mathrm a}(\mathcal Q_{64,\C}).
\]
Write $\Phi_0^{(64)}=\Phi_0|_{\mathcal H_{64}}$ and
$\Phi_\varphi^{(64)}=\Phi_\varphi|_{\mathcal H_{64}}$.
Proposition~\ref{prop:trace-cubic-stabilizer} and the decomposition
following it identify the trace-cubic kernel as $\mathbf7\oplus\mathbf{14}$, so $\Phi_0^{(64)}$ is injective.
Using the Koecher lower bound, we will show that
$\Phi_\varphi^{(64)}(\mathcal H_{64})\subseteq \Phi_0^{(64)}(\mathcal H_{64})$ and hence
define a relative operator. We then compute its matrix with respect
to the parity embeddings. Set
\[
\mathfrak m_\pm=\ker(\Phi_0^{(64)}+\lambda_\pm\Phi_\varphi^{(64)})
=\ker\Phi_{1,\lambda_\pm}|_{\mathcal H_{64}}.
\]

\begin{lemma}
\label{lem:64-block}
For every $Q\in\mathcal Q_{64,\C}$,
\begin{equation}
\begin{gathered}
\Phi_{C_\varphi}(J^{\mathrm s}_Q)=-\Phi_{C_0}(J^{\mathrm s}_Q)-\Phi_{C_0}(J^{\mathrm a}_Q),\\
\Phi_{C_\varphi}(J^{\mathrm a}_Q)=7\,\Phi_{C_0}(J^{\mathrm s}_Q)-\Phi_{C_0}(J^{\mathrm a}_Q).
\end{gathered}
\label{eq:two-phi-identities}
\end{equation}
Consequently, there is a unique endomorphism $T$ of
$\mathcal H_{64}$ such that $\Phi_\varphi^{(64)}=\Phi_0^{(64)}\circ T$.
In the ordered pair of embeddings $(J^{\mathrm s},J^{\mathrm a})$, with images written
as columns, the matrix of $T$ is
\[
T=\begin{pmatrix}-1&7\\-1&-1\end{pmatrix}.
\]
In particular, $\Phi_{a,b}=\Phi_{C_0}\circ(a\,\Id+bT)$ on this component, where, for the matrix of $T$,
\begin{equation}
\det(a\,\Id+bT)=(a-b)^2+7b^2=a^2-2ab+8b^2 .
\label{eq:relative-determinant}
\end{equation}
For $(a,b)\ne0$, the kernel of $\Phi_{a,b}$ on this component is
nonzero exactly when $a^2-2ab+8b^2=0$, equivalently when
$\lambda=b/a=\lambda_\pm$. At these parameters, we have
\[
\ker\Phi_{a,b}|_{\mathcal H_{64}}=\mathfrak m_\pm
=\{J^{\mathrm s}_Q+\zeta_\pm J^{\mathrm a}_Q:Q\in\mathcal Q_{64,\C}\},
\qquad \zeta_\pm=\pm\frac{\ii}{\sqrt7}.
\]
Here $\zeta=b/(a-b)=\lambda/(1-\lambda)$ is the ratio of the coefficient
of $J^{\mathrm a}_Q$ to that of $J^{\mathrm s}_Q$.
\end{lemma}

\begin{proof}
By Proposition~\ref{prop:koecher-jump}, each critical stabilizer
contains a copy of $V(1,1)$, which lies in $\mathcal H_{64}$ by the
definition of $\mathcal H_{64}$ as the $V(1,1)$-isotypic component. Thus $\dim \mathfrak m_\pm\ge64$. If $X$ lies in both kernels,
subtraction gives
$(\lambda_+-\lambda_-)\Phi_\varphi^{(64)}(X)=0$, hence $\Phi_\varphi^{(64)}(X)=\Phi_0^{(64)}(X)=0$.
Injectivity of $\Phi_0^{(64)}$ gives $X=0$. Since $\dim \mathcal H_{64}=128$, it follows that
\[
\mathcal H_{64}=\mathfrak m_+\oplus \mathfrak m_-,\qquad \dim \mathfrak m_+=\dim \mathfrak m_-=64.
\]
In particular each kernel is the irreducible copy just obtained.
For $X_\pm\in \mathfrak m_\pm$ we have
$\Phi_\varphi^{(64)}(X_\pm)=-\lambda_\pm^{-1}\Phi_0^{(64)}(X_\pm)$. Decomposing any $X\in \mathcal H_{64}$ as
$X_++X_-$ therefore proves $\Phi_\varphi^{(64)}(\mathcal H_{64})\subseteq \Phi_0^{(64)}(\mathcal H_{64})$. There is a unique
equivariant operator $T=(\Phi_0^{(64)})^{-1}\Phi_\varphi^{(64)}$ on $\mathcal H_{64}$, and
\[
T|_{\mathfrak m_\pm}=-\lambda_\pm^{-1}\Id=(-1\pm\ii\sqrt7)\Id.
\]
These two eigenvalues already give
$\det(a\,\Id+bT)=a^2-2ab+8b^2$ on the two-dimensional multiplicity space of $V(1,1)$ in $\mathcal H_{64}$.
It remains to identify the eigenlines in the chosen parity basis.

Use the standard $\varphi$ of \eqref{eq:standard-positive-form}, and let
$E_{ij}$ denote the matrix unit with its sole nonzero entry equal to $1$
in row $i$, column $j$.
For this computation, put $D_a=E_{aa}-E_{77}$ and $S_{ij}=E_{ij}+E_{ji}$.
Use $Q$ from Example~\ref{ex:q64-test-tensor} and choose $h=D_1+S_{37}$.
Direct multiplication and the definition of $\mathcal P$ give
\[
h^2=E_{11}+E_{33}+2E_{77}-S_{37},\qquad
\mathcal P(h,h)=-2E_{44}-2E_{66}+2S_{26}.
\]
For any symmetric $k$, the same four nonzero slice components give
\[
\varepsilon_Q(h,k)=k_{15}+k_{26}+k_{37},\qquad
\varepsilon_Q(k,h)=k_{11}-k_{22}-2k_{15}.
\]
Consequently
\begin{equation}
\bigl(\varepsilon_Q(h,h^2),\varepsilon_Q(h^2,h),
 \varepsilon_Q(h,\mathcal P(h,h)),
 \varepsilon_Q(\mathcal P(h,h),h)\bigr)=(-1,1,2,0).
\label{eq:q64-calibration}
\end{equation}
By \eqref{eq:quadratic-adjoint-pairings}, \eqref{eq:phi-cubic-contraction},
Lemma~\ref{lem:q64-model}\,(\ref{item:q64-model-traceless}), \eqref{eq:JQ-epsQ} and the half-sum and half-difference convention,
these values give
\[
\tfrac13\Phi_{1,\lambda}(J_Q^{\mathrm s})(h,h,h)=\lambda,\qquad
\tfrac13\Phi_{1,\lambda}(J_Q^{\mathrm a})(h,h,h)=\lambda-1.
\]
Since $\lambda_\pm-1\ne0$, neither critical kernel is the whole skew
copy $J^{\mathrm a}(\mathcal Q_{64,\C})$. The projection of each
$\mathfrak m_\pm$ onto the symmetric copy is therefore nonzero, hence an isomorphism by irreducibility and
Schur's lemma. Each $\mathfrak m_\pm$ is consequently a graph
$\{J_Q^{\mathrm s}+\zeta_\pm J_Q^{\mathrm a}:Q\in\mathcal Q_{64,\C}\}$.
Evaluating $\Phi_{1,\lambda_\pm}$ on $J_Q^{\mathrm s}+\zeta_\pm J_Q^{\mathrm a}$
at the displayed $Q$ and $h$ determines $\zeta_\pm$:
\[
\lambda_\pm+\zeta_\pm(\lambda_\pm-1)=0,\qquad
\zeta_\pm=\frac{\lambda_\pm}{1-\lambda_\pm}
 =\pm\frac{\ii}{\sqrt7}.
\]
In particular the plus critical parameter corresponds to the plus
sign of $\ii/\sqrt7$.

The two graph lines $(1,\zeta_\pm)$ and their eigenvalues
$-1\pm\ii\sqrt7$ determine the matrix of $T$ uniquely.
It is $\begin{pmatrix}-1&7\\-1&-1\end{pmatrix}$, since
$7\zeta_\pm=\pm\ii\sqrt7$ and $7\zeta_\pm^2=-1$.
Its columns give \eqref{eq:two-phi-identities}. Finally,
$\Phi_{a,b}=\Phi_0^{(64)}\circ(a\Id+bT)$ on $\mathcal H_{64}$ and $\Phi_0^{(64)}$ is injective, so the
determinant and the graph description give all the asserted kernels,
including the endpoints of the projective pencil.
\end{proof}

\begin{remark}
\label{rem:q64-literature}
The matrix of $T$ in Lemma~\ref{lem:64-block} gives $(\Id+T)^2=-7\,\Id$. In particular,
\[
(\Id+T)(J_Q^{\mathrm s}+\zeta_\pm J_Q^{\mathrm a})
=7\zeta_\pm(J_Q^{\mathrm s}+\zeta_\pm J_Q^{\mathrm a}),
\]
so these are the eigenspaces of $\Id+T$ with
eigenvalues $\pm\ii\sqrt7$, respectively.
We have not found the explicit parity embeddings, the identities
\eqref{eq:two-phi-identities}, or the resulting relative operator
$T$ with $(\Id+T)^2=-7\,\Id$ in the literature consulted.

The appearance of $\sqrt{-7}$ in the earlier $G_2\subset E_6$
constructions is discussed in Section~\ref{sec:introduction}. In particular,
Aschbacher's result, reported in \cite[p.~300]{testerman-1989-construction},
states that for a field $\mathbb F$ with $\operatorname{char}\mathbb F\neq2,3,7$,
the group $G_2(\mathbb F)$ preserves, up to scalar multiple, exactly two forms
similar to the $E_6(\mathbb F)$ form if and only if $\sqrt{-7}\in\mathbb F$.
\end{remark}

We can now assemble the kernels. For the critical kernels
$\mathfrak m_\pm$ of Lemma~\ref{lem:64-block}, define
\[
\mathcal E_\pm\colon\mathcal Q_{64,\C}\longrightarrow\mathfrak m_\pm,
\qquad \mathcal E_\pm(Q)=J^{\mathrm s}_Q+\zeta_\pm J^{\mathrm a}_Q.
\]
By Lemma~\ref{lem:q64-model}\,(\ref{item:q64-model-equivariance}) and its infinitesimal form above, $\mathcal E_\pm$ is a $\gtwoC$-equivariant map onto
$\mathfrak m_\pm$, and it is injective because the $W$-self-adjoint part of $\mathcal E_\pm(Q)$
is $J^{\mathrm s}_Q$ (Lemma~\ref{lem:parity-embeddings-injective}\,(\ref{item:parity-embeddings-injectivity})). Hence
$\mathfrak m_\pm\cong\mathcal Q_{64,\C}\cong V(1,1)$ by Corollary~\ref{cor:q64-is-64}.

\begin{proposition}
\label{prop:block-kernels}
Let $(a,b)\in\C^2\setminus\{0\}$. Then $\ker\Phi_{a,b}$ is the direct sum of its
intersections with the isotypic components $\End(\WC)_{V(j_1,j_2)}$ of $\End(\WC)$, and
\begin{equation}
\dim\ker\Phi_{a,b}=
\begin{cases}
14+7=21, & b=0,\\
14+64=78, & a^2-2ab+8b^2=0,\\
14, & \text{otherwise}.
\end{cases}
\label{eq:block-count}
\end{equation}
The contributions are
\[
\begin{array}{c|c|c}
\text{type}&\text{dimension in the kernel}&\text{condition}\\ \hline
V(0,0)&0&\text{all parameters}\\
V(0,1)&14&\text{all parameters}\\
V(1,0)&7&b=0\\
V(1,1)&64&a^2-2ab+8b^2=0\\
\text{five remaining types}&0&\text{all parameters}
\end{array}
\]
The contributions of $V(1,0)$ and $V(1,1)$ vanish outside the indicated
loci. At a critical parameter, the $64$-dimensional contribution is a
single copy of $V(1,1)$, and
\[
\mathfrak g_\pm:=\stab(N_\pm)=\ker\Phi_{1,\lambda_\pm}=\rho(\gtwoC)\oplus\mathfrak m_\pm .
\]
\end{proposition}

\begin{proof}
By Lemma~\ref{lem:common-g2}, $\Phi_{a,b}$ is
$\gtwoC$-equivariant, so its kernel decomposes
over the isotypic components of $\End(\WC)$ listed
in Table~\ref{tab:endomorphism-inventory}.

The component $V(0,1)=\rho(\gtwoC)$ lies in every
kernel and has dimension $14$. The scalar
component contributes nothing, since
$\Phi_{a,b}(\Id_W)=3(aC_0+bC_\varphi)\ne0$
by Lemma~\ref{lem:invariant-cubic-plane}.

The unique $V(1,0)$ component is
$\rho(\Lambda^2_7)_{\C}$. It lies in
$\ker\Phi_{C_0}$ by
Proposition~\ref{prop:trace-cubic-stabilizer},
but has trivial intersection with
$\ker\Phi_{C_\varphi}=\rho(\gtwoC)$ by
Proposition~\ref{prop:phi-cubic-stabilizer}.
Thus it contributes dimension $7$ when $b=0$
and zero otherwise.

On the $V(1,1)$-isotypic component,
Lemma~\ref{lem:64-block} gives a trivial kernel
unless $a^2-2ab+8b^2=0$. At the two critical
parameters $\lambda_\pm=b/a$, the kernel is
$\mathfrak m_\pm\cong V(1,1)$, of dimension $64$.

By Proposition~\ref{prop:five-blocks}, all
remaining components have trivial kernel.
These are the two copies of $V(2,0)$ and the
summands of types $V(0,2)$, $V(4,0)$,
$V(3,0)$, and $V(2,1)$.

Adding these contributions gives
\eqref{eq:block-count}. In particular, at the
critical parameters,
\[
\ker\Phi_{1,\lambda_\pm}
=\rho(\gtwoC)\oplus\mathfrak m_\pm
\]
has dimension $14+64=78$, attaining the bound
of Proposition~\ref{prop:koecher-jump}.
\end{proof}

Since $\stab(N_{a,b})=\ker\Phi_{a,b}$, the count \eqref{eq:block-count}
proves Theorem~\ref{thm:stabilizer-trichotomy}: at $b=0$ the quadratic
equals $a^2\neq0$, so the three cases are disjoint, and on the locus
$a^2-2ab+8b^2=0$ we have $a\neq0$, and $\lambda=b/a$ is a root of
$8\lambda^2-2\lambda+1$, one of the two non-real parameters $\lambda_\pm$ of \eqref{eq:critical-parameters}.

Proposition~\ref{prop:block-kernels} identifies the stabilizer itself
as well as its dimension: it is $\rho(\lieso(7,\mathbb C))$ for $C_0$,
$\rho(\gtwoC)\oplus\mathfrak m_\pm$ for $N_\pm$, and
$\rho(\gtwoC)$ for every other member, including $C_\varphi$.
Thus the trace endpoint has the orthogonal symmetry algebra, whereas
the Yukawa endpoint $C_\varphi$ retains exactly $\rho(\gtwoC)$.
Only the two mixed non-real members $N_\pm$ acquire the additional
$64$-dimensional summand. The exceptional symmetry therefore belongs
to this conjugate pair, rather than to the pointwise Yukawa cubic itself.
It remains to identify the Lie algebra structure at the critical parameters.

\subsection{The critical stabilizer}
\label{subsec:critical-type}

In this section we will identify the $78$-dimensional algebras
$\mathfrak g_\pm=\rho(\gtwoC)\oplus\mathfrak m_\pm$ of Proposition~\ref{prop:block-kernels},
with the isomorphism $\mathcal E_\pm\colon\mathcal Q_{64,\C}\to\mathfrak m_\pm$ fixed above. We will first recall some terminology from Lie theory.
An irreducible representation is \emph{minuscule}
if its weights form a single orbit under the
Weyl group, in which case every weight has
multiplicity one.

\begin{theorem}
\label{thm:critical-type}
For either sign, $\mathfrak g_\pm$ is a simple
complex Lie algebra of type $E_6$. Its action
on $\WC$ is one of the two $27$-dimensional
minuscule representations
$V(\omega_1^{E_6})$ and $V(\omega_6^{E_6})$.
These representations are dual to one another
and are exchanged by the nontrivial
automorphism of the $E_6$ Dynkin diagram.
\end{theorem}

\begin{proof}
\emph{Simplicity.}
For a real endomorphism $X$, write
$\|X\|_{\End}^2=\langle X,X\rangle_{\End}$ using
\eqref{eq:endomorphism-pairing}. We choose a sign and write
\[
\mathfrak a=\rho(\gtwoC),
\qquad
\mathfrak m=\mathfrak m_\pm,
\qquad
\mathfrak g=\mathfrak g_\pm=\mathfrak a\oplus\mathfrak m.
\]
The summands $\mathfrak a\cong V(0,1)$ and
$\mathfrak m\cong V(1,1)$ are nonisomorphic
irreducible $\gtwoC$-modules. Every ideal of $\mathfrak g_\pm$ is invariant
under the adjoint action of
$\mathfrak a=\rho(\gtwoC)$. Since
\(
\mathfrak g=\mathfrak a\oplus\mathfrak m
\)
is a sum of two nonisomorphic irreducible
$\gtwoC$-modules, each of multiplicity one,
the only possible ideals are
$0$, $\mathfrak a$, $\mathfrak m$,
and $\mathfrak g$.

The action of $\mathfrak a$ on $\mathfrak m$
is nontrivial and irreducible, so
\[
[\mathfrak a,\mathfrak m]=\mathfrak m.
\]
In particular, $\mathfrak a$ is not an ideal.

To exclude $\mathfrak m$, consider the
invariant symmetric bilinear form
\[
\mathfrak b(X,Y)=\tr_W(XY)
\]
on $\mathfrak g$. By Schur's lemma,
$\mathfrak b(\mathfrak a,\mathfrak m)=0$, since the
two summands are nonisomorphic irreducible
$\mathfrak a$-modules. Moreover, $\mathfrak b$ is nondegenerate on
$\mathfrak m$. Indeed, for any nonzero
real $Q\in\mathcal Q_{64}$, put
$X=\mathcal E_\pm(Q)=J_Q^{\mathrm s}+\zeta_\pm J_Q^{\mathrm a}$.
Since $\zeta_\pm^2=-1/7$ and $J_Q^{\mathrm s}$ and
$J_Q^{\mathrm a}$ are respectively $W$-self-adjoint and
$W$-skew-adjoint,
\[
\begin{aligned}
\mathfrak b(X,X)
&=\tr_W\bigl((J_Q^{\mathrm s})^2\bigr)
  -\tfrac17\tr_W\bigl((J_Q^{\mathrm a})^2\bigr)\\
&=\|J_Q^{\mathrm s}\|_{\End}^2
  +\tfrac17\|J_Q^{\mathrm a}\|_{\End}^2>0.
\end{aligned}
\]
Thus $\mathfrak b$ is positive definite on the real
span of these elements, and its
complex-bilinear extension is nondegenerate
on $\mathfrak m$. If $\mathfrak m$ were an ideal, then
$[\mathfrak m,\mathfrak m]\subseteq\mathfrak m$.
By invariance of $\mathfrak b$, this would give
\[
\mathfrak b([A,U],V)=\mathfrak b(A,[U,V])=0
\]
for all $A\in\mathfrak a$ and
$U,V\in\mathfrak m$, since
$\mathfrak b(\mathfrak a,\mathfrak m)=0$.
But $[\mathfrak a,\mathfrak m]=\mathfrak m$
and $\mathfrak b$ is nondegenerate on $\mathfrak m$,
a contradiction.

Hence $\mathfrak g$ has no nonzero
proper ideals and is simple.

\emph{Identification of the type.}
By the classification of complex simple Lie algebras,
the only types of dimension $78$ are $B_6$, $C_6$,
and $E_6$.

The action of $\mathfrak g_\pm$ on $\WC$ is
irreducible, since its restriction to $\gtwoC$
is already irreducible. Neither $B_6$ nor $C_6$
admits an irreducible representation of dimension $27$
\cite[proof of Lemma~4, Table~1, p.~4]{kramer-2016-cubic}.
Hence $\mathfrak g_\pm$ is of type $E_6$.

\emph{The $27$-dimensional representation.}
Finally, the only $27$-dimensional irreducible representations of $E_6$ are
$V(\omega_1^{E_6})$ and $V(\omega_6^{E_6})$; they are mutually dual, minuscule, and
exchanged by the diagram automorphism \cite[\S3, Table~1]{podder-2023-k-ring},
\cite[\S\S26.5, 30.4, 32.1]{etingof-2024-lie}. Since $\WC$ is irreducible of dimension $27$
under $\mathfrak g_\pm$, it is one of them.
\end{proof}

\begin{corollary}
\label{cor:stabilizer-types}
For $[a:b]\in\mathbb P^1(\C)$,
\[
\stab(N_{a,b})\cong
\begin{cases}
\lieso(7,\C),&[a:b]=[1:0],\\
\mathfrak e_6,&a^2-2ab+8b^2=0,\\
\gtwoC,&\text{otherwise}.
\end{cases}
\]
\end{corollary}

\begin{proof}
This follows from Proposition~\ref{prop:block-kernels} and
Theorem~\ref{thm:critical-type}.
\end{proof}

The classical embeddings and branching underlying this realization
are discussed in Section~\ref{sec:introduction}.

For either sign, let
$S_\pm(h_1,h_2)=S_{\lambda_\pm,\lambda_\mp,2/3}(h_1,h_2)$
be the Koecher operators of Definition~\ref{def:koecher}, in the notation of Corollary~\ref{cor:koecher}, where $h_1,h_2\in\WC$. Then, for any $h\in\WC$ they are given by
\begin{equation}\label{eq:koecher-op-critical}
S_{\pm}(h_1,h_2)h=\langle h,h_2\rangle_W h_1-2F_{\mp}\bigl(F_{\pm}(h_1,h),h_2\bigr)+\frac13\langle h_1,h_2\rangle_W h.
\end{equation}

By Proposition~\ref{prop:double-adjoint} and Corollary~\ref{cor:koecher}, all these operators lie in $\stab(N_\pm)$.
Proposition~\ref{prop:koecher-jump} shows that their span contains a
$64$-dimensional submodule. In Corollary~\ref{cor:koecher-span} we show that the
span is an ideal of $\mathfrak g_\pm$, so the simplicity established in
Theorem~\ref{thm:critical-type} makes it all of $\mathfrak g_\pm$.

\begin{corollary}
\label{cor:koecher-span}
For either sign, let
$S_{\pm}(h_1,h_2)$ be the Koecher operators at the critical parameters,
as in \eqref{eq:koecher-op-critical}. Then
\[
\mathfrak g_\pm
=
\Span_\C\{S_\pm(h_1,h_2):h_1,h_2\in\WC\}
\cong\mathfrak e_6.
\]
Thus the quadratic adjoints $F_\pm$ and $F_\mp$
give an explicit realization of $\mathfrak e_6$
on $\WC$.
\end{corollary}

\begin{proof}
Let $G$ be the identity component of the
stabilizer of $N_\pm$ in $GL(\WC)$, and write
$\gamma^\vee=(\gamma^\tadj)^{-1}$ for $\gamma\in G$.
Invariance of $N_\pm$ gives
\[
F_\pm(\gamma h)=\gamma^\vee F_\pm(h).
\]
The double-adjoint identity
\[
F_\mp(F_\pm(h))=\tfrac23N_\pm(h)h
\]
then implies
\[
F_\mp(\gamma^\vee u)=\gamma F_\mp(u).
\]
Indeed, this follows first for $u=F_\pm(h)$.
By Lemma~\ref{lem:critical-adjoint-dominance}, the identity then
extends to all of $\WC$, since both sides are polynomial in $u$.

Polarizing these identities and substituting
into the formula for $S_\pm$ gives
\[
S_\pm(\gamma h_1,\gamma^\vee h_2)
=\gamma S_\pm(h_1,h_2)\gamma^{-1}.
\]
Consequently, the linear span of the Koecher
operators is invariant under the adjoint action
of $G$, hence is an ideal of $\mathfrak g_\pm$.
It is nonzero by Proposition~\ref{prop:koecher-jump}.
Since $\mathfrak g_\pm$ is simple, this span
equals $\mathfrak g_\pm$.
\end{proof}

\begin{remark}
\label{rem:koecher-tkk}
The two critical stabilizers realize dual actions under the $W$-pairing:
\[
\mathfrak g_\mp=\{-X^{\tadj}:X\in\mathfrak g_\pm\}.
\]
Indeed, $X\mapsto-X^{\tadj}$ is the Lie algebra isomorphism for which
$\langle Xh,k\rangle_W+\langle h,-X^{\tadj}k\rangle_W=0$.
Thus the action of $\mathfrak g_\mp$ corresponds to the contragredient
of the action of $\mathfrak g_\pm$.

The $W$-pairing identifies $\WC$ with its dual
as $G_2$-modules, but not as
$\mathfrak e_6$-modules. Accordingly, although
the Koecher operators are expressed as a
$G_2$-equivariant map
\[
\WC\otimes\WC\longrightarrow \mathfrak g_\pm,
\]
their natural $\mathfrak e_6$-equivariant
interpretation is
\[
\WC\otimes\WC^*
\longrightarrow \mathfrak e_6,
\qquad
h_1\otimes h_2^\flat\longmapsto S_\pm(h_1,h_2),
\]
where $h_2^\flat=\langle h_2,\cdot\rangle_W$.
In terms of the minuscule representations
identified in Theorem~\ref{thm:critical-type},
this map realizes, up to a nonzero scalar,
the projection onto the adjoint \(78\)-dimensional summand in the decomposition into irreducible \(\esix\)-modules:
\[
\mathbf{27}\otimes\mathbf{27}^*
=\mathbf1\oplus\mathbf{78}\oplus\mathbf{650}.
\]

These operators are related to the Tits--Kantor--Koecher construction
in the theory of Jordan algebras
\cite[2.25--2.26 and 14.31]{springer-1973-jordan}.
\end{remark}

We have determined the stabilizer of every member of the projective
pencil and characterized the two critical parameters in the affine
family by the double-adjoint identity.
In Section~\ref{sec:e6-albert-norm} we identify the full stabilizer groups,
prove uniqueness of their invariant cubics, and identify the
critical cubics with the complex Albert norm.
\section{$E_6$, the minuscule $\mathbf{27}$, and the Albert norm}
\label{sec:e6-albert-norm}

In Section~\ref{sec:exceptional-pencil}, we identified the critical
cubics $N_\pm=C_0+\lambda_\pm C_\varphi$ and showed that their
infinitesimal stabilizers $\mathfrak g_\pm=\stab(N_\pm)$ are simple
Lie algebras of type $E_6$, each acting on $\WC$ through one of the two
minuscule representations $V(\omega_1^{E_6})$ and
$V(\omega_6^{E_6})$ (Theorem~\ref{thm:critical-type}).
We now determine the full stabilizer groups, establish the
uniqueness of their invariant cubics, and identify $N_\pm$
with complex Albert norms.

The distinction between a cubic norm and a unital algebra
is important here. The $G_2$-structure canonically determines
the two critical norms, but a corresponding Jordan product
requires a choice of unit of nonzero norm. We develop the
resulting Albert algebra structures, describe their $F_4$
stabilizers and residual $G_2$ symmetries, and construct
the real forms needed in Section~\ref{sec:reductions}.

Our approach combines direct arguments in the $G_2$ model
with classical results about the Albert algebra.
Section~\ref{subsec:critical-group} determines the full
stabilizer groups using connected algebraic groups and outer
automorphisms of $\mathfrak e_6$ and $\gtwoC$.
Sections~\ref{subsec:albert-norm} and~\ref{subsec:basepoints}
establish the norm identification and develop the corresponding
unital algebras, including transitivity on nonzero norm levels
and their $F_4$ stabilizers.
Section~\ref{subsec:g2-at-basepoint} uses the closed-orbit
theorem and the classical stabilizers of a unit vector and
an associative plane to analyze the residual $G_2$ symmetry.
Section~\ref{subsec:yukawa-comparison} distinguishes the pointwise Yukawa
cubic from the critical Albert norms.
Finally, Sections~\ref{subsec:real-form} and
\ref{subsec:euclidean-albert} construct compact $E_6$ forms
using Cartan decomposition and identify the Euclidean Albert
algebra associated with an $SU(3)$-type basepoint using the
Faraut--Kor\'anyi classification. No further classification
of simple Lie algebras is required.

\subsection{The stabilizer group and the invariant cubic}
\label{subsec:critical-group}
The Lie algebra determines the identity component of the stabilizer.
To determine the full group, we use the minuscule representation to exclude
additional components, then describe the common stabilizer of two distinct members of the pencil to establish
uniqueness of their invariant cubics.

Let
\[
G_2^{\C}
=\{O\in SO(\VC,g_{\C}):O^*\varphi_{\C}=\varphi_{\C}\}
\]
be the connected complex algebraic group of type $G_2$, with
Lie algebra $\gtwoC$.

Regard $\WC=S_0^2\VC^*$ as the space of traceless symmetric
$2$-tensors, with the complex-bilinear extensions of the conventions of
Section~\ref{sec:deformation-representation}. The natural left action is
\[
\rho(O)h=(O^{-1})^*h,
\]
as in~\eqref{eq:rho-group}.

Under the metric identification~\eqref{eq:metric-identification}, this is conjugation, $g^{-1}(\rho(O)h)=O(g^{-1}h)O^{-1}$, and its differential is the action $\rho(M)$ of~\eqref{eq:rho-action}. This action preserves the $W$-pairing and both cubics $C_0$ and $C_\varphi$. Since its differential is injective (Section~\ref{sec:deformation-representation}) and $G_2^{\C}$ has trivial center, the representation is faithful. Its image is a closed connected algebraic subgroup of $GL(\WC)$ with Lie algebra $\rho(\gtwoC)$.

As in Section~\ref{sec:introduction}, $E_6(\C)$ denotes the simply connected complex group of type $E_6$. The undecorated $\mathfrak e_6$ and $\mathfrak f_4$ denote the complex Lie algebras, while $\gtwo$ is the compact real Lie algebra and $\gtwoC$ its complexification. To pass from $\mathfrak g_\pm$ to the full stabilizer group, we use four standard facts over $\C$ \cite{garibaldi-2023-albert,milne-2017-algebraic,podder-2023-k-ring}:
\begin{enumerate}
\item[(i)] a connected algebraic subgroup of $GL(\WC)$ is determined by its Lie algebra;
\item[(ii)] a connected semisimple group with Lie algebra $\mathfrak e_6$ is a central quotient of $E_6(\C)$;
\item[(iii)] the center of $E_6(\C)$ is the group $\mu_3$ of cube roots of unity, and it acts nontrivially, by scalars, on a minuscule $\mathbf{27}$;
\item[(iv)] $\operatorname{Out}(\mathfrak e_6)\cong\mathbb Z/2$, and its nontrivial element exchanges $V(\omega_1^{E_6})$ and $V(\omega_6^{E_6})$.
\end{enumerate}

\begin{theorem}
  \label{thm:critical-group}
  For either sign, the stabilizer of the critical cubic is
  \begin{equation*}
    G_\pm:=\Stab_{GL(\WC)}(N_\pm)\cong E_6(\C).
  \end{equation*}
  This is the faithful minuscule realization of the simply connected group.
  It is connected, and its center $\mu_3$ acts by scalar cube roots of unity.
  The stabilizer of the cubic line $\C N_\pm$ is
  $\widehat G_\pm=G_\pm\cdot\C^\times$, where $\C^\times$ acts by scalars;
  its Lie algebra is $\mathfrak g_\pm\oplus\C\Id_W$ and its dimension is $79$.
\end{theorem}

A classical construction of $E_6$ as the isometry group of a trilinear
form on its $27$-dimensional module is given by Aschbacher
\cite[\S5]{aschbacher-1987-27}.
Under the identification of Theorem~\ref{thm:albert-norm},
$G_\pm$ and $\widehat G_\pm$ become, respectively, the
norm-isometry and norm-similarity groups of the complex
Albert algebra
\cite[12.3(iii), 14.20]{springer-1973-jordan},
\cite[\S15, p.~29]{garibaldi-2023-albert}.
The following proof derives their structure directly from
Theorem~\ref{thm:critical-type}.

\begin{proof}
Write $G_\pm^\circ$ for the identity component.
Since its Lie algebra is $\mathfrak g_\pm\cong\mathfrak e_6$,
the group $G_\pm^\circ$ is connected and almost simple of type
$E_6$, hence a central quotient of $E_6(\C)$.
The induced representation of $E_6(\C)$ on
$\WC$ is the $27$-dimensional minuscule representation.
Its kernel is contained in the center $\mu_3$, which acts
nontrivially and therefore faithfully.
Thus $G_\pm^\circ\cong E_6(\C)$, with center
$\mu_3\Id_W$.

Now consider the full group. Let $\gamma\in G_\pm$. Conjugation by $\gamma$ preserves
$\mathfrak g_\pm$ and therefore induces an automorphism
$\alpha$ of $\mathfrak g_\pm\cong\mathfrak e_6$.
Since $\gamma$ acts on $\WC$, this automorphism must preserve
the isomorphism class of the $27$-dimensional
$\mathfrak e_6$-module $\WC$.
However, the nontrivial outer automorphism of $\mathfrak e_6$
exchanges the two nonisomorphic minuscule representations
$V(\omega_1^{E_6})$ and $V(\omega_6^{E_6})$.
Thus $\alpha$ must be inner.

Choose $\delta\in G_\pm^\circ$ inducing the same automorphism.
Then $\delta^{-1}\gamma$ commutes with $\mathfrak g_\pm$.
Since $\WC$ is irreducible under $\mathfrak g_\pm$
(it is already irreducible under $\rho(\gtwoC)$), Schur's lemma gives
$\delta^{-1}\gamma=c\Id_W$ for some $c\in\C^\times$.
Both $\delta$ and $\gamma$ preserve the nonzero cubic $N_\pm$,
so $c^3=1$.
But $\mu_3\Id_W$ is the center of $G_\pm^\circ$, as noted above.
Hence $\delta^{-1}\gamma\in G_\pm^\circ$, and consequently
$\gamma\in G_\pm^\circ$.
Therefore $G_\pm=G_\pm^\circ$.

  \smallskip
  \emph{The cubic line.}
  If $L\in GL(\WC)$ satisfies $N_\pm(Lh)=cN_\pm(h)$ for all $h$, choose $t\in\C^\times$ with $t^3=c$.
  Then $t^{-1}L\in G_\pm$, proving
  $\widehat G_\pm=G_\pm\cdot\C^\times$.
  At the Lie-algebra level, $\Phi_{N_\pm}(\Id_W)=3N_\pm$ gives
  \begin{equation*}
    \Phi_{N_\pm}(X)=cN_\pm
    \quad\Longleftrightarrow\quad
    X-\tfrac c3\Id_W\in\mathfrak g_\pm.
  \end{equation*}
  The sum $\mathfrak g_\pm+\C\Id_W$ is direct because $\Phi_{N_\pm}(\Id_W)=3N_\pm\ne0$, and its dimension is $78+1=79$.
\end{proof}

The uniqueness of the invariant cubic will follow from a description of the
common stabilizer of two distinct members of the pencil.
As before, let $\mathbb F\in\{\R,\C\}$ and $N_{a,b}=aC_0+bC_\varphi$ on
$W_{\mathbb F}$.
All stabilizers below fix the indicated cubics exactly.

\begin{proposition}
  \label{prop:common-stabilizer}
  \renewcommand{\labelenumi}{(\arabic{enumi})}
  \renewcommand{\theenumi}{\arabic{enumi}}
  \leavevmode
  \begin{enumerate}
    \item\label{item:common-stabilizer-lie} If $[a_1:b_1]\ne[a_2:b_2]$ in $\mathbb P^1(\mathbb F)$, that is, $a_1b_2-a_2b_1\ne0$, then $\stab(N_{a_1,b_1})\cap\stab(N_{a_2,b_2})=\rho(\mathfrak g_2^{\mathbb F})$.
    \item\label{item:common-stabilizer-group} The intersection of stabilizers of the two critical cubics is given by
    \begin{equation}
      G_+\cap G_-=\mu_3\times\rho(G_2^{\C}),
      \label{eq:common-stabilizer}
    \end{equation}
    a direct product of algebraic groups, in which $\mu_3$ acts by scalar cube roots of unity and meets $\rho(G_2^{\C})$ trivially.
    Relative to the original real form $W\subset\WC$, for either sign,
    \[
      G_\pm\cap GL(W)=\rho(G_2),\qquad
      \mathfrak g_\pm\cap\End(W)=\rho(\gtwo),
    \]
    where real maps are identified with their complex-linear extensions.
  \end{enumerate}
\end{proposition}

\begin{proof}
  \emph{(\ref{item:common-stabilizer-lie}) Lie algebras.}
  The inclusion of $\rho(\mathfrak g_2^{\mathbb F})$ is
  Lemma~\ref{lem:common-g2}.
  Conversely, linearity of $\Phi_N(X)$ in $N$ and invertibility of the
  coefficient matrix imply that an $X$ annihilating both members annihilates
  $C_0$ and $C_\varphi$.
  By Propositions~\ref{prop:trace-cubic-stabilizer}
  and~\ref{prop:phi-cubic-stabilizer},
  \begin{equation*}
    \ker\Phi_{C_0}\cap\ker\Phi_{C_\varphi}
    =\rho(\lieso(7,\mathbb F))\cap\rho(\mathfrak g_2^{\mathbb F})
    =\rho(\mathfrak g_2^{\mathbb F}).
  \end{equation*}

  \emph{(\ref{item:common-stabilizer-group}) Groups.} Let $H=G_+\cap G_-$. By (\ref{item:common-stabilizer-lie}), $\Lie(H)=\rho(\gtwoC).$ Since $\rho(G_2^\C)\subset H$ is connected and has the same Lie algebra, it is the identity component $H^\circ$. \par

  Every $L\in H$ normalizes $H^\circ$ and thus induces an
automorphism of $\gtwoC$. Since $G_2$ has no nontrivial diagram
automorphisms, every such automorphism is inner. Hence there
is an $O\in G_2^\C$ such that $\rho(O)^{-1}L$ commutes with
$\rho(\gtwoC)$. By irreducibility of $\WC$ and Schur's lemma,
\[
L=\omega\rho(O),\qquad\omega\in\C^{\times}.
\]
Since $L$ and $\rho(O)$ preserve $N_+$, we have $\omega^3=1$.
Conversely, every such product preserves both critical cubics.
The product is direct: if $\omega\Id_{W}\in\rho(G_2^{\C})$,
preservation of the $W$-pairing gives $\omega^2=1$, and
hence $\omega=1$. Therefore
\[
G_+\cap G_-=\mu_3\times\rho(G_2^\C).
\]

Finally, both factors are preserved by complex conjugation.
Since $\mu_3\cap\R=\{1\}$ and the real points of
$\rho(G_2^\C)$ form $\rho(G_2)$, the elements preserving
$W$ are precisely $\rho(G_2)$.
Since $\overline{N_+}=N_-$, a real element of either critical
group lies in their intersection, giving the asserted real points of
each group. Likewise, conjugation exchanges $\mathfrak g_+$ and
$\mathfrak g_-$. A real element of either lies in
$\rho(\gtwoC)$ by part~(\ref{item:common-stabilizer-lie}), and the real
elements of $\rho(\gtwoC)$ are exactly $\rho(\gtwo)$.
\end{proof}

\begin{remark}
  The group $\mu_3\times\rho(G_2^{\C})$ is the normalizer of the irreducible
  subgroup $\rho(G_2^{\C})$ in either $G_\pm$; for the compact form, see Wolf's classification
  \cite[Theorem~11.1 and its table]{wolf-1968-geometry}.
  Proposition~\ref{prop:common-stabilizer} identifies that known group as the
  common stabilizer of the two critical cubics, using the standard normalizer
  argument \cite[p.~34]{garibaldi-2015-simple}.
\end{remark}

Part~(\ref{item:common-stabilizer-lie}) of Proposition~\ref{prop:common-stabilizer} now gives uniqueness
of the invariant cubic.

\begin{proposition}
  \label{prop:unique-cubic}
  For either sign, $(S^3\WdualC)^{\mathfrak g_\pm}=\C N_\pm$; that is, the cubics invariant under the critical $\mathfrak e_6$ form a complex line.
\end{proposition}

\begin{proof}
  Every $\mathfrak g_\pm$-invariant cubic is also
  $\rho(\gtwoC)$-invariant, and hence lies in the pencil by
  Lemma~\ref{lem:invariant-cubic-plane}.
  If such a cubic $P$ were independent of $N_\pm$, then
  Proposition~\ref{prop:common-stabilizer}\,(\ref{item:common-stabilizer-lie}) would give
  \begin{equation*}
    \mathfrak g_\pm
      \subseteq\stab(N_\pm)\cap\stab(P)=\rho(\gtwoC).
  \end{equation*}
  This is impossible, since $\dim\mathfrak g_\pm=78>14$.
\end{proof}

For the classical minuscule $\mathbf{27}$ of $E_6$,
uniqueness of the invariant cubic is also known from
invariant theory
\cite[\S4, Propositions~3 and~12; \S5, Proposition~47]
{sato-1977-classification},
\cite[12.4(ii)]{springer-1973-jordan}.
Proposition~\ref{prop:unique-cubic} establishes this
uniqueness directly in the present $G_2$ model.

\subsection{The Albert norm}
\label{subsec:albert-norm}

The complex Albert algebra carries a minuscule $\mathbf{27}$ of $E_6$ and an
invariant cubic norm.
Using the uniqueness just proved, we will identify that norm with $N_\pm$.
We first recall how a norm, together with a chosen unit, determines the
Jordan algebra.

Recall that a \emph{Jordan algebra} over $\R$ or $\C$ has a commutative
bilinear product $\circ$ satisfying the Jordan identity
\begin{equation}
  \label{eq:jordan-identity}
  x\circ\bigl((x\circ x)\circ y\bigr)
    =(x\circ x)\circ(x\circ y).
\end{equation}
Jordan algebras over these fields are power-associative: the subalgebra
generated by one element is associative, so the powers of that element are
unambiguous \cite[\S3.10]{petersson-2019-survey}.
A Jordan algebra is \emph{special} if it embeds in an associative algebra endowed with
the symmetrized product $x\circ y=(xy+yx)/2$, and \emph{exceptional}
otherwise \cite[\S\S1.1--1.3]{petersson-2019-survey}.
Let $\O_{\C}=\O\otimes_{\R}\C$ be the complexified octonion algebra.
Octonion conjugation extends complex-linearly, with norm
$n_{\O}(a)=a\bar a$ and trace $t_{\O}(a)=a+\bar a$.
A particular case of a Jordan algebra is the complex Albert algebra $\mathbb A=H_3(\O_{\C})$ that consists of matrices
\begin{equation*}
x=\begin{pmatrix}
\alpha_1&a_3&\bar a_2\\
\bar a_3&\alpha_2&a_1\\
a_2&\bar a_1&\alpha_3
\end{pmatrix},\qquad \alpha_i\in\C,\quad a_i\in\O_{\C},
\end{equation*}
with unit $1$ and product $x\circ y=(xy+yx)/2$
\cite[\S1.5, p.~224]{petersson-2019-survey}; see also Baez
\cite[\S3]{BaezOcto}. Here Hermitian symmetry uses octonion conjugation
and does not conjugate the complex coefficients. Its cubic norm is
\begin{equation*}
N_{\mathbb A}(x)=\alpha_1\alpha_2\alpha_3
 -\sum_{i=1}^3\alpha_i n_{\O}(a_i)+t_{\O}((a_1a_2)a_3)
\end{equation*}
\cite[Theorem~15, p.~244]{petersson-2019-survey}.
The trace makes the last term independent of its parenthesization.
The algebra is simple \cite[Theorem~1, p.~224]{petersson-2019-survey}.
The Albert algebra $\mathbb A$ is exceptional
\cite[Theorem~1 and Corollary~18(c)]{petersson-2019-survey}.

The algebra product is recovered from $(N_{\mathbb A},1)$ by the classical construction
of Springer and McCrimmon \cite[Theorem~5, p.~504]{mccrimmon-1969-freudenthal}.
We recall the conventions needed below, following Petersson
\cite[\S5.2, Theorem~13 and \S5.3]{petersson-2019-survey}.

For a cubic form $N$ on a finite-dimensional complex vector space and a point
$e$ with $N(e)=1$, define the \emph{based trace form} by
\begin{equation}
  \label{eq:based-trace-form}
  \begin{aligned}
    \langle x,y\rangle_e
      &=-\partial_x\partial_y\log N\big|_e\\
      &=dN|_e(x)\,dN|_e(y)-d^2N|_e(x,y)\\
      &=9N(e,e,x)N(e,e,y)-6N(e,x,y).
  \end{aligned}
\end{equation}
The last expression uses the symmetric trilinear polarization of $N$,
normalized as in~\eqref{eq:cubic-polarization}, and makes the
definition independent of any choice of logarithm.
Substituting $y=e$ in~\eqref{eq:based-trace-form} gives the associated linear trace:
\begin{equation}
  \label{eq:based-linear-trace}
  \tr_e(x):=\langle x,e\rangle_e=3N(e,e,x),
  \qquad \text{with } \tr_e(e)=3.
\end{equation}
To express the based algebra operations and construct the real structure in
Section~\ref{subsec:euclidean-albert}, we represent the based trace form
using the fixed $W$-pairing. On $\WC$, the \emph{based trace operator} $T_e$ is defined by
\begin{equation}
  \label{eq:based-trace-operator}
  \langle T_ex,y\rangle_W=\langle x,y\rangle_e.
\end{equation}
Explicit formulas for $T_e$ and its inverse at a critical cubic are given in
Lemma~\ref{lem:based-trace-operator}.

The pair $(N,e)$ is \emph{nonsingular} when the based trace form is
nondegenerate \cite[\S5.3]{petersson-2019-survey}.
In that case there is a unique quadratic map $x\mapsto x^\sharp$,
the \emph{based adjoint}, such that
$\langle x^\sharp,y\rangle_e=dN|_x(y)$.
We use the unhalved polarization
\begin{equation*}
  x\mathbin{\sharp}y=(x+y)^\sharp-x^\sharp-y^\sharp,
  \qquad x\mathbin{\sharp}x=2x^\sharp.
\end{equation*}
A nonsingular pair is a \emph{cubic norm structure} when it satisfies the
following identities:
\begin{equation}
  \begin{gathered}
    N(e)=1,\qquad e^\sharp=e,\qquad x^{\sharp\sharp}=N(x)\,x,\\
    dN|_x(y)=\langle x^\sharp,y\rangle_e,\qquad e\mathbin{\sharp}x=\tr_e(x)\,e-x.
  \end{gathered}
  \label{eq:cubic-norm-structure}
\end{equation}

The first and fourth identities hold by construction here; the remaining
conditions are the basepoint, adjoint and unit identities.
These are polynomial identities, and are required to hold after scalar extension as well.  Over $\C$, we can also define a bilinear product:

\begin{equation}
\label{eq:circ-e-Jordan-prod}
  \begin{aligned}
    x\circ_e y=\tfrac12\bigl(&x\mathbin{\sharp}y+\tr_e(x)y+\tr_e(y)x\\
      &-\bigl(\tr_e(x)\tr_e(y)-\langle x,y\rangle_e\bigr)e\bigr),
  \end{aligned}
\end{equation}
as in \cite[Proposition~3]{krutelevich-2007-jordan}.
By the classical construction theorem, this product satisfies the Jordan
identity~\eqref{eq:jordan-identity} and has unit $e$
\cite[Theorem~13]{petersson-2019-survey}.
The resulting algebra $\mathfrak J(N,e)$, equipped with its norm $N$,
is a \emph{cubic Jordan algebra}; its based trace form is nondegenerate.
We write $x^{\circ2}=x\circ_e x$ and, more generally, $x^{\circ n}$ for its
Jordan powers, with $x^{\circ0}=e$.
The same construction theorem gives the characteristic identity
\begin{equation}
  \label{eq:cubic-characteristic-identity}
  x^{\circ3}-\tr_e(x)x^{\circ2}
    +\tr_e(x^\sharp)x-N(x)e=0.
\end{equation}
Here $\tr_e(x^\sharp)$ is the quadratic trace.
An element $x$ is invertible exactly when $N(x)\ne0$, with Jordan inverse
$x^{-1}=N(x)^{-1}x^\sharp$
\cite[Theorem~13, p.~241]{petersson-2019-survey}.

\phantomsection
\label{par:cubic-degree-determinant}
For a finite-dimensional unital Jordan algebra, its \emph{rank}, also
called its degree here, is the largest degree of an element's minimal
polynomial in the Jordan powers. On the open dense set of elements
attaining this degree $r$, the minimal polynomial is the
\emph{generic characteristic polynomial}
$t^r-a_1(x)t^{r-1}+\cdots+(-1)^r a_r(x)$.
Its coefficients extend polynomially to the whole algebra; the
\emph{Jordan determinant} is $a_r(x)$
\cite[\S II.2, pp.~27--30]{faraut-1994-analysis}.
Thus \eqref{eq:cubic-characteristic-identity} bounds the rank by three;
when the rank is three, the determinant equals $N$ by comparison of
constant coefficients. In $\mathbb A$, diagonal matrices with three distinct
entries have minimal polynomial of degree three. Thus $\mathbb A$ has
rank three, and $N_{\mathbb A}$ is its Jordan determinant.

The nonsingularity and adjoint identity above are also the admissibility
conditions in \cite[Definition~2]{krutelevich-2007-jordan}.
Albert algebras are nonsingular cubic Jordan algebras
\cite[Corollary~18(a)]{petersson-2019-survey}.

Associated with this product is the quadratic operator
$\mathcal U_x$, defined by
\begin{equation}\label{eq:springer-u-operator}
  \mathcal U_x y
  =2x\circ_e(x\circ_e y)-x^{\circ2}\circ_e y.
\end{equation}
Equivalently, in terms of the cubic norm and its adjoint,
\[
  \mathcal U_x y
  =\langle x,y\rangle_e\,x
   -(x^\sharp)\mathbin{\sharp}y.
\]
The equivalence of these expressions is part of the
classical cubic norm construction
\cite[Theorem~13]{petersson-2019-survey}.
\begin{lemma}[Quadratic operators]
\label{lem:quadratic-operators}
Let $(N,e)$ be a nonsingular cubic norm structure,
with associated Jordan algebra $\mathfrak J(N,e)$.
For every $x,y$,
\begin{equation}\label{eq:quadratic-operator-properties}
  \mathcal U_x e=x^{\circ2},
  \qquad
  N(\mathcal U_x y)=N(x)^2N(y).
\end{equation}
If $N(x)\ne0$, then $\mathcal U_x$ is invertible,
with inverse $\mathcal U_{x^{-1}}$.

Consequently, if $N(x)^2=1$, then $\mathcal U_x$
preserves $N$ and carries the unit $e$ to $x^{\circ2}$.
\end{lemma}

\begin{proof}
Since $e$ is the unit of the Jordan algebra,
\eqref{eq:springer-u-operator} gives
\[
  \begin{aligned}
  \mathcal U_x e
    &=2x\circ_e(x\circ_e e)-x^{\circ2}\circ_e e\\
    &=2x^{\circ2}-x^{\circ2}=x^{\circ2}.
  \end{aligned}
\]
The second identity in
\eqref{eq:quadratic-operator-properties} is the
standard norm-composition identity
\cite[Theorem~13]{petersson-2019-survey}.
For invertible $x$, the inverse identity
$\mathcal U_x^{-1}=\mathcal U_{x^{-1}}$
is the standard identity for quadratic operators
\cite[33.11]{garibaldi-2024-albert}.
The final assertion follows immediately.
\end{proof}

These constructions are preserved by linear changes of variables
that preserve the cubic norm and basepoint. In the next lemma, the
operations on the left are defined by $(N',e')$, and those on the right
by $(N,e)$.

\begin{lemma}[Change of variables]
\label{lem:change-of-variables}
Let $N$ and $N'$ be cubic forms with points $e$ and $e'$, where $N(e)=1$,
and let $L$ be a linear isomorphism satisfying
\[
N'(Lx)=N(x),\qquad L(e)=e'.
\]
Then $N'(e')=1$ and
\begin{equation}\label{eq:L-based-trace}
\langle Lx,Ly\rangle_{e^{\prime}}=\langle x,y\rangle_{e}\qquad\text{and}\qquad\tr_{e^{\prime}}(Lx)=\tr_{e}(x),
\end{equation}
so $(N',e')$ is nonsingular if and only if $(N,e)$ is. In that case
\begin{equation}\label{eq:L-sharp}
(Lx)^{\sharp'}=L(x^\sharp),
\end{equation}
the same identity holds for the polarized adjoints, and $(N',e')$
is a cubic norm structure if and only if $(N,e)$ is. If it is, then
\begin{equation}\label{eq:L-prod}
L(x\circ_e y)=Lx\circ_{e'}Ly,
\end{equation}
and $L$ is an isomorphism of the corresponding Jordan
algebras $\mathfrak J(N,e)$ and $\mathfrak J(N',e')$.
\end{lemma}

\begin{proof}
First, $N'(e')=N'(Le)=N(e)=1$, so both based trace forms are defined.
Taking the second derivative of the identity $N'(Lx)=N(x)$ at $e$,
and using $L(e)=e'$, gives \eqref{eq:L-based-trace}.
Since $L$ is an isomorphism, one based trace form is nondegenerate
exactly when the other is.
In that case, the defining property of the based adjoint and
nondegeneracy of the based trace form $\langle\cdot,\cdot\rangle_{e'}$ give \eqref{eq:L-sharp}.
The polarized version follows by polarization.
The identities~\eqref{eq:cubic-norm-structure} are polynomial identities
in the norm, the basepoint, the based trace, the based trace form and the adjoint,
and $L$ intertwines each of these by \eqref{eq:L-based-trace} and
\eqref{eq:L-sharp}; this remains true after scalar extension. Hence they hold for
$(N',e')$ exactly when they hold for $(N,e)$.
Substituting these identities into the product formula~\eqref{eq:circ-e-Jordan-prod}
gives \eqref{eq:L-prod}.
Since $L$ is a linear isomorphism, it is therefore an isomorphism of Jordan algebras
\cite[Corollary~18(b)]{petersson-2019-survey}.
\end{proof}

\begin{theorem}
  \label{thm:albert-norm}
  For either sign there is a complex-linear isomorphism $\Upsilon\colon\mathbb A\to\WC$ with $\Upsilon^*N_\pm=N_{\mathbb A}$.
  Thus each critical cubic is linearly equivalent over $\C$ to the determinant
  of the complex Albert algebra, with any nonzero scalar factor absorbed into
  the linear map.
\end{theorem}

This equivalence is over $\C$; the critical members themselves are not real
(Theorem~\ref{thm:stabilizer-trichotomy}). Note that this theorem establishes linear equivalence but does not specify a preferred coordinate map $\Upsilon$.

\begin{proof}
  By a classical theorem, the exact linear stabilizer of $N_{\mathbb A}$
is the simply connected complex group of type $E_6$
\cite[\S15, p.~29]{garibaldi-2023-albert}.
As a representation of this group, the $27$-dimensional vector
space $\mathbb A$ is faithful. It is irreducible because $E_6(\C)$ has no
nontrivial irreducible representation of smaller dimension, and is therefore
one of the two nonisomorphic minuscule modules
$V(\omega_1^{E_6})$ and $V(\omega_6^{E_6})$, as in
Theorem~\ref{thm:critical-type}.

By Theorem~\ref{thm:critical-type}, $\WC$ is likewise one of
these two modules for $\mathfrak g_\pm\cong\mathfrak e_6$.
Since the diagram automorphism of $\mathfrak e_6$ exchanges
the two modules, we may choose an identification of the
stabilizer Lie algebras under which the representations
on $\mathbb A$ and $\WC$ are equivalent.
This gives an intertwining linear isomorphism
$\Upsilon'\colon\mathbb A\to\WC$.

By Proposition~\ref{prop:unique-cubic}, the space of invariant
cubics on $\WC$ is one-dimensional. Hence
\[
\Upsilon'^*N_\pm=cN_{\mathbb A}
\]
for some $c\in\C^\times$.
Choose $t\in\C^\times$ with $t^3c=1$ and set
$\Upsilon=t\Upsilon'$. Then $\Upsilon^*N_\pm=t^3cN_{\mathbb A}=N_{\mathbb A}$.

\end{proof}

\subsection{Basepoints and Jordan products}
\label{subsec:basepoints}

Fix either sign and choose a \emph{basepoint}, an element $e\in\WC$ with
$n:=N_\pm(e)\ne0$.
Write $\lambda=\lambda_\pm$ and $\bar\lambda=\lambda_\mp$.
Normalize the cubic at $e$ by
\begin{equation*}
  N_e=N_\pm/n,\qquad N_e(e)=1.
\end{equation*}
The following formulas express $T_e$ and its inverse in terms
of the quadratic adjoints of Definition~\ref{def:quadratic-adjoint}.
\begin{lemma}
  \label{lem:based-trace-operator}
  For $n=N_\lambda(e)\ne0$, the based trace operator for $N_\lambda/n$ is
  \begin{equation}
    \label{eq:based-trace-response}
    T_ex=\frac3n\tr_e(x)F_\lambda(e)-\frac6nF_\lambda(e,x).
  \end{equation}
  It is invertible, with
  \begin{equation}
    \label{eq:based-trace-inverse}
    T_e^{-1}y=\langle e,y\rangle_W e
      -F_{\bar\lambda}(F_\lambda(e),y).
  \end{equation}
\end{lemma}

\begin{proof}
  Differentiating $N_\lambda/n$ gives
  \begin{align*}
    \tr_e(x)&=\frac3n\langle F_\lambda(e),x\rangle_W,\\
    d^2(N_\lambda/n)|_e(x,y)&=\frac6n\langle F_\lambda(e,x),y\rangle_W.
  \end{align*}
  Substitute these in~\eqref{eq:based-trace-form} and use
  \eqref{eq:based-trace-operator} to obtain~\eqref{eq:based-trace-response}.
  Put $f=F_\lambda(e)$. Differentiating
  $F_{\bar\lambda}(F_\lambda(x))=\tfrac23N_\lambda(x)x$
  from Proposition~\ref{prop:double-adjoint} at $x=e$ gives
  \[
    4F_{\bar\lambda}(f,F_\lambda(e,x))
      =2\langle f,x\rangle_W e+\tfrac23nx.
  \]
  Since $F_{\bar\lambda}(f)=\tfrac23ne$ and
  $\langle e,f\rangle_W=n$, substitution of
  \eqref{eq:based-trace-response} yields
  \[
    F_{\bar\lambda}(f,T_ex)=\tfrac3n\langle f,x\rangle_W e-x,
    \qquad
    \langle e,T_ex\rangle_W=\tfrac3n\langle f,x\rangle_W.
  \]
  Multiplying the second identity by $e$ and subtracting the first shows that the right side of
  \eqref{eq:based-trace-inverse} is a left inverse of $T_e$.
  Since $\WC$ is finite-dimensional, it is the inverse.
\end{proof}

There is no rank or projection condition on $e$. In particular, the
geometric basepoints $E_u$ and $E_\Sigma$ of~\eqref{eq:test-tensors}, built
from orthogonal projection tensors and studied in Section~\ref{subsec:g2-at-basepoint}, are convenient examples of possible Jordan units.
The next theorem proves that $(N_e,e)$ is a nonsingular cubic norm structure and
we denote its Jordan algebra by $\mathfrak J_e=\mathfrak J(N_e,e)$.
Table~\ref{tab:based-conventions} collects the pairings and normalizations used
in the explicit operations below. The based adjoint uses the normalized
cubic and the based trace form, whereas $F_\lambda$ uses the fixed
$W$-pairing.

\begin{table}[htbp]
\centering
\small
\renewcommand{\arraystretch}{1.25}
\begin{tabular}{p{0.25\linewidth}p{0.67\linewidth}}
\hline
Object & Convention \\
\hline
Fixed $W$-pairing & $\langle x,y\rangle_W=\tr(xy)$, complex bilinear. \\
Hermitian form & $\langle x,\bar y\rangle_W$, with conjugation fixing $W$. \\
Normalized cubic & $N_e=N_\lambda/n$, where $n=N_\lambda(e)\ne0$. \\
Based trace form & $\langle x,y\rangle_e=-d^2\log N_e|_e(x,y)=\langle T_ex,y\rangle_W$. \\
Based linear trace & $\tr_e(x)=\langle x,e\rangle_e=\frac3n\langle F_\lambda(e),x\rangle_W$. \\
Quadratic adjoint & $F_\lambda=\frac13\nabla N_\lambda$; its polarization $F_\lambda(x,y)$ is halved. \\
Based adjoint & $x^\sharp=\frac3nT_e^{-1}F_\lambda(x)$; its polarization $x\mathbin\sharp y$ is unhalved. \\
\hline
\end{tabular}
\caption{Pairings and adjoints at a basepoint of a critical cubic.}
\label{tab:based-conventions}
\end{table}

\begin{theorem}
  \label{thm:based-algebras}
  \renewcommand{\labelenumi}{(\arabic{enumi})}
  \renewcommand{\theenumi}{\arabic{enumi}}
  \leavevmode
  \begin{enumerate}
    \item\label{item:based-algebras-albert} For every basepoint $e$, $(N_e,e)$ is a nonsingular cubic norm
    structure and $\mathfrak J_e\cong\mathbb A$, so $\mathfrak J_e$ is simple, exceptional and of degree exactly $3$.
    Its adjoint and Jordan product are
    \begin{equation}
      \label{eq:based-albert-operations}
      \begin{aligned}
        h^\sharp
          &=\frac3n\bigl(\langle e,F_\lambda(h)\rangle_W e
              -F_{\bar\lambda}(F_\lambda(e),F_\lambda(h))\bigr),\\
        h\circ_e k
          &=-\frac3nF_{\bar\lambda}(F_\lambda(e),F_\lambda(h,k))\\
          &\quad+\frac12\bigl(\tr_e(h)k+\tr_e(k)h\bigr).
      \end{aligned}
    \end{equation}
    For $\Upsilon$ as in Theorem~\ref{thm:albert-norm}, the point $e_{\mathbb A}=\Upsilon(1)$ has $N_\pm(e_{\mathbb A})=1$, and $\Upsilon$ is an isomorphism $\mathbb A\to\mathfrak J(N_\pm,e_{\mathbb A})$.
    \item\label{item:based-algebras-orbits} For any basepoints $e,e'\in\WC$,
    \begin{equation*}
      e'\in G_\pm\cdot e
      \quad\Longleftrightarrow\quad
      N_\pm(e')=N_\pm(e).
    \end{equation*}
    Equivalently, each nonzero norm level
    $\{h\in\WC:N_\pm(h)=c\}$ is a single $G_\pm\cong E_6(\C)$-orbit.
    Allowing scalar multiplication, all basepoints
    lie in one orbit of $\widehat G_\pm=G_\pm\cdot\C^\times$.
  \end{enumerate}
\end{theorem}

\begin{proof}
  We first construct the algebra at one basepoint, prove transitivity using that algebra, and then transport the construction to an arbitrary basepoint.

  \smallskip \textit{The Albert basepoint.} Let $e_{\mathbb A}=\Upsilon(1)$, where $\Upsilon$ is the norm equivalence of Theorem~\ref{thm:albert-norm}. Then

  \[
  N_\pm(e_{\mathbb A})=N_{\mathbb A}(1)=1.
  \]

  Since $\Upsilon$ preserves the cubic norm and basepoint, \eqref{eq:L-based-trace} identifies the based trace forms, and \eqref{eq:L-sharp} identifies the based adjoints. Consequently, $(N_\pm,e_{\mathbb A})$ is a nonsingular cubic norm structure satisfying \eqref{eq:cubic-norm-structure}. By \eqref{eq:L-prod}, $\Upsilon$ also preserves the Jordan products, giving

  \[
  \mathfrak J(N_\pm,e_{\mathbb A})\cong\mathbb A.
  \]

  \smallskip \textit{Transitivity on nonzero norm levels.} We first prove transitivity on $N_\pm^{-1}(1)$. Let $h\in\WC$ satisfy $N_\pm(h)=1$, and regard $h$ as an element of $\mathfrak J(N_\pm,e_{\mathbb A})$. By the definition of invertibility, $h$ is invertible since $N_{\pm}(h)\neq0$.

  Since Jordan algebras are power-associative, the Jordan powers $h^{\circ j}$ of $h$ span a commutative associative subalgebra $\C[h]$. Here the powers are taken with respect to the product $\circ_{e_{\mathbb A}}$, with $h^{\circ0}=e_{\mathbb A}$, and a polynomial $p(t)=\sum_ja_jt^j$ is evaluated in $\C[h]$ as $p(h)=\sum_ja_jh^{\circ j}$. All polynomial evaluations, powers and minimal polynomials below are taken in this subalgebra, and $\mathcal U_k$ is the quadratic operator of $\mathfrak J(N_\pm,e_{\mathbb A})$.

  Moreover, $h$ admits a square root $k\in\C[h]$. Indeed, by
  \eqref{eq:cubic-characteristic-identity}, its minimal polynomial divides
  $t^3-\tr_{e_{\mathbb A}}(h)t^2+\tr_{e_{\mathbb A}}(h^\sharp)t-1$,
  whose constant term is nonzero. A polynomial square root can therefore be constructed modulo each primary factor by a truncated binomial expansion, and by the Chinese remainder theorem these combine to give $k=p(h)$ with $k^{\circ2}=h$. This is a standard square-root argument used in
the theory of Albert isotopes
\cite[Corollary~31.25]{garibaldi-2024-albert}.

  By Lemma~\ref{lem:quadratic-operators},

  \[
  \mathcal U_k e_{\mathbb A}=k^{\circ2}=h.
  \]

  Applying the norm identity in \eqref{eq:quadratic-operator-properties} gives

  \[
  N_\pm(k)^2=N_\pm(h)=1.
  \]

  Thus $\mathcal U_k$ preserves $N_\pm$, belongs to $G_\pm$,
  and takes $e_{\mathbb A}$ to $h$. This proves transitivity on
  $N_\pm^{-1}(1)$, which is the classical transitivity theorem for
  elements of norm one \cite[Theorem~7]{jacobson-1961-groups}.

  For any $c\in\C^\times$, multiplication by a cube root of $c$ identifies $N_\pm^{-1}(1)$ with $N_\pm^{-1}(c)$. Since scalar multiplication commutes with $G_\pm$, the group acts transitively on every nonzero norm level,
as in the classical Albert setting
\cite[Proposition~12(i)]{krutelevich-2007-jordan}.

  \smallskip \textit{Arbitrary basepoints.} Let $e\in\WC$ with $n=N_\pm(e)\ne0$, and choose $t\in\C^\times$ with $t^3=n$. Then $N_\pm(e/t)=1$, so by transitivity there is $\gamma\in G_\pm$ such that

  \[
  \gamma(e_{\mathbb A})=e/t.
  \]

  Define $L:\WC\to\WC$ by $L(x)=\gamma^{-1}(x/t)$. Since $\gamma$ preserves $N_\pm$, we have

  \[
  L(e)=e_{\mathbb A},
      \qquad
      N_\pm(Lx)=\frac{N_\pm(x)}{n}=N_e(x).
  \]

  Thus $L$ preserves the normalized cubic norm and basepoint. Lemma~\ref{lem:change-of-variables} shows that $(N_e,e)$ is a
  nonsingular cubic norm structure and identifies the Jordan products: by~\eqref{eq:L-prod}, $L(x\circ_e y)=Lx\circ_{e_{\mathbb A}}Ly$. Hence

  \[
  \mathfrak J_e
        \cong\mathfrak J(N_\pm,e_{\mathbb A})
        \cong\mathbb A.
  \]

  In particular, every $\mathfrak J_e$ is simple, exceptional and of degree $3$.

  \smallskip \textit{The algebra operations.}
  The defining adjoint relation and~\eqref{eq:based-trace-operator} give
  $h^\sharp=3T_e^{-1}F_\lambda(h)/n$. Its unhalved polarization satisfies
  \[
    \tfrac12h\mathbin{\sharp}k=\tfrac3nT_e^{-1}F_\lambda(h,k).
  \]
  The based trace form~\eqref{eq:based-trace-form} gives
  \[
    \tr_e(h)\tr_e(k)-\langle h,k\rangle_e
      =\tfrac6n\langle F_\lambda(e,h),k\rangle_W.
  \]
  Substituting these in~\eqref{eq:circ-e-Jordan-prod}, we obtain
  \[
    h\circ_e k=\tfrac3nT_e^{-1}F_\lambda(h,k)
      +\tfrac12\bigl(\tr_e(h)k+\tr_e(k)h\bigr)
      -\tfrac3n\langle F_\lambda(e,h),k\rangle_W e.
  \]
  Now use~\eqref{eq:based-trace-inverse}. Symmetry of the polarized cubic gives
  $\langle e,F_\lambda(h,k)\rangle_W=\langle F_\lambda(e,h),k\rangle_W$,
  so the two scalar multiples of $e$ cancel in the product.
  The same inverse formula applied to the adjoint proves both identities
  in~\eqref{eq:based-albert-operations}.

  Finally, since $G_\pm$ preserves $N_\pm$ and acts transitively on each nonzero norm level, the value of $N_\pm$ completely determines the $G_\pm$-orbits on $\{N_\pm\ne0\}$. Allowing scalar multiplication, $\widehat G_\pm=G_\pm\cdot\C^\times$ acts transitively on the full set of basepoints.
\end{proof}

Fix any basepoint $e_0$ with $n_0=N_\pm(e_0)\ne0$.
For another basepoint $e$, put $n=N_\pm(e)$ and choose
$t\in\C^\times$ with $t^3=n/n_0$.
For $\gamma\in G_\pm$, write $L=t^{-1}\gamma$.

\begin{corollary}[Transport to a model basepoint]
\label{cor:basepoint-transport}
\label{par:basepoint-transport}
There is $\gamma\in G_\pm$ with $\gamma e=t e_0$ whose associated norm
similarity $L$ satisfies
\begin{equation}
  \label{eq:basepoint-transport}
  L(e)=e_0,\qquad
  N_{e_0}(Lx)=N_e(x),\qquad
  L(x\circ_e y)=Lx\circ_{e_0}Ly.
\end{equation}
\end{corollary}
\begin{proof}
Theorem~\ref{thm:based-algebras}\,(\ref{item:based-algebras-orbits})
applies because $N_\pm(e)=N_\pm(te_0)$.
Lemma~\ref{lem:change-of-variables} gives the product identity.
\end{proof}
Thus any chosen model basepoint yields, after a scalar adjustment, the
same abstract unital Albert algebra as every other basepoint.
The scalar adjustment is essential unless $n=n_0$: the exact group
$G_\pm$ itself preserves the value of $N_\pm$.
Neither $t$ nor $\gamma$ is canonically determined.
The transport is pointwise over $\C$; it need not preserve the original
real $G_2$-module or assemble into a global gauge transformation on a manifold.
In particular, $\gamma$ need not normalize the embedded $\rho(G_2)$,
so transport need not preserve the residual $G_2$-stabilizer.

\begin{remark}
\label{rem:rank-one-normal-form}
In~\eqref{eq:basepoint-transport}, take $e_0=E_u$ for a real unit vector
$u$ (a basepoint by Proposition~\ref{ex:structured-basepoints}), so $\gamma e=tE_u$. With $v=su$ and $s^2=t$, the model is
\begin{equation*}
tE_u=E_v=7(P_v)_0=7P_v-g(v,v)g,\qquad g(v,v)\ne0.
\end{equation*}
Thus every basepoint is $G_\pm$-equivalent to the traceless part of a
rank-one symmetric tensor. The vector $v$ need not
be real or of unit length. Requiring the target to be $E_u$ itself forces
$N_\pm(e)=N_\pm(E_u)$.
\end{remark}

\paragraph{Changing the unit: isotopy.}
For an invertible $p\in\mathbb A$, the $p$-isotope $\mathbb A^{(p)}$
has unit $p^{-1}$ and quadratic operators
$\mathcal U_x^{(p)}=\mathcal U_x\mathcal U_p$
\cite[\S3.11, p.~232]{petersson-2019-survey}.
Over $\C$, its product is recovered from these data: first
$x^{\circ2}_{(p)}=\mathcal U_x^{(p)}p^{-1}$, then polarize this square.
Its norm is $N_{\mathbb A^{(p)}}=N_{\mathbb A}(p)N_{\mathbb A}$,
and it is the Jordan algebra~\eqref{eq:circ-e-Jordan-prod} of the cubic norm structure
$(N_{\mathbb A}(p)N_{\mathbb A},p^{-1})$
\cite[Theorem~2, pp.~500--501]{mccrimmon-1969-freudenthal},
\cite[33.11--33.12, pp.~326--327]{garibaldi-2024-albert}.
With $a=\Upsilon^{-1}(e)$ and $p=a^{-1}$, these data are
$(N_{\mathbb A}/N_{\mathbb A}(a),a)$, so
Lemma~\ref{lem:change-of-variables} identifies
$\mathbb A^{(a^{-1})}$ with $\mathfrak J_e$ via $\Upsilon$.
Thus Theorem~\ref{thm:based-algebras}\,(\ref{item:based-algebras-albert})
realizes the classical isomorphism of complex Albert isotopes
\cite[Corollary~31.25]{garibaldi-2024-albert} in the geometric cubic model.

Fixing the norm and a basepoint is equivalent to fixing the resulting
unital Albert algebra. Its symmetry group is of type $F_4$. Recall the based linear trace
$\tr_e$ of~\eqref{eq:based-linear-trace}.
Since $\tr_e(e)=3$, the projection onto the $26$-dimensional summand is
$h\mapsto h-\tfrac13\tr_e(h)e$. For a basepoint $e$, put
$F_{4,e}=\Stab_{G_\pm}(e)$ and
$\mathfrak f_e=\{X\in\mathfrak g_\pm:Xe=0\}$.

\begin{proposition}
  \label{prop:f4-stabilizer}
  For every basepoint $e$, $\mathfrak f_e=\operatorname{Der}(\mathfrak J_e)\cong\mathfrak f_4$,
  a simple Lie algebra of dimension $52$ and rank $4$.
  The $\mathfrak f_e$-module $\WC$ splits as $\WC=\C e\oplus\mathbf{26}$, with $\mathbf{26}=\{h:\tr_e(h)=0\}$ irreducible, and $F_{4,e}$ is connected and isomorphic to the complex group $F_4(\C)$ of type $F_4$.
\end{proposition}

\begin{proof}
  A derivation $D$ of $\mathfrak J_e$ kills the unit, and $\exp(tD)$ is an
  algebra automorphism for every $t$.
  Such automorphisms preserve the cubic norm, so $D$ infinitesimally
  annihilates $N_e$, and hence $N_\pm$. Thus $D\in\mathfrak f_e$.

  Conversely, let $X\in\mathfrak g_\pm$ satisfy $Xe=0$. Then $\exp(tX)$ preserves
  both $N_e$ and $e$, and Lemma \ref{lem:change-of-variables} shows that it preserves the Jordan
  product. Hence $X$ is a derivation.

  The remaining claims are transported through the isomorphism
  $\mathfrak J_e\cong\mathbb A$ of Theorem~\ref{thm:based-algebras}\,(\ref{item:based-algebras-albert}). This
  isomorphism takes $e$ to $1$ and $\tr_e$ to the trace of $\mathbb A$. It
  takes $F_{4,e}$ to the stabilizer of $1$ in the norm-isometry group of
  $\mathbb A$, and $\operatorname{Der}(\mathfrak J_e)$ to $\operatorname{Der}(\mathbb A)$.
  For the complex Albert algebra, $\operatorname{Der}(\mathbb A)\cong\mathfrak f_4$ \cite[Theorem, p.~137]{chevalley-1950-exceptional}. Its trace-zero part is an irreducible $\mathbf{26}$ \cite[Proposition~2.4.6(1)]{yokota-2009-exceptional}. The stabilizer of $1$ is the connected group $F_4(\C)$ \cite[4.6 and 14.24]{springer-1973-jordan}, \cite[\S15, p.~29]{garibaldi-2023-albert}, and at group level also \cite[Theorem~5]{petersson-2019-survey}.
\end{proof}

The transitivity in Theorem~\ref{thm:based-algebras}\,(\ref{item:based-algebras-orbits}) and the
$F_4$ stabilizer give, for every $c\in\C^\times$ and every
$e_c$ with $N_\pm(e_c)=c$,
\begin{equation*}
  N_\pm^{-1}(c)=G_\pm\cdot e_c
  \cong E_6(\C)/F_4(\C).
\end{equation*}
Each nonzero norm level is therefore a $26$-dimensional homogeneous
hypersurface in $\WC$, and the values $c\in\C^\times$ distinguish
these orbits of the exact group $G_\pm$.

\subsection{The $G_2$ symmetry retained by a basepoint}
\label{subsec:g2-at-basepoint}

A basepoint is needed to turn a critical norm into a Jordan product.
We ask how much of the original $G_2$ symmetry survives this choice.
No nonzero tensor is $G_2$-invariant, and the bound below shows that
at most an $8$-dimensional subgroup survives, attained at the line basepoints $E_u$ of
Proposition~\ref{ex:structured-basepoints}. This motivates the line-field construction of
Section~\ref{sec:reductions}.
The original group $\rho(G_2)$ preserves both critical norms.
After a basepoint $e$ is chosen, the subgroup that also preserves the Jordan
product is $\rho(G_2)\cap F_{4,e}$.
For $h\in\WC$, set
\begin{equation*}
  \mathfrak g_{2,h}=\{M\in\gtwo:\rho(M)h=0\},
  \qquad
  \mathfrak g^{\C}_{2,h}=\{M\in\gtwoC:\rho(M)h=0\}.
\end{equation*}
These are the \emph{annihilators}, or infinitesimal stabilizers, of $h$ in
$\gtwo$ and in $\gtwoC$. For real $h$, $\mathfrak g^{\C}_{2,h}$ is the
complexification of $\mathfrak g_{2,h}$, since $\rho$ is real.

Since $F_{4,e}=\Stab_{G_\pm}(e)$ and $\rho(G_2)\subset G_\pm$,
the first identity below is immediate, and taking Lie
algebras gives the second:
\begin{equation}
  \label{eq:g2-basepoint-intersection}
  \rho(G_2)\cap F_{4,e}=\rho(\Stab_{G_2}(e)),
  \qquad
  \rho(\gtwo)\cap\mathfrak f_e=\rho(\mathfrak g_{2,e}).
\end{equation}
We compute this intersection for two geometric families of basepoints,
constructed from line and associative-plane projections, and bound its
dimension for every choice. These families are not claimed to classify
all basepoints up to the original $G_2$-action.
The transport~\eqref{eq:basepoint-transport} uses the larger complex
group $G_\pm$; its relation to the original $G_2$ is discussed after
Corollary~\ref{cor:basepoint-transport}.

For a real unit vector $u$ and a $3$-plane $\Sigma\subset V$, recall the
projection tensors $P_u=u^\flat\otimes u^\flat$ and $P_\Sigma$
of~\eqref{eq:projection-tensors} and the corresponding test tensors $E_u$
and $E_\Sigma$ of~\eqref{eq:test-tensors}.
Here $g^{-1}P_\Sigma=\pr_\Sigma$ and $P_{\Sigma^\perp}=g-P_\Sigma$.
Associative planes are as in Definition~\ref{def:associative-plane}.
\begin{proposition}[Line and associative-plane basepoints]
  \label{ex:structured-basepoints}
For a real unit vector $u$ and an associative $3$-plane $\Sigma$, the test
tensors are
\begin{equation}
  E_u=7P_u-g,
  \qquad
  E_\Sigma=7P_\Sigma-3g=4P_\Sigma-3P_{\Sigma^\perp}.
  \label{eq:structured-basepoints}
\end{equation}
Then, for the compact $G_2$, the stabilizer subgroups are
\begin{equation*}
  \Stab_{G_2}(E_u)\cong SU(3)\rtimes\mathbb{Z}/2\qquad \text{and } \Stab_{G_2}(E_\Sigma)\cong SO(4).
\end{equation*}
The corresponding annihilators $\mathfrak g_{2,E_u}\cong\mathfrak{su}(3)$ and $\mathfrak g_{2,E_\Sigma}\cong\lieso(4)$ have dimensions $8$ and $6$, respectively.
In each case the fixed subspace of the annihilator acting on $W$ is
exactly the line spanned by the indicated basepoint.
Their cubic values are
\begin{equation*}
  \begin{aligned}
    (C_0(E_u),C_\varphi(E_u))&=(210,84),\\
    (C_0(E_\Sigma),C_\varphi(E_\Sigma))&=(84,1680).
  \end{aligned}
\end{equation*}
In particular, both tensors have nonzero critical norm and are basepoints.
\end{proposition}

\begin{proof}
The natural tensor action~\eqref{eq:rho-group} gives
\[
  \rho(O)P_u=P_{Ou},\qquad \rho(O)P_\Sigma=P_{O\Sigma}
  \qquad (O\in G_2).
\]
As $G_2$ preserves $g$, and the projection tensors determine the
corresponding unoriented line and plane, respectively,
\[
  \Stab_{G_2}(E_u)=\Stab_{G_2}(\R u),
  \qquad
  \Stab_{G_2}(E_\Sigma)=\Stab_{G_2}(\Sigma).
\]
The classical stabilizer results give
$\Stab_{G_2}(u)\cong SU(3)$ and
$\Stab_{G_2}(\Sigma)\cong SO(4)$
\cite[Theorems~1.9.1 and~1.10.1]{yokota-2009-exceptional},
\cite[Theorem~IV.1.8]{HarveyLawson82}.
The line stabilizer has an additional component, interchanging $u$
and $-u$. The extension by $\mathbb Z/2$ splits: for $u=e_1$,
the involution
$\diag(-1,1,-1,1,-1,1,-1)\in G_2$ sends $e_1$ to $-e_1$,
and the general case follows by $G_2$-transitivity on the unit sphere
\cite[Theorem~1.9.2]{yokota-2009-exceptional}.
This proves the group statements and, on differentiation, the
assertions about the annihilators.

For the fixed-line claims, the two relevant real representations
decompose as
\[
  V=\R u\oplus u^\perp\quad (SU(3)),
  \qquad
  V=\Sigma\oplus\Sigma^\perp\quad (SO(4)).
\]
In each case the summands are inequivalent irreducible real modules.
On $u^\perp\cong\C^3$, the real commutant of $SU(3)$ is spanned by
the identity and the skew-adjoint complex structure, so its self-adjoint
elements are real scalars. On the other summands, real Schur's lemma gives
the same conclusion. A self-adjoint endomorphism commuting with the
stabilizer is therefore scalar on each summand. Taking its traceless part
leaves exactly $\R E_u$ or $\R E_\Sigma$, respectively.

Both tensors in~\eqref{eq:structured-basepoints} have nonzero critical norm.
Indeed, $C_0(h)$ and $C_\varphi(h)$ are real for real $h$, whereas
$\lambda_\pm\notin\R$.
Thus $N_\pm(h)=0$ would force both real cubics to vanish.
The eigenvalues of $g^{-1}E_u$ and $g^{-1}E_\Sigma$ give
\begin{equation*}
  C_0(E_u)=6^3-6=210,
  \qquad
  C_0(E_\Sigma)=3\cdot4^3-4\cdot3^3=84.
\end{equation*}
For $u=e_1$ and the associative plane $\Sigma=\Span(e_1,e_2,e_3)$, these
are the test tensors $E_1$ and $E_\Sigma$ of the proof of
Proposition~\ref{prop:double-adjoint}.
The values above agree with~\eqref{eq:line-plane-values}, which gives
$N_\lambda(E_1)=42(5+2\lambda)$ and $N_\lambda(E_\Sigma)=84(1+20\lambda)$.
By $G_2$-equivariance, these cubic values hold for every $u$ and $\Sigma$.
\end{proof}

\begin{remark}
Either tensor can therefore serve as the model basepoint $e_0$
in~\eqref{eq:basepoint-transport}.
For that transport, $F_{4,e}=\gamma^{-1}F_{4,e_0}\gamma$,
since multiplication of $e_0$ by $t$ does not change its stabilizer.

The transport by $G_\pm$ does not identify the two models under
$G_2$, even projectively: the ratios of the original invariant cubics are
\begin{equation*}
  \frac{C_\varphi(E_u)}{C_0(E_u)}=\frac25,
  \qquad
  \frac{C_\varphi(E_\Sigma)}{C_0(E_\Sigma)}=20.
\end{equation*}
These ratios are unchanged by $G_2$ and by nonzero scalar multiplication.
For example, choose $s_\pm^3=N_\pm(E_u)/N_\pm(E_\Sigma)$.
Then $E_u$ and $s_\pm E_\Sigma$ have equal critical norm
and are $G_\pm$-equivalent, but they are not projectively $G_2$-equivalent.
\end{remark}

The next proposition bounds the residual symmetry for every nonzero complex
tensor, including tensors of zero critical norm.
\begin{proposition}
  \label{prop:g2-at-basepoint}
For every nonzero $h\in\WC$, $\dim_{\C}\mathfrak g^{\C}_{2,h}\le8$.
Equality is attained both by the nonzero real multiples of $E_u$, for
every real unit vector $u$, and by the rank-one null tensors $P_v=v^\flat\otimes v^\flat$, where $0\ne v\in\VC$
and $g(v,v)=0$.
No $G_2^{\C}$-orbit meets both families.
The latter satisfy $N_\pm(P_v)=0$ and hence are not basepoints.
\end{proposition}

\begin{proof}
  The irreducible $G_2^\C$-module $\WC=V(2,0)$ has a unique closed
  projective orbit $Y$. It is the image of the null quadric
  \[
    Q^5=\{[v]\in\mathbb P(\VC):g(v,v)=0\}
  \]
  under the map $[v]\mapsto[P_v]$, and has dimension $5$.
  It lies in the closure of every projective orbit
  \cite[Theorem~III.1.2]{zak-1993-tangents}.
  If $[h]\notin Y$, the projective orbit of $[h]$ has dimension
  greater than $5$, since $Y$ lies in the boundary of its orbit closure
  and that boundary has strictly smaller dimension than the orbit.
  If $[h]\in Y$, a one-parameter subgroup of a maximal torus fixing
  $[h]$ acts on $\C h$ by a nontrivial character. It therefore rescales
  $h$ through all of $\C^\times$, so the affine orbit has dimension $6$.
  Thus in either case
  \[
    \dim_\C\mathfrak g^\C_{2,h}
       =14-\dim(G_2^\C\cdot h)\le8.
  \]
  Equality holds for $E_u$ by Proposition~\ref{ex:structured-basepoints} and for every nonzero null
  $P_v$ by the affine-orbit calculation. Finally,
  $C_0(P_v)=C_\varphi(P_v)=0$, since $g(v,v)=0$ and $\mathcal P(P_v,P_v)=0$ by~\eqref{eq:p-rank-one-rule}, so $N_\pm(P_v)=0$;
  these tensors are not basepoints and cannot lie in the orbit of a
  nonzero multiple of $E_u$, whose critical norm is nonzero.
\end{proof}

For comparison, Darp\"o classifies the $G_2$-orbits of real symmetric $2$-tensors,
written there as symmetric matrices;
the orbit of $E_u$ is the case treated in
\cite[Proposition~3.1(2)]{darpo-2007-normal}.
At a general complex tensor the stabilizer is much smaller: it is trivial
by \cite[Corollary~6]{guralnick-2024-generic}, which applies because
$\dim V(2,0)=27>14=\dim G_2^{\C}$.
Proposition~\ref{prop:g2-at-basepoint} instead gives a uniform bound at every
nonzero complex tensor, attained at the geometric basepoints $E_u$.
We next compare these constructions with the Yukawa cubic.

\subsection{Comparison with the Yukawa cubic}
\label{subsec:yukawa-comparison}
\label{par:yukawa-comparison}
Lee and Leung describe an exceptional Jordan interpretation of the cubic on
$\Lambda^3_{27}$, with an $E_6$ symmetry and an $F_4$ stabilizer after an
$SU(3)$ reduction \cite[\S1.2.1, Remark]{Lee:2002fa}.
This interpretation was given without proof or an explicit comparison with a
normalized Albert norm.
Proposition~\ref{prop:f4-stabilizer} and Proposition~\ref{ex:structured-basepoints} establish
such a construction, with a normalized norm and a basepoint, for the
critical cubics.
We now compare it with their cubic.

Recall from~\eqref{eq:yukawa-integrand} that the pointwise Yukawa
integrand is a nonzero multiple of $C_\varphi$.

An Albert algebra over $\R$ is a real form of $\mathbb A$, that is, a
real Jordan algebra $J$ with $J\otimes_{\R}\C\cong\mathbb A$
\cite[\S1.6, p.~225]{petersson-2019-survey}.
Its norm $N_J$ complexifies to $N_{\mathbb A}$.

\begin{proposition}
  \label{prop:yukawa-not-albert}
  Let $\mathbb F=\R$ or $\C$.
  For every Albert algebra $J$ over $\mathbb F$, with norm $N_J$, there is no $\mathbb F$-linear isomorphism $L\colon W_{\mathbb F}\to J$ and no $s\in\mathbb F^\times$ such that $C_\varphi(h,h,h)=sN_J(Lh)$ for every $h\in W_{\mathbb F}$.
  Thus $C_\varphi$ is not linearly equivalent, even up to a nonzero scalar, to any real or complex Albert norm.
\end{proposition}

\begin{proof}
  A nonzero rescaling leaves the exact infinitesimal stabilizer of a cubic
  unchanged, and pullback by a linear isomorphism conjugates it:
  \begin{equation*}
    \stab(sN_J\circ L)=L^{-1}\stab(N_J)L.
  \end{equation*}
  An equivalence would therefore give stabilizers of equal dimensions.
  By Proposition~\ref{prop:phi-cubic-stabilizer},
  $\dim_{\mathbb F}\stab(C_\varphi)=14$.
  The exact linear stabilizer of an Albert norm is an algebraic group of
  type $E_6$ \cite[\S15, p.~29]{garibaldi-2023-albert}, whose Lie algebra has
  dimension $6+72=78$.
  Over $\R$ it is a real form of the same complex Lie algebra and has real
  dimension $78$.
  Since $14\ne78$, the proposed equivalence is impossible over either field.
\end{proof}

Thus the pointwise Yukawa integrand on $\Lambda^3_{27}$ is not equivalent to an Albert norm.
The $E_6$ symmetry and the $F_4$ stabilizer of a basepoint only occur for the two
critical members of the same pencil.
At $E_u$, the compact $G_2$ subgroup preserving the based algebra is
$SU(3)\rtimes\mathbb Z/2$, and choosing the unit vector $u$ as additional data
selects its $SU(3)$ identity component.
This discussion only concerns the pointwise cubic, rather than the integrated Yukawa coupling on the moduli space.
We return to the comparison with Lee and Leung in
Section~\ref{sec:reductions}, using parallel pointings.

Definition of a Jordan product requires the additional choice of a basepoint of nonzero norm and every such choice gives an algebra isomorphic to $\mathbb A$, with
symmetry group $F_{4,e}$.
However, we have shown that there is no nonzero basepoint that is invariant under the original $G_2$:
Proposition~\ref{prop:g2-at-basepoint} bounds the dimension of the
annihilator $\mathfrak g^{\C}_{2,h}$ of every nonzero $h\in\WC$ by $8<14$.
Consequently no basepoint can be selected equivariantly from $(V,\varphi)$
alone.
Any choice retains a compact subgroup of dimension at most $8$, with this
maximum attained at the basepoints $E_u$ and their nonzero real multiples.
In Section~\ref{sec:reductions} we organize the choice of pointing by allowing
real tangent lines to vary.

\subsection{Real points and compact $E_6$ forms}
\label{subsec:real-form}

The critical cubic lines are non-real and are exchanged by the
complex conjugation $h\mapsto\bar h$ associated with the original real module
$W$, since $\lambda_-=\bar\lambda_+$ in~\eqref{eq:critical-parameters}
and $C_0,C_\varphi$ are real. We distinguish three
real objects attached to these cubics. The real endomorphisms in
either critical stabilizer form only $\rho(\gtwo)$, of real
dimension $14$; they do not form a real form of
$\mathfrak g_\pm$. The Hermitian form induced by the metric,
on the other hand, selects a compact real form of
$\mathfrak g_\pm$, of real dimension $78$, and no basepoint is
needed for this. Finally, a real line in $V$ determines a basepoint and
a Euclidean real form of the corresponding complex Albert algebra,
of real dimension $27$. Here \emph{Euclidean} means admitting a
positive-definite associative inner product, equivalently being
formally real; see Section~\ref{subsubsec:euclidean-albert-form}.
We construct the compact forms here and
treat the based real algebras separately in
Section~\ref{subsec:euclidean-albert}.

For a complex-linear or antilinear map $X$ of $\WC$, write
$\overline X(h)=\overline{X\bar h}$. Thus $\overline X=X$ precisely
when $X$ preserves $W$.

The $W$-pairing on $W$ defines the Hermitian form
$\langle h,\bar k\rangle_W$ on $\WC$. For an endomorphism $X$,
$X^{\mathsf H}$ denotes its adjoint for this form, whereas $X^{\tadj}$
is its $W$-adjoint for the bilinear pairing~\eqref{eq:deformation-module}.
The induced pairing on endomorphisms is~\eqref{eq:endomorphism-pairing}.
We use the parity
embeddings $J^{\mathrm s}_Q$ and $J^{\mathrm a}_Q$ from~\eqref{eq:q64-parity-action} to describe
the compact real form selected by this pairing.
Recall the real tensor space $\mathcal Q_{64}$ of
Definition~\ref{def:q64} and the map $\mathcal E_\pm$ onto
$\mathfrak m_\pm$ introduced before Proposition~\ref{prop:block-kernels}.
For the compact normalization, define
\[
  \mathcal E^c_\pm=\pm\ii\sqrt7\,\mathcal E_\pm,
  \qquad
  \mathcal E^c_\pm(Q)=-J^{\mathrm a}_Q\pm\ii\sqrt7\,J^{\mathrm s}_Q.
\]
The superscript $c$ stands for compact. This rescaling does not change
the complex image $\mathfrak m_\pm$. On $\mathcal Q_{64}$ the operators
$\mathcal E_\pm(Q)$ are Hermitian, so the imaginary factor makes
$\mathcal E^c_\pm(Q)$ skew-Hermitian, as the next lemma shows.

\begin{lemma}
  \label{lem:hermitian-W}
  \renewcommand{\labelenumi}{(\roman{enumi})}
  \renewcommand{\theenumi}{\roman{enumi}}
  \leavevmode
  \begin{enumerate}
    \item\label{item:hermitian-w-positive}
    The Hermitian form $\langle h,\bar k\rangle_W$ on $\WC$
    is positive definite and invariant under the natural action
    of $O(7)$.

    \item\label{item:hermitian-w-so}
    For every real $M\in\lieso(7)$, the endomorphism
    $\rho(M)$ of $\WC$ is skew-Hermitian.

    \item\label{item:hermitian-w-parity}
    For every $Q\in\mathcal Q_{64}$, the operator $\mathcal E^c_\pm(Q)$
    lies in $\mathfrak g_\pm$ and is skew-Hermitian.
  \end{enumerate}
\end{lemma}

\begin{proof}
  The Hermitian form is positive definite because the $W$-pairing
  is positive definite on $W$. For $O\in O(7)$, the
  natural action is $\rho(O)h=OhO^{-1}$. Since $O$ is real,
  cyclicity of trace gives
  \[
    \langle\rho(O)h,\overline{\rho(O)k}\rangle_W
      =\tr(OhO^{-1}O\bar kO^{-1})
      =\langle h,\bar k\rangle_W.
  \]
  This proves~(\ref{item:hermitian-w-positive}). Differentiating
  along $O(t)=\exp(tM)$ gives
  \[
    \langle\rho(M)h,\bar k\rangle_W
      +\langle h,\overline{\rho(M)k}\rangle_W=0,
  \]
  which proves~(\ref{item:hermitian-w-so}).

  For~(\ref{item:hermitian-w-parity}), Proposition~\ref{prop:block-kernels}
  and the rescaling above give
  $\mathcal E^c_\pm(Q)\in\mathfrak m_\pm\subset\mathfrak g_\pm$.
  By Lemma~\ref{lem:q64-model}\,(\ref{item:q64-model-equivariance}), the operators $J^{\mathrm s}_Q$ and
  $J^{\mathrm a}_Q$ are real and respectively self-adjoint and
  skew-adjoint for the $W$-pairing. Their Hermitian
  adjoints therefore satisfy
  \[
    (J^{\mathrm s}_Q)^{\mathsf H}=J^{\mathrm s}_Q,
    \qquad
    (J^{\mathrm a}_Q)^{\mathsf H}=-J^{\mathrm a}_Q.
  \]
  Consequently,
  \[
    \bigl(\mathcal E^c_\pm(Q)\bigr)^{\mathsf H}
      =J^{\mathrm a}_Q\mp\ii\sqrt7\,J^{\mathrm s}_Q
      =-\mathcal E^c_\pm(Q).
  \]
\end{proof}

Let $U(\WC)$ be the unitary group of this Hermitian form,
and let $\mathfrak u(\WC)$ be its real Lie algebra of
skew-Hermitian endomorphisms. Define
\begin{equation}
  \label{eq:K-compactform}
  K_\pm:=G_\pm\cap U(\WC)
\end{equation}
and
\begin{equation}
  \label{eq:frak-k-def}
  \mathfrak k_\pm
    :=\rho(\gtwo)\oplus\mathcal E^c_\pm(\mathcal Q_{64}).
\end{equation}

\begin{proposition}[The compact real forms]
  \label{prop:compact-form}
  \renewcommand{\labelenumi}{(\arabic{enumi})}
  \renewcommand{\theenumi}{\arabic{enumi}}
  \leavevmode
  \begin{enumerate}
    \item\label{item:compact-form-lie}
    The space $\mathfrak k_\pm$ is a real Lie algebra and
    \[
      \mathfrak k_\pm
        =\mathfrak g_\pm\cap\mathfrak u(\WC).
    \]
    It is a compact real form of
    $\mathfrak g_\pm\cong\mathfrak e_6$. In particular,
    \[
      \mathfrak g_\pm=\mathfrak k_\pm\oplus\ii\mathfrak k_\pm,
      \qquad
      \dim_\R\mathfrak k_\pm=78.
    \]

    \item\label{item:compact-form-group}
    The group $K_\pm$ is the compact, simply connected Lie
    group of type $E_6$, with Lie algebra $\mathfrak k_\pm$.
    It is a maximal compact subgroup of
    $G_\pm\cong E_6(\C)$, and
    \[
      G_\pm=K_\pm\exp(\ii\mathfrak k_\pm),
      \qquad
      \rho(G_2)\subset K_\pm.
    \]

    \item\label{item:compact-form-conjugation}
    For $Q\in\mathcal Q_{64}$,
    \[
      \overline{\mathcal E^c_+(Q)}=\mathcal E^c_-(Q).
    \]
    Complex conjugation exchanges $\mathfrak k_+$ and
    $\mathfrak k_-$, and likewise exchanges $K_+$ and $K_-$.
  \end{enumerate}
\end{proposition}

\begin{proof}
  By Proposition~\ref{prop:block-kernels} and the rescaling above, the complex Lie
  algebra $\mathfrak g_\pm$ is the direct sum of
  $\rho(\gtwoC)$ and
  $\mathfrak m_\pm=\mathcal E^c_\pm(\mathcal Q_{64,\C})$. Both maps $\rho$ and $\mathcal E^c_\pm$ are
  complex-linear extensions of their restrictions to real
  arguments. Hence every $X\in\mathfrak g_\pm$ has a unique
  expression
  \[
    X=Y+\ii Z,
    \qquad Y,Z\in\mathfrak k_\pm.
  \]
  By Lemma~\ref{lem:hermitian-W}, both $Y$ and $Z$ are
  skew-Hermitian, so
  \[
    X^{\mathsf H}=-Y+\ii Z.
  \]
  Thus $X$ is skew-Hermitian precisely when $Z=0$. This gives
  \[
    \mathfrak g_\pm\cap\mathfrak u(\WC)=\mathfrak k_\pm,
    \qquad
    \mathfrak g_\pm=\mathfrak k_\pm\oplus\ii\mathfrak k_\pm.
  \]
  As the intersection of $\mathfrak g_\pm$ and
  $\mathfrak u(\WC)$, $\mathfrak k_\pm$ is a real Lie
  subalgebra. Moreover, $-\tr_W(XY)$ is
  invariant and positive definite on it, since
  $-\tr_W(X^2)=\tr_W(X^{\mathsf H}X)>0$ for nonzero
  skew-Hermitian $X$. It follows that $\mathfrak k_\pm$ is
  compact. Since its complexification is
  $\mathfrak e_6$, it has real dimension $78$.

  The same decomposition shows that $\mathfrak g_\pm$ is
  invariant under $X\mapsto-X^{\mathsf H}$. Since $G_\pm$
  is connected (Theorem~\ref{thm:critical-group}), it is invariant under the corresponding
  group involution $\gamma\mapsto(\gamma^{\mathsf H})^{-1}$. The Cartan
  decomposition for self-adjoint complex reductive groups
  therefore gives
  \[
    G_\pm=K_\pm\exp(\ii\mathfrak k_\pm),
  \]
  with $K_\pm$ maximal compact and a deformation retract of
  $G_\pm$. By Theorem~\ref{thm:critical-group}, $G_\pm$ is
  connected and simply connected. The deformation retraction
  gives the same properties for $K_\pm$. Together with the
  Lie algebra identification, this proves that $K_\pm$ is
  the compact, simply connected group of type $E_6$.
  The inclusion $\rho(G_2)\subset K_\pm$ follows because
  $\rho(G_2)$ preserves both the critical cubic and the
  Hermitian form.

  Finally, the operators $J^{\mathrm s}_Q$ and $J^{\mathrm a}_Q$ are real, so
  $\overline{\mathcal E^c_+(Q)}=\mathcal E^c_-(Q)$.
  Conjugation exchanges the critical stabilizers and
  preserves $\mathfrak u(\WC)$ and $U(\WC)$. It therefore
  exchanges the compact real forms $\mathfrak k_\pm$ and the
  compact subgroups $K_\pm$.
\end{proof}

The compact $E_6$ group is classical; see
\cite[p.~68]{yokota-2009-exceptional} and
\cite[Theorem~3.1.1, p.~80]{yokota-2025-exceptional}.
The description above realizes its Lie algebra inside
$\End(\WC)$ through the parity embeddings. This compact
form is selected by the critical cubic together with the
fixed Hermitian form, not by the abstract complex cubic
alone.

\subsection{The Euclidean Albert algebra at an $SU(3)$-type basepoint}
\label{subsec:euclidean-albert}

We now choose a real line in $V$. A unit vector $u$ spanning
it determines
\[
  E_u=7P_u-g,
\]
as in~\eqref{eq:test-tensors}. Since $E_{-u}=E_u$, this
basepoint depends only on the unoriented line. Its stabilizer in $G_2$
is $SU(3)\rtimes\mathbb Z/2$, with identity component the stabilizer
$SU(3)$ of $u$ (Proposition~\ref{ex:structured-basepoints}); this is what
we mean by an \emph{$SU(3)$-type basepoint}.
The complex-bilinear based trace form and the positive Hermitian form
are related by an antilinear map, the \emph{Euclidean real structure}
$\varsigma_u$ constructed below. We show that this map is an
involutive Jordan automorphism. On its fixed space the based trace
form is $1/14$ times the Hermitian form, hence positive definite.

Fix either critical parameter. Where convenient, we suppress its sign on $N$, $F$,
$T_e$, $T_{E_u}$, $\varsigma_u$ and $\mathcal J_u$ (defined below), restoring
$\pm$ when displaying the root-dependent coefficients or
comparing the two constructions.

\subsubsection{The based trace form and $SU(3)$ sectors}
\label{subsubsec:su3-trace-sectors}

The $E_1$ row of~\eqref{eq:line-plane-values}, transported from
$e_1$ to $u$ by $G_2$, gives $F_\pm(E_u)=(5+2\lambda_\pm)E_u$.
Thus the norm and based linear trace are
\begin{equation}
  \label{eq:su3-norm-trace}
  n_\pm:=N_\pm(E_u)=42(5+2\lambda_\pm)\ne0,
  \qquad
  \tr_{E_u}(h)=\tfrac12h(u,u)=\tfrac1{14}\langle E_u,h\rangle_W.
\end{equation}
In particular, Lemma~\ref{lem:based-trace-operator} gives
\begin{equation}
  \label{eq:su3-trace-response}
  T_{E_u}h=\tfrac1{14}\tr_{E_u}(h)E_u
    -\tfrac6{n_\pm}F_\pm(E_u,h).
\end{equation}

Write $n=n_\pm$ when suppressing the chosen sign. The based trace form
and trace operator are those of~\eqref{eq:based-trace-form} and
\eqref{eq:based-trace-operator} for $N/n$; their pairing conventions
are collected in Table~\ref{tab:based-conventions}.

Recall from Section~\ref{subsec:koecher-jump} the operator
$\Theta_u(v)=v\times u$. Like $\Theta_1$, it satisfies
$\Theta_u^2=-I$ on $u^\perp$, by $G_2$-transitivity on unit
vectors. Recall also the real sectors $W_6(u)$,
$W_8(u)$ and $W_{12}(u)$, of real dimensions
$6$, $8$ and $12$, with complexifications $W_r(u)_\C$.
They give the orthogonal decompositions
\begin{equation}
  \label{eq:su3-sectors}
  \begin{aligned}
    W&=\R E_u\oplus W_6(u)
         \oplus W_8(u)\oplus W_{12}(u),\\
    \WC&=\C E_u\oplus W_6(u)_\C
         \oplus W_8(u)_\C\oplus W_{12}(u)_\C.
  \end{aligned}
\end{equation}
We use the same notation for an orthogonal projector
on $W$ and its complex-linear extension.
The sectors are those of Table~\ref{tab:koecher-sector-response},
transported from $e_1$ to $u$. Products below use the metric
identification, so $g$ acts as $I$. For $h\in\WC$, put
$a=h(u,u)/6$ and
$H=(g-P_u)h(g-P_u)+a(g-P_u)$, the traceless horizontal block.
We compute the trace spectrum because its common absolute value
will make the Euclidean real structure involutive.

\begin{lemma}[Sector projectors and the based trace operator]
  \label{lem:su3-sectors}
  \renewcommand{\labelenumi}{(\roman{enumi})}
  \renewcommand{\theenumi}{\roman{enumi}}
  With the splitting~\eqref{eq:su3-sectors}, the following hold.
  \begin{enumerate}
    \item\label{item:su3-sectors-projectors}
    The four orthogonal projectors are
    \begin{equation}
      \label{eq:su3-projectors}
      \begin{aligned}
        \Pi_1h&=aE_u,
        &\Pi_6h&=P_uh(g-P_u)+(g-P_u)hP_u,\\
        \Pi_8h&=\tfrac12(H-\Theta_uH\Theta_u),
        &\Pi_{12}h&=\tfrac12(H+\Theta_uH\Theta_u).
      \end{aligned}
    \end{equation}
    These projectors are polynomial in
    $(g,\varphi,P_u)$ and are invariant under $u\mapsto-u$.

    \item\label{item:su3-sectors-trace}
    The based trace operator and its inverse are
    \begin{equation}
      \label{eq:su3-trace-projectors}
      \begin{aligned}
        T_{E_u}
          &=\tfrac1{14}\Pi_1-\tfrac1{14}\Pi_6
            -\tfrac{1\pm3\ii\sqrt7}{112}\Pi_8
            +\tfrac{3\pm\ii\sqrt7}{56}\Pi_{12},\\
        T_{E_u}^{-1}
          &=14\Pi_1-14\Pi_6
            -\tfrac74(1\mp3\ii\sqrt7)\Pi_8
            +\tfrac72(3\mp\ii\sqrt7)\Pi_{12}.
      \end{aligned}
    \end{equation}
    Every eigenvalue of $T_{E_u}$ has absolute value $1/14$.
  \end{enumerate}
\end{lemma}

\begin{proof}
  Split a symmetric endomorphism into its mixed block,
  scalar block and horizontal traceless block. This gives
  the expressions for $\Pi_1$, $\Pi_6$ and $H$ in
  \eqref{eq:su3-projectors}. On the horizontal traceless
  block, the map $h\mapsto-\Theta_uh\Theta_u$ is an
  orthogonal involution. Its $+1$ and $-1$ eigenspaces are
  $W_8(u)$ and $W_{12}(u)$, respectively,
  giving the remaining projectors. The two occurrences of
  $u$ in $\Theta_uH\Theta_u$ enter quadratically, through
  $P_u$. Changing $u$ to $-u$ therefore leaves every
  projector unchanged.

  By $G_2$-equivariance, the sector values in
  Table~\ref{tab:koecher-sector-response} hold at every real
  unit vector $u$. Thus $F_\pm(E_u,\cdot)$ acts on the four
  sectors by
  \[
    5+2\lambda_\pm,\qquad
    \tfrac52+\lambda_\pm,\qquad
    -1+8\lambda_\pm,\qquad
    -1-6\lambda_\pm.
  \]
  Substituting in~\eqref{eq:su3-trace-response}, the scalar
  sector gives
  \[
    T_{E_u}E_u
      =\left(\frac3{14}
        -\frac{6(5+2\lambda_\pm)}{42(5+2\lambda_\pm)}\right)E_u
      =\frac1{14}E_u.
  \]
  The linear trace $\tr_{E_u}$ vanishes on the other three
  sectors. Their eigenvalues are therefore
  $-6f_r(\lambda_\pm)/[42(5+2\lambda_\pm)]$, where $f_r$ is the
  corresponding response scalar above. For the mixed sector,
  $f_6=(5+2\lambda_\pm)/2$ gives $-1/14$.
  With $\lambda_\pm=(1\pm\ii\sqrt7)/8$ from~\eqref{eq:critical-parameters},
  the other two become
  \[
    -\frac{-1+8\lambda_\pm}{7(5+2\lambda_\pm)}
      =-\frac{1\pm3\ii\sqrt7}{112},\qquad
    \frac{1+6\lambda_\pm}{7(5+2\lambda_\pm)}
      =\frac{3\pm\ii\sqrt7}{56}.
  \]
  These give~\eqref{eq:su3-trace-projectors}. Since
  \[
    |1\pm3\ii\sqrt7|^2=64,
    \qquad
    |3\pm\ii\sqrt7|^2=16,
  \]
  all four eigenvalues have absolute value $1/14$.
  Taking their reciprocals gives $T_{E_u}^{-1}$.
\end{proof}

\subsubsection{The Euclidean real structure}
\label{subsubsec:compatible-real-involution}

To show that the Euclidean real structure defined
in~\eqref{eq:euclidean-real-structure} below also
preserves the based Jordan product, we first relate the two
critical cubics through the based trace operator. The next
identity holds at any basepoint of nonzero norm.

\begin{lemma}[Opposite-cubic identity]
  \label{lem:opposite-cubic}
  Let $e\in\WC$ satisfy $n=N_\pm(e)\ne0$. The based trace
  operator for the normalized cubic $N_\pm/n$ satisfies
  \begin{equation}
    \label{eq:opposite-cubic}
    N_\mp(T_ex)=\frac{18}{n^2}N_\pm(x)
    \qquad\text{for all }x\in\WC.
  \end{equation}
\end{lemma}

\begin{proof}
  We prove the identity at $\lambda_+$. The based trace
  form is symmetric by~\eqref{eq:based-trace-form} and nondegenerate
  by Theorem~\ref{thm:based-algebras}\,(\ref{item:based-algebras-albert}), so $T_e$ is
  $W$-self-adjoint and invertible. Define
  \[
    Q(y)=n^{-2}N_+(T_e^{-1}y).
  \]
  By~\eqref{eq:F-adjoint}, its $W$-gradient and the based
  adjoint satisfy
  \[
    \nabla Q(y)=\frac3{n^2}T_e^{-1}F_+(T_e^{-1}y),
    \qquad
    x^\sharp=\frac3nT_e^{-1}F_+(x).
  \]
  Quadratic homogeneity and the based adjoint identity in
  \eqref{eq:cubic-norm-structure} give
  \[
    \nabla Q\bigl(3F_+(x)\bigr)
      =n(x^\sharp)^\sharp
      =N_+(x)x.
  \]
  By Proposition~\ref{prop:double-adjoint}, the last
  expression is also $\tfrac16F_-(3F_+(x))$.

  By Lemma~\ref{lem:critical-adjoint-dominance}, $F_+$ has dense
  image, so the two polynomial maps agree everywhere:
  \[
    \nabla Q=\tfrac16F_-=\nabla(N_-/18).
  \]
  Both cubics are homogeneous, so $Q=N_-/18$.
  Substituting $y=T_ex$ gives~\eqref{eq:opposite-cubic}.
  Complex conjugation gives the identity at $\lambda_-$.
\end{proof}

Return to $e=E_u$. We seek an antilinear map satisfying
$\langle x,\bar y\rangle_W=14\langle x,\varsigma_u y\rangle_{E_u}$.
The trace-operator definition~\eqref{eq:based-trace-operator} forces
$14T_{E_u}\varsigma_u y=\bar y$. Thus we first conjugate and then
apply the inverse trace operator, defining the Euclidean real structure by
\begin{equation}
  \label{eq:euclidean-real-structure}
  \varsigma_u(x)=\tfrac1{14}T_{E_u}^{-1}(\bar x).
\end{equation}
The factor $1/14$ is fixed by the common absolute value of
the eigenvalues of $T_{E_u}$. Define the fixed space of
$\varsigma_u$ by
\begin{equation}
  \label{eq:J_u-fixed-def}
  \mathcal J_u=\{x\in\WC:\varsigma_u(x)=x\}.
\end{equation}
From~\eqref{eq:su3-norm-trace} and
$5+2\lambda_\pm=(21\pm\ii\sqrt7)/4$, we have
$|5+2\lambda_\pm|^2=(441+7)/16=28$. Consequently,
\[
  n_\pm=N_\pm(E_u)=42(5+2\lambda_\pm),
  \qquad
  n_-=\bar n_+,
  \qquad
  n_+n_-=14^3\cdot18.
\]
Thus $N_\pm/n_\pm$ takes the value $1$ at $E_u$.

\begin{lemma}[The Euclidean real structure]
  \label{lem:euclidean-real-structure}
  \renewcommand{\labelenumi}{(\roman{enumi})}
  \renewcommand{\theenumi}{\roman{enumi}}
  For either critical parameter and every real unit vector $u$,
  the following hold.
  \begin{enumerate}
    \item\label{item:euclidean-real-structure-automorphism}
    The map $\varsigma_u$ is an antilinear involutive
    automorphism of $\mathfrak J_{E_u}$ fixing $E_u$.
    For every $x\in\WC$,
    \[
      \frac{N_\pm(\varsigma_u x)}{n_\pm}
        =\overline{\frac{N_\pm(x)}{n_\pm}}.
    \]

    \item\label{item:euclidean-real-structure-hermitian}
    For all $x,y\in\WC$,
    \begin{equation}
      \label{eq:euclidean-hermitian-pairing}
      \langle x,\bar y\rangle_W
        =14\langle x,\varsigma_u y\rangle_{E_u}.
    \end{equation}
    Moreover, $\varsigma_u$ is the unique antilinear
    involution $\eta$ for which
    \[
      \langle x,\bar y\rangle_W
        =s\langle x,\eta y\rangle_{E_u}
        \qquad\text{for all }x,y\in\WC
    \]
    with some $s>0$. Necessarily $s=14$.
  \end{enumerate}
\end{lemma}

\begin{proof}
  By Lemma~\ref{lem:su3-sectors}, the operator $T_{E_u}$
  acts by scalars on the complexifications of mutually
  orthogonal real subspaces. With
  $\overline{T_{E_u}}(x)=\overline{T_{E_u}\bar x}$, the spectral
  decomposition of $T_{E_u}$ gives
  \[
    T_{E_u}\overline{T_{E_u}}=\tfrac1{14^2}\Id.
  \]
  Consequently,
  \[
    \varsigma_u^2
      =\tfrac1{14^2}
        T_{E_u}^{-1}\overline{T_{E_u}}^{-1}
      =\Id.
  \]
  Since $E_u$ is real and $T_{E_u}E_u=E_u/14$, we also
  have $\varsigma_u E_u=E_u$.

  We next check the normalized cubic identity at
  $\lambda_+$. By Lemma~\ref{lem:opposite-cubic},
  \[
    N_-\circ T_{E_u}=\frac{18}{n_+^2}N_+.
  \]
  Applying this identity to $T_{E_u}^{-1}\bar x$ and using
  homogeneity gives
  \[
    \begin{aligned}
      N_+(\varsigma_u x)
        &=\frac1{14^3}N_+(T_{E_u}^{-1}\bar x)\\
        &=\frac{n_+^2}{14^3\cdot18}N_-(\bar x)\\
        &=\frac{n_+}{\bar n_+}\,
          \overline{N_+(x)}.
    \end{aligned}
  \]
  Dividing by $n_+$ proves the identity. The identity at the
  other root follows by conjugation.

  Since $\varsigma_u$ fixes $E_u$, polarization of the
  normalized cubic identity at $E_u$ yields
  \[
    \tr_{E_u}(\varsigma_u x)
      =\overline{\tr_{E_u}(x)}
  \]
  and
  \[
    \langle\varsigma_u x,\varsigma_u y\rangle_{E_u}
      =\overline{\langle x,y\rangle_{E_u}}.
  \]
  By the adjoint characterization in~\eqref{eq:cubic-norm-structure},
  the normalized cubic and the nondegenerate based trace
  form determine the based adjoint, as in~\eqref{eq:L-sharp}. Hence
  \[
    (\varsigma_u x)^\sharp=\varsigma_u(x^\sharp).
  \]
  Substitution in~\eqref{eq:circ-e-Jordan-prod} now gives
  \[
    \varsigma_u(x\circ_{E_u}y)
      =\varsigma_u x\circ_{E_u}\varsigma_u y.
  \]
  This proves~(\ref{item:euclidean-real-structure-automorphism}).

  For~(\ref{item:euclidean-real-structure-hermitian}),
  the definition of $\varsigma_u$ gives
  $14T_{E_u}\varsigma_u y=\bar y$. Therefore
  \[
    14\langle x,\varsigma_u y\rangle_{E_u}
      =14\langle x,T_{E_u}\varsigma_u y\rangle_W
      =\langle x,\bar y\rangle_W.
  \]
  Finally, suppose that another antilinear involution
  $\eta$ satisfies the stated identity with $s>0$.
  Nondegeneracy of the based trace form gives
  \[
    \eta=\frac{14}{s}\varsigma_u.
  \]
  Squaring yields $(14/s)^2=1$. Positivity forces $s=14$,
  and hence $\eta=\varsigma_u$.
\end{proof}

\subsubsection{The Euclidean real algebra and its compact stabilizer}
\label{subsubsec:euclidean-albert-form}

The fixed space of $\varsigma_u$ can be described separately
in each sector. Let $S$ be one of the real summands in
\eqref{eq:su3-sectors}, and let $\vartheta$ be the eigenvalue
of $T_{E_u}$ on $S_\C$. For nonzero $v\in S$,
\begin{equation}
  \label{eq:sector-real-line}
  \varsigma_u(av)=av
    \quad\Longleftrightarrow\quad
    \bar a=14\vartheta a.
\end{equation}
Since $|14\vartheta|=1$, this equation defines a real line in
$\C$. A nonzero solution $a$ determines the real subspace
$aS$ inside $S_\C$. Replacing
$a$ by any nonzero real multiple leaves this subspace
unchanged. Thus the coefficients below specify real lines
in $\C$; they need not be unit complex numbers. Such a
rescaling changes neither the real algebra nor its
inherited trace form.
For the last two sectors we choose generators
\[
  a_{8,\pm}=-\tfrac{7\pm3\ii\sqrt7}{4},\qquad
  a_{12,\pm}=\sqrt7\mp\ii.
\]

A real Jordan algebra $J$ is called \emph{formally real}
if
\[
  x_1\circ x_1+\cdots+x_k\circ x_k=0
    \quad\Longrightarrow\quad
    x_1=\cdots=x_k=0.
\]
For finite-dimensional unital Jordan algebras, formal reality
is equivalent to the existence of a positive-definite inner
product associative with respect to the Jordan product
\cite[Proposition~VIII.4.2]{faraut-1994-analysis}.
Such an algebra is called a Euclidean Jordan algebra.
Let $H_3(\O)$ be the real Jordan algebra of $3\times3$ Hermitian
matrices over the real octonions $\O$, with
$x\circ y=\tfrac12(xy+yx)$. It is Euclidean of rank $3$
\cite[Corollary~V.2.6]{faraut-1994-analysis}, and
$H_3(\O_\C)=H_3(\O)\otimes_\R\C$ is the complex Albert algebra
introduced in Section~\ref{subsec:albert-norm}.
Recall that the Jordan determinant is the signed constant coefficient
of the generic characteristic polynomial
(see the discussion of rank and determinant following~\eqref{eq:cubic-characteristic-identity}).

\begin{proposition}[The Euclidean Albert real form]
  \label{prop:euclidean-form}
  \renewcommand{\labelenumi}{(\arabic{enumi})}
  \renewcommand{\theenumi}{\arabic{enumi}}
  For either critical parameter and every real unit vector $u$,
  the following hold.
  \begin{enumerate}
    \item\label{item:euclidean-form-albert}
    With the real sectors of~\eqref{eq:su3-sectors},
    \begin{equation}
      \label{eq:euclidean-sector-decomposition}
      \begin{aligned}
        \mathcal J_u
          &=\{x\in\WC:\bar x=14T_{E_u}x\}\\
          &=\R E_u\oplus\ii W_6(u)
            \oplus a_{8,\pm}W_8(u)
            \oplus a_{12,\pm}W_{12}(u).
      \end{aligned}
    \end{equation}
    Each coefficient is a nonzero generator of the
    corresponding real scalar line and may be replaced
    independently by a nonzero real multiple.

    The real Jordan algebra $\mathcal J_u$, with unit $E_u$,
    is simple and Euclidean, and is isomorphic to $H_3(\O)$.
    Under this isomorphism, $N_\pm/n_\pm$ is the Jordan
    determinant. Its based trace form is positive definite:
    for every nonzero $x\in\mathcal J_u$,
    \[
      \langle x,x\rangle_{E_u}
        =\tfrac1{14}\langle x,\bar x\rangle_W>0.
    \]

    \item\label{item:euclidean-form-intersection}
    The intersection with the original real module is
    \[
      \mathcal J_u\cap W=\R E_u.
    \]
  \end{enumerate}
\end{proposition}

\begin{proof}
  The fixed-space equation follows directly from
  \eqref{eq:euclidean-real-structure}. Applying
  \eqref{eq:sector-real-line} to the four eigenvalues in
  Lemma~\ref{lem:su3-sectors} gives the respective real
  scalar lines
  \[
    \R,\qquad\ii\R,\qquad
    a_{8,\pm}\R,\qquad a_{12,\pm}\R.
  \]
  This proves~\eqref{eq:euclidean-sector-decomposition}
  and explains the freedom to rescale its coefficients
  by nonzero real numbers.

  By Lemma~\ref{lem:euclidean-real-structure}, the fixed
  space is a real Jordan algebra with unit $E_u$.
  The same lemma shows that its based trace form is real
  and that
  \[
    \langle x,x\rangle_{E_u}
      =\tfrac1{14}\langle x,\bar x\rangle_W
      \qquad(x\in\mathcal J_u).
  \]
  It is therefore positive definite.

  To verify formal reality directly, set $e=E_u$ and $N_e=N/n$, and polarize
$dN_e|_x(y)=\langle x^\sharp,y\rangle_e$ in~\eqref{eq:cubic-norm-structure}.
  The resulting symmetry of $\langle x\mathbin\sharp y,z\rangle_e$,
  together with $e\mathbin\sharp x=\tr_e(x)e-x$, gives
  \[
    \tr_e(x\mathbin\sharp y)
      =\langle x\mathbin\sharp y,e\rangle_e
      =\langle e\mathbin\sharp x,y\rangle_e
      =\tr_e(x)\tr_e(y)-\langle x,y\rangle_e.
  \]
  Apply $\tr_{E_u}$ to~\eqref{eq:circ-e-Jordan-prod}, using this identity
  and $\tr_{E_u}(E_u)=3$. This gives
  \[
    \tr_{E_u}(x\circ_{E_u}y)
      =\langle x,y\rangle_{E_u}.
  \]
  If $\sum_jx_j\circ_{E_u}x_j=0$ in $\mathcal J_u$, then
  \[
    0=\sum_j\tr_{E_u}(x_j\circ_{E_u}x_j)
      =\sum_j\langle x_j,x_j\rangle_{E_u}.
  \]
  Positivity forces every $x_j$ to vanish. Thus
  $\mathcal J_u$ is formally real and hence Euclidean.

  Since $\varsigma_u$ is an antilinear involution, the
  complexification of $\mathcal J_u$ is
  $\mathfrak J_{E_u}$, which is the simple complex Albert algebra of
  degree three by Theorem~\ref{thm:based-algebras}\,(\ref{item:based-algebras-albert}).
  A proper nonzero real ideal of $\mathcal J_u$ would
  complexify to a proper nonzero ideal of
  $\mathfrak J_{E_u}$. The latter algebra is simple, so
  $\mathcal J_u$ is simple as well.

  We must still determine its rank: the classification of simple Euclidean Jordan algebras, recalled below, also includes
  rank-two spin factors, so dimension and simplicity alone do not suffice.
  The discussion of rank and determinant following~\eqref{eq:cubic-characteristic-identity}
  bounds the rank by three.
  Suppose that every element had minimal polynomial of degree at most two.
  Then $E_u$, $x$ and $x^{\circ2}$ would be linearly
  dependent for every $x\in\mathcal J_u$. The corresponding
  coordinate minors are polynomial functions. If they
  vanished on the real form $\mathcal J_u$, they would
  vanish on its complexification. This contradicts the
  degree-three assertion of
  Theorem~\ref{thm:based-algebras}\,(\ref{item:based-algebras-albert}).
  Hence $\mathcal J_u$ has rank three.

  By the classification of simple Euclidean Jordan algebras
  \cite[Theorem~V.3.7]{faraut-1994-analysis}, a simple algebra
  of rank three has dimension $3+3d$, where
  $d\in\{1,2,4,8\}$. Since $\dim_\R\mathcal J_u=27$,
  we obtain $d=8$ and
  \[
    \mathcal J_u\cong H_3(\O).
  \]

  The rank-three determinant argument following~\eqref{eq:cubic-characteristic-identity}, applied to the real
  characteristic identity~\eqref{eq:cubic-characteristic-identity},
  identifies the Jordan determinant with $N_\pm/n_\pm$.
  This completes~(\ref{item:euclidean-form-albert}).

  Finally, if $x\in W$, then $\bar x=x$. Thus
  \[
    x\in\mathcal J_u\cap W
      \quad\Longleftrightarrow\quad
      x\in W\ \text{and}\ T_{E_u}x=\tfrac1{14}x.
  \]
  By Lemma~\ref{lem:su3-sectors}, the $1/14$-eigenspace
  is precisely $\C E_u$. Its intersection with $W$ is
  $\R E_u$, proving~(\ref{item:euclidean-form-intersection}).
\end{proof}

We now compare the based algebra with the compact group of
Section~\ref{subsec:real-form}. For this comparison,
restore the signs in $\varsigma_{u,\pm}$ and
$\mathcal J_{u,\pm}$. By~\eqref{eq:euclidean-hermitian-pairing},
\begin{equation}
  \label{eq:compact-via-based-pairing}
  \begin{split}
    K_\pm=\{\gamma\in G_\pm:\;&
      \langle\gamma x,\varsigma_{u,\pm}\gamma y\rangle_{E_u}
        =\langle x,\varsigma_{u,\pm}y\rangle_{E_u}\\
      &\text{for all }x,y\in\WC\}.
  \end{split}
\end{equation}
The based trace form here is taken at the corresponding
critical parameter. Although this expression uses $E_u$, the
group $K_\pm$ itself does not depend on $u$. The dependence
on the basepoint enters through the stabilizer of $E_u$ in
$K_\pm$.
Define $K_{\pm,u}=\Stab_{K_\pm}(E_u)$.

\begin{corollary}[The compact stabilizer of the basepoint]
  \label{cor:compact-basepoint}
  For either critical parameter,
  \begin{equation}
    \label{eq:compact-basepoint-stabilizer}
    K_{\pm,u}
      =\{\gamma\in G_\pm:
        \gamma E_u=E_u,\quad
        \gamma\varsigma_{u,\pm}
          =\varsigma_{u,\pm}\gamma\}.
  \end{equation}
  Restriction to $\mathcal J_{u,\pm}$ gives an isomorphism
  \[
    K_{\pm,u}
      \cong\operatorname{Aut}_{\R}(\mathcal J_{u,\pm}).
  \]
  In particular, $K_{\pm,u}$ is the compact Lie group of
  type $F_4$. Moreover,
  \[
    \rho(\Stab_{G_2}(E_u))\subset K_{\pm,u}\subset K_\pm,
    \qquad
    \Stab_{G_2}(E_u)\cong SU(3)\rtimes\mathbb Z/2.
  \]
\end{corollary}

\begin{proof}
  Fix a sign and suppress it. Suppose that $\gamma\in G$
  fixes $E_u$. Since it preserves the normalized cubic
  and the basepoint, it preserves the based trace form
  and the based Jordan product by
  Lemma~\ref{lem:change-of-variables}, equations~\eqref{eq:L-based-trace}
  and~\eqref{eq:L-prod}. Therefore
  \[
    \langle\gamma x,\varsigma_u\gamma y\rangle_{E_u}
      =\langle x,\gamma^{-1}\varsigma_u\gamma y\rangle_{E_u}.
  \]
  By nondegeneracy, this equals
  $\langle x,\varsigma_u y\rangle_{E_u}$ for every $x,y$
  precisely when $\gamma\varsigma_u=\varsigma_u\gamma$.
  Together with~\eqref{eq:compact-via-based-pairing}, this
  proves~\eqref{eq:compact-basepoint-stabilizer}.

  Such a $\gamma$ restricts to a real Jordan automorphism
  of $\mathcal J_u$. Restriction is injective because
  $\mathcal J_u$ complexifies to $\WC$. Conversely, a real
  Jordan automorphism of $\mathcal J_u$ fixes its unit
  $E_u$ and preserves Jordan powers, hence minimal polynomials.
  It therefore preserves the generic characteristic polynomial and its
  signed constant coefficient, the Jordan determinant
  (see the discussion following~\eqref{eq:cubic-characteristic-identity}). Its
  complex-linear extension therefore preserves $N/n$,
  and hence $N$, and commutes with $\varsigma_u$.
  By~\eqref{eq:compact-basepoint-stabilizer}, this extension
  belongs to $K_{\pm,u}$. Restriction is thus an isomorphism.

  By Proposition~\ref{prop:euclidean-form},
  $\mathcal J_u\cong H_3(\O)$. The automorphism group of
  this real algebra is compact $F_4$; see
  \cite[\S~2.2 and Theorem~3.7.1]{yokota-2009-exceptional}.
  Finally, $\rho(G_2)\subset K_\pm$ by Proposition~\ref{prop:compact-form}\,(\ref{item:compact-form-group}), and
  Proposition~\ref{ex:structured-basepoints} identifies the
  stabilizer of $E_u$ in $G_2$ with
  $SU(3)\rtimes\mathbb Z/2$. This proves the inclusions in
  the last display of the corollary.
\end{proof}

Thus $K_{\pm,u}=K_\pm\cap F_{4,E_u}$ is a compact form of the
complex stabilizer $F_{4,E_u}$ of Proposition~\ref{prop:f4-stabilizer}.
Indeed, the Lie algebra of $K_{\pm,u}$ consists of skew-Hermitian operators, so it
intersects its imaginary multiple trivially. This real Lie algebra and
$\mathfrak f_{E_u}$ have dimension $52$ over their respective fields;
the complexification of the former is therefore $\mathfrak f_{E_u}$.

The based construction depends only on the unoriented line
$\R u$. Indeed, $E_{-u}=E_u$, and the sector projectors,
trace operator, involution and fixed algebra are all
invariant under $u\mapsto-u$. The construction is also
$G_2$-equivariant: for $O\in G_2$,
\[
  \varsigma_{Ou,\pm}\rho(O)
    =\rho(O)\varsigma_{u,\pm},
  \qquad
  \mathcal J_{Ou,\pm}=\rho(O)\mathcal J_{u,\pm}.
\]
These identities follow from the equivariance of the
normalized cubic and based trace form, together with the
reality of $\rho(O)$. Conjugation exchanges the two based
constructions as well. Since
$\overline{T_{E_u,+}}=T_{E_u,-}$,
\[
  \overline{\varsigma_{u,+}}=\varsigma_{u,-},
  \qquad
  \overline{\mathcal J_{u,+}}=\mathcal J_{u,-}.
\]

Thus the metric-induced Hermitian form determines the
compact $E_6$ groups without a basepoint, whereas a tangent
line determines the Euclidean Albert algebras and their
compact $F_4$ stabilizers. On a manifold, a smooth line
subbundle of the tangent bundle gives these based objects
fiberwise: local unit sections differ only by sign, so the
constructions agree on overlaps. In
Section~\ref{sec:reductions}, we globalize these pointwise
constructions. By themselves, the resulting bundle
reductions do not imply parallelism or a holonomy reduction.
\section{The exceptional reductions on a \texorpdfstring{$7$}{7}-manifold}
\label{sec:reductions}

In the preceding sections we established the exceptional nature of the
two critical cubics on a single copy of the deformation representation
$\WC$.
We now return to a $7$-manifold $(\mathcal M,\varphi)$ and ask what these pointwise
cubics define globally.

The exceptional reductions $\mathcal F_\pm$ and their compact forms
$\mathcal F^c_\pm$ are canonical on every $7$-manifold with a $G_2$-structure
(Sections~\ref{subsec:exceptional-reductions}
and~\ref{subsec:compact-reductions}). An Albert algebra bundle requires
an additional choice of a section with nonzero critical norm, called a
\emph{pointing}. A tangent line field supplies one; such fields exist
on every $7$-manifold without boundary, and the resulting reductions
are homotopic (Section~\ref{subsec:albert-bundles}). Finally, among
metric connections on $T\mathcal M$, exactly the $G_2$-connections
preserve the critical cubics. Orthogonal projection of the Levi-Civita connection
selects the minimal connection, while parallel pointings require
reduced holonomy (Section~\ref{subsec:connections-pointings}).
As in Section~\ref{sec:introduction},
no torsion assumption is made, except for the statements on Levi-Civita
parallelism and the compact non-flat example in
Remark~\ref{rem:geometric-octonions}, which assume torsion-freeness.
Compactness is required only for that example.

It is important to distinguish the projective and the exact statements.
A line of cubic forms determines a reduction to the group preserving that
line; elements of this group may rescale a cubic.
To obtain a reduction to the \emph{exact}
stabilizer $E_6(\C)$, one must select a nowhere-vanishing
cubic in the line. In our setting this selection is automatic: the fixed
normalizations of $C_0$ and $C_\varphi$ give the global critical fields
$N_\pm=C_0+\lambda_\pm C_\varphi$. Thus the exceptional reductions depend
only on $\varphi$ and the labeling of the two conjugate roots.

\subsection{The exceptional reductions}
\label{subsec:exceptional-reductions}

Recall from Definition~\ref{def:deformation-bundle} the real deformation
bundle $\mathcal W_\varphi=S^2_0T^*\mathcal M$, its complexification
$\mathcal W_{\varphi,\C}$, and the principal $G_2$-bundle
$\mathrm{Fr}_\varphi$ of adapted coframes, together with its
complexification $\mathrm{Fr}_{\varphi,\C}$. Every adapted coframe $p$ over
$x\in\mathcal M$ induces a complex-linear frame
\[
  f_p:\WC\longrightarrow\mathcal W_{\varphi,\C,x},
  \qquad f_{pO}=f_p\rho(O)\quad (O\in G_2).
\]
Here $\rho$ is the representation of Section~\ref{sec:deformation-representation},
extended to $G_2^{\C}$ when required. By Theorem~\ref{thm:critical-group},
the exact stabilizer of the model critical cubic $N_\pm$ is
\[
  G_\pm=\Stab_{GL(\WC)}(N_\pm)\cong E_6(\C),
\]
whereas the stabilizer of its line is
$\widehat G_\pm=\Stab_{GL(\WC)}(\C N_\pm)=G_\pm\cdot\C^\times$. Since both $C_0$ and $C_\varphi$ are $G_2$-invariant,
the representation $\rho$ induces embeddings
\[
  \jmath_\pm:G_2\hookrightarrow G_\pm,
  \qquad \jmath_\pm(O)=\rho(O).
\]
We use the same symbol $N_\pm$ for the fixed model cubic on $\WC$ and
its global realization; we write $N_{\pm,x}$ for the cubic on the fiber
over $x$.
Conjugation of frames means $\bar f(h)=\overline{f(\bar h)}$.

Define the corresponding cubic-preserving frame bundles by
\begin{equation}
  \mathcal F_\pm
  =\Bigl\{f:\WC\xrightarrow{\sim}
      \mathcal W_{\varphi,\C,x}\ \Bigm|\
      x\in\mathcal M,\ f^*N_{\pm,x}=N_\pm\Bigr\}.
  \label{eq:explicit-reductions}
\end{equation}

\begin{theorem}[Global critical cubics and exceptional reductions]
  \label{thm:exceptional-reductions}
  \leavevmode
  \begin{enumerate}
    \renewcommand{\labelenumi}{(\arabic{enumi})}
    \renewcommand{\theenumi}{\arabic{enumi}}
\item\label{item:exceptional-reductions-cubics}
The invariant cubics $C_0$ and $C_\varphi$ frame the trivial rank-two
subbundle
\[
  \mathcal C_\varphi
    =\mathbb R C_0\oplus\mathbb R C_\varphi
    \subset S^3\mathcal W_\varphi^*.
\]
The critical equation $a^2-2ab+8b^2=0$ defines two constant sections of
$\mathbb P(\mathcal C_{\varphi,\C})\cong\mathcal M\times\mathbb P^1(\C)$,
corresponding to the trivial complex line subbundles
\[
  \mathcal N_\pm=\C N_\pm,\qquad
  N_\pm=C_0+\lambda_\pm C_\varphi,\qquad
  \lambda_\pm=\frac{1\pm\ii\sqrt7}{8}.
\]
These lines are complex conjugate, with $N_-=\overline{N_+}$.
\item\label{item:exceptional-reductions-frames}
Each $\mathcal F_\pm$ is a canonical reduction of the complex frame
bundle of $\mathcal W_{\varphi,\C}$ to $G_\pm\cong E_6(\C)$, with
\[
  \mathrm{Fr}_\varphi\times_{\jmath_\pm}G_\pm
  \xrightarrow{\ \cong\ }\mathcal F_\pm,
  \qquad [p,L]\longmapsto f_pL,
\]
and
$\mathcal W_{\varphi,\C}
  \cong\mathcal F_\pm\times_{G_\pm}\WC$.
The lines $\mathcal N_\pm$ also define reductions to
their respective line stabilizers $\widehat G_\pm$.
\item\label{item:exceptional-reductions-conjugation}
The bundle $\mathrm{Fr}_{\varphi,\C}$ is a common $G_2^\C$-reduction
of $\mathcal F_+$ and $\mathcal F_-$, and conjugation of frames
exchanges the two exceptional reductions.
\item\label{item:exceptional-reductions-adjoint}
The infinitesimal stabilizers form adjoint bundles
\[
  \operatorname{ad}\mathcal F_\pm
    \cong\mathrm{Fr}_\varphi\times_{G_2}\mathfrak g_\pm,
  \qquad
  \mathfrak g_\pm=\rho(\mathfrak g_2^\C)\oplus\mathfrak m_\pm
    \cong\mathbf{14}\oplus\mathbf{64}.
\]
The $\mathbf{64}$ summand gives the infinitesimal symmetries of
$N_\pm$ beyond those of $\varphi$.
  \end{enumerate}
\end{theorem}

\begin{proof}
  The universal contractions defining $C_0$ and $C_\varphi$ patch
  under changes of adapted coframe and are pointwise independent by
  Lemma~\ref{lem:invariant-cubic-plane}. Thus they frame
  $\mathcal C_\varphi$, and Theorem~\ref{thm:stabilizer-trichotomy}
  gives the constant critical equation and its conjugate solutions.

  Let $p$ be a local $\varphi$-adapted coframe over
  $U\subset\mathcal M$. By invariance,
  $f_p^*N_{\pm,x}=N_\pm$, so
  \[
    \mathcal F_\pm|_U
    =\{f_pL:L\in G_\pm\}
    \cong U\times G_\pm.
  \]
  Since $f_{pO}=f_p\rho(O)$ for $O\in G_2$,
  these local identifications give the stated extension
  of structure group and associated-bundle isomorphism.
  The exact stabilizer is identified with $E_6(\C)$ in
  Theorem~\ref{thm:critical-group}.
  The line-stabilizer reduction follows by extension
  of structure group.

  Since $\rho(G_2^\C)\subset G_\pm$, the maps
  $[p,O]\mapsto f_p\rho(O)$ identify
  $\mathrm{Fr}_{\varphi,\C}$ with a reduction common to
  $\mathcal F_+$ and $\mathcal F_-$.
  Finally, the adapted frames $f_p$ are real and
  $N_-=\overline{N_+}$, so conjugation sends
  $f_pL\in\mathcal F_+$ to
  $f_p\overline L\in\mathcal F_-$.
  The adjoint-bundle identification follows from the same extension of
  structure group, with its splitting given by
  Proposition~\ref{prop:block-kernels}.
\end{proof}

Together with Proposition~\ref{prop:double-adjoint},
Theorems~\ref{thm:stabilizer-trichotomy} and~\ref{thm:critical-group},
and Theorem~\ref{thm:based-algebras}, this proves
Theorem~\ref{thm:intro-exceptional-pencil}.

The fixed normalization gives the canonical reductions above.
Changing it can change the subbundle of frames.
\begin{remark}[Normalization and rescaling]
  Any nowhere-vanishing section
  $s_\pm=aN_\pm$, where $a:\mathcal M\to\C^\times$ is smooth,
  defines a principal $G_\pm$-bundle
  \begin{equation}
    \mathcal F_\pm(s_\pm)
    =\Bigl\{f:\WC\xrightarrow{\sim}
        \mathcal W_{\varphi,\C,x}\ \Bigm|\
        x\in\mathcal M,\ f^*s_{\pm,x}=N_\pm\Bigr\}.
    \label{eq:normalized-section-reductions}
  \end{equation}
  This need not coincide with the canonical reduction
  $\mathcal F_\pm$. Indeed, if we choose $t$ locally with
  $t^3=a^{-1}$, the frames $f_pt$ trivialize
  $\mathcal F_\pm(s_\pm)$.
  Different choices of cube root differ by elements of
  the center $\mu_3\subset G_\pm$, so no global cube
  root is required.
\end{remark}

\paragraph{Geometric meaning and common $G_2$ origin.}
The two exceptional reductions realize the pointwise Albert norms
as smooth cubic fields on the same complexified deformation bundle.
They arise by extension of structure group from $\mathrm{Fr}_\varphi$,
so their existence imposes no additional topological obstruction
beyond the given $G_2$-structure.
Here we use \emph{extension of structure group} to describe
the passage from $G_2$ to the larger group $G_\pm \cong E_6(\C)$.
Viewed inside the full $GL(27,\C)$ frame bundle
of $\mathcal W_{\varphi,\C}$, however, the resulting
$G_\pm$-structure is a \emph{reduction} of structure group.
Since $\rho(G_2)\subset G_\pm\subset GL(\WC)$,
this reduction exists automatically.
 The new content is the
identification of the exact stabilizers of the two critical cubics as
$G_\pm\cong E_6(\C)$, while the passage from invariant
tensors to reductions is standard
(cf.\ \cite[\S3.1.1]{bryant-2003} for the standard approach in the \(G_2\) context).
The common $G_2^\C$-reduction of
Theorem~\ref{thm:exceptional-reductions}\,(\ref{item:exceptional-reductions-conjugation})
is not the full intersection of $\mathcal F_+$ and $\mathcal F_-$.
By Proposition~\ref{prop:common-stabilizer},
\[
  G_+\cap G_-=\mu_3\times\rho(G_2^\C),
\]
so $\mathcal F_+\cap\mathcal F_-$ has this structure group,
with $\mathrm{Fr}_{\varphi,\C}$ as its distinguished
identity-component reduction. Thus the central $\mu_3$ survives in
the intersection of the two exceptional reductions.
The real adapted frames also lie in both compact reductions
of Section~\ref{subsec:compact-reductions}, which do not
contain $\mathrm{Fr}_{\varphi,\C}$.

\subsection{The Hermitian form and the compact reductions}
\label{subsec:compact-reductions}

Apply the Hermitian form of Lemma~\ref{lem:hermitian-W}
fiberwise to $\mathcal W_{\varphi,\C}$. Its value on $h,k$ is
$\langle h,\bar k\rangle_W$, with conjugation fixing
$\mathcal W_\varphi$. The unitary frames inside each exceptional
reduction give its compact form.

\begin{proposition}[The compact exceptional reductions]
  \label{prop:compact-reductions}
Recall $K_\pm$ from~\eqref{eq:K-compactform}, which contains
$\jmath_\pm(G_2)=\rho(G_2)$ by
Proposition~\ref{prop:compact-form}\,(\ref{item:compact-form-group}), and put
  \[
    \mathcal F^c_\pm
    =\bigl\{f\in\mathcal F_\pm:
      \langle fh,\overline{fk}\rangle_W
      =\langle h,\bar k\rangle_W
      \text{ for all }h,k\in\WC\bigr\}.
  \]
  Then $\mathcal F^c_\pm$ is a principal $K_\pm$-subbundle of
  $\mathcal F_\pm$, and the map $[p,L]\mapsto f_pL$ gives
  \[
    \mathrm{Fr}_\varphi\times_{\jmath_\pm}K_\pm
      \cong\mathcal F^c_\pm,
    \qquad
    \mathcal W_{\varphi,\C}
      \cong\mathcal F^c_\pm\times_{K_\pm}\WC.
  \]
  Each $\mathcal F^c_\pm$ contains $\mathrm{Fr}_\varphi$ but does
  not contain $\mathrm{Fr}_{\varphi,\C}$. Conjugation of frames exchanges
  $\mathcal F^c_+$ and $\mathcal F^c_-$, covering conjugation
  $K_+\to K_-$.
\end{proposition}

\begin{proof}
Adapted coframes are Hermitian isometries, so a local adapted coframe over an
open set $U\subset\mathcal M$ identifies the unitary frames in
$\mathcal F_\pm|_U$ with $U\times K_\pm$. Conjugation preserves the
unitarity condition and exchanges the critical cubics, while
$G_2^\C$ contains nonunitary elements, so $\mathcal F^c_\pm$ does not
contain $\mathrm{Fr}_{\varphi,\C}$.
\end{proof}

The compact real form is classical
(cf.\ \cite[Theorem~3.9.2]{yokota-2009-exceptional}); here it is
identified inside the deformation module, and the passage to a bundle
uses the same standard reduction correspondence as
Theorem~\ref{thm:exceptional-reductions}.

Thus every $G_2$-structure gives two canonical compact-$E_6$ structures
on its complexified deformation bundle, exchanged by conjugation.

\subsection{Pointings by line fields and Albert algebra bundles}
\label{subsec:albert-bundles}

A global Albert algebra structure on
$\mathcal W_{\varphi,\C}$ can be obtained from a
\emph{pointing}, namely a smooth section
$e\in\Gamma(\mathcal W_{\varphi,\C})$ satisfying
$N_\pm(e)\ne0$ everywhere. By
Theorem~\ref{thm:based-algebras}, the pair
$(N_\pm/N_\pm(e),e)$ defines a smooth bundle of
complex Albert algebras, and hence an
$F_4(\C)$-reduction of the $GL(27,\C)$-frame bundle of
$\mathcal W_{\varphi,\C}$ preserving the normalized cubic.
To obtain local trivializations, choose a smooth cube root of $N_\pm(e)$
locally. Transitivity in Theorem~\ref{thm:based-algebras}\,(\ref{item:based-algebras-orbits})
and local sections of the orbit projection
$G_\pm\to G_\pm/F_4(\C)$, with stabilizer taken at a fixed model
norm-one basepoint, then give smooth Albert frames.
The $F_4(\C)$-reduction obtained in this way need not be a reduction
of the original exceptional bundle $\mathcal F_\pm$,
whose frames preserve $N_\pm$ itself.

More generally, consider a smooth complex line
subbundle $\mathcal{L}\subset\mathcal W_{\varphi,\C}$ on
which $N_\pm$ is nowhere zero. Fix a model basepoint
$e_0\in\WC$ with $n_0=N_\pm(e_0)\ne0$.
The stabilizer of its line in $G_\pm$ is
\[
  \operatorname{Stab}_{G_\pm}(\C e_0)
  \cong F_4(\C)\times\mu_3,
\]
where $\mu_3$ is the group of cube roots of unity,
acting by scalar multiplication. Thus $\mathcal{L}$ determines
an $F_4(\C)\times\mu_3$-reduction of $\mathcal F_\pm$.

The restriction of $N_\pm$ to $\mathcal L$ is
nowhere zero, and therefore trivializes
$(\mathcal L^*)^{\otimes3}$. Consequently,
\[
  \mathcal L^{\otimes3}\cong\underline{\C}, \qquad \text{with}\qquad
  \qquad 3c_1(\mathcal L)=0.
\]
The local sections of $\mathcal L$ satisfying
$N_\pm(e)=n_0$ form a principal $\mu_3$-bundle.
Its class
\[
  \alpha_\pm(\mathcal L)\in H^1(\mathcal M;\mathbb Z/3)
\]
is the obstruction to a compatible
$F_4(\C)$-reduction of $\mathcal F_\pm$.
This class detects both possible three-torsion
in $c_1(\mathcal L)$ and nontrivial monodromy
of the cubic normalization when $\mathcal L$
is trivial. In the latter case, if
$\mathcal L=\C e$ for a global nowhere-vanishing
section $e$, the obstruction is the mod-three
winding class of $N_\pm(e)$. It vanishes
precisely when $N_\pm(e)/n_0$ admits a global
smooth cube root. Thus a constant critical norm
is sufficient, but not necessary, for an
$F_4(\C)$-reduction of the original exceptional
bundle compatible with $\mathcal L$.

A real pointing $e\in\Gamma(\mathcal W_\varphi)$
with nonzero critical norm defines conjugate Albert
algebra structures for both signs, since
$N_-(e)=\overline{N_+(e)}$. The $\mu_3$ ambiguity may
persist for real pointings, as the critical norms
may have nontrivial winding even on real sections.
The corresponding obstruction classes satisfy
\[
  \alpha_-(\C e)=-\alpha_+(\C e).
\]
Consequently, a compatible $F_4(\C)$-reduction
exists for one sign if and only if it exists
for the other.

Table~\ref{tab:pointing-choices} distinguishes the data fixed in these
constructions from the unrestricted reduction problem.
\begin{table}[ht]
  \centering
  \begin{tabular}{@{}p{0.29\linewidth}p{0.67\linewidth}@{}}
    \hline
    Fixed data & Result or obstruction \\
    \hline
    A section $e$ with $N_\pm(e)\ne0$
      & The pair $(N_\pm/N_\pm(e),e)$ defines an Albert algebra bundle.
        Its frames preserve the normalized cubic. \\[3pt]
    A line $\mathcal L$ on which $N_\pm$ is nowhere zero, and a nonzero constant $n_0$
      & A section $e\in\Gamma(\mathcal L)$ with $N_\pm(e)=n_0$ exists
        exactly when the principal $\mu_3$-bundle of normalized sections has
        $\alpha_\pm(\mathcal L)=0$. It gives a compatible
        $F_4(\C)$-reduction of $\mathcal F_\pm$. \\[3pt]
    The exceptional bundle $\mathcal F_\pm$, with no prescribed line
      & A norm-one pointing exists and is unique up to homotopy
        over a $7$-manifold, by Proposition~\ref{prop:no-obstruction}.
        Equivalently, the unrestricted $F_4(\C)$-reductions are
        homotopy-unique. \\
    \hline
  \end{tabular}
  \caption{Three global pointing questions.}
  \label{tab:pointing-choices}
\end{table}

A natural geometric source of real pointings is a
\emph{line field}, by which we mean a smooth real
line subbundle $\ell\subset T\mathcal M$. Let $P_\ell$
denote the orthogonal projection onto $\ell$,
viewed as a symmetric covariant tensor using
$g_\varphi$, and define
\begin{equation}\label{eq:El-def}
E_\ell=7P_\ell-g_\varphi
  \in\Gamma(\mathcal W_\varphi).
\end{equation}
Locally, $P_\ell=P_u=u^\flat\otimes u^\flat$ and
$E_\ell=E_u$ for a unit section $u$ of $\ell$, as in
\eqref{eq:projection-tensors} and \eqref{eq:test-tensors}.
At each point, $E_\ell$ is a basepoint of $SU(3)$ type.
The stabilizer of an unoriented line in the compact
$G_2$ is $SU(3)\rtimes\mathbb Z/2$, with identity
component $SU(3)$ (Proposition~\ref{ex:structured-basepoints}). Accordingly, the construction
requires no orientation of $\ell$, although a
corresponding global $SU(3)$-reduction of the
$G_2$-bundle would require such an orientation.

Finally, by the transport~\eqref{eq:basepoint-transport}
and Remark~\ref{rem:rank-one-normal-form},
every basepoint is $G_\pm$-equivalent, up to a
nonzero scalar, to one of the form
$E_u=7P_u-g$. This transport is pointwise and uses the
complex group $G_\pm$, so it need not assemble into a global
gauge transformation (see the discussion after
Corollary~\ref{cor:basepoint-transport}). Real line fields
therefore provide a canonical geometric construction
of constant-norm pointings, but do not exhaust
all possible pointings.

Fix the model unit vector $e_1\in V$ and recall
$E_1=7P_1-g$ from \eqref{eq:test-tensors}.
Write $H_{E_1}=\Stab_{G_2}(E_1)\cong SU(3)\rtimes\mathbb Z/2$,
as in Proposition~\ref{ex:structured-basepoints}.
Define the bundle of coframes adapted to $\ell$ by
\begin{equation}
  \label{eq:line-adapted-frames}
  \mathrm{Fr}_\ell
    =\{p\in\mathrm{Fr}_\varphi|_x:x\in\mathcal M,\ f_p(E_1)=E_\ell(x)\}.
\end{equation}
Equivalently, the dual frame of $p$ sends the model line
$\mathbb R e_1$ to $\ell_x$. The transitive action of $G_2$ on
unit vectors ensures that such coframes exist locally
\cite[Theorem~1.9.2, p.~17]{yokota-2009-exceptional}.
If an orientation of $\ell$ is fixed by a unit section $u$,
requiring instead that the dual frame send $e_1$ to $u_x$ gives
the principal $SU(3)$-subbundle of $\mathrm{Fr}_\ell$.

\begin{lemma}[The basepoint of a line field]
  \label{lem:line-field-basepoint}
  \leavevmode
  \begin{enumerate}
    \renewcommand{\labelenumi}{(\roman{enumi})}
    \renewcommand{\theenumi}{\roman{enumi}}
    \item\label{item:line-field-existence} Every $7$-manifold without boundary carries a smooth line field.
    \item\label{item:line-field-pointing} For every line field $\ell$,
      $E_\ell=7P_\ell-g_\varphi$ is a smooth global section of
      $\mathcal W_\varphi$, independent of the orientation of $\ell$, and
      \[
        N_\pm(E_\ell)=n_\pm=42(5+2\lambda_\pm)\ne0.
      \]
    \item\label{item:line-field-frames} For the model basepoint $E_1=7P_1-g$, the bundle $\mathrm{Fr}_\ell$ is a principal
      $H_{E_1}$-subbundle
      of $\mathrm{Fr}_\varphi$.
  \end{enumerate}
\end{lemma}

\begin{proof}
The first assertion follows componentwise from the classical
existence of a nowhere-vanishing section of $T\mathcal M$; see
\cite[Corollary~14.7 and the following paragraph, p.~145]{kreck-2010-differential}
for closed odd-dimensional manifolds, including nonorientable ones,
and \cite[Chapter~5, Section~2, Exercise~6(c), p.~139]{hirsch-1976-differential}
for the connected noncompact case.
Since $u^\flat\otimes u^\flat$ is invariant under $u\mapsto-u$,
$E_\ell$ is globally defined even when $\ell$ is nonorientable.
Its constant norm follows from~\eqref{eq:su3-norm-trace}.
For the final assertion, a local adapted coframe identifies the
frames taking $E_1$ to $E_\ell$ with the stabilizer $H_{E_1}$.
\end{proof}

Thus every real line field determines global
Albert algebra structures and
$F_4(\C)$-reductions of both original exceptional
bundles $\mathcal F_\pm$. In Table~\ref{tab:pointing-choices},
this is an explicit constant-norm solution for the line
$\mathcal L=\C E_\ell$, with $n_0=n_\pm$.

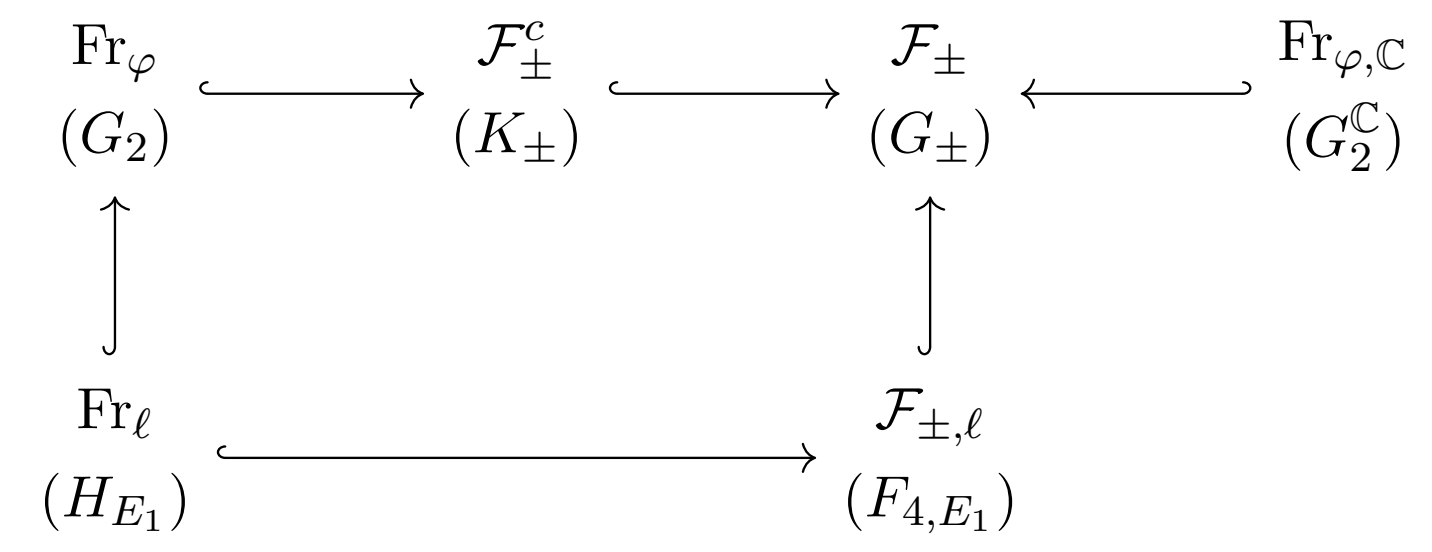
\begin{figure}[ht]
  \centering
  \begin{tikzcd}[column sep=large, row sep=large]
    \begin{gathered}\mathrm{Fr}_\varphi\\(G_2)\end{gathered}
      \arrow[r, hook]
    & \begin{gathered}\mathcal F^c_\pm\\(K_\pm)\end{gathered}
      \arrow[r, hook]
    & \begin{gathered}\mathcal F_\pm\\(G_\pm)\end{gathered}
    & \begin{gathered}\mathrm{Fr}_{\varphi,\C}\\(G_2^\C)\end{gathered}
      \arrow[l, hook'] \\
    \begin{gathered}\mathrm{Fr}_\ell\\(H_{E_1})\end{gathered}
      \arrow[u, hook] \arrow[rr, hook]
    & {}
    & \begin{gathered}\mathcal F_{\pm,\ell}\\(F_{4,E_1})\end{gathered}
      \arrow[u, hook]
    & {}
  \end{tikzcd}
  \caption{Inclusions of principal bundles, with their structure
groups. The lower row requires a line subbundle \(\ell\) of \(T\mathcal M\) and is constructed in
  Lemma~\ref{lem:line-field-basepoint} and Theorem~\ref{thm:albert-bundle}
  (Section~\ref{subsec:albert-bundles}).}
  \label{fig:structure-groups}
\end{figure}

Figure~\ref{fig:structure-groups} shows these reductions together with
the pointed reductions constructed below.

For each tangent line we write the classical based Albert operations
explicitly on the complex deformation bundle. We use the real sector
bundles $\mathcal W_r(\ell)\subset\mathcal W_\varphi$, $r=6,8,12$,
obtained from the splitting~\eqref{eq:su3-sectors}; their
complexifications will be denoted by $\mathcal W_r(\ell)_\C$.

\paragraph{The decomposition and trace operators.}
Fix a line field $\ell\subset T\mathcal M$, and let $E_\ell$ and $n_\pm$
be as above. We suppress the chosen root on $T_{E_\ell}$ throughout,
restoring signs on real structures when both occur.
The pointwise decomposition and projectors of
Lemma~\ref{lem:su3-sectors} apply fiberwise to $\ell$.
For a local unit section $u$ of $\ell$, use the projectors
$\Pi_1,\Pi_6,\Pi_8,\Pi_{12}$ of~\eqref{eq:su3-projectors}, with
$P_u=P_\ell$ and $E_u=E_\ell$. Their formulas depend smoothly on
$g_\varphi$, $\varphi$ and $P_\ell$ and are invariant under
$u\mapsto-u$, so their images define smooth global sector bundles.
The based trace and trace form are given by
\eqref{eq:su3-norm-trace} and~\eqref{eq:based-trace-form};
the trace operator $T_{E_\ell}$ and its inverse have the sector
expressions~\eqref{eq:su3-trace-projectors}.
Define the pointed frame subbundle and based algebra bundle by
\[
  \mathcal F_{\pm,\ell}
    =\{f\in\mathcal F_\pm|_x:x\in\mathcal M,\ f(E_1)=E_\ell(x)\},\qquad
  \mathfrak J_{E_\ell}
    =\bigl(\mathcal W_{\varphi,\C},N_\pm/n_\pm,
          E_\ell,\circ_{E_\ell}\bigr).
\]
For the explicit operations, write $\lambda=\lambda_\pm$ and
$\bar\lambda=\lambda_\mp$. We use $F_\nu=A+\nu B$ also on symmetric
tensors with nonzero trace, with $A$ and $B$ as defined
in~\eqref{eq:trace-product}--\eqref{eq:phi-product}.

\begin{theorem}[The Albert algebra bundle of a line field]
  \label{thm:albert-bundle}
 Let $\ell\subset T\mathcal M$ be a line field. Then
  \[
    \mathcal F_{\pm,\ell}
      \cong\mathrm{Fr}_\ell\times_{H_{E_1}}F_{4,E_1}
  \]
  is a principal $F_{4,E_1}\cong F_4(\C)$-subbundle
  of $\mathcal F_\pm$.
  The based algebras $\mathfrak J_{E_\ell}$
  form a smooth bundle of complex Albert algebras,
  locally trivialized by sections of $\mathrm{Fr}_\ell$.

  Their adjoint and Jordan product are given explicitly by
  \begin{equation}
    \label{eq:line-field-albert-operations}
    \begin{aligned}
      h^\sharp
        &=\frac{3+4\lambda}{56}F_\lambda(h)
          -\frac12F_{\bar\lambda}(P_\ell,F_\lambda(h))\\
        &\quad+\frac3{n_\pm}\langle E_\ell,F_\lambda(h)\rangle_WE_\ell,\\
      h\circ_{E_\ell}k
        &=\frac{3+4\lambda}{56}F_\lambda(h,k)
          -\frac12F_{\bar\lambda}(P_\ell,F_\lambda(h,k))\\
        &\quad+\frac14\bigl(h(u,u)k+k(u,u)h\bigr).
    \end{aligned}
  \end{equation}
  Here $u$ is any local unit section of $\ell$.
These operations are polynomial contractions of
  $g_\varphi$, $\varphi$ and $P_\ell$, and are defined
  also at noninvertible elements.
\end{theorem}

\begin{proof}
  By~\eqref{eq:g2-basepoint-intersection},
  $H_{E_1}\subset F_{4,E_1}$, and by
  Proposition~\ref{prop:f4-stabilizer}, we have
  $F_{4,E_1}\cong F_4(\C)$.
  A local section $p$ of $\mathrm{Fr}_\ell$
  induces a frame $f_p\in\mathcal F_\pm$ with
  $f_p(E_1)=E_\ell$.
  Every other such frame is of the form
  $f_pL$, with $L\in F_{4,E_1}$.
  This proves the principal-bundle isomorphism.
  Theorem~\ref{thm:based-algebras} then gives
  the smooth Albert algebra bundle.

  For the Jordan operations, recall
  $F_\lambda(E_\ell)=(5+2\lambda)E_\ell$, and
  $n_\pm=42(5+2\lambda)$ from \eqref{eq:su3-norm-trace}.
  For traceless $\mathsf X$, the definition~\eqref{eq:trace-product} and the identity \eqref{eq:p-identity-rule} give
  $A(g_\varphi,\mathsf X)=\mathsf X$ and $B(g_\varphi,\mathsf X)=-\mathsf X$, hence
  $F_{\bar\lambda}(g_\varphi,\mathsf X)=(1-\bar\lambda)\mathsf X$.
  Since $E_\ell=7P_\ell-g_\varphi$ and
  $\bar\lambda=\tfrac14-\lambda$, it follows that
  \[
    -\frac3{n_\pm}F_{\bar\lambda}(F_\lambda(E_\ell),\mathsf X)
      =-\frac1{14}F_{\bar\lambda}(E_\ell,\mathsf X)
      =\frac{3+4\lambda}{56}\mathsf X-\frac12F_{\bar\lambda}(P_\ell,\mathsf X).
  \]
  Apply this identity with $\mathsf X=F_\lambda(h)$ and $\mathsf X=F_\lambda(h,k)$
  in~\eqref{eq:based-albert-operations}, and use
  $\tr_{E_\ell}(h)=\tfrac12h(u,u)$. This proves
  \eqref{eq:line-field-albert-operations} directly.
  The formulas depend on the local unit section $u$ only through
  $P_\ell$ and evaluations such as $h(u,u)$, so they are globally defined.
\end{proof}

The reconstruction of a Jordan product from a normalized
cubic and a basepoint is classical, due to Springer
\cite[p.~261, equations~(9)--(12)]{springer-1962-characterization},
and was recalled in Section~\ref{subsec:albert-norm}.
Here we realize the classical construction globally,
with explicit tensor operations determined by the line field.

On $(\mathcal M,\varphi,\ell)$ the exceptional Jordan product is therefore a
natural tensor operation.

The Euclidean real forms of Section~\ref{sec:e6-albert-norm} also
transfer to a manifold. We keep the root as a subscript on the real structure
$\varsigma_{\ell,\pm}$ and its fixed bundle $\mathcal J_{\ell,\pm}$,
since both signs occur below.

\paragraph{The Euclidean real structure.}
Write $s\mapsto\bar s$ for the natural complex conjugation on
$\mathcal W_{\varphi,\C}=\mathcal W_\varphi\otimes_\R\C$.
For either critical parameter, define the antilinear map
\begin{equation}
  \begin{aligned}
    \varsigma_{\ell,\pm}:
      \Gamma(\mathcal W_{\varphi,\C})
        &\longrightarrow\Gamma(\mathcal W_{\varphi,\C}),\\
    s&\longmapsto
      \tfrac1{14}T_{E_\ell}^{-1}(\bar s).
  \end{aligned}
  \label{eq:global-euclidean-real-structure}
\end{equation}
Here $T_{E_\ell}$ is the trace operator for the chosen
root, and the bar denotes fiberwise conjugation.
Thus, for $a\in C^\infty(\mathcal M,\C)$,
\[
  \varsigma_{\ell,\pm}(as)
    =\bar a\,\varsigma_{\ell,\pm}(s).
\]
Equivalently, using the fiberwise projectors of~\eqref{eq:su3-projectors},
\[
  \varsigma_{\ell,\pm}(s)
    =\Pi_1\bar s-\Pi_6\bar s
      +\frac{-1\pm3\ii\sqrt7}{8}\Pi_8\bar s
      +\frac{3\mp\ii\sqrt7}{4}\Pi_{12}\bar s.
\]
The next proposition shows that this map is an
involutive automorphism of the based Albert algebra
bundle, with a Euclidean Albert algebra as its fixed
real form.

\begin{proposition}[The Euclidean real form]
  \label{prop:albert-real-form}
  Fix either sign and a line field $\ell$.
  The map $\varsigma_{\ell,\pm}$ is an antilinear
  involutive automorphism of $\mathfrak J_{E_\ell}$.
  Its fixed bundle is
  \[
    \mathcal J_{\ell,\pm}
      =\R E_\ell\oplus\ii\mathcal W_6(\ell)
       \oplus a_{8,\pm}\mathcal W_8(\ell)
       \oplus a_{12,\pm}\mathcal W_{12}(\ell),
  \]
  where the coefficients $a_{8,\pm}$ and $a_{12,\pm}$
  are those of Proposition~\ref{prop:euclidean-form}.
  In particular,
  \[
    \Gamma(\mathcal J_{\ell,\pm})
      =\{s\in\Gamma(\mathcal W_{\varphi,\C}):
          \varsigma_{\ell,\pm}(s)=s\}.
  \]
  Each fiber is a Euclidean Albert algebra with
  determinant $N_\pm/n_\pm$, and complex conjugation
  exchanges $\mathcal J_{\ell,+}$ and $\mathcal J_{\ell,-}$.

  At $\lambda_+$, the compact reduction
  $\mathcal F^c_+$ consists precisely of the frames
  $f\in\mathcal F_+$ satisfying
  \[
    \langle fh,\varsigma_{\ell,+}(fk)\rangle_{E_\ell}
      =\langle h,\varsigma_{e_1,+}k\rangle_{E_1}
\qquad (h,k\in\WC),
\]
where $\langle\cdot,\cdot\rangle_{E_1}$ is the based trace form
at the model basepoint $E_1$.
\end{proposition}

\begin{proof}
  A local section of $\mathrm{Fr}_\ell$ transfers
  Lemma~\ref{lem:euclidean-real-structure} and
  Proposition~\ref{prop:euclidean-form} to each fiber.
  These transfers agree because
  $\varsigma_{Au,\pm}\rho(A)=\rho(A)\varsigma_{u,\pm}$ for $A\in G_2$,
  as recalled at the end of Section~\ref{subsec:euclidean-albert}.
  Since the projectors and trace operator
  are smooth and invariant under $u\mapsto-u$,
  the resulting involutions and their fixed spaces
  are globally defined. Complex conjugation exchanges
  the trace operators for the two critical parameters,
  and hence exchanges the fixed real forms
  $\mathcal J_{\ell,+}$ and $\mathcal J_{\ell,-}$.

  Finally, Lemma~\ref{lem:euclidean-real-structure}
  gives the Hermitian identity
  \[
    14\langle h,\varsigma_{\ell,\pm} k\rangle_{E_\ell}
      =\langle h,\bar k\rangle_W.
  \]
  Thus the displayed frame condition at $\lambda_+$
  is equivalent to preservation of the Hermitian
  form induced by the $W$-pairing, which defines
  $\mathcal F^c_+$.
\end{proof}

Proposition~\ref{prop:albert-real-form} applies the pointwise results of
Proposition~\ref{prop:euclidean-form}, and their classical sources,
fiberwise; it involves no new classification of real forms. As in
Proposition~\ref{prop:euclidean-form}, the Euclidean real form meets
the real deformation bundle only in $\R E_\ell$.

Together, Lemma~\ref{lem:line-field-basepoint},
Theorem~\ref{thm:albert-bundle} and
Proposition~\ref{prop:albert-real-form} prove
Theorem~\ref{thm:intro-albert-bundles}.

We now consider the unrestricted reduction problem in the last row of
Table~\ref{tab:pointing-choices}. The choice of a tangent line field
changes the pointing geometrically, but does not change the homotopy
class of the resulting $F_4(\C)$-reduction.

\begin{proposition}[Homotopy uniqueness]
  \label{prop:no-obstruction}
  For either sign,
  \[
    \{N_\pm=1\}\cong E_6(\C)/F_4(\C)\simeq E_6/F_4,
  \]
  where $E_6/F_4$ denotes the quotient of the compact simply connected groups of these types, a compact symmetric space that is $8$-connected
  and has $\pi_9\cong\mathbb Z$. On a $7$-manifold every
  $F_4(\C)$-reduction of $\mathcal F_\pm$ is homotopic to
  $\mathcal F_{\pm,\ell}$, and its principal-bundle isomorphism
  class is independent of $\ell$.
\end{proposition}

The homotopy equivalence follows from the polar decompositions
\cite[Theorems~2.13.2 and~3.14.1]{yokota-2009-exceptional}; the
connectivity and $\pi_9$ follow from Conlon's cell structure
\cite[Theorem~1.1, p.~411, the following cell-structure remarks on p.~412, and the report on p.~414]{conlon-1966-application}.
Obstruction theory for sections then gives homotopy uniqueness over
a $7$-manifold. Thus an unrestricted norm-one pointing exists
and is unique up to homotopy; compatibility with a prescribed
complex line remains subject to the obstruction described above.

We can also treat all tangent lines at once.

\begin{remark}[The tautological family]
  Let $\pi:\mathbb P(T\mathcal M)\to\mathcal M$ be the real
  projectivized tangent bundle, with tautological
  line subbundle
  $\ell^{\mathrm{taut}}\subset\pi^*T\mathcal M$.
  Applying the preceding constructions fiberwise
  to $\ell^{\mathrm{taut}}$ gives a canonical smooth
  family of based complex Albert algebras with
  Euclidean real forms on $\pi^*\mathcal W_{\varphi,\C}$.
  Every line field $\ell\subset T\mathcal M$ determines a
  section $s_\ell:\mathcal M\to\mathbb P(T\mathcal M)$, and pulling
  back this family by $s_\ell$ recovers the
  Albert algebra bundle associated with $\ell$.
\end{remark}

The $G_2$-structure also supplies a geometric octonion bundle.
We compare its entrywise Hermitian construction with the Albert bundle above.
\begin{remark}[The geometric octonion bundle]
  \label{rem:geometric-octonions}
  Given the geometric octonion bundle
  $\O_\varphi=\underline\R\oplus T\mathcal M$, with its $\varphi$-induced
  product \cite{grigorian-2017-octonion-bundles},
  one might consider the Hermitian-matrix bundle
  $H_3(\O_\varphi)$. However, its underlying real bundle
  is $6\underline\R\oplus3T\mathcal M$, of $G_2$-type
  $6\cdot\mathbf1\oplus3\cdot\mathbf7$, whereas
  $\mathcal W_\varphi$ has irreducible type $\mathbf{27}$. The Whitney sum formula and the splitting principle for
characteristic classes give the integral identities
\[
  p_1(6\underline\R\oplus3T\mathcal M)=3p_1(T\mathcal M),
  \qquad
  p_1(S^2_0T^*\mathcal M)=9p_1(T\mathcal M).
\]
The second identity is the rank-seven case of
$p_1(S^2_0E^*)=(\operatorname{rank}E+2)p_1(E)$
for an oriented real vector bundle $E$. Thus the bundles differ even after complexification whenever
  $6p_1(T\mathcal M)\ne0$. In particular, this holds on every compact non-flat torsion-free
  $G_2$-manifold: Joyce's formula
  \cite[Proposition~10.2.7, p.~246]{Joycebook} gives
  $(p_1(T\mathcal M)\cup[\varphi])[\mathcal M]<0$, so $p_1(T\mathcal M)$ is not torsion.
\end{remark}

\subsection{Connections and parallel structures}
\label{subsec:connections-pointings}

We now determine which connections preserve these exceptional
structures and show that the metric singles out one of them.
Write $D$ for a connection on $T\mathcal M$ and for the induced connections
on the deformation bundles and frame bundles. In particular,
$D^{LC}$ is the Levi-Civita connection. Fix the intrinsic torsion convention
\[
  D^{LC}_X\varphi=\iota_{T(X)}\psi,
  \qquad
  D^{LC}_m\varphi_{abc}=T(e_m)^r\psi_{rabc},
\]
as in Bryant's description of intrinsic torsion and its four classes
\cite[\S4]{bryant-2003}. Extend the notation
$\Theta_v(w)=w\times v$ to arbitrary tangent vectors $v$.
From the decomposition $\lieso(7)=\gtwo\oplus\Lambda^2_7$, write every
metric connection on $T\mathcal M$ uniquely as
\begin{equation}
  D=D^{LC}+\beta+\Theta_{w(\cdot)},\qquad
  \beta\in\Omega^1(\mathcal M,\mathfrak g_2),\quad
  w\in\Omega^1(\mathcal M,T\mathcal M),
  \label{eq:metric-connection-decomposition}
\end{equation}
where the Lie-algebra-valued forms are read in adapted frames.
The minimal connection is
\[
  D^{\min}=D^{LC}-\tfrac13\Theta_{T(\cdot)}.
\]
It is metric and preserves $\varphi$ and $\psi$  \cite[\S 3.1]{CleytonIvanovCurv}.

Among metric connections on $T\mathcal M$, preservation of an exceptional
cubic is precisely the $G_2$-compatibility condition. First we need a deformation property of \(C_\varphi\).
With the metric held fixed, define the infinitesimal variation by
\[
  L_\varphi:T\mathcal M\longrightarrow S^3\mathcal W_\varphi^*,
  \qquad L_\varphi(v)=\left.\frac{d}{dt}\right|_{t=0}
       C_{\varphi+t\iota_v\psi}.
\]

\begin{lemma}[Variation of the $\varphi$-cubic]
  \label{lem:phi-cubic-variation}
  The map $L_\varphi$ is $G_2$-equivariant and fiberwise injective,
  and its image is an
  irreducible $\mathbf 7$ component. In particular,
  \[
    \operatorname{im}L_\varphi
      \cap\operatorname{span}_{\C}\{C_0,C_\varphi\}
      =\{0\}.
  \]
\end{lemma}

\begin{proof}
  We work in a $G_2$-adapted orthonormal frame.
  This deformation keeps the metric fixed to first order, so differentiating the contraction \eqref{eq:natural-cubics} defining $C_\varphi$ with a constant metric gives
  \[
    L_\varphi(v)(h_1,h_2,h_3)
      =2(\iota_v\psi)_{a_1a_2a_3}\varphi_{b_1b_2b_3}
        h_1^{a_1b_1}h_2^{a_2b_2}h_3^{a_3b_3}.
  \]
  This expression is $G_2$-equivariant since it involves
  only the invariant tensors $g$, $\varphi$ and $\psi$.

  Recall that $\lieso(7)$ decomposes as
  $\mathfrak g_2\oplus\Lambda^2_7$, where $\Lambda^2_7$
  consists of the maps $\Theta_v$. The standard contraction
  identities give
  \[
    \theta_{\Theta_v}=-3\iota_v\psi,
  \]
  where $\theta_M$ is the $3$-form \eqref{eq:insert-M}. Indeed,
  \[
    (\theta_{\Theta_v})_{ijk}
      =v^b\bigl(
        \varphi_{iba}\varphi^a{}_{jk}
+\varphi_{jba}\varphi_i{}^a{}_k
+\varphi_{kba}\varphi_{ij}{}^a
      \bigr)
      =-3v^b\psi_{bijk}.
  \]

  Set $R_t=\exp(t\Theta_v)$. By naturality under
  orthogonal transformations,
  \[
    C_{R_t^*\varphi}(h_1,h_2,h_3)
      =C_\varphi\bigl(
        \rho(R_t)h_1,\rho(R_t)h_2,\rho(R_t)h_3
      \bigr).
  \]
  Recall the infinitesimal action map $\Phi_N$ of~\eqref{eq:action-map}.
  Differentiation at $t=0$ gives
  \[
    -3L_\varphi(v)
      =\Phi_{C_\varphi}(\rho(\Theta_v)).
  \]
  By Proposition~\ref{prop:phi-cubic-stabilizer},
  $\ker\Phi_{C_\varphi}=\rho(\mathfrak g_2^\C)$.
  Since $\rho$ is injective and
  $\Theta_v\in\Lambda^2_7$, the right-hand side
  vanishes only when $v=0$. Hence $L_\varphi$
  is injective.

  Its image is therefore isomorphic to the irreducible
  $\mathbf 7$ representation of $G_2$. The invariant
  plane spanned by $C_0$ and $C_\varphi$ consists of
  trivial representations, so its intersection with
  $\operatorname{im}L_\varphi$ is zero.
\end{proof}

The variation formula determines the critical-cubic derivative of
any connection induced by a metric connection on $T\mathcal M$.

\begin{proposition}[Metric connections and critical cubics]
  \label{prop:metric-critical-connections}
  Let $D$ be a metric connection on $T\mathcal M$, written as in
  \eqref{eq:metric-connection-decomposition}, and use the same
  symbol for its induced connection on $\mathcal W_{\varphi,\C}$. Then
      \begin{equation}
        \label{eq:metric-critical-derivative}
        D_XN_\pm
          =\lambda_\pm L_\varphi(T(X)+3w(X)).
      \end{equation}
\end{proposition}

\begin{proof}
Since $\beta_X$ annihilates $\varphi$ and
$\theta_{\Theta_v}=-3\iota_v\psi$, we have
$D_X\varphi=\iota_{T(X)+3w(X)}\psi$.
Metric compatibility gives $DC_0=0$, and since $Dg=0$, $D_XC_\varphi$ is the
variation of $C_\varphi$ in the direction $D_X\varphi$ with the metric fixed.
The definition of $L_\varphi$ therefore gives~\eqref{eq:metric-critical-derivative}.
\end{proof}

The injectivity and representation type of $L_\varphi$ turn this
derivative formula into preservation criteria.
\begin{corollary}[Critical-cubic preservation criteria]
  \label{cor:critical-connection-criteria}
  \begin{enumerate}
    \renewcommand{\labelenumi}{(\arabic{enumi})}
    \renewcommand{\theenumi}{\arabic{enumi}}
    \item\label{item:metric-critical-equivalences}
      For a metric connection $D$ on $T\mathcal M$, the following conditions are equivalent:
      \begin{enumerate}
        \item[(a)] $D\varphi=0$;
        \item[(b)] $DN_\pm=0$;
        \item[(c)] $D$ preserves the line $\C N_\pm$.
      \end{enumerate}
      In particular, a metric connection preserves
      either critical cubic if and only if it is a
      $G_2$-connection, in which case it preserves both
      critical cubics.

    \item\label{item:metric-critical-levi-civita} For the Levi-Civita connection,
      \[
        D^{LC}_XN_\pm
          =\lambda_\pm L_\varphi(T(X)).
      \]
      Thus $D^{LC}N_\pm=0$ if and only if $T=0$,
      that is, if and only if $\varphi$ is torsion-free.

    \item\label{item:metric-critical-line-connections} For an arbitrary complex connection on
      $\mathcal W_{\varphi,\C}$ preserving the line
      $\C N_\pm$, its local connection form in a
      $G_2$-adapted frame is $\omega+\alpha\Id_W$, where
      $\omega$ is $\mathfrak g_\pm$-valued and $\alpha$
      is a complex-valued one-form. Then
      \[
        DN_\pm=-3\alpha\otimes N_\pm.
      \]
      Hence preserving the line $\C N_\pm$ is in general
      weaker than preserving the cubic $N_\pm$ itself.
  \end{enumerate}
\end{corollary}

\begin{proof}
  By the preceding proof, $D\varphi=0$ if and only if $T+3w=0$.
  Since $L_\varphi$ is injective and
  $\lambda_\pm\ne0$, we have $DN_\pm=0$ if and
  only if $T+3w=0$, equivalently,
  $D\varphi=0$.

  Suppose instead that $D$ preserves the line
  $\C N_\pm$. Then $D_XN_\pm$ is proportional
  to $N_\pm$ for every $X$. However, by
  \eqref{eq:metric-critical-derivative}, $D_XN_\pm$ lies in
  $\operatorname{im}L_\varphi$, which has zero
  intersection with $\C N_\pm$ by
  Lemma~\ref{lem:phi-cubic-variation}.
  Hence $D_XN_\pm=0$, and the preceding argument
  gives $D\varphi=0$. Conversely, a
  $G_2$-connection preserves both invariant cubics,
  and therefore both critical cubics.

  Part~\ref{item:metric-critical-levi-civita} follows by taking $\beta=w=0$.
  The injectivity of $L_\varphi$ then gives
  $D^{LC}N_\pm=0$ precisely when $T=0$.

  For part~\ref{item:metric-critical-line-connections}, recall that the infinitesimal stabilizer
  of the line $\C N_\pm$ is
  $\mathfrak g_\pm\oplus\C\Id_W$.
  In a $G_2$-adapted frame, the components of
  $N_\pm$ are constant, and the induced action
  of a connection form on covariant cubic tensors
  is the negative of $\Phi_{N_\pm}$.
  Since
  \[
    \Phi_{N_\pm}(\omega)=0,
    \qquad
    \Phi_{N_\pm}(\Id_W)=3N_\pm,
  \]
  we obtain
  \[
    DN_\pm=-3\alpha\otimes N_\pm.
  \]
  Thus the exact cubic is parallel precisely when
  $\alpha=0$.
\end{proof}

Beyond connections induced from $T\mathcal M$, compatible exceptional
connections have additional freedom in the $\mathbf{64}$ directions.
Recall from
      Proposition~\ref{prop:block-kernels} that
      \[
        \mathfrak g_\pm
          =\rho(\mathfrak g_2^\C)\oplus\mathfrak m_\pm,
      \]
where $\mathfrak m_\pm$ is the irreducible $64$-dimensional $G_2^\C$-module.
The minimal $G_2$-connection gives a reference connection on each
$\mathcal F_\pm$. The difference from any other connection is an
adjoint-bundle-valued one-form. The displayed equivariant splitting
therefore gives a unique expression
\[
  D=D^{\min}+\rho(\beta)+\Psi,
\]
where $\beta$ is now allowed to be $\mathfrak g_2^\C$-valued:
\begin{align*}
  \beta&\in\Omega^1\bigl(\mathcal M,
       \mathrm{Fr}_\varphi\times_{G_2}\mathfrak g_2^\C\bigr),\\
  \Psi&\in\Omega^1\bigl(\mathcal M,
       \mathrm{Fr}_\varphi\times_{G_2}\mathfrak m_\pm\bigr).
\end{align*}
Conversely, both summands are $\mathfrak g_\pm$-valued, so every such
pair defines a connection preserving $N_\pm$.

\begin{proposition}[Connections on the exceptional reductions]
  \label{prop:exceptional-connections}
  For a connection $D$ on $\mathcal F_\pm$ written as above,
      the form $\Psi$ is the second fundamental form
      of the intermediate reduction
      $\mathrm{Fr}_{\varphi,\C}\subset\mathcal F_\pm$, with group $G_2^\C$.
      It vanishes if and only if $D$ preserves the
      intermediate reduction, equivalently, if $D$
      preserves both critical cubics.

      In particular, a connection preserving one critical
      cubic and the real structure of
      $\mathcal W_{\varphi,\C}$ preserves both critical
      cubics and has $\Psi=0$. Such a connection is
      induced by a real $G_2$-connection on $T\mathcal M$.
\end{proposition}

Here $\Psi$ represents the class of the connection form along
$\mathrm{Fr}_{\varphi,\C}$ in
$\mathfrak g_\pm/\rho(\mathfrak g_2^\C)$, identified with
$\mathfrak m_\pm$ by the reductive splitting.
This class is the covariant derivative of a principal reduction
\cite[authors' PDF, \S~2.10, Lemma~2.10.1 and
Definition~2.10.4]{piccione-2006-theory}, also called its
second fundamental form
\cite[arXiv~v4, pp.~34--35]{lubke-2006-universal}.

\begin{proof}
  The quotient
  $\mathfrak g_\pm/\rho(\mathfrak g_2^\C)$ is
  identified with $\mathfrak m_\pm$.
  The projection of the connection form to this
  quotient is the second fundamental form of
  $\mathrm{Fr}_{\varphi,\C}$ in $\mathcal F_\pm$.
  Since $D^{\min}$ and $\rho(\beta)$ preserve
  the intermediate reduction, this projection
  is precisely $\Psi$. Thus $D$ preserves
  $\mathrm{Fr}_{\varphi,\C}$ if and only if
  $\Psi=0$.

  By Proposition~\ref{prop:common-stabilizer}\,
  (\ref{item:common-stabilizer-lie}),
  \[
    \mathfrak g_+\cap\mathfrak g_-
      =\rho(\mathfrak g_2^\C).
  \]
  Hence a connection on $\mathcal F_\pm$
  preserves both critical cubics if and only if
  its connection form is
  $\rho(\mathfrak g_2^\C)$-valued, which is
  equivalent to $\Psi=0$.

 Finally, suppose that $D$ also preserves the canonical
real structure on $\mathcal W_{\varphi,\C}$, given by
complex conjugation with respect to $\mathcal W_\varphi$. Since
  $N_-=\overline{N_+}$, preservation of one
  critical cubic implies preservation of the
  other. We have just shown that this forces
  $\Psi=0$. In a real $G_2$-adapted frame,
  the remaining connection form is real and
  lies in $\rho(\mathfrak g_2^\C)$, hence in
  $\rho(\mathfrak g_2)$. It therefore comes
  from a real $G_2$-connection on $T\mathcal M$.
\end{proof}

An exceptional reduction does not by itself determine a compatible
connection. Using the metric $g_\varphi$, however, we can select a
canonical compatible connection by orthogonally projecting the
Levi-Civita connection. We first identify the relevant orthogonal
complement.

By Remark~\ref{rem:koecher-tkk}, we have the decomposition
\[
  \End(\WC)\cong\mathbf1\oplus\mathbf{78}\oplus\mathbf{650}
\]
into irreducible $\esix$-modules, where the summand $\mathbf{78}$
is the adjoint representation on $\mathfrak g_\pm$.
The invariant form $\tr_W(XY)$ is nondegenerate and preserved by
conjugation. Since these three irreducible modules are pairwise
nonisomorphic and have distinct dimensions, they are mutually
orthogonal. We therefore have the $\tr_W$-orthogonal decomposition
\[
  \End(\WC)=\C\Id_W\oplus\mathfrak g_\pm\oplus\mathfrak c_\pm,
  \qquad \text{so that} \qquad
  \mathfrak g_\pm^\perp=\C\Id_W\oplus\mathfrak c_\pm,
  \qquad \text{where} \qquad
  \mathfrak c_\pm\cong\mathbf{650}.
\]
Let $\pi_\pm:\End(\WC)\to\mathfrak g_\pm$ denote the
corresponding orthogonal projection. The \emph{projected connection} of a
connection on $\mathcal W_{\varphi,\C}$ is the connection on $\mathcal F_\pm$
whose connection form, in a local frame of $\mathcal F_\pm$, is the
$\pi_\pm$-component of the original one (see Remark~\ref{rem:connection-defect}).
Write $\pi_{\mathfrak g_2}$ for the projection associated with
$\lieso(7)=\mathfrak g_2\oplus\{\Theta_v:v\in V\}$.

\begin{lemma}[Orthogonality of the $\mathbf7$-directions]
  \label{lem:exceptional-torsion-orthogonality}
  For either sign and every $v\in V$,
  \[
    \rho(\Theta_v)\in\mathfrak c_\pm.
  \]
  In particular, on $\lieso(7)$,
  \[
    \pi_\pm\circ\rho=\rho\circ\pi_{\mathfrak g_2}.
\]
\end{lemma}

\begin{proof}
  Since $\rho(G_2)\subset G_\pm$ and the invariant form
  $\tr_W(XY)$ is preserved by conjugation, $\pi_\pm$ is
  $G_2$-equivariant. Hence the complexified composition
  \[
    V_\C\xrightarrow{\ v\mapsto\Theta_v\ }\lieso(7,\C)
    \xrightarrow{\ \rho\ }\End(\WC)
    \xrightarrow{\ \pi_\pm\ }\mathfrak g_\pm
  \]
  is $G_2$-equivariant. The decomposition
  \[
    \mathfrak g_\pm
      =\rho(\mathfrak g_2^\C)\oplus\mathfrak m_\pm
      \cong\mathbf{14}\oplus\mathbf{64}
  \]
  contains no copy of $\mathbf7$, so this composition
  vanishes by Schur's lemma. Moreover, $\rho(\Theta_v)$
  is skew-adjoint on the real space $W$, hence traceless.
  Thus it belongs to $\mathfrak c_\pm$.
  The projection identity follows because $\pi_\pm$
  is the identity on $\rho(\mathfrak g_2)$.
\end{proof}

We apply Lemma~\ref{lem:exceptional-torsion-orthogonality} to the connections induced on
$\mathcal W_{\varphi,\C}$. As before, we use the same
notation for a metric connection on $T\mathcal M$ and its
induced connection on the deformation bundle.

\begin{proposition}[The minimal connection as projected Levi-Civita connection]
  \label{prop:projected-levi-civita}
  \leavevmode
  \begin{enumerate}
    \renewcommand{\labelenumi}{(\arabic{enumi})}
    \renewcommand{\theenumi}{\arabic{enumi}}

    \item\label{item:projected-connection-metric}
      Let $D$ be the connection on $\mathcal W_{\varphi,\C}$
      induced by a metric connection on $T\mathcal M$, written as in
      \eqref{eq:metric-connection-decomposition}. Then
      \[
        D=\bigl(D^{\min}+\rho(\beta)\bigr)
          +\rho\bigl(\Theta_{w(\cdot)+T(\cdot)/3}\bigr).
      \]
      For either sign, the first term preserves
      $\mathcal F_\pm$ and the second is
      $\mathfrak c_\pm$-valued. Thus both exceptional
      reductions give the same projected connection,
      $D^{\min}+\rho(\beta)$.

    \item\label{item:projected-connection-minimal}
      In particular, both signs give $D^{\min}$ as the
      projected Levi-Civita connection. For either sign,
      it is the unique connection $\widetilde D$ on
      $\mathcal F_\pm$ for which
      $D^{LC}-\widetilde D$ is
      $\mathfrak g_\pm^\perp$-valued.
      It is also the unique pointwise closest connection
      to $D^{LC}$ among all connections on $\mathcal F_\pm$,
      in the norm induced by $g_\varphi$ and the Hermitian
      form on $\mathcal W_{\varphi,\C}$.
      Moreover, $D^{\min}$ preserves the compact reduction
      $\mathcal F^c_\pm$.
  \end{enumerate}
\end{proposition}

\begin{proof}
  Substituting
  \[
    D^{LC}=D^{\min}
      +\tfrac13\rho(\Theta_{T(\cdot)})
  \]
  gives the decomposition in
  part~\ref{item:projected-connection-metric}.
  Its first term is induced by a $G_2$-connection,
  so it preserves both exceptional reductions and
  their compact reductions. By
  Lemma~\ref{lem:exceptional-torsion-orthogonality}, the
  second term lies in $\mathfrak c_\pm$.

  In particular, $D^{LC}-D^{\min}$ is
  $\mathfrak g_\pm^\perp$-valued.
  Every connection on $\mathcal F_\pm$ has the form
  $D^{\min}+\omega$, with $\omega$ a $\mathfrak g_\pm$-valued
  one-form, so the $\tr_W$-orthogonal splitting gives
  the asserted uniqueness.

  By Proposition~\ref{prop:compact-form},
  $\mathfrak g_\pm$ is the complexification of a
  compact real form consisting of skew-Hermitian
  operators. It is therefore stable under Hermitian
  adjunction, and, writing $\langle\cdot,\cdot\rangle_{\mathsf H}$
  for the Hermitian extension of the pairing
  $\langle\cdot,\cdot\rangle_{\End}$ of
  Section~\ref{sec:exceptional-pencil},
  \[
    \langle X,Y\rangle_{\mathsf H}
      =\tr_W(XY^{\mathsf H})=0
    \qquad
    (X\in\mathfrak g_\pm^\perp,\;
     Y\in\mathfrak g_\pm).
  \]
  Thus the $\tr_W$-orthogonal splitting is also
  Hermitian-orthogonal. For every such $\omega$,
  \[
    |D^{LC}-(D^{\min}+\omega)|^2
      =|D^{LC}-D^{\min}|^2+|\omega|^2
  \]
  pointwise, which proves that $D^{\min}$ is the unique
  pointwise closest connection to $D^{LC}$ among all
  connections on $\mathcal F_\pm$.
\end{proof}

Together with Lemma~\ref{lem:phi-cubic-variation},
Corollary~\ref{cor:critical-connection-criteria} and
Lemma~\ref{lem:exceptional-torsion-orthogonality},
Proposition~\ref{prop:projected-levi-civita} proves
Theorem~\ref{thm:intro-compatible-connections}.

For the compact reduction $\mathcal F^c_\pm$,
Proposition~\ref{prop:projected-levi-civita} reduces to the
orthogonal-projection characterization of the minimal
connection given by Cleyton and Swann
\cite[Definition~3.1 and Proposition~3.2]{cleyton-2004-einstein}.
The proposition shows that the same connection
is also the unique closest connection among
all connections on the complex exceptional
reduction $\mathcal F_\pm$.

\begin{remark}[Projection and the defect of a connection]
  \label{rem:connection-defect}
  The projection construction applies to any
  connection $D^0$ on $\mathcal W_{\varphi,\C}$, although its two
  projections need not coincide. Fix a sign and choose a local
  frame in $\mathcal F_\pm$. Decompose the connection form of
  $D^0$ according to the $\tr_W$-orthogonal splitting
  \[
    \End(\WC)=\mathfrak g_\pm\oplus\mathfrak g_\pm^\perp.
  \]
  The $\mathfrak g_\pm$-component defines a connection $D^\pm$
  preserving $\mathcal F_\pm$. Since the projection is
  $G_\pm$-equivariant, this construction is independent of the
  chosen frame. For this fixed sign, we may therefore write
  \[
    D^0=D^\pm+\varpi,
    \qquad
    \varpi=\varpi_1\Id_W+\varpi_{\mathfrak c},
  \]
  using the decomposition
  \[
    \mathfrak g_\pm^\perp
      =\C\Id_W\oplus\mathfrak c_\pm.
  \]
  Here $\varpi_1$ is a complex-valued one-form and
  $\varpi_{\mathfrak c}$ is a $\mathfrak c_\pm$-valued one-form.
  We call $\varpi$ the \emph{defect} of $D^0$.

  Since $D^\pm$ preserves $N_\pm$, the covariant derivative of
  the critical cubic is
  \[
    D^0N_\pm
      =-3\varpi_1\otimes N_\pm
       -\Phi_{N_\pm}(\varpi_{\mathfrak c}).
  \]
  Recall that $\ker\Phi_{N_\pm}=\mathfrak g_\pm$ and
  $\Phi_{N_\pm}(\Id_W)=3N_\pm$. Hence $\Phi_{N_\pm}$ is
  injective on $\mathfrak g_\pm^\perp$, and
  \[
    \Phi_{N_\pm}(\mathfrak c_\pm)\cap\C N_\pm=0.
  \]
  It follows that $D^0$ preserves the line $\C N_\pm$ if and
  only if $\varpi_{\mathfrak c}=0$, and preserves $N_\pm$
  itself if and only if $\varpi=0$.
\end{remark}

We return to the Albert algebra bundle of a line field and ask when
its product is parallel under a compatible connection.

\begin{proposition}[Parallel pointings]
  \label{prop:parallel-pointings}
  Let $D$ preserve $g$ and $\varphi$, and use the same symbol for
  its induced connections.
  \begin{enumerate}
    \renewcommand{\labelenumi}{(\arabic{enumi})}
    \renewcommand{\theenumi}{\arabic{enumi}}
    \item\label{item:parallel-pointings-product}
      The product $\circ_{E_\ell}$ is $D$-parallel if and
      only if $D\ell\subset\ell$. When this holds, $D$ also
      preserves $N_\pm/n_\pm$, the based trace form, the adjoint,
      and $\varsigma_{\ell,\pm}$.

    \item\label{item:parallel-pointings-holonomy}
      If $\varphi$ is torsion-free, then
      $D^{LC}E_\ell=0$ if and only if $\ell$ is parallel.
      On each connected component this forces
      $\operatorname{Hol}(g)\subseteq H_{E_1}$ up to conjugacy,
      with $H_{E_1}$ as in Proposition~\ref{ex:structured-basepoints};
      conversely, such holonomy gives a parallel $E_\ell$.

      If $\ell$ is oriented, then $D\ell\subset\ell$ if
      and only if its globally defined unit section $u$
      is $D$-parallel. In this case,
      \[
        \operatorname{Hol}(D)\subseteq
        \operatorname{Stab}_{G_2}(u)\cong SU(3).
      \]
      In particular, for torsion-free $\varphi$ and
      $D=D^{LC}$, an oriented parallel line field forces
      $\operatorname{Hol}(g)\subseteq SU(3)$ up to conjugacy.
      Conversely, such holonomy gives an oriented parallel
      line field.

  \end{enumerate}
\end{proposition}

\begin{proof}
  If $D$ preserves $\ell$, it preserves $P_\ell$ and every tensor
  contraction defining the operations, which proves one
  implication of~(\ref{item:parallel-pointings-product}) and the last sentence
  of~(\ref{item:parallel-pointings-product}).
  Conversely, suppose that $D$ preserves the product.
  Differentiating $E_\ell\circ_{E_\ell}h=h$ gives
  \[
    (D_XE_\ell)\circ_{E_\ell}h
       +E_\ell\circ_{E_\ell}D_Xh=D_Xh.
  \]
  The second term equals $D_Xh$, so
  $(D_XE_\ell)\circ_{E_\ell}h=0$.
  Taking $h=E_\ell$ and using the unit identity gives
  $D_XE_\ell=0$. Since $Dg=0$ and
  $P_\ell=(E_\ell+g)/7$, it follows that
  $DP_\ell=0$, equivalently $D\ell\subset\ell$.

  In the torsion-free case $D^{LC}$ preserves $\varphi$.
  By the holonomy principle, parallel lines correspond to
  lines preserved by holonomy, and the stabilizer in $G_2$
  of such a line is $H_{E_1}$. This proves the first assertion of~(\ref{item:parallel-pointings-holonomy}),
  including its converse.

  For the oriented assertion in~(\ref{item:parallel-pointings-holonomy}), orientation gives a global unit section $u$
  of $\ell$. If $D$ preserves $\ell$, then
  $D_Xu=a(X)u$ for some one-form $a$.
  Metric compatibility implies
  \[
    0=X(g(u,u))=2g(D_Xu,u)=2a(X),
  \]
  so $Du=0$. The converse is immediate.
  Since the stabilizer of a unit vector in $G_2$ is
  $SU(3)$, the holonomy inclusion follows.
  Conversely, if $\operatorname{Hol}(g)$ lies in an
  $SU(3)$ subgroup of $G_2$, it fixes a unit vector,
  which extends by parallel transport to a global
  parallel unit vector field.

\end{proof}

The holonomy principle and the absence of fixed elements in the
nontrivial irreducible $\mathbf{27}$ immediately give the following.
\begin{corollary}[Full restricted holonomy]
  \label{cor:full-holonomy-pointings}
  If $\operatorname{Hol}^0(g)=G_2$ up to conjugacy, then
  $\mathcal W_{\varphi,\C}$ has no nonzero $D^{LC}$-parallel section.
\end{corollary}

The case $\operatorname{Hol}(g)=H_{E_1}\cong SU(3)\rtimes\mathbb Z/2$ of
Proposition~\ref{prop:parallel-pointings}\,(\ref{item:parallel-pointings-holonomy}),
with disconnected holonomy, occurs in the compact barely $G_2$-manifolds
constructed as quotients of $\mathcal Y\times S^1$ by
$(y,t)\mapsto(\sigma(y),-t)$, where $\mathcal Y$ is a Calabi--Yau
threefold and $\sigma$ is a free antiholomorphic isometric involution.
Their holonomy is stated to be $SU(3)\rtimes\mathbb Z/2$
in~\cite{GrigorianBarelyG2, Joycebook}.
The parallel-pointing criterion of
Proposition~\ref{prop:parallel-pointings} can be seen directly:
the parallel line spanned by $\partial_t$ on the product
descends to the quotient, since the involution sends
$\partial_t$ to $-\partial_t$. This gives an unoriented parallel
line, explaining the disconnected stabilizer rather than $SU(3)$.

These results give a geometric realization of the exceptional
algebraic structures suggested by Lee and Leung
\cite{Lee:2002fa}. The \(E_6\)-invariant cubics themselves are distinguished complex linear combinations of the two geometrically relevant real cubics $C_0$ and $C_\varphi$. On any $7$-manifold with a $G_2$-structure, these critical cubics define canonical
conjugate $E_6(\C)$-reductions of the rank-$27$ deformation
bundle. A tangent line field supplies the additional choice of a unit,
yielding Albert algebra bundles and $F_4(\C)$-reductions.
The reduced-holonomy setting considered by Lee and Leung
appears when this choice is parallel.

The differential geometry of these exceptional structures is
governed by the intrinsic torsion of $\varphi$.
The Levi-Civita derivative of either critical cubic determines
the full intrinsic torsion, while orthogonal projection onto
either exceptional reduction gives the connection induced by
the minimal $G_2$-connection.
Thus the exceptional cubic structures are intrinsic to the
$G_2$-geometry, whereas parallel Albert products require
additional holonomy reduction.
The constructions provide a starting point for investigating
the curvature and differential invariants of exceptional
connections, and their relationship with global $G_2$
deformation theory.

\appendix
\section{Proofs of the contraction and quartic identities}
\label{app:contraction-proofs}

This appendix contains the contraction derivations deferred from
Section~\ref{sec:exceptional-pencil}.
In Section~\ref{subsec:cubic-tensor-contractions}, we derive the one-index
c-contraction identities \eqref{eq:one-index-c-contraction-mixed} and
\eqref{eq:one-index-c-contraction-phiphi} and the two-index constants
\eqref{eq:two-index-c-contraction}.
Section~\ref{subsec:stabilizer-contractions} uses these identities to
evaluate the Gram matrix and normal-operator scalars for stabilizer exclusion.
In Section~\ref{subsec:quartic-identities}, we prove
Lemma~\ref{lem:universal-quartic-identities}.
Section~\ref{subsec:notation} collects the notation used across the paper.
These derivations use the conventions and $G_2$-contraction identities of
Section~\ref{subsec:standard-form-contractions},
Lemma~\ref{lem:product-contractions}, and
\eqref{eq:projection-identity}.
We retain the up/down contraction convention of Section~\ref{sec:introduction}.
Metric-dual basis indices are written $k^\alpha$, $\Theta^i$ and $L^i$;
in an adapted orthonormal frame their numerical entries agree with those
of $k_\alpha$, $\Theta_i$ and $L_i$, respectively.

\subsection{The c-contraction identities}
\label{subsec:cubic-tensor-contractions}

We use the cubic components $c^s_{\alpha\beta\gamma}$ and the ordered
one-index c-contractions $\mathsf T^{st}$ of
Section~\ref{subsec:cubic-contraction-calculus}. Greek indices refer to
$W$ and Latin indices to $V$. We first evaluate the one-index
c-contractions, then sum them to obtain the two-index constants.
To evaluate the one-index contractions without choosing coordinates
on $W$, we write the pairings $\mathsf T^{st}(h_1,h_2;h_3,h_4)$ of
Section~\ref{subsec:cubic-contraction-calculus} as
\begin{equation*}
\mathsf T^{st}(h_1,h_2;h_3,h_4)
=\sum_\sigma C_s(h_1,h_2,k_\sigma)C_t(h_3,h_4,k^\sigma).
\end{equation*}
This sum pairs the corresponding traceless products $A$ and $B$,
so $\mathsf T^{st}_{\alpha\beta\gamma\delta}
=\mathsf T^{st}(k_\alpha,k_\beta;k_\gamma,k_\delta)$.
The projection identity \eqref{eq:projection-identity} gives, for
arbitrary symmetric matrices $h_1,h_2$,
\begin{equation*}
\sum_\sigma\tr(h_1k_\sigma)\tr(h_2k^\sigma)
=\tr(h_1h_2)-\tfrac17\tr(h_1)\tr(h_2).
\end{equation*}
For traceless inputs, the trace corrections are determined by
$\tr\tfrac12(h_1h_2+h_2h_1)=\langle h_1,h_2\rangle_W$ and
$\tr\mathcal P(h_1,h_2)=-\langle h_1,h_2\rangle_W$. Substituting these
into the projection identity evaluates the first three one-index
c-contractions and yields \eqref{eq:one-index-c-contraction-mixed}.
Exchanging the two cubic types while swapping the two index pairs gives the symmetry
\[\mathsf T^{st}(h_1,h_2;h_3,h_4)=\mathsf T^{ts}(h_3,h_4;h_1,h_2).\]
For the remaining one-index c-contraction, two applications of
\eqref{eq:phi-phi-contractions} give
\begin{equation}
\begin{aligned}
\mathsf T^{\varphi\varphi}(h_1,h_2;h_3,h_4)
 &=(\delta_{a_1a_3}\delta_{a_2a_4}-\delta_{a_1a_4}\delta_{a_2a_3}
     +\psi_{a_1a_2a_3a_4})\\
 &\qquad\qquad(\delta_{b_1b_3}\delta_{b_2b_4}-\delta_{b_1b_4}\delta_{b_2b_3}
     +\psi_{b_1b_2b_3b_4})\\
 &\qquad\qquad\cdot h_1^{a_1b_1}h_2^{a_2b_2}h_3^{a_3b_3}h_4^{a_4b_4}
   -\tfrac17\langle h_1,h_2\rangle_W\langle h_3,h_4\rangle_W.
\end{aligned}
\label{eq:one-index-c-contraction-phiphi-raw}
\end{equation}
To evaluate the $4$-form product in this expression, we first derive
a local triple-form contraction and take its trace. The same
calculation will also be used in the proof of
Lemma~\ref{lem:universal-quartic-identities}, and it does not use any of
the later stabilizer arguments.

For $h\in W$, put $Y=\mathcal P(h,h)$, $r=\tr(h^2)$ and
$\gamma_\varphi=C_\varphi(h)$, and write
$\{\mathsf X,\mathsf Y\}=\mathsf X\mathsf Y+\mathsf Y\mathsf X$
for the anticommutator.

\begin{lemma}[Triple-form reduction]
\label{lem:triple-form-reduction}
For every $h\in W$, we have
\begin{equation}
\mathcal P(Y,Y)
=4\gamma_\varphi h+4hYh+4rY-8\mathcal P(h,h^3)-8\{h^2,Y\},
\label{eq:pyy-reduction}
\end{equation}
\begin{equation}
\tr(Y^2)=5r^2-8\tr(h^4)+12C_\varphi(h,h,h^2).
\label{eq:tr-y-squared}
\end{equation}
\end{lemma}

\begin{proof}
The reduction below uses the following one index $\varphi\psi$
$G_2$-contraction identity:
\begin{equation}
\begin{aligned}
  \varphi_{ijm}\psi^{m}{}_{abc}
    ={}&\delta_{ia}\varphi_{jbc}+\delta_{ib}\varphi_{jca}
          +\delta_{ic}\varphi_{jab}\\
       &-\delta_{ja}\varphi_{ibc}-\delta_{jb}\varphi_{ica}
          -\delta_{jc}\varphi_{iab}.
\end{aligned}
\label{eq:six-term-contraction}
\end{equation}
In Bryant's notation, after relabeling, this is \cite[(2.9)]{bryant-2003}.

We reduce
$\mathcal T_{i;abcd}=\varphi_i{}^{pq}\varphi_{pab}\varphi_{qcd}$ in two
steps: contract its first two factors using
\eqref{eq:phi-phi-contractions}, then contract the resulting $4$-form
using \eqref{eq:six-term-contraction}. This gives
\begin{equation}
\begin{aligned}
\mathcal T_{i;abcd}
 &=-\delta_{ia}\varphi_{bcd}+\delta_{ib}\varphi_{acd}
   +\delta_{ic}\varphi_{abd}-\delta_{id}\varphi_{abc}\\
 &\quad-\delta_{ac}\varphi_{ibd}+\delta_{bc}\varphi_{iad}
   +\delta_{ad}\varphi_{ibc}-\delta_{bd}\varphi_{iac}.
\end{aligned}
\label{eq:triple-contraction}
\end{equation}
In particular, the first contraction is
$\varphi_i{}^p{}_q\varphi_{pab}
=-\delta_{ia}\delta_{qb}+\delta_{ib}\delta_{qa}-\psi_{iqab}$.
Denote the four signed summands in the first line of
\eqref{eq:triple-contraction} by $\mathcal T^{(1)}_1,\ldots,
\mathcal T^{(1)}_4$ and those in its second line by
$\mathcal T^{(2)}_1,\ldots,\mathcal T^{(2)}_4$.
Pair two such summands by
\begin{equation*}
[\mathcal T^{(p)}_{u,i},\mathcal T^{(q)}_{v,j}]_h
=\mathcal T^{(p)}_{u,i;abcd}\mathcal T^{(q)}_{v,j;a'b'c'd'}
  h^{aa'}h^{bb'}h^{cc'}h^{dd'}.
\end{equation*}
Here $i,j$ are the free indices and $p,q=1,2$; the subscripts $1,\ldots,4$ below
enumerate summands. To display the pairing arrays, abbreviate
\begin{equation*}
\mathsf C=\gamma_\varphi h,\quad \mathsf S=hYh,\quad \mathsf R=h^2Y,\quad
\mathsf D=rY,\quad \mathsf W=\mathcal P(h,h^3).
\end{equation*}
We pair the summands within each group and between the two groups.
Substituting the Kronecker deltas gives the following three arrays,
with rows and columns in the order of the summands above:
\begin{equation}
\begin{gathered}
\bigl([\mathcal T^{(1)}_u,\mathcal T^{(1)}_v]_h\bigr)_{u,v}
=\begin{pmatrix}
\mathsf C&-\mathsf S&\mathsf S&\mathsf S\\
-\mathsf S&\mathsf C&\mathsf S&\mathsf S\\
\mathsf S&\mathsf S&\mathsf C&-\mathsf S\\
\mathsf S&\mathsf S&-\mathsf S&\mathsf C
\end{pmatrix},\\[2pt]
\bigl([\mathcal T^{(1)}_u,\mathcal T^{(2)}_v]_h\bigr)_{u,v}
=-\begin{pmatrix}
\mathsf R&0&\mathsf R&0\\
0&\mathsf R&0&\mathsf R\\
\mathsf R&\mathsf R&0&0\\
0&0&\mathsf R&\mathsf R
\end{pmatrix},\qquad
\bigl([\mathcal T^{(2)}_u,\mathcal T^{(2)}_v]_h\bigr)_{u,v}
=\begin{pmatrix}
\mathsf D&-\mathsf W&-\mathsf W&0\\
-\mathsf W&\mathsf D&0&-\mathsf W\\
-\mathsf W&0&\mathsf D&-\mathsf W\\
0&-\mathsf W&-\mathsf W&\mathsf D
\end{pmatrix}.
\end{gathered}
\label{eq:triple-contraction-pairings}
\end{equation}
The reverse mixed array, pairing group~$(2)$ with group~$(1)$,
is obtained from the displayed $(1)$--$(2)$ array by
transposing both its $4\times4$ matrix and the free
indices $i,j$ in each entry. Representative entries account for every tensor type and every sign:
\begin{align*}
[\mathcal T^{(1)}_{1,i},\mathcal T^{(1)}_{1,j}]_h
 &=h_{ij}\varphi_{bcd}\varphi_{b'c'd'}
   h^{bb'}h^{cc'}h^{dd'}=\mathsf C_{ij},\\
[\mathcal T^{(1)}_{1,i},\mathcal T^{(1)}_{2,j}]_h
 &=-h_i{}^{a'}h_j{}^b\varphi_{bcd}\varphi_{a'c'd'}
   h^{cc'}h^{dd'}=-\mathsf S_{ij},\\
[\mathcal T^{(1)}_{1,i},\mathcal T^{(2)}_{1,j}]_h
 &=(h^2)_i{}^c\varphi_{bcd}\varphi_{jb'd'}
   h^{bb'}h^{dd'}=-\mathsf R_{ij},\\
[\mathcal T^{(1)}_{1,i},\mathcal T^{(2)}_{2,j}]_h
 &=-h_i{}^{a'}(h^2)^{bc}\varphi_{bcd}\varphi_{ja'd'}h^{dd'}=0,\\
[\mathcal T^{(2)}_{1,i},\mathcal T^{(2)}_{1,j}]_h
 &=r\varphi_{ibd}\varphi_{jb'd'}h^{bb'}h^{dd'}=\mathsf D_{ij},\\
[\mathcal T^{(2)}_{1,i},\mathcal T^{(2)}_{2,j}]_h
 &=-\varphi_{ibd}\varphi_{ja'd'}(h^3)^{ba'}h^{dd'}=-\mathsf W_{ij}.
\end{align*}
In the fourth entry, the symmetric tensor $(h^2)_{bc}$
is contracted against alternating slots, so the entry vanishes.
The other entries follow by permuting the four paired dummy slots.
To recover the full contraction, we sum the entries of each array:
the first, mixed and last give $4\mathsf C+4\mathsf S$, $-8\mathsf R$ and
$4\mathsf D-8\mathsf W$, respectively. The reverse mixed array contributes
$-8\mathsf R^{\transp}$. Since
$[\mathcal T_i,\mathcal T_j]_h=\mathcal P(Y,Y)_{ij}$, the sum is
\eqref{eq:pyy-reduction}.
Taking its trace, using \eqref{eq:p-trace-rule} and $\tr Y=-r$, gives
\eqref{eq:tr-y-squared}.

The derivation of \eqref{eq:pyy-reduction} and
\eqref{eq:tr-y-squared} does not use the quartic identities
of Lemma~\ref{lem:universal-quartic-identities}.
\end{proof}

\begin{proposition}
\label{prop:one-index-c-contractions}
The one-index c-contractions satisfy
\eqref{eq:one-index-c-contraction-mixed} and
\eqref{eq:one-index-c-contraction-phiphi}.
\end{proposition}

\begin{proof}
The first three contractions were evaluated above. For the remaining
contraction, we use Lemma~\ref{lem:triple-form-reduction}.
We use the alternating vector $\eta(h_1,h_2)$ of
Section~\ref{subsec:cubic-contraction-calculus}, extended to
arbitrary symmetric inputs by the same formula.
Extend $C_\varphi$ to arbitrary symmetric inputs by the same ordered
contraction as in its definition, and set
\begin{equation*}
\begin{aligned}
\mathcal C(h_1,h_2,h_3,h_4)={}&C_\varphi\bigl(h_1,h_2,\tfrac12(h_3h_4+h_4h_3)\bigr)
 +C_\varphi\bigl(h_1,h_3,\tfrac12(h_2h_4+h_4h_2)\bigr)\\
 &+C_\varphi\bigl(h_1,h_4,\tfrac12(h_2h_3+h_3h_2)\bigr)
 +C_\varphi\bigl(h_2,h_3,\tfrac12(h_1h_4+h_4h_1)\bigr)\\
 &+C_\varphi\bigl(h_2,h_4,\tfrac12(h_1h_3+h_3h_1)\bigr)
 +C_\varphi\bigl(h_3,h_4,\tfrac12(h_1h_2+h_2h_1)\bigr).
\end{aligned}
\end{equation*}
It differs from the sum $\Xi$ of
Section~\ref{subsec:cubic-contraction-calculus} by trace terms.
We derive \eqref{eq:one-index-c-contraction-phiphi} by expanding
\eqref{eq:one-index-c-contraction-phiphi-raw} into the delta, mixed and
$4$-form contributions. Write each parenthesis in
\eqref{eq:one-index-c-contraction-phiphi-raw} as $D+\psi$, with
$D_{abef}=\delta_{ae}\delta_{bf}-\delta_{af}\delta_{be}$, and abbreviate
\begin{equation*}
\begin{gathered}
p_0=\langle h_1,h_2\rangle_W\langle h_3,h_4\rangle_W,\quad
p_1=\langle h_1,h_3\rangle_W\langle h_2,h_4\rangle_W,\quad
p_2=\langle h_1,h_4\rangle_W\langle h_2,h_3\rangle_W,\\
t_1=\tr(h_1h_2h_3h_4),\quad t_2=\tr(h_1h_2h_4h_3),\quad
t_3=\tr(h_1h_3h_2h_4),\\
e=\eta(h_1,h_3)\cdot\eta(h_2,h_4)+\eta(h_1,h_4)\cdot\eta(h_2,h_3).
\end{gathered}
\end{equation*}
The two delta factors contribute $p_1+p_2-2t_3$.
For the mixed factors, we rearrange \eqref{eq:phi-phi-contractions} as
\begin{equation*}
\psi_{abef}=\delta_{ab}\delta_{ef}-\delta_{af}\delta_{eb}
             -\varphi^i{}_{ae}\varphi_{ibf}.
\end{equation*}
Write $\mathsf X\lrcorner\varphi$ for the vector with components
$(\mathsf X\lrcorner\varphi)_k=\mathsf X^{ij}\varphi_{ijk}$.
Thus for arbitrary matrices $\mathsf X,\mathsf Y$,
\begin{equation*}
\psi_{abef}\mathsf X^{ae}\mathsf Y^{bf}
=\tr(\mathsf X^{\transp}\mathsf Y)-\tr(\mathsf X\mathsf Y)
 -(\mathsf X\lrcorner\varphi)\cdot(\mathsf Y\lrcorner\varphi).
\end{equation*}
We apply this identity to $(\mathsf X,\mathsf Y)=(h_1h_3,h_2h_4)$ and
$(h_1h_4,h_2h_3)$,
then collect the traces using cyclicity and reversal; by the definition
of $\eta$, $(h_1h_3)\lrcorner\varphi=\eta(h_1,h_3)$, and similarly for the
other three products. The sum of the
two mixed contributions is $2t_1+2t_2-4t_3-2e$.

For the product of the $4$-forms, contracting the two $\varphi$ factors
sharing the middle index of $Y^2$ gives
\begin{equation*}
Y^2-2\mathcal P(h^2,h^2)
=\bigl(\varphi_{ia_1a_2}\varphi_{ja_3a_4}\psi_{b_1b_2b_3b_4}
        h^{a_1b_1}h^{a_2b_2}h^{a_3b_3}h^{a_4b_4}\bigr)_{ij}.
\end{equation*}
Taking the trace gives the left-hand side of the scalar identity below,
because the two metric terms in the $i,j$ contraction vanish against
$\psi$. We can therefore
use the trace calculation \eqref{eq:tr-y-squared} to obtain the scalar identity
\begin{equation*}
\psi_{a_1a_2a_3a_4}\psi_{b_1b_2b_3b_4}
h^{a_1b_1}h^{a_2b_2}h^{a_3b_3}h^{a_4b_4}
=3r^2-6\tr(h^4)+12C_\varphi(h,h,h^2).
\end{equation*}
We recover the four-input expression by normalized polarization.
The contraction of two copies of $\psi$ with four matrix inputs is
fully symmetric in those inputs, because interchanging two inputs
interchanges two slots in each alternating form.
Polarizing the scalar identity in the ordered inputs
$h_1,h_2,h_3,h_4$ therefore gives the explicit contraction
\begin{equation*}
\psi_{a_1a_2a_3a_4}\psi_{b_1b_2b_3b_4}
h_1^{a_1b_1}h_2^{a_2b_2}h_3^{a_3b_3}h_4^{a_4b_4}
  =p_0+p_1+p_2-2(t_1+t_2+t_3)+2\mathcal C.
\end{equation*}
The trace term splits into three classes under cyclicity and reversal.
For the last term, there are six choices of the two linear inputs,
with the Jordan product taken in the other pair. We now add the delta,
mixed and $4$-form contributions, and subtract the trace correction
$p_0/7$:
\begin{equation}
\begin{aligned}
\mathsf T^{\varphi\varphi}
 &=(p_1+p_2-2t_3)+(2t_1+2t_2-4t_3-2e)\\
 &\quad+p_0+p_1+p_2-2(t_1+t_2+t_3)+2\mathcal C-\tfrac17p_0.
\end{aligned}
\label{eq:one-index-c-contraction-phiphi-j}
\end{equation}
This is the intermediate expression with symmetrized products. If every symmetrized
product in $\mathcal C$ is
replaced by its traceless part $A(\cdot,\cdot)$, the new sum $\Xi$ satisfies
$\mathcal C=\Xi-2(p_0+p_1+p_2)/7$, because
$C_\varphi(h_1,h_2,g)=-\langle h_1,h_2\rangle_W$. The coefficients of
$p_0,p_1,p_2$ in \eqref{eq:one-index-c-contraction-phiphi-j} then become
$2/7,10/7,10/7$, respectively,
which gives \eqref{eq:one-index-c-contraction-phiphi}.
Both versions keep the alternating-vector terms and the order of the
inputs.
\end{proof}

\paragraph{The two-index c-contraction constants.}
To derive the entries of $\mathsf a$ in \eqref{eq:two-index-c-contraction},
we use the projection sums, a mixed cubic contraction, and the
alternating-vector contraction. We evaluate these before
summing the one-index c-contraction identities. For $h\in W$,
the projection identity gives
\begin{equation}
\sum_\alpha k_\alpha k^\alpha=\tfrac{27}{7}g,\qquad
\sum_\alpha k_\alpha h k^\alpha=\tfrac5{14}h,\qquad
\sum_\alpha\mathcal P(k_\alpha,k^\alpha)=-\tfrac{27}{7}g.
\label{eq:projection-sums}
\end{equation}
For the last identity, inserting the projection identity in
$\varphi_{iab}\varphi_{jcd}\sum_\alpha
(k_\alpha)^{ac}(k^\alpha)^{bd}$ gives
$0-3\delta_{ij}-6\delta_{ij}/7$.
The mixed cubic contraction needed below is
\begin{equation*}
u(h_1,h_2):=\sum_\alpha
C_\varphi\bigl(h_1,k_\alpha,\tfrac12(h_2k^\alpha+k^\alpha h_2)\bigr)
=\tfrac9{14}\langle h_1,h_2\rangle_W.
\end{equation*}
The two terms of the symmetrized product agree by transposition.
After inserting the
projection identity, their sum reduces to
$\varphi_{a_1a_2a_3}\varphi_{b_1b_2b_3}h_1^{a_1b_1}h_2^{a_3p}
\Pi^{a_2b_2}{}_{p}{}^{b_3}$,
where $\Pi$ denotes the right-hand side of \eqref{eq:projection-identity}.
The three summands contribute $0$, $\langle h_1,h_2\rangle_W/2$ and
$\langle h_1,h_2\rangle_W/7$,
respectively. To evaluate the latter two, we apply
\eqref{eq:phi-phi-contractions} to the alternating tensors.
The $4$-form part vanishes on contraction with a symmetric input,
and the pair-trace term vanishes by tracelessness.

For the vector terms of \eqref{eq:one-index-c-contraction-phiphi-j}, put
$\Theta_i:=\Theta_{e_i}$, with $(\Theta_i)_{jk}=\varphi_{ijk}$,
as in Section~\ref{subsec:koecher-jump}. The $G_2$-contraction identities imply
$\sum_i\Theta_i\Theta^i=-6g$ and
$\sum_i\Theta_i h\Theta^i=h$ on $W$.
Also $\eta(h_1,h_2)_i=\langle h_1,[\Theta_i,h_2]/2\rangle_W$.
The projection identity and cyclicity of trace therefore give
\begin{equation*}
\begin{aligned}
\sum_\alpha\eta(h_1,k_\alpha)\cdot\eta(h_2,k^\alpha)
 &=\tfrac14\sum_i\tr([\Theta_i,h_1][\Theta^i,h_2])\\
 &=\tfrac12\sum_i\tr(\Theta_i h_1\Theta^i h_2
                        -\Theta_i\Theta^i h_1h_2)
 =\tfrac72\langle h_1,h_2\rangle_W.
\end{aligned}
\end{equation*}
We can now sum the one-index contractions using
$\mathsf a_{st}\langle h_1,h_2\rangle_W=\sum_\alpha
\mathsf T^{st}(h_1,k_\alpha;h_2,k^\alpha)$. For the first two entries, substituting
the sums above into \eqref{eq:one-index-c-contraction-mixed} gives
\begin{equation*}
\mathsf a_{00}=\tfrac12(\tfrac5{14}+\tfrac{27}{7})-\tfrac17
=\tfrac{55}{28},\qquad
\mathsf a_{0\varphi}=\tfrac9{14}+\tfrac17=\tfrac{11}{14}.
\end{equation*}
For the remaining entry $\mathsf a_{\varphi\varphi}$, we substitute
$(h_1,k_\alpha;h_2,k^\alpha)$ in \eqref{eq:one-index-c-contraction-phiphi-j}
and sum each type of term separately.
The three metric terms contribute $(6/7+54+2)\langle h_1,h_2\rangle_W$,
and the crossed trace contributes $-8(27/7)\langle h_1,h_2\rangle_W$.
The six terms of $\mathcal C$ sum to
$4u(h_1,h_2)-54\langle h_1,h_2\rangle_W/7
=-36\langle h_1,h_2\rangle_W/7$:
four are copies of $u(h_1,h_2)$, while the other two use
$\sum_\alpha k_\alpha k^\alpha=27g/7$ and
$\sum\mathcal P(k_\alpha,k^\alpha)=-27g/7$.
Of the two vector terms, one vanishes since $\eta(k_\alpha,k^\alpha)=0$.
The other sums to $-7\langle h_1,h_2\rangle_W/2$ before multiplication
by $-2$.
Combining these contributions gives
\begin{equation}
\begin{gathered}
\mathsf a_{\varphi\varphi}
=\left(\tfrac67+54+2\right)-\tfrac{216}{7}-\tfrac{72}{7}+7
=\tfrac{159}{7}.
\end{gathered}
\label{eq:two-index-c-contraction-constants}
\end{equation}
These entries give the matrix and determinant in
\eqref{eq:two-index-c-contraction}.
All sums use the complete ordered cubic contraction; introducing a
factor $1/3!$ would change these constants. The projection identity
represents the identity operator on $W$, so this derivation and the resulting
contraction are independent of the choice of orthonormal basis. These
contractions extend complex-bilinearly.

\subsection{Contractions for stabilizer exclusion}
\label{subsec:stabilizer-contractions}

We prove the matrix evaluations of Lemmas~\ref{lem:four-map-gram}
and~\ref{lem:four-block-scalars}. We use the action maps $\Phi_s$,
their algebraic transposes and pairings from
\eqref{eq:action-map-ortho}--\eqref{eq:action-map-adjoint-ortho},
the multiplication operators $m_h^{(s)}$ of \eqref{eq:mn-operators},
and the operators $\mathcal Z_{st}$ of \eqref{eq:back-contraction}.
The coefficients $\mathsf a_{st}$ and one-index c-contractions
$\mathsf T^{st}$ have been evaluated above.
We first describe the four multiplicity-one blocks as tensors and
evaluate the seven-generator contraction used in their normal operators.

Represent $X\in \End \WC$ by its four-tensor $X_{ij,kl}$, symmetric and traceless in
each pair, with $(Xh)_{ij}=X_{ij,kl}h^{kl}$. Here the comma is only used to separate groups of indices, and does not denote differentiation. Then
$(X^{\tadj})_{ij,kl}=X_{kl,ij}$. Define the partial trace $r(X)$ and the
$\varphi$-contraction $q(X)$ of the four-tensor of $X$,
\[
r(X)_{ik}=X_{ij,k}{}^j,\qquad q(X)_{ikm}=X_{ij,kl}\varphi^{jl}{}_{m}.
\] From \eqref{eq:phi-phi-contractions}, we obtain
\begin{align}
q(X)_{ikm}\varphi^k{}_{n}{}^m&=X_{ij,kl}\varphi^{jl}{}_{m}\varphi^k{}_{n}{}^m \notag\\
&=X_{ij,kl}(\delta^{jk}\delta^l{}_{n}-\delta^j{}_{n}\delta^{lk}+\psi^{jlk}{}_{n}) \notag\\
&=r(X)_{in}, \label{eq:q-r-contraction}
\end{align}
since each pair of indices of $X_{ij,kl}$ is traceless and symmetric. Therefore, $q(X)=0$ forces $r(X)=0$. Moreover, if $X^{\tadj}=\epsilon X$ with
$\epsilon=\pm1$, then $q(X)_{kim}=-\epsilon\,q(X)_{ikm}$. In particular, for $W$-self-adjoint $X$, $q(X)$ is skew in the indices $i,k$, and for $W$-skew-adjoint $X$, $q(X)$ is instead symmetric on these indices.
Let $\mathcal S(X)$ denote the orthogonal projection of the full symmetrization of $X_{ij,kl}$ in its four $V$-indices onto $S^4_0V^*$, the space of fully symmetric traceless four-tensors, regarded again as an endomorphism through its two pairs of indices. For $W$-skew-adjoint $X$, the full symmetrization vanishes, so $\mathcal S(X)=0$.
For $W$-self-adjoint $X$ the full symmetrization is the average of the three pairings below, and its trace over one pair is $\tfrac23r(X)$. Hence on a $W$-self-adjoint $X$ with $r(X)=0$, no trace correction is needed, and
\begin{equation}
\label{eq:EndWC-symmetrization}
\mathcal S(X)_{ij,kl}
=\frac13\bigl(X_{ij,kl}+X_{ik,jl}+X_{il,jk}\bigr).
\end{equation}
Here $\mathcal S(X)$ is a fully symmetric traceless four-tensor, $\mathcal S^2(X)=\mathcal S(X)$, and $q(\mathcal S(X))=0$.

We will also use the equivariant map
$\jmath_3\colon S^3_0V^*\to\Lambda^2\WC$ defined by
\begin{equation*}
  \jmath_3(\mathsf t)_{ij,kl}=\mathsf t_{ikm}\varphi_{jl}{}^m+\mathsf t_{ilm}\varphi_{jk}{}^m
  +\mathsf t_{jkm}\varphi_{il}{}^m+\mathsf t_{jlm}\varphi_{ik}{}^m.
\end{equation*}
The transpose $q^\top$ is taken with respect to the full ordered pairing
$\langle u,v\rangle_3=u_{ikm}v^{ikm}$ on three-tensors and the pairing
$\langle X,Y\rangle_{\End}=\tr_W(X^{\tadj}Y)$ on $\End(\WC)$
in~\eqref{eq:endomorphism-pairing}.

\begin{lemma}
\label{lem:four-block-tensors}
The four higher-dimensional multiplicity-one blocks in \(\End \WC\) (from Table \ref{tab:endomorphism-inventory}) have the following properties:
\[
\begin{array}{c|c|cccc}
\text{Component}&\dim&X^{\tadj}&q(X)&r(X)&\mathcal S(X)\\ \hline
V(4,0)=S^4_0V^*&182&X&0&0&X\\
V(0,2)&77&X&0&0&0\\
V(3,0)=\jmath_3(S^3_0V^*)&77&-X&8\mathsf t&0&0\\
V(2,1)&189&-X&0&0&0
\end{array}
\]
In the third row, $X=\jmath_3(\mathsf t)$ for a unique
$\mathsf t\in S^3_0V^*$; moreover, $q^\top q(X)=2X$ on this block.
\end{lemma}

\begin{proof}
From \eqref{eq:q-r-contraction}, $q(X)=0$ implies $r(X)=0$, hence to prove \(r(X)=0\), it is sufficient to show \(q(X)=0\). The identification $V(4,0)=S^4_0V^*$ gives the first row: full symmetry implies $q(X)=0$, and $\mathcal S(X)=X$.

For $V(0,2)$, the vanishing of $q$ and $\mathcal S$ follows by equivariance.
The map $q$ sends $W$-self-adjoint $X$ into
$$\Lambda^2V^*\otimes V^*\cong\mathbf1\oplus2\cdot\mathbf7\oplus\mathbf{14}\oplus2\cdot\mathbf{27}\oplus\mathbf{64},$$
where the decomposition follows from
$\Lambda^2V=\Lambda^2_7\oplus\Lambda^2_{14}$, $\Lambda^2_7\cong V$,
$V\otimes V\cong\mathbf1\oplus\mathbf7\oplus\mathbf{14}\oplus\mathbf{27}$, and
\eqref{eq:seven-adjoint-product}. On the other hand, $\mathcal S$ maps into
$S^4V^*=S^4_0V^*\oplus(g\cdot S^2_0V^*)\oplus\C\,g\cdot g$, where $\cdot$ denotes the product in the symmetric algebra $S^\bullet V^*$, and neither target contains $V(0,2)$. Therefore the images must be $0$ by Schur's lemma.

To compute $q(\jmath_3(\mathsf t))$, we apply \eqref{eq:phi-phi-contractions}. The four terms then contribute $6\mathsf t,\mathsf t,\mathsf t,0$ to
$q(\jmath_3(\mathsf t))=8\mathsf t$. Hence, $\jmath_3$ is injective, and since $S^3_0V^*=V(3,0)$ occurs once in
$\Lambda^2\WC$, its image is precisely that component. Note that since $\mathsf t$ is symmetric, we have $r(\jmath_3(\mathsf t))=0$.

The $189$-dimensional component $V(2,1)$ is $W$-skew-adjoint and satisfies $q(X)=0$, and hence $r(X)=0$.
Indeed, $q$ maps it into the invariant direct sum
\[
S^2V^*\otimes V^*=S^3V^*\oplus\ker(\Sym),
\]
whose summands have dimensions $84$ and $112$. Each projection of $q$
to these summands must vanish: a nonzero equivariant map from the
irreducible $189$-dimensional module would be injective, which is
impossible for either target.

For the two skew-adjoint components, $\mathcal S(X)=0$ by antisymmetry
under exchange of the two index pairs.
Full ordered pairings also give
$q^\top(\mathsf t)=\jmath_3(\mathsf t)/4$, hence
$q^\top q(X)=2X$ on $V(3,0)$.
\end{proof}

To evaluate $\mathcal Z_{\varphi\varphi}$ on these blocks, we must
contract the two alternating-vector terms in
\eqref{eq:one-index-c-contraction-phiphi}. One of them gives
$\sum_iL_iXL^i$. Here $\Theta_i$ is the generator of
$\Lambda^2_7$ dual to $e_i$, and $L_i=\tfrac12\rho(\Theta_i)$ is one half of
its action on $W$, given by
\[
(L_i)_{\alpha\beta}=\eta(k_\alpha,k_\beta)_i,
\qquad L_i(h)=\tfrac12[\Theta_i,h].
\]
The next lemma expresses this sum through $q$ and $\mathcal S$, whose
values on the four blocks are known from Lemma~\ref{lem:four-block-tensors}.
This lets us compute the required contraction using only the
two-$\varphi$ contraction identity \eqref{eq:phi-phi-contractions}.

\begin{lemma}
\label{lem:seven-generator}
For a $W$-self-adjoint or $W$-skew-adjoint endomorphism
$X\in\End(\WC)$ with $r(X)=0$, the
following identity holds:
\begin{equation}
\sum_iL_iXL^i=
\begin{cases}
q^\top q(X)+3\mathcal S(X)-2X,&X^{\tadj}=X,\ r(X)=0,\\
q^\top q(X)-X,&X^{\tadj}=-X,\ r(X)=0.
\end{cases}
\label{eq:seven-generator-sandwich}
\end{equation}
\end{lemma}

\begin{proof}
For $r(X)=0$, put
\[
\mathcal D(X)_{ab,cd}=\tfrac14\bigl(
 X_{cb,ad}+X_{db,ca}+X_{ac,bd}+X_{ad,cb}\bigr).
\]
Expanding $L_i(h)=\tfrac12[\Theta_i,h]$ on both sides of $X$ gives
\begin{equation*}
\begin{aligned}
\left(\sum_i L_iXL^i\right)_{ab,cd}
=-\tfrac14\bigl(
 &\varphi_{ia}{}^j\varphi^i{}_c{}^lX_{bj,dl}
  +\varphi_{ia}{}^j\varphi^i{}_d{}^lX_{bj,cl}\\
 &+\varphi_{ib}{}^j\varphi^i{}_c{}^lX_{aj,dl}
  +\varphi_{ib}{}^j\varphi^i{}_d{}^lX_{aj,cl}\bigr).
\end{aligned}
\end{equation*}
For example, the first product contracts by
\eqref{eq:phi-phi-contractions} as
\[
\varphi_{iaj}\varphi^i{}_{cl}
 =\delta_{ac}\delta_{jl}-\delta_{al}\delta_{jc}+\psi_{ajcl}.
\]
The first delta product contributes $-\tfrac14\delta_{ac}r(X)_{bd}=0$,
and the second contributes $\tfrac14X_{bc,da}$.
Applying the same identity to the other three products gives
\begin{equation*}
\begin{aligned}
\left(\sum_i L_iXL^i\right)_{ab,cd}
=\mathcal D(X)_{ab,cd}-\tfrac14\bigl(
 &\psi_a{}^j{}_c{}^lX_{bj,dl}+\psi_a{}^j{}_d{}^lX_{bj,cl}\\
 &+\psi_b{}^j{}_c{}^lX_{aj,dl}+\psi_b{}^j{}_d{}^lX_{aj,cl}\bigr).
\end{aligned}
\end{equation*}
Thus the surviving delta terms sum to $\mathcal D(X)$.
On the other hand, pair symmetrization in the definition of $q^\top$ gives
\begin{equation*}
(q^\top q(X))_{ab,cd}
 =\tfrac14\bigl(
 \varphi_{bd}{}^m q(X)_{acm}+\varphi_{bc}{}^m q(X)_{adm}
 +\varphi_{ad}{}^m q(X)_{bcm}+\varphi_{ac}{}^m q(X)_{bdm}\bigr).
\end{equation*}
Its pair traces vanish because
$q(X)_{ikm}\varphi^k{}_n{}^m=r(X)_{in}=0$.
After substituting $q(X)_{ikm}=X_{ij,kl}\varphi^{jl}{}_{m}$, the
two-$\varphi$ products appear explicitly:
\begin{equation*}
\begin{aligned}
(q^\top q(X))_{ab,cd}
=\tfrac14\bigl(
 &\varphi_{bd}{}^m\varphi^{jl}{}_{m}X_{aj,cl}
  +\varphi_{bc}{}^m\varphi^{jl}{}_{m}X_{aj,dl}\\
 &+\varphi_{ad}{}^m\varphi^{jl}{}_{m}X_{bj,cl}
  +\varphi_{ac}{}^m\varphi^{jl}{}_{m}X_{bj,dl}\bigr).
\end{aligned}
\end{equation*}
For the first product, \eqref{eq:phi-phi-contractions} gives
\[
\varphi_{bd}{}^m\varphi^{jl}{}_{m}X_{aj,cl}
 =X_{ab,cd}-X_{ad,cb}+\psi_{bd}{}^{jl}X_{aj,cl}.
\]
The four delta contributions therefore sum to $X-\mathcal D(X)$.
The four $\psi$ contributions agree with those in $\sum_iL_iXL^i$:
for example, $\psi_{bd}{}^{jl}X_{aj,cl}=-\psi_b{}^j{}_d{}^lX_{aj,cl}$,
and the other three terms agree by the same interchange of indices.
Consequently,
\[
\sum_iL_iXL^i=q^\top q(X)-X+2\mathcal D(X).
\]
For self-adjoint $X$, $\mathcal D(X)=(3\mathcal S(X)-X)/2$;
for skew-adjoint $X$, $\mathcal D(X)=0$. This proves
\eqref{eq:seven-generator-sandwich}.
\end{proof}

The four maps $\mathcal B_p$ and their ordering are fixed before
Lemma~\ref{lem:four-map-gram}. We now evaluate their Gram matrix.

\begin{proof}[Proof of Lemma~\ref{lem:four-map-gram}]
\label{proof:four-map-gram}
For $s,t,u,v\in\{0,\varphi\}$, the definition of $\Phi_s^\top$,
\eqref{eq:back-contraction} and \eqref{eq:m-trace-pairing} give
\[
\begin{aligned}
\langle\mathcal B_{st}(h_1),\mathcal B_{uv}(h_2)\rangle_3
&=\langle m_{h_1}^{(t)},
   \Phi_s^\top\Phi_u(m_{h_2}^{(v)})\rangle_{\End}\\
&=3\mathsf a_{su}\mathsf a_{tv}\langle h_1,h_2\rangle_W
  +6\tr_W\bigl(
    m_{h_1}^{(t)}\mathcal Z_{su}(m_{h_2}^{(v)})
  \bigr).
\end{aligned}
\]
By equivariance and Schur's lemma,
$\mathcal B_p^\top\mathcal B_q=\mathcal G_{pq}\Id_W$.
Taking the trace over $W$ in the preceding pairing formula gives
\begin{equation}
\mathcal G_{(s,t),(u,v)}
 =3\mathsf a_{su}\mathsf a_{tv}
  +\frac29\mathsf T^{ts}_{\delta\alpha\beta\gamma}
                (\mathsf T^{uv})^{\alpha\gamma\delta\beta}.
\label{eq:four-map-gram-contraction}
\end{equation}
To see the index order, note that the second trace is the sum
\[
\sum_\delta\tr_W\bigl(
m_{k_\delta}^{(t)}\mathcal Z_{su}(m_{k^\delta}^{(v)})\bigr)
=c^t_{\delta\alpha\beta}(c^s)^{\beta\gamma\epsilon}
  (c^u)^\alpha{}_{\sigma\epsilon}(c^v)^\sigma{}_{\gamma}{}^\delta.
\]
Renaming the contracted slots pairs the first two cubics into
$\mathsf T^{ts}$ and the last two into $\mathsf T^{uv}$ as in
\eqref{eq:four-map-gram-contraction}. Division by $\dim W=27$
changes the coefficient $6$ into $2/9$.

For a representative entry, expand the two $\mathsf T^{00}$ factors
using \eqref{eq:one-index-c-contraction-mixed}. Successive contractions with the
projection identity and the sums \eqref{eq:projection-sums} give
\begin{equation*}
\begin{aligned}
\mathsf T^{00}_{\delta\alpha\beta\gamma}
 (\mathsf T^{00})^{\alpha\gamma\delta\beta}
 &=\frac14\left(3\frac{4995}{392}+\frac{11043}{49}\right)
    -\frac17\left(\frac{135}{14}+\frac{729}{7}\right)
    +\frac{27}{49}\\
 &=\frac{78705}{1568}.
\end{aligned}
\end{equation*}
The first bracket consists of the four trace-times-trace terms; three
agree by cyclicity and reversal of trace. The next bracket contains the
two trace corrections, and the last term comes from contracting the
two metric products. Thus
$\mathcal G_{11}=3(55/28)^2+(2/9)(78705/1568)=4455/196$.
The full table of contractions in \eqref{eq:four-map-gram-contraction}, in the
ordered pairs $(0,0),(\varphi,0),(0,\varphi),(\varphi,\varphi)$, is
\begin{equation*}
\bigl(\mathsf T^{ts}_{\delta\alpha\beta\gamma}
      (\mathsf T^{uv})^{\alpha\gamma\delta\beta}\bigr)_{(s,t),(u,v)}
 =\frac1{1568}
 \begin{pmatrix}
 78705&31482&31482&-359964\\
 31482&-359964&-359964&601128\\
 31482&-359964&-359964&601128\\
 -359964&601128&601128&5523984
 \end{pmatrix}.
\end{equation*}
For every other entry we use the same contraction order in
\eqref{eq:four-map-gram-contraction}. We first expand the two
$\mathsf T$ factors, using \eqref{eq:one-index-c-contraction-mixed}
for the $00$ and mixed factors and
\eqref{eq:one-index-c-contraction-phiphi} for a $\varphi\varphi$
factor. We then sum
each repeated basis index with the projection identity and contract
the resulting $V$-indices with \eqref{eq:phi-phi-contractions}.
In particular, the two alternating-vector terms of the
$\varphi\varphi$ identity retain the displayed index order throughout
this expansion. This gives the contraction table above, including
the off-diagonal entries. Substitution in
\eqref{eq:four-map-gram-contraction} gives
the matrix in \eqref{eq:four-map-gram-matrix}. Its determinant, $179397504/7$, is nonzero.
\end{proof}

On each multiplicity-one block, $\Gamma_{st}$ is the scalar of
$\Phi_s^\top\Phi_t$ as in Lemma~\ref{lem:four-block-scalars}.

\begin{proof}[Proof of Lemma~\ref{lem:four-block-scalars}]
\label{proof:four-block-scalars}
By Lemma~\ref{lem:four-block-tensors}, $r(X)=0$ on each of these
components, and the same lemma gives the values of $q$ and
$\mathcal S$ there.
First contract the three identities in
\eqref{eq:one-index-c-contraction-mixed} with $X$ in the index order
specified by \eqref{eq:back-contraction}. A single insertion of the projection identity gives
\begin{equation}
\mathcal Z_{00}(X)=
\begin{cases}
-\tfrac{11}{28}X+\tfrac34\mathcal S(X),&X^{\tadj}=X,\ r(X)=0,\\
\tfrac17X,&X^{\tadj}=-X,\ r(X)=0.
\end{cases}
\label{eq:trace-sandwich}
\end{equation}
More explicitly, before the normalized trace tensor $g/\sqrt{7}$ is
removed, the four-tensor of $\mathcal Z_{00}(X)$ is
\[
\tfrac18\bigl(
(X^{\tadj})_{lj,ki}+(X^{\tadj})_{kj,li}
 +(X^{\tadj})_{li,kj}+(X^{\tadj})_{ki,lj}\bigr),
\]
plus terms involving $r(X)$. Removing $g/\sqrt7$ subtracts
$X^{\tadj}/7$. The bracket equals
$2(3\mathcal S(X)-X)$ in the self-adjoint case and zero in the
skew-adjoint case, giving \eqref{eq:trace-sandwich}.

For $X^{\tadj}=\epsilon X$ with $r(X)=0$, the same contraction
of the mixed identity gives
\begin{equation*}
\begin{aligned}
\mathcal Z_{0\varphi}(X)_{ij,kl}
=-\frac{\epsilon}{4}\bigl(
 &\varphi^c{}_{ik}q(X)_{jcl}+\varphi^c{}_{il}q(X)_{jck}\\
 &+\varphi^c{}_{jk}q(X)_{icl}+\varphi^c{}_{jl}q(X)_{ick}\bigr)
 +\frac{\epsilon}{7}X_{ij,kl}.
\end{aligned}
\end{equation*}
The last term is the correction from the normalized trace tensor
$g/\sqrt7$, since
$\mathcal P(g,h)=-h$ on $W$. The reversed-order operator is obtained
by $W$-adjoint transposition:
$\mathcal Z_{\varphi0}(X)
=\mathcal Z_{0\varphi}(X^{\tadj})^{\tadj}$.
Consequently both mixed operators equal $\epsilon X/7$ when $q(X)=0$.
For $V(3,0)$, Lemma~\ref{lem:four-block-tensors} gives
$q(\jmath_3(\mathsf t))=8\mathsf t$, so the bracket is
$8\jmath_3(\mathsf t)$. Thus both mixed operators equal
$13X/7$.

Finally evaluate \eqref{eq:one-index-c-contraction-phiphi} at
$(h_1,h_2,h_3,h_4)=(k_\alpha,k^\gamma,k_\beta,k_\sigma)$ and contract
with $X^\sigma{}_\gamma$. The crossed-trace term $-8\tr(h_1h_3h_2h_4)$ vanishes because its
contraction
contains the transpose of $r(X)$.
Write $m^\top,n^\top$ for the algebraic transposes of
$h\mapsto m_h,n_h$. Then the six mixed cubic terms contract to
\begin{equation*}
\mathcal Z_{\varphi0}(X+X^{\tadj})
 +\mathcal Z_{0\varphi}(X+X^{\tadj})
 +n_{m^\top X}+m_{n^\top X}.
\end{equation*}
The last two terms vanish for $W$-skew-adjoint $X$. For self-adjoint $X$ with
$q(X)=0$ they vanish as well: $m^\top X=r(X)_0$, and $n^\top X$
is the traceless part of
$\varphi^i{}_{a}{}^b\varphi^j{}_{c}{}^dX_{ij,bd}$, a contraction of $q(X)$.
Of the two $\eta$-terms, the first contracts to zero by $r(X)=0$;
the second is $\sum_iL_iXL^i$. The metric contribution is
$12X/7$ for $W$-self-adjoint $X$ and $8X/7$ for $W$-skew-adjoint $X$.
For the former, the mixed contribution is
$2(4X/7)$ when $q(X)=0$; for the latter it is zero.
Lemma~\ref{lem:seven-generator} therefore gives
\begin{equation}
\begin{aligned}
\mathcal Z_{\varphi\varphi}(X)
 &=\tfrac{48}{7}X-6\mathcal S(X)
 &&(X^{\tadj}=X,\ q(X)=0),\\
\mathcal Z_{\varphi\varphi}(X)
 &=\tfrac{22}{7}X-2q^\top q(X)
 &&(X^{\tadj}=-X,\ r(X)=0).
\end{aligned}
\label{eq:crossed-sandwich-scalars}
\end{equation}
By Lemma~\ref{lem:four-block-tensors}, $\mathcal S$ acts as $1$ on
$V(4,0)$ and as $0$ on $V(0,2)$, and $q^\top q$ acts as $2$ on $V(3,0)$
and as $0$ on $V(2,1)$. Substituting these values gives:

{\renewcommand{\arraystretch}{1.35}
\begin{center}
\begin{tabular}{lccc}
  Component & $\mathcal Z_{00}$
    & $\mathcal Z_{0\varphi}=\mathcal Z_{\varphi0}$
    & $\mathcal Z_{\varphi\varphi}$\\ \hline
  $V(4,0)$ & $\tfrac{5}{14}$ & $\tfrac17$ & $\tfrac67$\\
  $V(0,2)$ & $-\tfrac{11}{28}$ & $\tfrac17$ & $\tfrac{48}{7}$\\
  $V(3,0)$ & $\tfrac17$ & $\tfrac{13}{7}$ & $-\tfrac67$\\
  $V(2,1)$ & $\tfrac17$ & $-\tfrac17$ & $\tfrac{22}{7}$
\end{tabular}
\end{center}
}

Finally, by \eqref{eq:back-contraction}, $\Gamma_{st}$ is
$3\mathsf a_{st}$ plus six times the scalar of $\mathcal Z_{st}$ in the
table above,
which yields the stated normal-operator scalars and determinants.
\end{proof}

\subsection{Quartic identities}
\label{subsec:quartic-identities}

In this section, we prove Lemma~\ref{lem:universal-quartic-identities} of
Section~\ref{sec:exceptional-pencil}. The proof uses one contraction of
$\varphi$ with the Hodge dual of a $3$-form, which we record first.

\begin{lemma}
\label{lem:hodge-contraction}
For $u,v\in V$,
\begin{equation}
  \ast\bigl((\iota_u\varphi)\wedge(\iota_v\varphi)\bigr)
    =u^\flat\wedge\iota_v\varphi+v^\flat\wedge\iota_u\varphi.
      \label{eq:insertions-hodge}
\end{equation}
For any $3$-form $\xi$,
\begin{equation}
  \varphi_{iab}\varphi_{jcd}(\ast\xi)^{abcd}
    =2\bigl(\xi_{iab}\varphi_{j}{}^{ab}+\xi_{jab}\varphi_{i}{}^{ab}\bigr).
  \label{eq:three-form-hodge-contraction}
\end{equation}
\end{lemma}
\begin{proof}
Since $\psi=\ast\varphi$, the mixed contraction
\eqref{eq:mixed-contraction} implies
\begin{equation*}
  \bigl[\ast((\iota_u\varphi)\wedge\varphi)\bigr]_{ij}
    =\tfrac12\psi_{ijab}(\iota_u\varphi)^{ab}
    =2(\iota_u\varphi)_{ij}.
\end{equation*}
Together with $\ast(\iota_u\varphi)=u^\flat\wedge\psi$, this gives
$(\iota_u\varphi)\wedge\varphi=2u^\flat\wedge\psi$. Write $w=u\times v$
for the cross product \eqref{eq:cross-product}. Contracting that
identity with $v$ gives
\begin{equation*}
  (\iota_u\varphi)\wedge(\iota_v\varphi)
    =2g(u,v)\psi-2u^\flat\wedge\iota_v\psi
        +\varphi\wedge w^\flat.
\end{equation*}
We interchange $u,v$ in this last equation and average the two
equations:
\begin{equation*}
  (\iota_u\varphi)\wedge(\iota_v\varphi)
     =2g(u,v)\psi-u^\flat\wedge\iota_v\psi
         -v^\flat\wedge\iota_u\psi.
\end{equation*}
The standard Hodge identity for these degrees in dimension seven is
\begin{equation*}
  \ast(u^\flat\wedge\iota_v\psi)
      =g(u,v)\varphi-v^\flat\wedge\iota_u\varphi.
\end{equation*}
Here $\ast$ takes $4$-forms to $3$-forms. Applying this identity to the
two insertion terms
of the averaged equation, with $\ast\psi=\varphi$, proves
\eqref{eq:insertions-hodge}. Finally,
\eqref{eq:insertions-hodge}, paired with an arbitrary $3$-form $\xi$,
gives \eqref{eq:three-form-hodge-contraction} in components.
Converting the normalized pairing of forms to a sum over all ordered
indices gives the factor $2$.
\end{proof}

For $h\in W$, recall $Y=\mathcal P(h,h)$ and $r=\tr(h^2)$, and put
$\gamma_0=C_0(h)$. Also recall that $\{\mathsf X,\mathsf Y\}=\mathsf X\mathsf Y+\mathsf Y\mathsf X$ denotes the anticommutator.
Introduce a formal variable $t$ and set
$L(t)^i{}_j=\delta^i{}_j+th^i{}_j$, using the metric identification
\eqref{eq:metric-identification}. Its inverse $M(t)=L(t)^{-1}$ is given
by the Neumann series
\begin{equation}
\begin{aligned}
M(t)^i{}_j
 &=\sum_{m\geq0}(-t)^m(h^m)^i{}_j\\
 &=\delta^i{}_j-th^i{}_j+t^2(h^2)^i{}_j-t^3(h^3)^i{}_j
     +t^4(h^4)^i{}_j+\mathcal O(t^5).
\end{aligned}
\label{eq:formal-inverse-series}
\end{equation}
Since $\tr(h)=0$, the identity
$\det L=\exp(\tr\log L)$ gives the determinant expansion
\begin{equation}
\det L(t)=1-\tfrac r2t^2+\tfrac{\gamma_0}3t^3
 +\left(\tfrac{r^2}{8}-\tfrac{\tr(h^4)}{4}\right)t^4+\mathcal O(t^5).
\label{eq:formal-determinant-series}
\end{equation}
For the inverse series, the logarithm series and the determinant--exponential
identity, see \cite[\S3.4, eq.~(187), and \S10.1.5, eqs.~(497), (502)]{petersen-2012-matrix}.
Here we use these identities in $\End(V)[[t]]$ and $\mathbb R[[t]]$;
the constant term of $L$ is $I$, so no invertibility assumption on $h$
is needed.

Both $L$ and $M$ are self-adjoint with respect to $g$. In tensor products
such as $\mathcal P(M,M)$ we write $M$ for the symmetric covariant
$2$-tensor $M_{ij}=g_{ia}M^a{}_j$, as in
\eqref{eq:metric-identification}. Pullbacks insert the endomorphism in
each covariant slot; for example,
$(L^*\psi)_{abcd}=L^p{}_aL^q{}_bL^r{}_cL^s{}_d\psi_{pqrs}$.
Define the symmetric polynomial tensor
\begin{equation}
\mathsf Z(t)_{ij}
 =\varphi_i{}^{ab}\varphi_j{}^{cd}(L^*\psi)_{abcd}.
\label{eq:formal-z-definition}
\end{equation}

\begin{lemma}[Formal resolvent]
\label{lem:formal-resolvent}
The pullback of $\psi$ satisfies
\begin{equation}
\ast(L^*\psi)=\det(L)M^*\varphi.
\label{eq:formal-hodge-resolvent}
\end{equation}
Consequently,
\begin{equation}
\mathsf Z(t)=2\det(L)\{M,\mathcal P(M,M)\}.
\label{eq:formal-z-resolvent}
\end{equation}
\end{lemma}

\begin{proof}
Let $\varepsilon$ be the oriented volume tensor. The cofactor identity is
\begin{equation*}
\varepsilon_{b_1\cdots b_7}L^{b_1}{}_{a_1}\cdots L^{b_7}{}_{a_7}
 =\det(L)\varepsilon_{a_1\cdots a_7}.
\end{equation*}
Contracting three slots with $M$ and raising indices with $g$, using
self-adjointness, gives
\begin{gather*}
\varepsilon_{ijk}{}^{abcd}L^p{}_aL^q{}_bL^r{}_cL^s{}_d
 =\det(L)M^a{}_iM^b{}_jM^c{}_k\varepsilon_{abc}{}^{pqrs}.
\end{gather*}
Contracting with $\psi_{pqrs}/4!$ and using $\ast\psi=\varphi$ proves
\eqref{eq:formal-hodge-resolvent}.

Take $\xi=\det(L)M^*\varphi$ in \eqref{eq:three-form-hodge-contraction}.
By \eqref{eq:formal-hodge-resolvent}, $\ast\xi=L^*\psi$, so the
left-hand side is $\mathsf Z(t)_{ij}$. On the right-hand side,
\begin{equation*}
\xi_{iab}\varphi_j{}^{ab}
 =\det(L)M^p{}_iM^q{}_aM^r{}_b\varphi_{pqr}\varphi_j{}^{ab}
 =\det(L)M_i{}^p\mathcal P(M,M)_{pj}.
\end{equation*}
The other term is its transpose. Their sum gives
\eqref{eq:formal-z-resolvent}, with the contracted products of
\eqref{eq:contracted-product}.
\end{proof}

We now compare fourth-order coefficients in
Lemma~\ref{lem:formal-resolvent} to compute $Y^2$.

\begin{lemma}[Formal-resolvent quartic identity]
\label{lem:formal-resolvent-quartic}
For every $h\in W$, we have
\begin{equation}
\begin{aligned}
Y^2={}&6\mathcal P(h^2,h^2)+8\mathcal P(h,h^3)
 +4\{h,\mathcal P(h,h^2)\}+2\{h^2,Y\}-2rY\\
 &-8h^4+4rh^2+\tfrac83\gamma_0h+(2\tr(h^4)-r^2)g.
\end{aligned}
\label{eq:formal-y-squared}
\end{equation}
\end{lemma}

\begin{proof}
We compute the fourth-order coefficient in \eqref{eq:formal-z-resolvent},
then compute it directly in \eqref{eq:formal-z-definition}. Write
$\mathsf P_1=\mathcal P(h,h^2)$ and
$\mathsf P_2=\mathcal P(h^2,h^2)$, and recall that
$\mathsf W=\mathcal P(h,h^3)$ is as in
Section~\ref{subsec:cubic-tensor-contractions}.
The inverse series \eqref{eq:formal-inverse-series} and the elementary
identity $\mathcal P(g,h_1)=\tr(h_1)g-h_1$ \eqref{eq:p-identity-rule} give each
coefficient of the quadratic expression:
\begin{equation}
\begin{aligned}
\mathcal P(M,M)={}&6g+2th+t^2(2rg-2h^2+Y)\\
 &+t^3(-2\gamma_0g+2h^3-2\mathsf P_1)\\
 &+t^4(2\tr(h^4)g-2h^4+2\mathsf W+\mathsf P_2)+\mathcal O(t^5).
\end{aligned}
\label{eq:formal-p-resolvent-coefficients}
\end{equation}
Before multiplying by the determinant, we multiply by $M$ using
\eqref{eq:formal-inverse-series}. The
coefficients of orders $0,1,2$ are $6g$, $-4h$, $2rg+2h^2+Y$,
respectively, and the coefficient of degree four is
\begin{equation*}
2\tr(h^4)g+2\gamma_0h+2rh^2-2h^4+2\mathsf W+\mathsf P_2
 +2h\mathsf P_1+h^2Y.
\end{equation*}
The determinant series \eqref{eq:formal-determinant-series} contributes
to degree four through its terms of degrees
$0,2,3,4$. Thus the required coefficient is
\begin{equation}
\begin{aligned}
[t^4]\det(L)M\mathcal P(M,M)
={}&\left(\tfrac{\tr(h^4)}{2}-\tfrac{r^2}{4}\right)g
 +\tfrac{2\gamma_0}{3}h+rh^2-2h^4\\
 &+2\mathsf W+\mathsf P_2+2h\mathsf P_1+h^2Y-\tfrac r2Y.
\end{aligned}
\label{eq:formal-resolvent-fourth-coefficient}
\end{equation}
For the direct calculation, we return to the defining four-slot
contraction. The coefficient $[t^4]L^*\psi$ is $\psi$ with $h$ inserted
in each of its four slots. We identify this term in $Y^2$ by contracting the two
$\varphi$ factors that share the middle index using
\eqref{eq:phi-phi-contractions}. The two metric terms each equal
$\mathcal P(h^2,h^2)$, and the $4$-form term is exactly the
coefficient being computed. Thus
\begin{equation}
[t^4]\mathsf Z(t)=Y^2-2\mathsf P_2.
\label{eq:formal-z-direct-coefficient}
\end{equation}
Equating twice the sum of \eqref{eq:formal-resolvent-fourth-coefficient}
and its transpose to \eqref{eq:formal-z-direct-coefficient} gives \eqref{eq:formal-y-squared}.
\end{proof}

\begin{proof}[Proof of Lemma~\ref{lem:universal-quartic-identities}]
We work over $\mathbb R$. We combine the differentiated scalar
c-contraction with Lemma~\ref{lem:formal-resolvent-quartic} to prove the
first identity, then use Lemma~\ref{lem:triple-form-reduction} for the
second. Throughout the proof put
\begin{equation*}
a=A(h),\quad b=B(h),\quad r=\tr(h^2),\quad
\gamma_0=C_0(h),\quad \gamma_\varphi=C_\varphi(h),\quad Y=\mathcal P(h,h).
\end{equation*}
The operators $m_h=A(h,\cdot)$ and $n_h=B(h,\cdot)$
are those defined in \eqref{eq:mn-operators}; in this proof we write
$m,n$ for them. The six quartic products and eight iterated products
that occur below will be denoted by
\begin{equation}
\begin{gathered}
\begin{array}{c|cccccc}
j&1&2&3&4&5&6\\\hline
\mathcal A_j&A(a)&A(a,b)&A(b)&B(a)&B(a,b)&B(b)
\end{array},\\[2pt]
\begin{array}{c|cccccccc}
j&1&2&3&4&5&6&7&8\\\hline
\mathcal V_j&m^2a&m^2b&mna&mnb&nma&nmb&n^2a&n^2b
\end{array}.
\end{gathered}
\label{eq:quartic-notation-table}
\end{equation}
These are elements of $W$, with subscripts enumerating products rather
than tensor components. We keep the mixed iterated products distinct
throughout the calculation, so that their cancellations in the two
eliminations can be followed explicitly.

\emph{Differentiated scalar c-contraction.}
We use the scalar identity~\eqref{eq:scalar-cubic-contraction},
$8\langle a,a\rangle_W-12\langle a,b\rangle_W+\langle b,b\rangle_W=2r^2$.
The factors of two in the differentiation come from the quadratic
maps and their halved polarizations:
$D_ha[k]=2m_hk$ and $D_hb[k]=2n_hk$.
We differentiate the scalar identity, use self-adjointness of $m,n$,
and divide by four to obtain
\begin{equation}
8ma-6mb-6na+nb=2rh.
\label{eq:quartic-e1-first-derivative}
\end{equation}
Applying $m$ and $n$ to this equation gives two relations among the
iterated products:
\begin{equation}
\begin{aligned}
R_m&:=8\mathcal V_1-6\mathcal V_2-6\mathcal V_3
          +\mathcal V_4-2ra=0,\\
R_n&:=8\mathcal V_5-6\mathcal V_6-6\mathcal V_7
          +\mathcal V_8-2rb=0.
\end{aligned}
\label{eq:quartic-chain-relations}
\end{equation}
We also differentiate \eqref{eq:quartic-e1-first-derivative} in a
direction $k\in W$, held fixed during differentiation and set equal to
$a$ or $b$ only afterwards. On the left we obtain
\begin{equation*}
8A(k,a)-6A(k,b)-6B(k,a)+B(k,b)
+16m^2k-12(mn+nm)k+2n^2k.
\end{equation*}
On the right, the derivative of $2rh$ is
$4\langle h,k\rangle_Wh+2rk$. We now take $k=a,b$ and use
$\langle h,a\rangle_W=\gamma_0$ and $\langle h,b\rangle_W=\gamma_\varphi$
to obtain
\begin{equation}
\begin{aligned}
H_a:={}&8\mathcal A_1-6\mathcal A_2-6\mathcal A_4+\mathcal A_5
 +16\mathcal V_1-12(\mathcal V_3+\mathcal V_5)
 +2\mathcal V_7-4\gamma_0h-2ra=0,\\
H_b:={}&8\mathcal A_2-6\mathcal A_3-6\mathcal A_5+\mathcal A_6
 +16\mathcal V_2-12(\mathcal V_4+\mathcal V_6)
 +2\mathcal V_8-4\gamma_\varphi h-2rb=0.
\end{aligned}
\label{eq:quartic-fork-chain-relations}
\end{equation}
To reduce these relations, we express two of the iterated terms in
the product notation:
\begin{equation}
\mathcal V_1=\mathcal A_1+\tfrac17ra-\tfrac17\gamma_0h,
\qquad
\mathcal V_7=\mathcal A_5+\tfrac87ra-\tfrac87\gamma_0h.
\label{eq:quartic-crossed-e1}
\end{equation}
The first identity is an elementary Jordan-associator calculation
starting from $ma=h^3-rh/7-\gamma_0g/7$. Applying $m$ once more and removing
the trace gives the claimed expression. For the second, we subtract
\eqref{eq:one-index-c-contraction-phiphi-j} evaluated at $(h,a;h,k)$
from its value at $(h,h;a,k)$. Full symmetry of $C_\varphi$ and
self-adjointness of $n$ identify the two contractions as
\begin{equation*}
\begin{aligned}
\mathsf T^{\varphi\varphi}(h,h;a,k)
 &=\langle b,B(a,k)\rangle_W=\langle\mathcal A_5,k\rangle_W,\\
\mathsf T^{\varphi\varphi}(h,a;h,k)
 &=\langle na,nk\rangle_W=\langle\mathcal V_7,k\rangle_W.
\end{aligned}
\end{equation*}
The six cubic terms cancel
by their full symmetry. The crossed traces cancel because $ha=ah$, and the
alternating-vector terms vanish because $\eta(h,a)=0$.
Only the metric terms remain, so
$\langle\mathcal A_5-\mathcal V_7,k\rangle_W
=-8(r\langle a,k\rangle_W-\gamma_0\langle h,k\rangle_W)/7$.
The second identity in \eqref{eq:quartic-crossed-e1} follows because $k$
is arbitrary.
Substituting \eqref{eq:quartic-crossed-e1} into $H_a=0$ yields
\begin{equation}
8(\mathcal V_3+\mathcal V_5)
=16\mathcal A_1-4\mathcal A_2-4\mathcal A_4+2\mathcal A_5
 +\tfrac{12}{7}ra-\tfrac{40}{7}\gamma_0h.
\label{eq:quartic-mixed-chain-reduction}
\end{equation}

\emph{Trace removal and the first quartic identity.}
First we take traceless parts, expressing the matrix terms in the
notation of \eqref{eq:quartic-notation-table}:
\begin{equation}
\begin{aligned}
(Y^2)_0&=\mathcal A_3-\tfrac27rb,&
(\mathcal P(h^2,h^2))_0&=\mathcal A_4-\tfrac27ra,\\
(\mathcal P(h,h^3))_0&=\mathcal V_5-\tfrac17\gamma_0h+\tfrac17rb,&
\{h,\mathcal P(h,h^2)\}_0
 &=2\mathcal V_3-\tfrac27\gamma_0h-\tfrac27ra,\\
\{h^2,Y\}_0&=2\mathcal A_2+\tfrac27rb-\tfrac27ra,&
(h^4)_0&=\mathcal A_1+\tfrac27ra.
\end{aligned}
\label{eq:projection-table}
\end{equation}
For example,
$\mathcal P(h,h^2)=na-\gamma_0g/7-rh/7$ and
$h^3=ma+\gamma_0g/7+rh/7$; these follow from $a=h^2-rg/7$,
$\tr\mathcal P(h,a)=-\tr(ha)=-\gamma_0$ and
$\mathcal P(h,g)=-h$. The other entries follow by expanding
$Y=b-rg/7$ and applying the trace projection. In particular the
factors $1/7$ enter before any elimination.
Using the table to take the traceless part of
\eqref{eq:formal-y-squared}, we obtain
\begin{equation}
\mathcal A_3=6\mathcal A_4+8(\mathcal V_5+\mathcal V_3)
 +4\mathcal A_2-8\mathcal A_1-\tfrac{12}{7}ra+\tfrac8{21}\gamma_0h.
\label{eq:quartic-first-reduction}
\end{equation}
Collecting the trace corrections cancels the $rb$ terms and gives
the displayed coefficients of $ra$ and $\gamma_0h$.
We eliminate the mixed iterated products by substituting
\eqref{eq:quartic-mixed-chain-reduction}. The mixed product
$\mathcal A_2$ and the $ra$ terms cancel as follows:
\begin{equation*}
\begin{aligned}
\mathcal A_3
 &=6\mathcal A_4+16\mathcal A_1-4\mathcal A_2-4\mathcal A_4
   +2\mathcal A_5+\tfrac{12}{7}ra-\tfrac{40}{7}\gamma_0h\\
 &\qquad+4\mathcal A_2-8\mathcal A_1-\tfrac{12}{7}ra+\tfrac8{21}\gamma_0h\\
 &=8\mathcal A_1+2\mathcal A_4+2\mathcal A_5-\tfrac{16}{3}\gamma_0h.
\end{aligned}
\end{equation*}
Thus $8\mathcal A_1-\mathcal A_3+2\mathcal A_4+2\mathcal A_5
=16\gamma_0h/3$, the first identity of the lemma.

\emph{The second quartic identity.}
For the second identity, we return to \eqref{eq:pyy-reduction}, proved
above using the eight-term contraction \eqref{eq:triple-contraction}.
We first express that relation in traceless form, using
\begin{equation*}
\begin{aligned}
\mathcal P(Y,Y)_0&=\mathcal A_6+\tfrac27rb,\\
(hYh)_0&=2\mathcal V_2-\mathcal A_2-\tfrac17rb
                         +\tfrac27\gamma_\varphi h-\tfrac17ra,\\
\mathcal P(h,h^3)_0&=\mathcal V_5-\tfrac17\gamma_0h+\tfrac17rb,
\qquad
\{h^2,Y\}_0=2\mathcal A_2+\tfrac27rb-\tfrac27ra.
\end{aligned}
\end{equation*}
For the second equality, write $J(h_1,h_2)=(h_1h_2+h_2h_1)/2$.
The matrix identity and its traceless conversion are
\begin{align*}
hbh&=2J(h,J(h,b))-J(h^2,b),\\
(hbh)_0&=2m^2b-A(a,b)-\tfrac17rb+\tfrac27\gamma_\varphi h.
\end{align*}
Indeed, $J(h,b)=mb+\gamma_\varphi g/7$ and
$J(h^2,b)=A(a,b)+rb/7+\langle a,b\rangle_Wg/7$.
Substituting $Y=b-rg/7$ then gives the displayed expression for
$(hYh)_0$. Substitution in \eqref{eq:pyy-reduction} gives
\begin{equation}
R_3:=\mathcal A_6-8\mathcal V_2+20\mathcal A_2+8\mathcal V_5
 -\tfrac{36}{7}\gamma_\varphi h-\tfrac87\gamma_0h-\tfrac{12}{7}ra+\tfrac27rb=0.
\label{eq:quartic-triple-residual}
\end{equation}
We complete the elimination by writing the second quartic expression
as a combination of the relations already proved. Let
$Q_1=8\mathcal A_1-\mathcal A_3+2\mathcal A_4+2\mathcal A_5
-16\gamma_0h/3$, which vanishes by the first identity. Using $R_m,R_n,H_a,H_b$
from \eqref{eq:quartic-chain-relations} and
\eqref{eq:quartic-fork-chain-relations}, we collect terms to obtain
\begin{equation}
\begin{aligned}
16\mathcal A_2+2\mathcal A_3-8\mathcal A_4+\mathcal A_6
 -\tfrac{16}{3}\gamma_\varphi h
={}&-2R_m+\tfrac13R_n+H_a-\tfrac16H_b-Q_1+\tfrac76R_3.
\end{aligned}
\label{eq:quartic-second-residual}
\end{equation}
Collecting coefficients on the right cancels every iterated product
and every $ra,rb,\gamma_0h$ term, leaving the displayed quartic expression.
All residuals on its right-hand side vanish, proving the second identity.

All identities are polynomial in $h$, so they extend to $W_{\mathbb C}$
by complexification. Equivalently, their normalized polarizations give
the symmetric four-input tensor identities among the components
$\Sym_{\beta\gamma\delta\epsilon}
(c^s_{\alpha\sigma\tau}(c^t)^\sigma{}_{\beta\gamma}
 (c^u)^\tau{}_{\delta\epsilon})$ and
$\Sym_{\beta\gamma\delta\epsilon}
(c^s_{\beta\gamma\delta}\delta_{\alpha\epsilon})$.
The proof has used only the $G_2$-contraction identities, the evaluated
one-index c-contraction identity, the formal inverse of $I+th$, and the local triple-form
reduction; no critical parameter relation or stabilizer result enters.
\end{proof}

\subsection{Notation}
\label{subsec:notation}

The following table collects the notation used across sections. Its final
column points to the definition or convention. Pairings are complex
bilinear unless a conjugate or a Hermitian adjoint is written explicitly.

\begingroup
\small
\setlength{\tabcolsep}{4pt}
\renewcommand{\arraystretch}{1.15}
\begin{longtable}{@{}p{0.25\linewidth}p{0.46\linewidth}p{0.25\linewidth}@{}}
\caption{Notation and its defining locations.}\label{tab:notation}\\
\hline
Symbol & Meaning & Location \\
\hline
\endfirsthead
\multicolumn{3}{l}{\tablename~\thetable\ (continued)}\\
\hline
Symbol & Meaning & Location \\
\hline
\endhead
\hline
\endfoot
\multicolumn{3}{l}{\textbf{Linear algebra on $V$ and $W$}}\\*
$\mathcal M$, $V$ & The $7$-manifold and its model tangent space $\mathbb R^7$. & Section~\ref{subsec:standard-form-contractions} \\
$\varphi$, $g$, $\psi$, $\vol$ & The positive $3$-form, its metric, $\psi=\ast\varphi$, and its volume form. & \eqref{eq:standard-positive-form}--\eqref{eq:standard-dual-form} \\
$e_i$, $e^i$; $i,j,k$ & Adapted frame, dual coframe, and tensor indices from $1$ to $7$. & Section~\ref{subsec:standard-form-contractions} \\
$u\times v$ & The cross product determined by $\varphi$ and $g$. & \eqref{eq:cross-product} \\
$W$, $\WC$ & $S^2_0V^*$ and its complexification. & Definition~\ref{def:deformation-module} \\
$hk$, $h^m$, $h_0$ & Contracted tensor product, its powers, and traceless part. & \eqref{eq:contracted-product}, \eqref{eq:traceless} \\
$\langle h,k\rangle_W$ & $\tr(hk)$, the fixed $W$-pairing. & Definition~\ref{def:deformation-module} \\
$P_u$, $P_\Sigma$, $E_u$, $E_\Sigma$ & Rank-one tensor, plane projection tensor, and their traceless test tensors; $E_i=E_{e_i}$. & \eqref{eq:projection-tensors}, \eqref{eq:test-tensors} \\
$\mathcal I_\varphi$, $\chi_h$ & Insertion map and $\chi_h=\mathcal I_\varphi(h)$. & \eqref{eq:model-map} \\
$\rho(O)$, $\rho(M)$ & Group and infinitesimal actions on $W$. & \eqref{eq:rho-group}, \eqref{eq:rho-action} \\
$\theta_M$ & Insertion of an endomorphism into covariant slots. & \eqref{eq:insert-M} \\
$V(j_1,j_2)$, $\mathbf n$ & Irreducible $\mathfrak g_2^{\mathbb C}$-modules, indexed by highest weight or dimension. & Section~\ref{subsec:deformation-module} \\
$X^{\tadj}$ & Adjoint for the fixed bilinear $W$-pairing. & Section~\ref{subsec:deformation-module} \\
$\langle h,\bar k\rangle_W$, $X^{\mathsf H}$ & Hermitian form and its adjoint; conjugation fixes $W$. & Section~\ref{subsec:real-form} \\
\multicolumn{3}{l}{\textbf{Cubics, products and adjoints}}\\*
$C_0$, $C_\varphi$ & Trace and $\varphi$-cubics; the same symbols denote their normalized symmetric trilinear polarizations. & \eqref{eq:natural-cubics}, \eqref{eq:cubic-polarization} \\
$A(h,k)$, $\mathcal P(h,k)$, $B(h,k)$ & Symmetric traceless matrix product, full $\varphi$-product, and its traceless part. & \eqref{eq:trace-product}, \eqref{eq:phi-product} \\
$N_{a,b}$, $N_\lambda$ & Homogeneous pencil $aC_0+bC_\varphi$ and affine pencil $C_0+\lambda C_\varphi$. & \eqref{eq:invariant-pencil} \\
$\lambda_\pm$, $N_\pm$ & Critical parameters $(1\pm\ii\sqrt7)/8$ and their cubics. & \eqref{eq:critical-parameters} \\
$F_\lambda$, $F_\pm$ & Quadratic adjoints for the fixed $W$-pairing; $F_\lambda(h,k)$ is their halved polarization. & \eqref{eq:F-lambda-def}, Definition~\ref{def:quadratic-adjoint} \\
$k_\alpha$, $k^\alpha$ & Orthonormal basis of $W$ and its metric dual; Greek indices range from $1$ to $27$. & Section~\ref{subsec:cubic-contraction-calculus} \\
$c^s$, $\mathsf T^{st}$, $\mathsf a_{st}$ & Cubic components and their one-index and two-index contractions, $s,t\in\{0,\varphi\}$. & \eqref{eq:cubic-tensor-coordinates}, \eqref{eq:two-index-c-contraction} \\
$\eta(h_1,h_2)$, $\Xi$ & Alternating vector and sum of mixed cubic-product contractions. & Section~\ref{subsec:cubic-contraction-calculus} \\
\multicolumn{3}{l}{\textbf{Stabilizers and the $\mathbf{64}$}}\\*
$\stab(N)$, $\Phi_N$ & Exact infinitesimal stabilizer and cubic action map. & Definition~\ref{def:exact-stabilizer} \\
$\Phi_{a,b}$, $\Phi_s$ & Action maps of $N_{a,b}$ and $C_s$. & \eqref{eq:action-map}, \eqref{eq:action-map-ortho} \\
$\langle U,U'\rangle_3$, $\langle X,Y\rangle_{\End}$ & Full ordered tensor pairing and $\tr_W(X^{\tadj}Y)$, respectively. & \eqref{eq:endomorphism-pairing} \\
$\Phi_s^\top$ & Algebraic transpose between these two paired spaces. & \eqref{eq:action-map-adjoint-ortho} \\
$m_h$, $n_h$, $m_h^{(s)}$ & Multiplication operators for $A$, $B$, and the indexed choice of either product. & \eqref{eq:mn-operators} \\
$\mathcal Z_{st}$ & Operators occurring in $\Phi_s^\top\Phi_t$. & \eqref{eq:back-contraction} \\
$\mathcal B_{st}$, $\mathcal G$, $\Gamma_{st}$ & Four cubic maps, their Gram matrix, and normal-operator scalars on the multiplicity-one blocks. & Lemmas~\ref{lem:four-map-gram}, \ref{lem:four-block-scalars} \\
$S(h_1,h_2)$, $D(h_1,h_2)$, $S_\pm$ & Koecher operators, their associated operators $D$, and the two critical families. & Definition~\ref{def:koecher}, Corollary~\ref{cor:koecher} \\
$\Theta_u$, $\Theta_i$ & $\Theta_u(v)=v\times u$ and $\Theta_i=\Theta_{e_i}$. & Section~\ref{subsec:koecher-jump} \\
$W_6(u)$, $W_8(u)$, $W_{12}(u)$ & Real $SU(3)$ sectors of the pointwise space $W$. & Section~\ref{subsec:koecher-jump} \\
$\beta(u,v)$ & Symmetric bilinear map used in rank-one recognition. & Section~\ref{subsec:explicit-64} \\
$\mathcal Q_{64}$, $\varepsilon_Q$ & Tensor model of the $\mathbf{64}$ and its induced bilinear form. & Definition~\ref{def:q64}, \eqref{eq:JQ-epsQ} \\
$J_Q^{\mathrm s}$, $J_Q^{\mathrm a}$ & The two parity embeddings of $\mathcal Q_{64}$. & \eqref{eq:JQ-epsQ} \\
$\mathcal H_{64}$, $T$, $\zeta_\pm$ & The $64$-isotypic block, relative operator, and its critical eigenvalues. & Section~\ref{subsec:nothing-else-jumps} \\
$\mathfrak m_\pm$, $\mathcal E_\pm$, $\mathfrak g_\pm$ & Critical $64$-dimensional summands, their tensor parametrizations, and $\stab(N_\pm)$. & Section~\ref{subsec:nothing-else-jumps} \\
$G_\pm$, $\widehat G_\pm$, $\mu_3$ & Exact and line stabilizers of the critical cubics and cube roots of unity. & Section~\ref{subsec:critical-group} \\
\multicolumn{3}{l}{\textbf{Albert algebra at a basepoint}}\\*
$\mathbb A$, $N_{\mathbb A}$ & Complex Albert algebra and its cubic norm. & Section~\ref{sec:e6-albert-norm} \\
$e$, $n$, $N_e$ & Basepoint, its nonzero value $n=N_\pm(e)$, and normalized cubic $N_e=N_\pm/n$. & Table~\ref{tab:based-conventions} \\
$\langle x,y\rangle_e$, $\tr_e$, $T_e$ & Based trace form, linear trace, and operator relative to the $W$-pairing. & Table~\ref{tab:based-conventions} \\
$x^\sharp$, $x\mathbin\sharp y$, $\circ_e$ & Based adjoint, its unhalved polarization, and Jordan product. & Theorem~\ref{thm:based-algebras} \\
$x^{\circ m}$ & Jordan powers at the basepoint, with $x^{\circ0}=e$. & Section~\ref{sec:e6-albert-norm} \\
$v^\flat$, $\xi^\sharp$ & Metric dual of a vector and the inverse musical map on covectors; for an Albert algebra element, $x^\sharp$ has the based-adjoint meaning above. & Section~\ref{subsec:standard-form-contractions} \\
$\mathfrak J(N,e)$, $\mathfrak J_e$, $\mathcal U_x$ & Algebra of a cubic norm structure, normalized critical based algebra, and quadratic representation. & Section~\ref{sec:e6-albert-norm} \\
$F_{4,e}$, $\mathfrak f_e$, $\mathfrak g_{2,h}$ & Basepoint stabilizer, its Lie algebra, and the $\mathfrak g_2$ annihilator of $h$. & Section~\ref{subsec:g2-at-basepoint} \\
$K_\pm$, $\mathfrak k_\pm$, $\mathcal E^c_\pm$ & Compact critical stabilizers, their Lie algebras, and compact $64$-dimensional parametrizations. & Section~\ref{subsec:real-form} \\
$n_\pm$, $\Pi_1$, $\Pi_6$, $\Pi_8$, $\Pi_{12}$ & Values $N_\pm(E_u)$ and projections onto the $SU(3)$ sectors. & Section~\ref{subsubsec:su3-trace-sectors} \\
$\varsigma_u$, $\mathcal J_u$ & Euclidean real structure and its fixed real Albert algebra. & Section~\ref{subsubsec:euclidean-albert-form} \\
\multicolumn{3}{l}{\textbf{Bundles}}\\*
$\mathcal W_\varphi$, $\mathcal W_{\varphi,\mathbb C}$ & Deformation bundle and its complexification. & Definition~\ref{def:deformation-bundle} \\
$\mathrm{Fr}_\varphi$, $\mathrm{Fr}_{\varphi,\mathbb C}$, $f_p$ & Adapted frame bundles and the induced maps on deformation tensors. & Section~\ref{subsec:exceptional-reductions} \\
$\mathcal C_\varphi$, $\mathcal N_\pm$ & Bundle of invariant cubic planes and its critical line subbundles. & Theorem~\ref{thm:exceptional-reductions} \\
$\mathcal F_\pm$, $\mathcal F^c_\pm$ & Complex exceptional reductions and compact reductions. & Sections~\ref{subsec:exceptional-reductions}, \ref{subsec:compact-reductions} \\
$\ell$, $P_\ell$, $E_\ell$ & Tangent line field, its projection tensor, and $7P_\ell-g$. & \eqref{eq:El-def} \\
$\mathrm{Fr}_\ell$, $H_{E_1}$ & Line-adapted frames and stabilizer of the model basepoint. & \eqref{eq:line-adapted-frames} \\
$\mathcal F_{\pm,\ell}$, $\mathfrak J_{E_\ell}$ & Based $F_4$ reductions and Albert algebra bundles. & Theorem~\ref{thm:albert-bundle} \\
$\mathcal W_r(\ell)$ & Real sector subbundles, $r=6,8,12$. & Section~\ref{subsec:albert-bundles} \\
$\varsigma_{\ell,\pm}$, $\mathcal J_{\ell,\pm}$ & Global real structures and their Euclidean fixed bundles. & \eqref{eq:global-euclidean-real-structure} \\
\multicolumn{3}{l}{\textbf{Connections}}\\*
$D^{LC}$, $D^{\min}$, $T$ & Levi-Civita connection, minimal $G_2$-connection, and intrinsic torsion. & Section~\ref{subsec:connections-pointings} \\
$\beta$, $\omega$ & $\mathfrak g_2$-valued metric connection correction and exceptional connection form. & Section~\ref{subsec:connections-pointings} \\
$L_\varphi$, $\pi_\pm$, $\mathfrak c_\pm$ & Cubic variation map, exceptional projections, and $650$-dimensional complementary modules in $\mathfrak{sl}(\WC)$. & Section~\ref{subsec:connections-pointings} \\
$\tr_W(XY)$ & Invariant bilinear form used for the connection projection. & Section~\ref{subsec:connections-pointings} \\
\end{longtable}
\endgroup

\end{document}